%% file: thesis.tex
\documentclass[12pt]{puthesis}
\usepackage{amsmath,amsfonts,amsthm,amssymb,amscd,epigraph,ibycus4}
\newcommand{\st}{s.t.}
\newcommand{\quadrec}[2]{\left(\frac{#1}{#2}\right)}
\newcommand{\njac}[2]{\lbrack {#1},{#2}\rbrack}
\DeclareMathOperator{\sgn}{sgn}

\DeclareMathOperator{\ord}{ord}
\DeclareMathOperator{\av}{av}

\DeclareMathOperator{\num}{num}
\DeclareMathOperator{\den}{den}
\DeclareMathOperator{\lcm}{lcm}
\DeclareMathOperator{\Lcm}{Lcm}
\DeclareMathOperator{\rad}{rad}
\DeclareMathOperator{\mo}{\,mod}
\DeclareMathOperator{\sq}{sq}

\DeclareMathOperator{\vol}{vol}

\DeclareMathOperator{\SL}{SL}

\DeclareMathOperator{\Area}{Area}

\DeclareMathOperator{\Vol}{Vol}
\DeclareMathOperator{\Res}{Res}
\DeclareMathOperator{\Disc}{Disc}
\DeclareMathOperator{\Sp}{Sp}
\DeclareMathOperator{\charac}{char}
\DeclareMathOperator{\irr}{irr}
\DeclareMathOperator{\rnk}{rank}
\DeclareMathOperator{\BSD}{\mathit{BSD}}

\DeclareMathOperator{\Frob}{Frob}
\DeclareMathOperator{\Gal}{Gal}

\DeclareMathOperator{\Tr}{Tr}

\DeclareMathOperator{\cl}{cl}
\newtheorem{prop}{Proposition}[section]
\newtheorem{thm}[prop]{Theorem}

\newtheorem{cor}[prop]{Corollary}
\newtheorem*{conja1}{Conjecture $\mathfrak{A}_1$}
\newtheorem*{conja2}{Conjecture $\mathfrak{A}_2$}
\newtheorem*{conjb1}{Conjecture $\mathfrak{B}_1$}
\newtheorem*{conjb2}{Conjecture $\mathfrak{B}_2$}
\newtheorem*{conja1d}{Conjecture $\mathfrak{A}_1$}
\newtheorem*{conja2d}{Conjecture $\mathfrak{A}_2$}
\newtheorem*{conjb1d}{Hypothesis $\mathfrak{B}_1$}
\newtheorem*{conjb1e}{Hypothesis $\mathfrak{B}_1(K,P,\eta(N),\epsilon(N))$}
\newtheorem*{conjb2d}{Hypothesis $\mathfrak{B}_2$}
\newtheorem*{conjb2e}{Hypothesis $\mathfrak{B}_2(K,P,\eta(N),\epsilon(N))$}
\newtheorem*{conja1q}{Conjecture $\mathfrak{A}_1(K,P,\delta(N))$}
\newtheorem*{conja2q}{Conjecture $\mathfrak{A}_2(K,P,\delta(N))$}
\newtheorem*{thm11t}{Theorem 1.1$^\prime$}

\newtheorem*{thm13t}{Theorem 1.3$^\prime$}
\newtheorem*{thm14t}{Theorem 1.4$^\prime$}

\newtheorem{lem}[prop]{Lemma}
\newtheorem{defn}{Definition}
\newtheorem*{defn*}{Definition}
\newtheorem*{bilinc}{Bilinear condition}
\newenvironment{Rem}{{\bf Remark.}}{}
\numberwithin{equation}{section}
\author{Harald A. Helfgott}
\title{Root Numbers and the Parity Problem}
\submitted{June, 2003}
\abstract{\input{abstract}}
\acknowledgements{\input{acknow}}
\begin{document}
\setlength{\epigraphwidth}{170pt}
\input{intro}
\setlength{\epigraphwidth}{220pt}
\input{ell}

\setlength{\epigraphwidth}{300pt}
\input{lambda}

\input{lambirr}
\setlength{\epigraphwidth}{265pt}
\input{square}
\appendix
\input{appa}

\input{appb}

\input{biblio}
\end{document}

%% file: abstract.tex
Let $\mathcal{E}$ be a one-parameter family of elliptic curves over a number field $K$. It is natural
to expect the average root number of the curves
in the family to be zero. All known counterexamples to this folk conjecture
occur for families obeying a certain degeneracy condition.
We prove that the average root number is zero for a large class of families
of elliptic curves of fairly general type. Furthermore, we show that any
non-degenerate family $\mathcal{E}$ has average root number $0$, provided that
two classical arithmetical conjectures hold 
for two homogeneous polynomials with integral coefficients
constructed explicitly in terms of $\mathcal{E}$. 

The first such conjecture -- commonly associated with Chowla --
asserts the equidistribution of the
 parity of the number of primes dividing the integers
represented by a polynomial. More precisely: given a homogeneous polynomial
$f\in \mathbb{Z}\lbrack x,y\rbrack$, it is believed that $\mu(f(x,y))$
averages to zero. This conjecture can be said to represent the parity
problem in its pure form, while covering the same notional ground as the
Bunyakovsky-Schinzel and Hardy-Littlewood conjectures taken together.

For $\deg f = 1$ and $\deg f = 2$, Chowla's conjecture is essentially 
equivalent to the prime number theorem. For $\deg f > 2$, the conjecture
has been unproven up to now; the traditional approaches by means of analysis and sieve theory fail. We prove the conjecture for $\deg f = 3$.

There remains to state
the second arithmetical conjecture referred to previously. It is
believed that any non-constant homogeneous polynomial 
$f\in \mathbb{Z}\lbrack x,y\rbrack$ yields to a square-free sieve. We sharpen
the existing bounds on the known cases by a sieve refinement and a new 
approach combining height functions, sphere packings and sieve methods.

%% file: acknow.tex
\epigraph{\greek{w(`s a)'ra fwnh'sas po're fa'rmakon a)rgei+fo'nths\\
e)k gai'hs e)ru'sas, kai' moi fu'sin au)tou= e)'deice.\\
r(i'zh| me`n me'lan e)'ske, ga'lakti de` ei)'kelon a)'nqos:\\
mw=lu de' min kale'ousi qeoi': xalepo`n de' t' o)ru'ssein\\
a)ndra'si ge qnhtoi=si, qeoi` de' te pa'nta du'nantai.}
}{Homer, {\em Odyssey}, 10.302--10.306}
As my own words do not suffice to express my gratitude to my advisor,
Henryk Iwaniec, the reader is referred to the passage above. The present work, 
however, is dedicated to those who authored the author, namely, Michel Helfgott
and Edith Seier. To them, then, for love and geometry.

        I am indebted to Gergely Harcos
for his careful reading of early versions of the manuscript and for having
prodded me to put my thesis in its present form. The second reader of the
thesis, Peter Sarnak, has been helpful throughout my stay at Princeton. Thanks
are due as well to Keith Conrad, Jordan Ellenberg,
Chris Hall and Emmanuel Kowalski
 for their useful advice and to Keith Ramsay for our discussions 
on his unpublished work. This listing is not meant to be exhaustive.

%% file: intro.tex
\chapter{Introduction}
\epigraph{
Por qu\'e los \'arboles esconden\\el esplendor de sus ra\'{\i}ces?}{
Neruda, {\em El libro de las preguntas}}
\section{Root numbers of elliptic curves}
Let $E$ be an elliptic curve over $\mathbb{Q}$. The reduction $E$ mod $p$
can
\begin{enumerate}
\item be an elliptic curve over $\mathbb{Z}/p$,
\item have a node, or
\item have a cusp.
\end{enumerate}

        We call the reduction {\em good} in the first case, {\em multiplicative} in the
second case and {\em additive} in the third case. If the reduction is 
not good, then, as might be expected, we call it 
{\em bad}. If the reduction at $p$
is multiplicative, we call it {\em split} if the slopes at the node lie
in $\mathbb{Z}/p$, and {\em non-split} if they do not. Additive reduction
becomes either good or multiplicative in some finite 
extension of $\mathbb{Q}$. Thus every $E$ must fall into one
of two categories: either it has good reduction over a finite extension
of $\mathbb{Q}$, possibly $\mathbb{Q}$ itself, or it has multiplicative
reduction over a finite extension of $\mathbb{Q}$, possibly
$\mathbb{Q}$ itself. 
We speak accordingly of {\em potential good reduction} and {\em potential multiplicative reduction}. 

        The $L$-function of $E$ is defined to be
        \[L(E,s) = \prod_{\text{$p$ good}} (1 - a_p p^{-s} + p^{1-2s})^{-1} \prod_{\text{$p$ bad}} (1 - a_p p^{-s})^{-1},\]
where $a_p$ is $p+1$ minus the number of points in $E$ mod $p$. As can
be seen, $L(E,s)$ encodes the local behaviour of $E$. It follows from the modularity theorem 
(\cite{Wi}, \cite{TW}, \cite{BCDT}) that $L(E,s)$
has analytic continuation to all of $\mathbb{C}$ and satisfies the
following functional equation:
\[\mathcal{N}_E^{(2-s)/2} (2\pi )^{s-2} \Gamma(2-s) L(E,2-s) = W(E) \mathcal{N}_E^{s/2} (2\pi)^{-s} \Gamma(s) L(E,s),\]
where $W(E)$, called the {\em root number} of $E$, equals $1$ or $-1$,
and $\mathcal{N}_E$ is the conductor of $E$.
The function $L(E,s)$ corresponds to a modular form $f_E$ of weight $2$ and
level $\mathcal{N}_E$. The canonical
involution $W_{\mathcal{N}}$ acting on modular
forms of level $\mathcal{N}$ has $f_E$ as an eigenvector with eigenvalue $W(E)$.

The set $E(\mathbb{Q})$ of points on $E$ with rational coordinates is an abelian group under
the standard operation + (see e.g. \cite{Si}, III.1). A classical theorem of
 Mordell's states that $E(\mathbb{Q})$ is finitely generated. We define the 
{\em algebraic rank} of $E$ to be the rank of $E(\mathbb{Q})$.
We denote the algebraic rank of $E$ by $\rnk(E)$.
The Birch-Swinnerton-Dyer conjecture asserts that 
$\rnk(E)$ equals the order of vanishing $\ord_{s=1} L(E,s)$ of
$L(E,s)$ at $s=1$. Since $W(E)$ is one if the order of
vanishing is even and minus one if it is odd, the root number gives us
the parity of the algebraic rank, conditionally on the conjecture. This
fact makes the root number even more interesting than it already is on its own.
Assuming the Birch-Swinnerton-Dyer conjecture for curves of algebraic
rank zero, 
it suffices to prove that the root number of an elliptic curve is minus
one to show that the rank is positive. 
If, on the other hand, we prove
that $W(E)=1$ and find by other means that 
there are infinitely points on $E$, we have that the rank is ``high", that is, 
at least two. (Rank 2 is considered high as it may already be atypical in
certain contexts.)

        It is a classical result (\cite{De} -- cf. \cite{Ro2}, \cite{Ta}) that the root
number can be expressed as a product of local factors,
\[W(E) = \prod_v W_v(E),\]
where each $W_v(E)$ can be expressed in terms of a canonical
representation $\sigma'_{E,v}$ of the Weil-Deligne group of 
$\mathbb{Q}_v$:
\[W_v(E) = \frac{\epsilon(\sigma'_{E,v},\psi, d x)}{
|\epsilon(\sigma'_{E,v},\psi, d x)|},\]
where $\psi$ is any nontrivial unitary character of $\mathbb{Q}_p$ and
$d x$ is any Haar measure on $\mathbb{Q}_p$. This expression has been made
explicit in terms of the coefficients of $E$ (\cite{Ro}, \cite{Con},
\cite{Ha}). Thus many questions about the distribution
of $W(E) = (-1)^{\ord_{s=1} L(E,s)}$ have become somewhat more approachable
than the corresponding questions about the distribution of
$\ord_{s=1} L(E,s)$.

The natural expectation is that $W(E)$ be $1$
as often as $-1$ when $E$ varies within a family of elliptic curves that is 
in some sense typical or naturally defined. This is consistent with what
is currently known about average ranks and seems to have become a folk
conjecture (see for example \cite{Si2}, section 5). As we will see below,
families in which this is known not to hold are in some sense degenerate.

\section{Families of elliptic curves and questions of distribution}\label{ssec:famdef}
        By a {\em family}\, $\mathcal{E}$ of elliptic curves over $\mathbb{Q}$
on one variable we mean an elliptic
surface over $\mathbb{Q}$, 
or, equivalently, an elliptic curve over $\mathbb{Q}(t)$.
In the latter formulation, a family is given by two rational functions
$c_4,c_6\in \mathbb{Q}(t)$ such that $\Delta = (c_4^3-c_6^2)/1728$ is not
identically zero, and its fiber $\mathcal{E}(t)$ at a point $t\in \mathbb{Q}$ 
is the curve given by the equation
\[y^2 = x^3 - \frac{c_4(t)}{48} x - \frac{c_6(t)}{864}.\]
For finitely many $t\in \mathbb{Q}$, the curve $\mathcal{E}(t)$ will
be singular. In such a case we set $W(\mathcal{E}(t)) = 1$. 

Every primitive irreducible polynomial $Q\in \mathbb{Z}[t]$ determines a valuation
(or {\em place}) of $\mathbb{Q}(t)$. An additional valuation is given
by $\deg(\den) - \deg(\num)$, that is, the map taking an element
of $\mathbb{Q}(t)$ to the degree of its denominator minus the degree of
its numerator.
 Given a valuation $v$ of $\mathbb{Q}(t)$ and 
an elliptic curve $\mathcal{E}$
over $\mathbb{Q}(t)$, we can examine the reduction
$\mathcal{E} \mo v$ and give it a type in exactly the same way we have
described for reductions $E \mo p$: the curve $\mathcal{E}$ will be said to
have {\em good} reduction if its reduction at $v$ is an elliptic curve over
the residue field, resp. {\em multiplicative} reduction if the reduction
at $v$ has a node, {\em additive} if it has a cusp, {\em split
multiplicative} if it has a node and the slopes at the node
are in the residue field, {\em non-split multiplicative} if it is has a node
but the slopes at the node are not in the residue field,
{\em potentially good} if $\mathcal{E}$ has good reduction at the
place lying over $v$ in some finite extension of $\mathbb{Q}(t)$,
{\em potentially multiplicative} if $\mathcal{E}$ has multiplicative 
reduction at the place lying over $v$ in some 
finite extension of $\mathbb{Q}(t)$. The type of reduction of 
$\mathcal{E}$ at a given place $v$ 
can be determined by the usual valuative criteria
(see e.g. \cite{Si}, 179--183).

Define
\begin{equation}\begin{aligned}\label{eq:M}
M_\mathcal{E}(x,y) &= \prod_{\text{$\mathcal{E}$ has mult.\,red. at $v$}} P_v(x,y),
\end{aligned}\end{equation}
where $P_v=x$ if $v$ is the place $\deg(\den)-\deg(\num)$,
$P_v=x^{\deg Q} Q\left(\frac{y}{x}\right)$ 
if $v$ is a valuation given by a primitive 
irreducible polynomial $Q\in \mathbb{Z}[t]$.

Given a function $f:\mathbb{Z}\to \mathbb{C}$ and an arithmetic progression 
$a + m \mathbb{Z}$, we define
\[\av_{a + m \mathbb{Z}} f = \lim_{N\to \infty}
\frac{1}{N/m} \mathop{\sum_{1\leq n\leq N}}_{n\equiv a \mo m}
 f(n) .\]
If $\av_{a + m \mathbb{Z}} f = 0$ for all 
$a\in \mathbb{Z}$, $m\in \mathbb{Z}^+$,
we say that $f$ {\em averages to zero over the integers.}
Given a function $f:\mathbb{Z}^2\to \mathbb{C}$, a lattice coset 
$L\subset \mathbb{Z}^2$ and a sector $S\subset \mathbb{R}^2$ (see section
\ref{sec:notpre}), we define
\[\av_{S\cap L} f = \lim_{N\to \infty}
 \frac{1}{\#(S\cap L\cap \lbrack -N,N\rbrack^2)}
 \sum_{(x,y)\in S\cap L\cap \lbrack -N,N\rbrack^2} f(x,y) .\]
We say that $f$ {\em averages to zero over $\mathbb{Z}^2$} if 
$\av_{S\cap L} f = 0$ for all choices of $S$ and $L$.
Given a function $f:\mathbb{Q}\to \mathbb{Z}$, a lattice
coset $L\subset \mathbb{Z}^2$ and a sector $S\subset \mathbb{R}^2$, we define
\[\av_{\mathbb{Q},S\cap L} f =
\lim_{N\to \infty}
 \frac{
 \sum_{(x,y)\in S\cap L\cap \lbrack -N,N\rbrack^2,\, \gcd(x,y)=1} f(y/x) 
}{\# 
\{(x,y)\in S\cap L\cap \lbrack -N,N\rbrack^2 :
\gcd(x,y)=1\}} .\]
 We say that $f$ {\em averages to zero over the rationals} if
$\av_{\mathbb{Q},S\cap L} f = 0$ for all choices of $S$ and $L$.
We are making our definition of zero average strict enough for it to be invariant under fractional linear transformations. Moreover, by
letting $S$ be arbitrary, we allow sampling to be restricted to any
open interval in $\mathbb{Q}$. Thus our results will not be imputable
to peculiarities in averaging order or to superficial cancellation.

In the literature, a family $\mathcal{E}(t)$ for which 
\[j(\mathcal{E}(t))=\frac{c_4(\mathcal{E}(t))^3}{\Delta(\mathcal{E}(t))}\]
is constant is sometimes called a {\em constant family}. 
When examples were found
(\cite{Ro}, \cite{Riz1}) of constant families of elliptic curves in which the 
root number did not average to zero, it seemed plausible that such behaviour might be a degeneracy peculiar to constant families. 
It was thus somewhat of a surprise when non-constant families of non-zero average root number were found to exist.

        All non-constant families considered in \cite{Ma} and \cite{Riz2}
had $M_\mathcal{E}(x,y)$ equal to a constant, or, what is the same,
$M_\mathcal{E}(x,y)=1$; in other words, they had no places
of multiplicative reduction as elliptic curves over $\mathbb{Q}(t)$.
Families with non-constant $M_\mathcal{E}$ were hardly touched upon,
as they were felt to present severe number-theoretical difficulties
(see, e. g., \cite{Ma}, p. 34, third paragraph). The subject of the present
work is precisely such families. 

        We will see how families with non-constant $M_\mathcal{E}$ are not
only heuristically different from families with constant $M_\mathcal{E}$
but also quite different in their behaviour. As we will prove -- 
in some cases conditionally on two standard conjectures in analytic
number theory, and in the other cases unconditionally -- 
$W(\mathcal{E}(t))$ averages to zero over the integers and over
the rationals for any family of elliptic
curves $\mathcal{E}$ with non-constant $M_\mathcal{E}$. 
All autocovariances of $W(\mathcal{E}(t))$ other than the variance
are zero as well. In other words, for any family $\mathcal{E}$ with
at least one place of multiplicative reduction over $\mathbb{Q}(t)$,
the function $t\mapsto W(\mathcal{E}(t))$ behaves essentially like white noise.

        We may thus see the constancy of $M_\mathcal{E}$ as the
proper criterion of degeneracy for our problem. The generic case is
 that of non-constant $M_\mathcal{E}$: for
a typical pair of polynomials or rational functions $c_4(t)$, $c_6(t)$, 
the numerator of the discriminant $\Delta = (c_4(t)^3 - c_6(t)^2)/1728$ does in general have polynomial factors not present in $c_4(t)$ or $c_6(t)$. Any such factor will be present in $M_\mathcal{E}$ as well,
making it non-constant.

\section{Issues and definitions}\label{sec:isdf}

        The main analytical difficulty in case (3) lies in the parity of the number of primes
dividing an integer represented by a polynomial. A precise discussion necessitates
some additional definitions.

 We will say that the reduction of $\mathcal{E}$ at $v$ is {\em quite bad} if there is no
non-zero rational function $d(t)$ for which the family
\[\mathcal{E}_{d} : d(t) y^2 = x^3 - \frac{c_4(t)}{48} - \frac{c_6(t)}{864}\]
has good reduction at $v$. If the reduction of $\mathcal{E}$ at $v$ is bad but not
quite bad, we say it is {\em half bad}.

We let
\begin{equation}\begin{aligned}\label{eq:bb'}
B_\mathcal{E}(x,y) &= \prod_{\text{$\mathcal{E}$ has bad red. at $v$}} P_v(x,y),\\
B'_\mathcal{E}(x,y) &= \prod_{\text{$\mathcal{E}$ has quite bad red. at $v$}} P_v(x,y) ,\end{aligned}\end{equation}
where $P_v$ is as in (\ref{eq:M}). It follows immediately from the 
definitions that $M_\mathcal{E}(x,y)$, $B_\mathcal{E}(x,y)$
and $B'_\mathcal{E}(x,y)$ are square-free and can be constant
only if identically equal to one. By saying that a polynomial $P$ 
is square-free we mean that no irreducible non-constant polynomial $P_i$ 
appears in the factorization $P = P_1  P_2  \dotsb P_n$ more than once.

Given a function $f:\mathbb{Z}\to \{-1,1\}$, 
a non-zero integer $k$ and an arithmetic
progression $a + m \mathbb{Z}$, we define
\[\gamma_{a + m \mathbb{Z}, k}(f) = 
 \lim_{N\to \infty} \frac{1}{N/m}
\mathop{\sum_{1\leq n\leq N}}_{n\equiv a \mo m} f(n) f(n + k) .
\]
If $\av_{\mathbb{Z}} f = 0$, then $\gamma_{\mathbb{Z}, k}(f)$
equals the $k$th autocorrelation and the $k$th autocovariance of 
the sequence $f(1), f(2), f(3),\dotsb$. (Note that, since $f(n) = \pm 1$ 
for all $n$, the concepts of autocorrelation
and autocovariance coincide when $\av_\mathbb{Z} f =0$.)
We say that $f$ is {\em white noise over the integers}
if $\av_{a + m \mathbb{Z}} f = 0$ and
$\av_{a + m \mathbb{Z}, k} f = 0$
 for all choices of $a + m \mathbb{Z}$ and $k$.
 
Given a function $f:\mathbb{Z}\to \{-1,1\}$, 
a lattice coset $L\subset \mathbb{Z}^2$, a sector $S\subset \mathbb{R}^2$ and
a non-zero rational $t_0$,
we define
\[\gamma_{L\cap S,t}(f) = 
\lim_{N\to \infty}
 \frac{\sum_{(x,y)\in S\cap L\cap \lbrack -N,N\rbrack^2,\, \gcd(x,y)=1}
f\left(\frac{y}{x}\right) f\left(\frac{y}{x} + t\right)}{\# 
\{(x,y)\in S\cap L\cap \lbrack -N,N\rbrack^2 :
\gcd(x,y)=1\}}
 .\]
We say that {\em $f$ is white noise over the rationals} if
$\av_{\mathbb{Q}, L\cap S} f = 0$ and 
$\gamma_{L\cap S,t}(f) = 0$ for all choices of $L$, $S$ and $t$.

We can now list all questions addressed here and in the previous
literature as follows:
\begin{enumerate}
\item are $\{t\in \mathbb{Q} : W(\mathcal{E}(t))=1\}$ and
$\{t\in \mathbb{Q} : W(\mathcal{E}(t)) = -1\}$ both infinite?
\item are $\{t\in \mathbb{Q} : W(\mathcal{E}(t))=1\}$ and
$\{t\in \mathbb{Q} : W(\mathcal{E}(t)) = -1\}$ both dense in $\mathbb{Q}$?
\item does $W(\mathcal{E}(t))$ average to zero over the integers?
\item is $W(\mathcal{E}(t))$ white noise over the integers?
\item does $W(\mathcal{E}(t))$ average to zero over 
the rationals?
\item is $W(\mathcal{E}(t))$ white noise over the rationals?
\end{enumerate}

Evidently, an affirmative answer to (2) implies one to (1). An affirmative
answer to (5) implies that the answers to (1) and (2) are ``yes'' as well.

\section{The square-free sieve}
Starting with \cite{GM} and \cite{Ro}, the square-free sieve has appeared
time and again in the course of nearly every endeavour
to answer any of the questions above.
It seems by now to be an analytic difficulty that cannot be avoided.

\begin{defn}
We say that a polynomial $P\in \mathbb{Z}\lbrack x \rbrack$ {\em yields to a square-free
sieve} if
\begin{equation}\label{eq:sqfr}\lim_{N\to \infty} \frac{1}{N} \#\{1\leq x\leq N:
\exists p>N^{1/2} \text{\,\st\,} p^2|P(x)\} = 0.\end{equation}
We say that a homogeneous polynomial $P\in \mathbb{Z}\lbrack x,y\rbrack$ 
{\em yields to a square-free sieve}
\begin{equation}\lim_{N\to \infty} \frac{1}{N^2} 
\# \{-N\leq x,y \leq N : \gcd(x,y)=1, 
 \exists p> N \text{\,\st\,} p^2|P(x,y)\} = 0 .
\end{equation}
\end{defn}

There is very little we can say unconditionally
about a family $\mathcal{E}$ unless we can prove that
$B'_\mathcal{E}(x,y)$ yields to a square-free sieve. 
\begin{conja1}
Every square-free polynomial $P\in \mathbb{Z}\lbrack x\rbrack$ 
yields to a square-free sieve.
\end{conja1}
\begin{conja2}
Every square-free homogeneous polynomial $P\in \mathbb{Z}\lbrack x,y\rbrack$ 
yields to a square-free sieve.
\end{conja2}

Conjecture $\mathfrak{A}_1(P)$ is clear for $P$ linear.\footnote{By $X(P)$ 
we denote the validity of a conjecture $X$ for a specific polynomial $P$.
Thus Conjecture $\mathfrak{A}_1(P)$ is the same as the statement ``$P$ yields to
a square-free sieve.''}
Estermann \cite{Es} proved it for $\deg(P)=2$.
Hooley (\cite{Ho}, Chapter 4) proved it for $\deg(P)=3$. By then it was
expected that $\mathfrak{A}_1$ would hold for any square-free polynomial;
in some sense $\mathfrak{A}_1$ and $\mathfrak{A}_2$
 are much weaker than the conjectures $\mathfrak{B}_1$, $\mathfrak{B}_2$
to be treated in section \ref{sec:chow}, though $\mathfrak{B}_i$
does not imply $\mathfrak{A}_i$.
Greaves \cite{Gr} proved $\mathfrak{A}_2(P)$ for $\deg(P)\leq 6$. 
Both Hooley's and Greaves' bounds on
the speed of convergence of (\ref{eq:sqfr}) will be strengthened in Chapter \ref{chap:square}.

Note that, if $P_1$ and $P_2$ have no factors in common 
and $\mathfrak{A_i}(P_1)$ and $\mathfrak{A_i}(P_2)$ both hold, then $\mathfrak{A_i}(P_1 P_2)$ holds. Let \[\deg_{\irr}(P) = \max_i\, \deg(Q_i),\] where 
$P = Q_1^{k_1} Q_2^{k_2}\dotsb Q_n^{k_n}$ is the decomposition of $P$ into irreducible factors. Given this notation,
 we can say that we know $\mathfrak{A}_1(P)$ for
$\deg_{\irr}(P)\leq 3$ and $\mathfrak{A}_2(P)$ for
$\deg_{\irr}(P)\leq 6$.

Granville has shown \cite{Gran} that Conjectures $\mathfrak{A}_1$ and $\mathfrak{A}_2$ follow
in general from the $abc$ conjecture. 
 Unlike the unconditional results
just mentioned, this general conditional result
does not give us any explicit bounds.

\section{Previous results}\label{subs:prevres}
We can now state what is known about the answers to the questions posed at the
end of section \ref{sec:isdf}. A family $\mathcal{E}$ will present one
of three very different kinds of behaviour depending on whether
$j(\mathcal{E}(t))$ or $M_\mathcal{E}(x,y)$ is constant. Notice
that, if $j(\mathcal{E}(t))$ is constant, then $M_{\mathcal{E}}(x,y)$
is constant.
\begin{enumerate}
\item $j$ constant\\In this case $\mathcal{E}$ consists of
quadratic twists
\[\mathcal{E}_{d}(t) : d(t) y^2 = x^3 - \frac{c_4}{48} x - \frac{c_6}{864}\]
of a fixed elliptic curve over $\mathbb{Q}$. Rohrlich
\cite{Ro} showed that, depending on the twisting function $d$, either 
(1) $\{t\in \mathbb{Q} : W(\mathcal{E}(t)) = 1\}$ and 
$\{t\in \mathbb{Q} : W(\mathcal{E}(t)) = - 1\}$ 
are both dense in $\mathbb{Q}$, or
(2) $W(\mathcal{E}(t))$ is constant on $\{t\in \mathbb{Q}: d(t)>0\}$ and
on $\{t\in \mathbb{Q}: d(t)<0\}$. Rizzo \cite{Riz1} pointed out
that in the latter case the set of values of
$\av_{\mathbb{Q}} W(\mathcal{E}(t))$
 for different functions $d$ is dense is $[-1,1]$.
\item $j$ non-constant, $M_\mathcal{E}$ constant\\ 
Here Manduchi showed \cite{Ma} that 
$\{t\in \mathbb{Q}: d(t)>0\}$ and
on $\{t\in \mathbb{Q}: d(t)<0\}$
are both dense 
provided that Conjecture $\mathfrak{A}_2(B'_\mathcal{E})$ holds.
 
Rizzo has given examples \cite{Riz2} of families $\mathcal{E}$
with non-constant $j$ and $M_\mathcal{E}=1$ such that 
$\av_{\mathbb{Z}} W(\mathcal{E}(n)) \ne 0$.
In section \ref{subs:fallada}
we will see an example of a familiy with non-constant $j$, $M_\mathcal{E}=1$
and $\av_{\mathbb{Q}} W(\mathcal{E}(t)) \ne 0$.

\item $j$ and $M_\mathcal{E}$ non-constant\\
Manduchi \cite{Ma} showed that, if $\deg(M_\mathcal{E})=1$, then both 
$\{t\in \mathbb{Q}: d(t)>0\}$ and
on $\{t\in \mathbb{Q}: d(t)<0\}$
are infinite. Nothing else has been known until now for this case.

\end{enumerate}

The main difference between cases (1) and (2), on the other hand, and
case (3), on the other, can be roughly outlined as follows.
Assume $M_\mathcal{E}$ is constant; in other
words, assume that $\mathcal{E}$ has no places of multiplicative reduction
when considered as an elliptic curve over $\mathbb{Q}(t)$.
Then for every $\epsilon$ there is a finite set $S$ of primes such
that, for any large $N$, the elliptic curve $\mathcal{E}(t)$ can have
multiplicative reduction at places $p$ not in $S$ only for a proportion
less than $\epsilon$ of all values of $t$. As will become clear later,
this eliminates what would
otherwise be the analytical heart of the matter, namely, the
estimation of 
\[\prod_{\text{$p$ mult.}} W_p(\mathcal{E}(t)) ,\]
that is, the product of the local root numbers at the places $p$
of multiplicative reduction.

\section{A conjecture of Chowla's. The parity problem}\label{sec:chow}
Our main purpose is to determine the behavior
of the root number in families $\mathcal{E}$ with $M_\mathcal{E}$
non-constant. We will see that in this case all issues raised in
Section \ref{sec:isdf} amount to a classical arithmetical question
 in disguise. Consider the {\em Liouville function}
\[\lambda(n) = \begin{cases}
\prod_{p|n} (-1)^{v_p(n)} &\text{if $n\ne 0$}\\
0 &\text{if $n=0$.}\end{cases}\]

\begin{conjb1}
Let $P\in \mathbb{Z}\lbrack x \rbrack$ be a polynomial not of the
form $c Q^2(x)$, $c\in \mathbb{Z}$, $Q\in \mathbb{Z}[x]$.
 Then $\lambda(P(n))$ averages to zero over the integers.
\end{conjb1}
\begin{conjb2}
Let $P\in \mathbb{Z}\lbrack x,y \rbrack$ be a homogeneous
polynomial not of the form $c Q^2(x,y)$, $c\in \mathbb{Z}$,
$Q\in \mathbb{Z}[x,y]$. Then $\lambda(P(x,y))$ averages to zero
over $\mathbb{Z}^2$.
\end{conjb2}

In the present form, Conjecture $\mathfrak{B}_1$ is credited to S. Chowla.
(Some cases of $\mathfrak{B}_1$ were already included in the 
Hardy-Ramanujan conjectures.) As
stated in \cite{Ch}, p. 96: 
\begin{quotation}
If \lbrack $P$ is linear, Conjecture $\mathfrak{B}_1(P)$\rbrack\, is equivalent to
the Prime Number Theorem.
If \lbrack the degree of $P$\rbrack\, is at least $2$,
this seems an extremely hard conjecture.
\end{quotation}
In fact $\mathfrak{B}_1(x(x+1))$ is commonly considered to be roughly as hard as the Twin 
Prime Number conjecture. 

Conjecture $\mathfrak{B}_2(P)$ is equivalent to the Prime Number Theorem when
$P$ is linear. In the case of $P$ quadratic, the main ideas needed for a proof
of $\mathfrak{B}_2(P)$ were supplied by
de la Vall\'ee-Poussin (\cite{DVP1}, \cite{DVP2}) and Hecke (\cite{Hec}).
(We provide a full treatment in section \ref{sec:parquad}.)
The attacks on $\mathfrak{B}_1(P)$ for $\deg(P)=1$ and on $\mathfrak{B}_2(P)$ for
$\deg(P)=1,2$ rely on the fact that one can reduce the problem to a
question about $L$-functions. This approach breaks down for
$\mathfrak{B}_1(P)$, $\deg(P)>1$ and $\mathfrak{B}_2(P)$, $\deg(P)>2$, as there seems to be no
analytic object corresponding to $P\in \mathbb{Z}[x]$, $\deg P>1$
or $P\in \mathbb{Z}[x,y]$, $\deg P>2$.

A classical sieve treatment of conjectures $\mathfrak{B}_1$ and $\mathfrak{B}_2$ is doomed to fail; they may be said to represent the parity problem in its purest form.
(The {\em parity problem} is the fact that,
as was pointed out by Selberg \cite{Se2}, 
a standard sieve framework cannot distinguish
between numbers with an even number of prime factors and
numbers with an odd number of prime factors.) Until recently, the
parity problem was seen as an unsurmountable difficulty whenever
the sets to be examined were sparser than the integers. The
sets in question here are $S_1(P)=\{P(n):n\in \mathbb{Z}\}$
and $S_2(P)=\{P(n,m):n,m\in \mathbb{Z}\}$. For a set $S\in \mathbb{Z}$,
define the {\em logarithmic density} $d(S)$ to be
\[d(S) = \lim_{N\to \infty} \frac{\log\,(\#\{x\in S: |x|<N\})}{\log N}\]
when defined. A set $S$ is said to be sparser than the integers if
$d(S)<1$. Since $d(S_1(P))=1/\deg(P)$ and $d(S_2(P))=2/\deg(P)$,
the set $S_1$ is sparser than the integers for $\deg(P)>1$ and
$S_2$ is sparser than the integers for $\deg(P)>2$.

We prove
conjecture $\mathfrak{B}_2(P)$ for $\deg(P)=3$.
For $P$ irreducible, the approach taken follows the same lines
as the novel results of the last few years 
(\cite{FI1}, \cite{FI2}, \cite{HB}, \cite{HBM}) 
on the number of primes represented by a polynomial. Friedlander and
Iwaniec (\cite{FI1}, \cite{FI2}) broke through the difficulties
imposed by the parity problem in proving that there are infinitely
many primes of the form $x^2+y^4$. While the specifics in their
extremely ingenious method do not seem to carry over simply to any
other polynomial, Heath-Brown (\cite{HB}) succeeded
in proving the existence of infinitely many primes of the form
$x^3+2 y^3$ while following akin general lineaments. In the same
way, while $\mathfrak{B}_2(P)$ for $\deg P=3$ demands a great deal of ad-hoc work,
it can be said to be a new instance of the general approach of
Friedlander and Iwaniec. Note that one cannot deduce $\mathfrak{B}_2(P)$, $\deg P=3$ from the corresponding result about the existence or number of primes
represented by $P$; such an implication exists only for $\deg P=1$.
For $\mathfrak{B}_2(P)$, $P$ reducible, there is not even a corresponding question on
prime numbers, and in fact the methods used then are quite different
from those for $P$ irreducible.

\section{Results}
By {\bf Theorem 0.0} ($\mathfrak{X}(P)$, $\mathfrak{Y}(Q)$)
we mean a theorem conditional on conjectures $\mathfrak{X}$ and $\mathfrak{Y}$
in so far as they concern the objects $P$ and $Q$, respectively. A result
whose statement does not contain parentheses after the numeration should
be understood to be unconditional. 

\begin{thm}[$\mathfrak{A}_1(B_\mathcal{E}'(1,t))$,
$\mathfrak{B}_1(M_\mathcal{E}(1,t))$]\label{thm:nav}
Let $\mathcal{E}$ be a family of elliptic curves over $\mathbb{Q}$
on one variable. Assume that 
$M_\mathcal{E}(1,t)$ is not constant. Then $W(\mathcal{E}(t))$ averages
to zero over the integers.
\end{thm}

\begin{thm}[$\mathfrak{A}_1(B_\mathcal{E}'(1,t))$,
$\mathfrak{B}_1(M_\mathcal{E}(1,t) M_\mathcal{E}(1,t+k))$ for all non-zero
$k\in \mathbb{Z}$]
Let $\mathcal{E}$ be a family of elliptic curves over $\mathbb{Q}$
on one variable. Let $k$ be an integer other than zero.
Assume that 
$M_\mathcal{E}(1,t)$ is not constant. 
Then $W(\mathcal{E}(t))$ is white noise over the integers.
\end{thm}

\begin{thm}[$\mathfrak{A}_2(B_\mathcal{E}')$,
$\mathfrak{B}_2(M_\mathcal{E})$]\label{thm:qav}
Let $\mathcal{E}$ be a family of elliptic curves over $\mathbb{Q}$
on one variable. Assume that 
$M_\mathcal{E}$ is not constant. Then $W(\mathcal{E}(t))$ averages
to zero over the rationals.
\end{thm}

\begin{thm}[$\mathfrak{A}_2(B_\mathcal{E}')$,
$\mathfrak{B}_2(M_\mathcal{E}(x,y) 
M_\mathcal{E}(k_0 x, k_0 y + k_1 x))$ for all non-zero $k_0\in \mathbb{Z}$ and all
$k_1\in \mathbb{Z}$]
Let $\mathcal{E}$ be a family of elliptic curves over $\mathbb{Q}$
on one variable. Let $k = k_1/k_0$ be a non-zero rational number,
$\gcd(k_0,k_1)=1$.
Assume that 
$M_\mathcal{E}$ is not constant.  Then $W(\mathcal{E}(t))$ is white noise
over the rationals.
\end{thm}

The unconditional cases of the theorems above 
can be stated as follows.

\begin{thm11t}
Let $\mathcal{E}$ be a family of elliptic curves over $\mathbb{Q}$
on one variable. Assume
$\deg_{\irr}(B'_\mathcal{E}(1,t))\leq 3$ and
$\deg(M_\mathcal{E}(1,t))=1$.
Then $W(\mathcal{E}(t))$ averages to zero over the integers.
Explicitly, for any arithmetic progression $a + m \mathbb{Z}$,
$m\leq (\log N)^{A_1}$,
\[ 
\av_{a + m \mathbb{Z}}
W(\mathcal{E}(n))\ll
\begin{cases} (\log N)^{-A_2}
&\text{if $\deg_{\irr}(B'_\mathcal{E}(1,t))=1,2$,}\\
(\log N)^{-0.5718\dotsc}
&\text{if $\deg_{\irr}(B'_\mathcal{E}(1,t))=3$,}\end{cases}\]
where $A_1$ and $A_2$ are arbitrarily large constants, and the implicit
constant depends only on $\mathcal{E}$, $A_1$ and $A_2$.
\end{thm11t}

\begin{thm13t}
Let $\mathcal{E}$ be a family of elliptic curves over $\mathbb{Q}$
on one variable. Assume that 
$M_\mathcal{E}$ is not constant. Suppose that
$\deg_{\irr}(B'_\mathcal{E})\leq 6$ and $\deg(M_\mathcal{E})\leq 3$. Then
$W(\mathcal{E}(t))$ averages to zero over the rationals.
Explicitly, for any sector $S\subset \mathbb{R}^2$ and every lattice
coset $L\subset \mathbb{Z}^2$ of index $\lbrack \mathbb{Z}^2 : L\rbrack
\leq (\log N)^{A_1}$, we have that
$\av_{\mathbb{Q}, S\cap L} (W(\mathcal{E}(t)))$ is bounded above by
\[\begin{aligned}
C \cdot (\log N)^{-A_2}
&\,\,\text{if $\deg_{\irr}(B'_\mathcal{E})\leq 5$,
$\deg(M_\mathcal{E})=1,2$,}\\
C \cdot \frac{\log \log N}{\log N}
&\text{
if $\deg_{\irr}(B'_\mathcal{E})\leq 5$,
$\deg(M_\mathcal{E})=3$, $M_{\mathcal{E}}$ red.,}\\
C \cdot 
\frac{(\log \log N)^5 (\log \log \log N)}{\log N}
&\text{
if $\deg_{\irr}(B'_\mathcal{E})\leq 5$,
$\deg(M_\mathcal{E})=3$, $M_{\mathcal{E}}$ irr.,}\\
 C \cdot (\log N)^{-1/2} &\text{\,\,if $\deg_{\irr}(B'_\mathcal{E})=6$,
$\deg(M_\mathcal{E})\leq 3$,}
\end{aligned}\]
where $A_1$ and $A_2$ are arbitrarily large constants, and $C$ 
depends only on $\mathcal{E}$, $S$, $A_1$ and $A_2$.
\end{thm13t}

\begin{thm14t}
Let $\mathcal{E}$ be a family of elliptic curves over $\mathbb{Q}$
on one variable.
 Suppose that
$\deg_{\irr}(B'_\mathcal{E})\leq 6$ and $\deg(M_\mathcal{E})=1$.
Then $W(\mathcal{E})$ is white noise over the rationals.
Explicitly, for any sector $S\subset \mathbb{R}^2$, any lattice
coset $L\subset \mathbb{Z}^2$ of index $\lbrack \mathbb{Z}^2 : L\rbrack
\leq (\log N)^{A_1}$, and any non-zero rational number $t_0$,
 we have that \[\gamma_{L\cap S, t_0}(W(\mathcal{E}(t)))
 \ll
 \begin{cases}
(\log N)^{-A_2} &\text{if $\deg_{\irr}(B'_\mathcal{E})\leq 5$,}\\
(\log N)^{-0.5718\dotsc} &\text{if $\deg_{\irr}(B'_\mathcal{E}) = 6$,}
\end{cases}
\]
where $A_1$ and $A_2$ are arbitrarily large constants, and the implied
constant depends only on $\mathcal{E}$, $S$, $A_1$ and $A_2$.
\end{thm14t}

By $\BSD(E)$ we denote
the validity of the Birch-Swinnerton-Dyer conjecture for the elliptic curve
$E$ over $\mathbb{Q}$.
As consequences of Theorems \ref{thm:nav} and \ref{thm:qav}, we have

\begin{cor}[$\mathfrak{A}_1(B_\mathcal{E}'(1,t))$,
$\mathfrak{B}_1(M_\mathcal{E}(1,t))$,
$\BSD(\mathcal{E}(t))$
for every $t\in \mathbb{Z}$]\label{lem:small1}
Let $\mathcal{E}$ be a family of elliptic curves over $\mathbb{Q}$
on one variable. Assume that $j(\mathcal{E}(t))$ and
$M_\mathcal{E}(1,t)$ are not constant. 
Then
\[\av_{\mathbb{Z}} \rnk(\mathcal{E}(t)) \geq \rnk(\mathcal{E}) + 1/2\]
for every interval $I\subset \mathbb{R}$.
\end{cor}
\begin{cor}[$\mathfrak{A}_2(B_\mathcal{E}')$,
$\mathfrak{B}_2(M_\mathcal{E})$, $\BSD(\mathcal{E}(t))$
for every $t\in \mathbb{Q}$]\label{lem:small2}
Let $\mathcal{E}$ be a family of elliptic curves over $\mathbb{Q}$
on one variable. Assume that $j(\mathcal{E}(t))$ and
$M_\mathcal{E}$ are not constant. 
Then
\[\av_I \rnk(\mathcal{E}(t)) \geq \rnk(\mathcal{E}) + 1/2\]
for every interval $I\subset \mathbb{R}$.
\end{cor}

For conditional upper bounds on $\av_{\mathbb{Z}} \rnk(\mathcal{E}(t))$
and a general discussion of what is currently believed about
the distribution of $\rnk(\mathcal{E}(t))$, see \cite{Si2}.

From Corollaries \ref{lem:small1} and \ref{lem:small2} we obtain the 
following two statements, which are far weaker than the preceding but,
in general, seem
to be still inaccessible otherwise.

\begin{cor}[$\mathfrak{A}_1(B_\mathcal{E}'(1,t))$,
$\mathfrak{B}_1(M_\mathcal{E}(1,t))$,
$\BSD(\mathcal{E}(t))$
for every $t\in \mathbb{Z}$]\label{lem:tiny1}
Let $\mathcal{E}$ be a family of elliptic curves over $\mathbb{Q}$
on one variable. Assume that $j(\mathcal{E}(t))$ and
$M_\mathcal{E}(1,t)$ are not constant. 
Then $\mathcal{E}(t)$ has infinitely many rational points for infinitely many
$t\in \mathbb{Z}$.
\end{cor}
\begin{cor}[$\mathfrak{A}_2(B_\mathcal{E}')$,
$\mathfrak{B}_2(M_\mathcal{E})$, $\BSD(\mathcal{E}(t))$
for every $t\in \mathbb{Q}$]\label{lem:tiny2}
Let $\mathcal{E}$ be a family of elliptic curves over $\mathbb{Q}$
on one variable. Assume that $j(\mathcal{E}(t))$ and
$M_\mathcal{E}$ are not constant. 
Then $\mathcal{E}(t)$ has infinitely many rational points for infinitely many
$t\in \mathbb{Q}$.
\end{cor}

The reader may wonder whether it is possible to dispense with conjectures
$\mathfrak{B}_1$, $\mathfrak{B}_2$ and still obtain results along the lines
of Theorems \ref{thm:nav} and \ref{thm:qav}. That this is not the case
is the import of the following two results.

\begin{prop}[$\mathfrak{A}_1(B'_\mathcal{E})$]\label{prop:forgot1}
Let $\mathcal{E}$ be a family of elliptic curves over $\mathbb{Q}$ on one
variable. Assume that $M_\mathcal{E}(1,t)$ is not
constant. Suppose that $W(\mathcal{E}(t))$ averages to zero over the integers.
Then $\mathfrak{B}_1(M_\mathcal{E}(1,t))$ holds.
\end{prop}

\begin{prop}[$\mathfrak{A}_2(B'_\mathcal{E})$]\label{prop:forgot2}
Let $\mathcal{E}$ be a family of elliptic curves over $\mathbb{Q}$ on one
variable. Assume that $M_\mathcal{E}$ is not
constant. Suppose that $W(\mathcal{E}(t))$ averages to zero over the rationals.
 Then $\mathfrak{B}_2(M_\mathcal{E})$ holds.
\end{prop}

Thus, if we assume $\mathfrak{A}_1$ and $\mathfrak{A}_2$, or the 
$abc$-conjecture, which implies them, we have that the problem of averaging
the root number is equivalent to the problem of averaging $\lambda$ over
the values taken by a polynomial.

\section{Families of curves over number fields}\label{subs:whenneedit}

We know considerably less about elliptic curves over arbitrary
number 
fields than we do about elliptic curves over $\mathbb{Q}$.
The $L$-function of an elliptic curve $E$ 
over a number field $K$ is known to have
a functional equation only for some special choices of $E$
over totally real number fields other than $\mathbb{Q}$ \cite{SW}. 
Nevertheless,
we know that, if the $L$-function of an elliptic curve over a number 
field $K$ does have a functional equation, its sign must be equal to
the product of the local root numbers \cite{De}. Thus we can simply define
the root number of an elliptic curve $E$ over $K$ as the product of the local
root numbers $W_\mathfrak{p}(E)$, knowing that the sign of a hypothetical
functional equation would have to equal such a product.

Let $E$ be an elliptic curve over a number field $K$.
The local root numbers $W_\mathfrak{p}(E)$ have been explicited by
Rohrlich \cite{Rog} for every prime $\mathfrak{p}$ not dividing $2$ or $3$.
To judge from Halberstadt's tables for $K=\mathbb{Q}$, $\mathfrak{p} = 2,3$
\cite{Ha}, a solution for $\mathfrak{p}|2,3$ and arbitrary $K$ is likely to
admit only an exceedingly unwieldy form. One of our results
(Proposition \ref{prop:corona}) will allow us to ignore $W_\mathfrak{p}$
for finitely many $\mathfrak{p}$, and, in particular, for all 
$\mathfrak{p}$ dividing $2$ or $3$. Due to this simplification, we will find
working with root numbers over number fields no 
harder than working with root numbers over the rationals. 

Averaging is a different matter. It is not immediately clear what kind of
average should be taken when the elliptic surface in question is defined over
$K\ne \mathbb{Q}$. Should one take the average root number of
the fibers lying over $\mathbb{Z}$ or $\mathbb{Q}$, as before? Or
should one take the average over all fibers, where the base $K$
is ordered by norm? (It is not clear what this would mean when $K$ has
real embeddings.) Or should one consider all elements of the base
inside a box in $K \otimes_{\mathbb{Q}} \mathbb{R}$? The basic descriptive
machinery presented in Section \ref{sec:plifun} is independent of the kind
of average settled upon. As our main purpose in generalizing our results
is to understand the root number better,
not to become involved in the difficulties inherent in applying analytic 
number theory to arithmetic over number fields, we choose to take
averages over $\mathbb{Q}$ and $\mathbb{Z}$. However, we work over number fields
whenever one can proceed in 
general without complicating matters; see subsections 
\ref{subs:defba}--\ref{subs:quadr} and section \ref{sec:cuacriba}.

By a {\em family}\, $\mathcal{E}$ 
of elliptic curves over a number field $K$ on one variable
we mean an elliptic curve over $K(t)$. 
Let $\mathfrak{O}_K$ be the ring of integers of $K$.
We can state conjectures $\mathfrak{A}_1$, $\mathfrak{A}_2$, $\mathfrak{B}_1$
and $\mathfrak{B}_2$ almost exactly as before.

\begin{defn}
Let $K$ be a number field.
We say that a polynomial $P\in \mathfrak{O}_K\lbrack x \rbrack$ {\em yields to a square-free
sieve} if
\[\lim_{N\to \infty} \frac{1}{N} \,\#\{1\leq x\leq N:
\exists \mathfrak{p} \text{\,\st\,} \rho(\mathfrak{p})>N^{1/2}, \mathfrak{p}^2|P(x)\} = 0,\]
where $\rho(\mathfrak{p})$ is the positive integer generating 
 $\mathfrak{p}\cap \mathbb{Z}$.
We say that a homogeneous polynomial $P\in \mathfrak{O}_K\lbrack x,y\rbrack$ 
{\em yields to a square-free sieve}
\[\lim_{N\to \infty} \frac{1}{N^2} 
\,\# \{-N\leq x,y \leq N : \exists \mathfrak{p} \text{\,\st\,} \rho(\mathfrak{p})>N,
\mathfrak{p}^2|P(x,y)\} = 0 .\]
\end{defn}

\begin{defn}
Let $K$ be a number field. We define the {\em generalized Liouville function} $\lambda_K$
on the set of ideals of $\mathfrak{O}_K$ as follows:
\[\lambda_K(\mathfrak{a}) = \begin{cases}
\prod_{\mathfrak{p}|\mathfrak{a}} (-1)^{v_\mathfrak{p}(\mathfrak{a})} & \text{if $\mathfrak{a}\ne 0$,}\\
0 &\text{if $\mathfrak{a}=0$.}\end{cases}\]
If $x\in \mathfrak{O}_K$, we take $\lambda_K(x)$ to mean $\lambda_K((x))$.
\end{defn}

\begin{conja1d}
Let $K$ be a number field.
Every square-free polynomial $P\in \mathfrak{O}_K\lbrack x\rbrack$ 
yields to a square-free sieve.
\end{conja1d}
\begin{conja2d}
Let $K$ be a number field.
Every square-free homogeneous polynomial $P\in \mathfrak{O}_K\lbrack x,y\rbrack$ factors
yields to a square-free sieve.
\end{conja2d}

\begin{conjb1d}
Let $P\in \mathfrak{O}_K\lbrack x \rbrack$ be a polynomial not of the
form $c Q^2(x)$, $c\in \mathfrak{O}_K$, $Q\in \mathfrak{O}_K[x]$.
 Then $\lambda_K(P(n))$ averages to zero over the (rational) integers.
\end{conjb1d}

\begin{conjb2d}
Let $P\in \mathfrak{O}_K\lbrack x,y \rbrack$ be a homogeneous
polynomial not of the form $c Q^2(x,y)$, $c\in \mathfrak{O}_K$,
$Q\in \mathfrak{O}_K[x,y]$. Then 
$\lambda_K(P(x,y))$ averages to zero over $\mathbb{Z}^2$.
\end{conjb2d}

Notice that we speak of Hypotheses $\mathfrak{B}_1$ and $\mathfrak{B}_2$, not of
Conjectures $\mathfrak{B}_1$ and $\mathfrak{B}_2$. This is so because $\mathfrak{B}_i$
fails to hold for some polynomials $P$ over number fields other than $\mathbb{Q}$. Take, for example, $K = \mathbb{Q}(i)$, $P(x) = x$. Then $\lambda_K(P(x))=1$ for all $x\in \mathbb{Z}$
with $x\equiv 1 \mo 4$. 
        
We can, however, reduce Hypothesis $\mathfrak{B}_i(K,P)$ 
to the case $K=\mathbb{Q}$ for which
it is thought always to hold, provided that $K$ and $P$ satisfy certain
conditions. (The counterexample $K = \mathbb{Q}(i)$, $P(x)=x$ does
not fulfill these criteria.) In particular, if $K/\mathbb{Q}$ is Galois,
the situation can be described fully (Corollaries \ref{cor:aga1} and
\ref{cor:aga2}).
Conjecture $\mathfrak{A}_i(K,P)$ can always be reduced to
$\mathfrak{A}_i(\mathbb{Q},P')$ for some polynomial $P'$ over $\mathbb{Q}$. See Appendix
 \ref{sec:app}.

Theorems 1.1--1.4 and Propositions \ref{prop:forgot1}, \ref{prop:forgot2} carry over word by word with $\mathbb{Q}$ replaced by $K$. Corollaries \ref{lem:small1} to \ref{lem:tiny2} carry over easily as well.

\section{Guide to the text}
                        
The main body of the present work is divided into three parts. They are independent from each other as 
far as notation and background are concerned. The first part (Chapter \ref{chap:ell}) 
applies the main results 
of the other two parts, which address
 the analytical side of the matter. The reader who is
interested only in the distribution of the root number may want to confine his attention to Chapter \ref{chap:ell} on a first reading.

In the second part, we prove that $\lambda$ and $\mu$ average to zero over the integers represented
by a homogeneous polynomial of degree at most $3$. In the third part, we strengthen the available
results on square-free sieves by using a mixture of techniques based in part on elliptic curves.
The appendices deal with several related topics of possible interest, including the behavior of
$\lambda(x^2+y^4)$, the relation between certain hypotheses for different number fields, and the average of the root number of cusp forms. 

%% file: ell.tex
\chapter{The distribution of root numbers\\in families of elliptic curves}\label{chap:ell}
\section{Outline}

We will start by describing the behavior of the local root number 
$W_\mathfrak{p}(\mathcal{E}(t))$ for fixed $\mathfrak{p}$ and varying $t$.
It will be necessary to introduce and explay a class of objects, 
{\em pliable functions}, which, among other properties, have desirable
qualities as multipliers.

The global root number $W(\mathcal{E}(t))$ can be written as the product of
a pliable function, a term of the form $\lambda(P(x,y))$ and a correction
factor reflecting the fact that square-free polynomials may adopt 
values that are not square-free. The last factor will be dealt with by
means of a square-free sieve.

\section{Notation and preliminaries}\label{sec:notpre}

Let $n$ be a non-zero integer.
We write $\tau (n)$ for the number of positive 
divisors of $n$, $\omega (n)$ for the
number of the prime divisors of $n$, and $\rad(n)$ for the product of the prime divisors
of $n$. For any $k\geq 2$, we write $\tau_k(n)$ for the number of 
$k$-tuples $(n_1,n_2,\dotsc, n_k)\in (\mathbb{Z}^+)^k$ such that
$n_1\cdot n_2\cdot \dotsb n_k = |n|$. Thus $\tau_2(n) = \tau(n)$.
We adopt the convention that $\tau_1(n) = 1$.
By $d|n^{\infty }$ we will mean that $p|n$ for every
prime $p$ dividing $d$. We let
\[\sq(n) = \prod_{p^2|n} p^{v_p(n)-1} .\]

We denote by $\mathfrak{O}_K$ the ring of integers of a global or local
field $K$. We let $I_K$ be the semigroup of non-zero ideals of $\mathfrak{O}_K$.
If $K$ is a global field and $v$ is a place of $K$, we
will write $\mathfrak{O}_v$ instead of $\mathfrak{O}_{K_v}$.
By a {\em $\mathfrak{p}$-adic field}
we mean a local field of characteristic zero and finite residue field.

Let $K$ be a number field. Let $\mathfrak{a}$ be a non-zero ideal
of $\mathfrak{O}_K$. 
We write $\tau_K(\mathfrak{a})$ for the number of ideals dividing
 $\mathfrak{a}$,
$\omega_K(\mathfrak{a})$ for the number of prime ideals dividing 
$\mathfrak{a}$,
and $\rad_K(\mathfrak{a})$ for the product of the prime ideals dividing
$\mathfrak{a}$.
Given a positive integer $k$, we write $\tau_{K,k}(\mathfrak{a})$ for the number of 
$k$-tuples $(\mathfrak{a}_1,\mathfrak{a}_2,\dotsc ,\mathfrak{a}_k)$ of
ideals of $\mathfrak{O}_K$ such that
$\mathfrak{a} = \mathfrak{a}_1 \mathfrak{a}_2 \dotsb \mathfrak{a}_k$. 
Thus $\tau_2(\mathfrak{a}) = \tau(\mathfrak{a})$.
We let
\[\begin{aligned}
\sq_K(\mathfrak{a}) &= \begin{cases}
\prod_{\mathfrak{p}^2 | \mathfrak{a}} 
\mathfrak{p}^{v_\mathfrak{p}(\mathfrak{a}) - 1} 
&\text{if $\mathfrak{a}\ne 0$,}\\0 &\text{if $\mathfrak{a}=0$,}\end{cases}\\
\mu_K(\mathfrak{a}) &= \begin{cases}
\prod_{\mathfrak{p}|\mathfrak{a}} (-1) &\text{if $\sq_K(\mathfrak{a})=1$,}\\
0 &\text{otherwise.}\end{cases}\end{aligned}\]
We define $\rho(\mathfrak{a})$ to be the positive integer generating $\mathfrak{a}\cap \mathbb{Z}$.

Let $\mathfrak{a}$, $\mathfrak{b}$ be ideals of $\mathfrak{O}_K$.
By $\mathfrak{a}|\mathfrak{b}^\infty$ we mean that 
$\mathfrak{p}|\mathfrak{b}$ for every prime ideal 
$\mathfrak{p}$ dividing $\mathfrak{a}$. 
We write
\[\begin{aligned}
\gcd(\mathfrak{a},\mathfrak{b}) &= \prod_{\mathfrak{p}|\mathfrak{a},
\mathfrak{b}}
\mathfrak{p}^{\min(v_{\mathfrak{p}}(\mathfrak{a}),
v_{\mathfrak{p}}(\mathfrak{b}))},\\
\lcm(\mathfrak{a},\mathfrak{b}) &= 
\mathfrak{a}\cdot \mathfrak{b} \cdot 
(\gcd(\mathfrak{a},\mathfrak{b}))^{-1}.\end{aligned}\]

Throughout, we will say that two polynomials
$f,g\in \mathfrak{O}_K\lbrack x\rbrack$ 
 {\em have no common factors} if they are 
coprime as elements of $K\lbrack x\rbrack$.
We will say that $f\in \mathfrak{O}_K\lbrack x\rbrack$ is {\em square-free}
if there are no polynomials $f_1, f_2\in K\lbrack x\rbrack$, $f_1\notin K$,
such that $f = f_1^2 \cdot f_2$. The same usage will hold for
polynomials in two variables: $f,g \in \mathfrak{O}_K\lbrack x,y\rbrack$
{\em have no common factors} if they are coprime in $K\lbrack x,y\rbrack$,
and $f\in \mathfrak{O}_K\lbrack x,y\rbrack$ is {\em square-free} if
it is not of the form $f_1^2 \cdot f_2$, $f_1, f_2\in K\lbrack x \rbrack$,
$f_1\in K$.

We define the resultant $\Res(f,g)$ of two polynomials
$f,g\in \mathfrak{O}_K\lbrack x\rbrack $ 
as the determinant of the corresponding Sylvester matrix:
\begin{equation}\label{eq:sylv}\left(\begin{array}{ccccccccc}
a_n & a_{n-1} & \dotsb & a_1 & a_0 & 0 & 0 & \dotsb & 0\\
0   & a_n & a_{n-1} & \dotsb & a_1 & a_0 & 0 & \dotsb & 0\\
\vdots & \vdots & \vdots & \vdots & \vdots & \vdots & \vdots & \vdots \\
0 & \dotsb & 0 & a_n & a_{n-1} & \dotsb & a_1 & a_0 & 0\\
0 & \dotsb & 0 & 0 & a_n & a_{n-1} & \dotsb & a_1 & a_0\\
b_m & b_{m-1} & \dotsb & b_1 & b_0 & 0 & 0 & \dotsb & 0\\
0 & b_m & b_{m-1} & \dotsb & b_1 & b_0 & 0 & \dotsb & 0\\
\vdots & \vdots & \vdots & \vdots & \vdots & \vdots & \vdots & \vdots \\
0 & \dotsb & 0 & b_m & b_{m-1} & \dotsb & b_1 & b_0 & 0\\
0 & \dotsb & 0 & 0 & b_m & b_{m-1} & \dotsb & b_1 & b_0
\end{array}\right)
\end{equation}
where we write out $f = \sum_{j=0}^n a_j x^j$,
$g = \sum_{j=0}^m b_i x^i$.

Assume $f$ and $g$ have no common factors. Then $\Res(f,g)$ is a
non-zero element
of $\mathfrak{O}_K$. Moreover,
 $\gcd(f(x),g(x)) | \Res(f,g)$ for any integer $x$. We adopt the
convention that the discriminant $\Disc(f)$ equals $\Res(f,f^{\prime})$.

The resultant of two homogeneous polynomials $f,g \in \mathfrak{O}_K\lbrack
x,y\rbrack $ is also defined as the determinant of the 
Sylvester matrix (\ref{eq:sylv}), where we write out
\[f = \sum_{j=0}^n a_j x^j y^{n-j},\;\;
g = \sum_{j=0}^m b_i x^i y^{m-i}.\]
 Assume $f$ and $g$ have no common factors. Then $\Res(f,g)$
is a non-zero element of $\mathfrak{O}_K$. Moreover, $\gcd(f(x,y),g(x,y)) | \Res(f,g)$ for any coprime integers 
$x$, $y$. 

For a homogeneous polynomial $f\in \mathfrak{O}_K\lbrack x,y\rbrack $ 
we define 
\begin{equation*}
\Disc(f) = \lcm \left(\Res\left(f(x,1),\frac{\partial f(x,1)}{\partial x}%
\right), \Res \left(f(1,y),\frac{\partial f(1,x)}{\partial x}\right)\right).
\end{equation*}

Note that a polynomial $f\in \mathfrak{O}_K\lbrack x\rbrack$ has a
factorization (in general not unique) into polynomials 
$f_1,\dotsb ,f_n\in \mathfrak{O}_K\lbrack x\rbrack$ irreducible in
$\mathfrak{O}_K\lbrack x\rbrack$. In any such factorization,
$f_1,\dotsb , f_n$ are in fact irreducible in $K\lbrack x\rbrack$. The
same is true for homogeneous polynomials $f\in \mathfrak{O}_K\lbrack x,y\rbrack$
and factorization into irreducibles in $\mathfrak{O}_K\lbrack x,y\rbrack$
and $K\lbrack x,y\rbrack$.

A \emph{lattice} is a subgroup of $\mathbb{Z}^n$ of finite index; a \emph{%
lattice coset} 
 is a coset of such a subgroup. By the {\em index} of a lattice
coset we mean the index of the lattice of which it is a coset.
For any lattice cosets $L_1$, $%
L_2$ with $\gcd(\lbrack \mathbb{Z}^n : L_1 \rbrack , \lbrack \mathbb{Z}^n :
L_2 \rbrack) = 1$, the intersection $L_1\cap L_2$ is a lattice coset with 
\begin{equation}
\lbrack \mathbb{Z}^n : L_1\cap L_2\rbrack = \lbrack \mathbb{Z}^n :
L_1\rbrack \lbrack \mathbb{Z}^n : L_2\rbrack .
\end{equation}
In general, if $L_1$, $L_2$ are lattice cosets, then $L_1\cap L_2$ is either
the empty set or a lattice coset such that 
\begin{equation}  \label{eq:inters} \begin{aligned}
\lcm(\lbrack \mathbb{Z}^n : L_1\rbrack,\lbrack \mathbb{Z}^n : L_2\rbrack) 
&\mid
\lbrack \mathbb{Z}^n : L_1\cap L_2\rbrack,\\ \lbrack \mathbb{Z}^n :
L_1\cap L_2\rbrack &\mid \lbrack \mathbb{Z}^n : L_1\rbrack \lbrack \mathbb{Z}^n
: L_2\rbrack .
\end{aligned} \end{equation}

Let $L_1$, $L_2$ be lattices in $\mathbb{Z}^2$. Let 
$m = \lcm(\lbrack \mathbb{Z}^2 : L_1\rbrack, \lbrack \mathbb{Z}^2 : L_2\rbrack)$. Let
$R = \{(x,y) \in \mathbb{Z}^2 : \gcd(m,\gcd(x,y)) = 1\}$. Then $(R\cap L_1)\cap (R\cap L_2)$
is either the empty set or the intersection of $R$ and a lattice $L_3$ of index $m$:
\begin{equation}\label{eq:hmhm}\lbrack \mathbb{Z}^2 : L_3\rbrack  = m =
\lcm(\lbrack \mathbb{Z}^2 : L_1\rbrack, \lbrack \mathbb{Z}^2 : L_2\rbrack) .\end{equation}

For $S\subset \lbrack -N,N\rbrack^n$ a convex set and $L\subset \mathbb{Z}^n$ a lattice
coset, 
\begin{equation}  \label{eq:wbr}
\#(S\cap L) = \frac{\Area(S)}{\lbrack \mathbb{Z}^n :L\rbrack} + O(N^{n-1}),
\end{equation}
where the implied constant depends only on $n$.

The following lemma will serve us better than (\ref{eq:wbr}) when
$L$ is a lattice of index greater than $N$.
\begin{lem}\label{lem:ramsay}
Let $L$ be a lattice of index $\lbrack \mathbb{Z}^2 :L\rbrack \leq N^2$.
Then
\[\#(\{-N\leq x,y\leq N : \gcd(x,y)=1\} \cap L) \ll 
\frac{N^2}{\lbrack \mathbb{Z}^2 : L\rbrack}.\]
\end{lem}
\noindent {\em Proof.}
Let \[M_0 = \min_{(x,y)\in L} \max(|x|,|y|).\]
By \cite{Gr}, Lemma 1, 
\[\#(\lbrack -N,N\rbrack^2 \cap L) \ll \frac{N^2}{\lbrack \mathbb{Z}^2 :
L\rbrack} + O\left(\frac{N}{M_0}\right).\]
If $M_0\geq 
\frac{\lbrack \mathbb{Z}^2 : L\rbrack}{2 N}$ we are done. Assume
$M_0< \frac{\lbrack \mathbb{Z}^2 : L \rbrack}{2 N}$. 
Suppose 
\[\#(\{-N\leq x,y\leq N : \gcd(x,y)=1\} \cap L) > 2.\]
Let $(x_0,y_0)$ be a point
such that $\max(|x_0|,|y_0|) = M_0$. Let $(x_1,y_1)$ be a point
in $\#(\{-N\leq x,y\leq N : \gcd(x,y)=1\} \cap L)$ other than
$(x_0,y_0)$ and $(-x_0,-y_0)$. Since $\gcd(x_0,y_0)=\gcd(x_1,y_1)=1$,
it cannot happen that $0$, $(x_0,y_0)$ and $(x_1,y_1)$ lie on the same
line. Therefore we have a non-degenerate parallelogram
$(0,(x_0,y_0),(x_1,y_1),(x_0+x_1,y_0+y_1))$ whose area has to be at least
$\lbrack \mathbb{Z}^2 : L\rbrack$. On the other hand, its area can
be at most $\sqrt{x_0^2+y_0^2} \cdot \sqrt{x_1^2+y_1^2} \leq 
\sqrt{2} M_0 \cdot \sqrt{2} N = 2 M_0 N$. Since we have assumed
$M_0< \frac{\lbrack \mathbb{Z}^2 : L \rbrack}{2 N}$ we arrive at a
contradiction. \qed

By a \emph{sector} we will mean a connected component of a set of the form $%
\mathbb{R}^{n}-(T_{1}\cap T_{2}\cap \dotsb \cap T_{n})$, where $T_{i}$ is a
hyperplane going through the origin. Every sector $S$ is convex.

Let $x\in \mathbb{R}$ be given. We write $\lfloor x\rfloor$ for the largest integer no
greater than $x$, $\lceil x\rceil$ for the smallest integer no smaller than $x$,
and $\{x\}$ for $x-\lfloor x\rfloor$.

We define $\lbrack \text{true}\rbrack$ to be $1$ and $\lbrack \text{false}\rbrack$ to be $0$.
Thus, for example, $x\mapsto \lbrack x \in S\rbrack$ is the characteristic function of a 
set $S$.
\section{Pliable Functions}\label{sec:plifun}

Since this section is devoted to a newly defined class of objects,
we might as well start by attempting to give an intuitive sense of their
meaning. Take a function $f:\mathbb{Z}_p\to \mathbb{C}$. For $f$
to be {\em affinely pliable}, it is necessary but not sufficient that 
$f$ be locally constant almost everywhere. We say that $f$ is {\em affinely
pliable at $0$} if there is an integer $k\geq 0$ such that the value of
$f(x)$ depends 
only on $v_p(x)$ and on $p^{-v_p(x)} x \mo p^k$. Thus, if, say, $p=3$
and $k=1$, each of the following values is uniquely defined:
\[\begin{array}{llllll}
f(\dotsc 01_3) &f(\dotsc 02_3) & f(\dotsc 11_3) &f(\dotsc 12_3) &
 f(\dotsc 21_3) &f(\dotsc 22_3) \\
f(\dotsc 010_3) &f(\dotsc 020_3) & f(\dotsc 110_3) &f(\dotsc 120_3) &
 f(\dotsc 210_3) &f(\dotsc 220_3) \\
f(\dotsc 0100_3) &f(\dotsc 0200_3) & f(\dotsc 1100_3) &f(\dotsc 1200_3) &
 f(\dotsc 2100_3) &f(\dotsc 2200_3) \\
\dotsc & \dotsc & \dotsc & \dotsc & \dotsc & \dotsc 
\end{array}\]
A function $f$ on $\mathbb{Z}_p$ is affinely pliable at $t_1, \dotsc, t_n$
if it displays the same behaviour near $t_1, t_2,\dotsc ,t_n$ as the
example above displays near $0$. A function $f$ on $\mathbb{R}$ is affinely
pliable at $t_1<\dotsb < t_n$ if it is constant on $(-\infty,t_1)$,
$(t_1,t_2)$, \dots , $(t_n,\infty)$. A function $f$ on $\mathbb{Q}$ is
affinely pliable if it is affinely pliable when seen at finitely many
places simultaneously, in a sense to be made precise now.
\subsection{Definition and basic properties}\label{subs:defba}
\begin{defn}\label{def:rpli}
Let $K$ be a number field or a
$\mathfrak{p}$-adic field. A function $f$ on a subset $S$ of $K$ is said to be
{\em affinely pliable}
if there are finitely many triples
\[(v_j,U_j,t_j)\]
with $v_j$ a place of $K$, $U_j$ an open subgroup of $K_{v_j}^*$
and $t_j$ an element of $K_{v_j}$ such that $f(t) = f(t')$
for all $t, t'\in S$ such that $t-t_j$ and $t'-t_j$ are non-zero and
equal in $K_{v_j}^*/U_j$ for all $j$.
\end{defn}
If $K$ is a $\mathfrak{p}$-adic field, then $v_j$ has no choice but to equal
the valuation $v_{\mathfrak{p}}$ of $K$.
When $f$ is affinely pliable with respect to $(v_1,U_1,t_1)$,\dots,$(v_n,U_n,t_n)$,
we say $f$ is {\em affinely pliable at $t_1$,\dots,$t_n$.}
\begin{defn}
\label{def:pli} 
Let $K$ be a number field or a $\mathfrak{p}$-adic field.
A function $f$ on a subset of
$K^n$ is said to be \emph{pliable}
if there are finitely many triples 
\begin{equation*}
(v_j,U_j,\vec{q}_j)
\end{equation*}
with $v_j$ a place of $K$, $U_j$ an open subgroup of 
$K_{v_j}^*$ and $\vec{q}_j$ an element of $K_{v_j}^n-\{(0,\dotsc,0)\}$ such that $%
f(x_1,x_2,\dotsb, x_n) = f(x_1^{\prime},x_2^{\prime},\dotsb , x_n^{\prime})$
whenever the scalar products $\vec{x} \cdot \vec{q}_j$ and $\vec{x}%
^{\prime}\cdot \vec{q}_j$ 
are non-zero and equal in $K_{v_j}^*/U_j$ for all $j$.
\end{defn}

The following are some typical examples of pliable and
affinely pliable functions. Let $K$ be a
number field or a $\mathfrak{p}$-adic field, $v$ a place of $K$. Then 
$t\mapsto v(t)$ is affinely pliable. So are 
$t\mapsto t \mo \mathfrak{p}_{v}$ (defined on $K\cap \mathfrak{O}_{K_{v}}$) and
$t\mapsto t \pi_{v}^{-v(t)} \mo \mathfrak{p}_{v}$ (defined on $K^*$), 
where $\mathfrak{p}_{v}$ is
the prime ideal of $\mathfrak{O}_{v}$ and $\pi_{v}$ is a generator of
$\mathfrak{p}_{v}$. If $K$ is a $\mathfrak{p}$-adic field, 
any continuous character $\chi:K^*\mapsto \mathbb{C}$ is affinely
pliable. For any ball $B = \{t\in K:|t-t_0|_{v}<r\}$, the characteristic function
$x\mapsto \lbrack x\in B\rbrack$ is affinely pliable. An example of a pliable
function would be $(x,y)\mapsto v_{\mathfrak{p}}(3 x + 5 y)$, or $(x,y,z)\mapsto
\chi(3 y - 2 x + z)$. A function is affinely pliable at $0$ if and only 
if it is a pliable function on one variable ($n=1$). Of the examples of affinely pliable
functions given above, all are affinely pliable at $0$, save for 
$x\mapsto \lbrack x\in B\rbrack$, which is affinely pliable at $t_0$.

It is clear that $g \circ (f_1\times f_2 \times ... \times f_n)$ is pliable
(resp. affinely pliable)
for $f_1, f_2, \dotsb f_n$ pliable (resp. affinely pliable)
and $g$ an arbitrary function whose domain is a subset of the range
of $f_1\times f_2\times \dotsb \times f_n$.
 Note, in particular, that $f_1 f_2 \dotsb f_n$ is pliable
(resp. affinely pliable) for $f_1$,\dots,$f_n$ pliable. We will now prove that,
under certain circumstances, pliability is preserved under composition in the other order:
not only is $t\mapsto \chi^3(t)+\chi(t)+5$ affinely pliable, but
$t\mapsto \chi(t^3+t+5)$ is affinely pliable as well.

\begin{lem}\label{lem:rinv}
Let $K$ be a number field or a
$\mathfrak{p}$-adic field. Let $v$ be a place of $K$, $f\in K_{v}(t)$ a rational
function and $U$ an open subgroup of $K_{v}^*$. 
Let $t_1, t_2, \dotsc, t_n\in K$
be the zeroes and poles of $f$ in $K_{v}$. Let $t_0=0$.
Then there is an open subgroup $U_{v}'$ of $K_{v}^*$ 
 such that
$f(t)$ is in the
same coset $r U_{v}$ of $U_{v}$ as $f(t')$ whenever\, $t-t_j$ and $t'-t_j$ lie
in the same coset $r_j U_{v}$ of $U_{v}$ for every $0\leq j\leq m$.
\end{lem}
\begin{proof}
We will choose $U_{v}\subset \mathfrak{O}_{K_{v}}^*$. If $t$ and $t'$ belong to the
same coset of $U_{v}$, then $t\in \mathfrak{O}_{K_{v}}$ implies 
$t'\in \mathfrak{O}_{K_{v}}$.
For any $t\in K_{v}$, either $t\in \mathfrak{O}_{K_{v}}$ or 
$1/t\in \mathfrak{O}_{K_{v}}$.
Let $\hat{f}\in K_{v}(t)$ be the rational function taking $t$ to $f(1/t)$.
If we prove the statement of the lemma for both $f$ and $\hat{f}$ under the
assumption that $t, t'\in \mathfrak{O}_{K_{v}}$, we will have proven it
for any $t, t'\in K_{v}$.  Thus we need consider only $t,t'\in \mathfrak{O}_{K_{v}}$.

As in Lemma \ref{lem:inv}, we can assume $f$ is an irreducible polynomial
with integer coefficients. If $f$ is linear, the statement is immediate.
Hence we can assume $f\in \mathfrak{O}_{K_{v}}\lbrack t \rbrack$, $f$ irreducible,
$\deg(f)\geq 2$.

Hensel's lemma implies that $v(f(t))\leq 2 v(\Disc(f))$ for every
$t\in \mathfrak{O}_{K_{v}}$, as the contrary would be enough for $f(x)=0$ to
have a non-trivial solution in $K_{v}$. Since $U_{v}$ is open, it contains
a set of the form $1+ \pi^k \mathfrak{O}_{K_{v}}$, where $\pi$ is a prime of
$\mathfrak{O}_{K_{v}}$. Set $U_{v}' = 1 + \pi^{k+ 2 v(\Disc(f))}$. 
Suppose $t,t'\in \mathfrak{O}_{K_{v}}$ lie in the same coset of $U_{v}'$. Then
$v(t-t')\geq k + 2 v(\Disc(f)) + v(t)$. Since $v$ is non-archimedean,
\[|f(t)-f(t')|\leq |t-t'|\leq |\pi|^{k+2\Disc(f) + v(t)}\leq
|\pi|^{k+ f(t)} = |\pi|^k |f(t)|.\]
Therefore $f(t)$ and $f(t')$ lie in the same coset of $U_{v}$.
\end{proof}
\begin{prop}\label{prop:simil}
Let $K$ be a number field or a $\mathfrak{p}$-adic field. Let $f\in K(t)$. 
Let a function $g$ on $S\subset K$ be affinely pliable. Then $g\circ f$ on $S'=
\{t\in K: f(t)\in S\}$ is affinely pliable.
\end{prop}
\begin{proof}
Immediate from Definition \ref{def:rpli} and Lemma \ref{lem:rinv}.
\end{proof}
\begin{prop}\label{prop:similtrans}
Let $K$ be a number field or a $\mathfrak{p}$-adic field. Let $f_1,\dotsb,f_n\in K(t)$.
Let a function $g$ on $S\subset K^n$ be pliable. Then the map
$t\mapsto g(f_1(t),\dotsb,f_n(t))$ on $S'=\{t\in K: (f_1(t),\dotsb,f_n)\in S\}$
is affinely pliable.
\end{prop}
\begin{proof}
Immediate from Definitions \ref{def:rpli} and \ref{def:pli} and Lemma \ref{lem:rinv}.
\end{proof}
\begin{lem}
\label{lem:inv} Let $K$ be a number field or a $\mathfrak{p}$-adic field.
Let $v$ be a place of $K$, 
$F\in K_v\lbrack x,y\rbrack $ a homogeneous polynomial and $U_v$ an open
subgroup of $K_{v}^*$. Then there is an open subgroup 
$U^{\prime}_v$ of $K_{v}^*$ and a finite subset $\{\vec{x}_j\}$
of $K_{v}^2$ such that $F(x,y)$ is in the same coset $r U_v$ of 
$U_v$ as $F(x^{\prime},y^{\prime})$ whenever $(x,y) \cdot \vec{x}_j$ and $(x^{\prime},y^{\prime}) \cdot \vec{x}_j$ lie in the same coset $r_j
U^{\prime}_v$ of $K_{v}^*$ for all $j$.
\end{lem}
\begin{proof}
Suppose that $F = F_1 F_2$ and that the lemma holds for $(F_1,v,U_v)$
and $(F_2,v,U_v)$ with conditions $(U^{\prime}_{v,1},\{\vec{x}_{i,1}\})
$ and $(U^{\prime}_{v,2},\{\vec{x}_{k,2}\})$, respectively. Set $%
U^{\prime}_v = U^{\prime}_{v,1} \cap U^{\prime}_{v,2}$ and $\{\vec{x}%
_j\} = \{\vec{x}_{i,1}\} \cup \{\vec{x}_{k,2}\}$. Assume that $(x,y) \cdot 
\vec{x}_j$ and $(x^{\prime},y^{\prime}) \cdot \vec{x}_j$ lie in the same
coset of $U^{\prime}_v$ for all $j$. Then $F_1(x,y)$ is in the same coset
of $U_v$ as $F_1(x^{\prime},y^{\prime})$ and $F_2(x,y)$ is in the same
coset as $F_2(x^{\prime},y^{\prime})$. Hence $F_1(x,y) F_2(x,y)$ is in the
same coset as $F_1(x^{\prime},y^{\prime}) F_2(x^{\prime},y^{\prime})$.

We can thus assume that $F$ is irreducible. Suppose $F$ is linear. Write $%
F(x,y) = a x + b y$. Then the lemma holds with $U_v^{\prime}= U_v$ and $%
\{\vec{x}_j\}=\{(a,b)\}$. We are left with the case when $F$ is irreducible
of degree greater than one.

Suppose $v$ is finite. We can assume $F\in \mathfrak{O}_v[x,y]$.
  Hensel's
Lemma implies that $v(F(x,y)) - (\deg F) \min(v(x),v(y)) \leq 2 v(\Disc(F))$
for all $x,y\in K^*$,
as the contrary would be enough for $F(x,y)=0$ to have a non-trivial
solution in $K_v^2$. Since $U_v$ is open, it contains a set of
the form $1 + \pi^k \mathfrak{O}_v$, 
where $\pi$ is a prime of $\mathfrak{O}_v$.
 Set $U_v^{\prime}= 1 + \pi^{k+2 v(\Disc(F))} \mathfrak{O}_v$, $\vec{x}_1=(1,0)$, $\vec{x}_2 = (0,1)$. 
Suppose that $(x,y)$ and $(x^{\prime},y^{\prime})$ 
satisfy the conditions in the lemma, that is, $x$
and $x^{\prime}$ lie in the same coset of $U_v^{\prime}$, and so do $y$
and $y^{\prime}$.
It follows that $v(x-x^{\prime}) \geq k + 2 v(\Disc(F)) + v(x)$
and $v(y-y^{\prime}) \geq k + 2 v(\Disc(F)) + v(y)$. 
Since $v$ is
non-archimedean, 
\begin{equation*}
|F(x,y) - F(x^{\prime},y^{\prime})|_v \leq
|\pi|^{(\deg(F) -1) \min(v(x),v(y))} \max(|x-x^{\prime}|_v,
|y-y^{\prime}|_v) .
\end{equation*}
Now 
\begin{equation*}
\begin{aligned}
\max(|x-x'|_v,|y-y'|_v) &= |\pi|_v^{\min(v(x-x'),v(y-y'))}\\
&\leq |\pi|_v^k |\pi|_v^{- (\deg(F) - 1) \min(v(x),v(y))}
|\pi|^{2 v(\Disc(F)) + 
\deg(F) \min(v(x), v(y))}  \\
&\leq |\pi|_v^k |\pi|_v^{- (\deg(F) - 1) \min(v(x),v(y))}
|F(x,y)|_v .\end{aligned}
\end{equation*}
Thus 
\begin{equation*}
|F(x,y) - F(x^{\prime},y^{\prime})|_v \leq |\pi|^k |F(x,y)|_v .
\end{equation*}
This means that $F(x,y)$ and $F(x^{\prime},y^{\prime})$ are in the same
coset of $U_v$. 

Suppose now that $v$ is infinite and $F(x,y)$ is irreducible and of degree
greater than one. Then the degree of $F$ must be two. We have either $U_v
= \mathbb{R}^*$ or $U_v = \mathbb{R}^+$. Since $F$ is either positive
definite or negative definite, $F(x,y)$ and $F(x^{\prime},y^{\prime})$ lie
in the same coset of $U_v$ for any $x$, $y$ not both zero. Since we are
given that $x$ and $y$ are coprime they cannot both be zero. Choose $%
U_v^{\prime}= \mathbb{R}^*$, $\{x_j\}$ empty.
\end{proof}

As usual, we write $\vec{e}_1 = (1,0,\dotsc,0), \vec{e}_2=(0,1,\dotsc,0),
\dotsc ,\vec{e}_n = (0,0,\dotsc,1)$.

\begin{prop}
\label{prop:plinvar} Let $K$ be a number field or a $\mathfrak{p}$-adic field.
Let $F_1, F_2, \dotsc, F_n \in K\lbrack
x,y\rbrack$ be homogeneous polynomials. Let a function $f$ on $S\subset K^n$ 
be pliable with respect to $\{(v_j,U_j,\vec{q}_j)\}_j$. Suppose $\vec{q}_j\in \{\vec{e}_1,\vec{e}_2,\dotsb,\vec{e}_n\}$ for every $j$. Then $(x,y)\mapsto f(F_1(x,y),F_2(x,y),\dotsb ,F_n(x,y))$ is a
pliable function on \[S'=\{(x,y)\in K^2:(F_1(x,y),\dotsb,F_n(x,y))\in S\}.\]
\end{prop}
\begin{proof}
Immediate from Definition \ref{def:pli} and Lemma \ref{lem:inv}.
\end{proof}

\begin{prop}
\label{prop:plisam} Let $K$ be a number field or a $\mathfrak{p}$-adic field.
Let $F_1, F_2, \dotsc, F_n \in K\lbrack
x,y\rbrack$ be homogeneous polynomials of the same degree. Let a function $f$ on $S\subset K^n$ 
be pliable with respect to $\{(v_j,U_j,\vec{q}_j)\}_j$. 
Then \[(x,y)\mapsto f(F_1(x,y),F_2(x,y),\dotsb ,F_n(x,y))\] is a
pliable function on $S'=\{(x,y)\in K^2:(F_1(x,y),\dotsb,F_n(x,y))\in S\}$.
\end{prop}
\begin{proof}
Immediate from Definition \ref{def:pli} and Lemma \ref{lem:inv}.
\end{proof}

\begin{lem}\label{lem:pingpong}
Let $K$ be a number field or a $\mathfrak{p}$-adic field. 
Let $f$ be an affinely pliable function on a subset $S$ of $K$. Then the map
\[(x,y)\to f(y/x)\]
on $S'=\{(x,y)\in K^2 : y/x \in S\}$ is pliable.
\end{lem}
\begin{proof}
Say $f$ is affinely pliable with respect to $\{(v_j, U_j, t_j)\}_j$.
Let $(x,y),(x',y')\in S'$, be such that
$t_j x - y$ and $t_j x' - y'$ belong to the same coset 
$r_j U_j\subset K_{v_k}^*$ of
$U_j$ for every $j$. Assume furthermore that $x$ and $x'$ belong
to the same coset of $U_j$. Then $y/x-t_j$ and $y'/x'-t_j$ belong
to the same coset of $U_j$ for every $j$.
Therefore $(x,y)\to f(y/x)$
is pliable with respect to $\{(v_j,U_j,(t_j,-1))\}_j \cup
\{(v_j, U_j, (1,0))\}_j$.
\end{proof}

\begin{lem}\label{lem:pongping}
Let $K$ be a number field or a $\mathfrak{p}$-adic field. Let $f$ be a pliable 
function on a subset $S$ of $K^2$. 
Then the map \[t\mapsto f(1,t)\] on $S'=\{t\in K:(1,t)\in S\}$ is
affinely pliable.
\end{lem}
\begin{proof}
Say $f$ is pliable with respect to $\{v_j,U_j,(q_{j 1},q_{j 2})\}_{j\in J}$. 
Let $t, t'\in S$ be
such that $t + q_{j 1}/q_{j 2}$ and $t' + q_{j 1}/q_{j 2}$ belong to the same coset
$r_j U_j \subset K_{v_j}^*$ of $U_j$ for every $j$ such that $q_{j 2}\ne 0$.
Then $q_{j 1} + q_{j 2} t$ and $q_{j 1} + q_{j 2} t'$ belong to the same coset
$r_j U_j \subset K_{v_j}^*$ of $U_j$ for every $j$. Therefore $t\mapsto f(1,t)$
is affinely pliable with respect to 
$\{(v_j,U_j,-q_{j 1}/q_{j 2})\}_{j\in J'}$, where $J'=\{j\in J: q_{j 2}\ne 0\}$.
\end{proof}

\begin{lem}\label{lem:opres}
If $K$ is a number field or a $\mathfrak{p}$-adic field, 
$L$ a finite extension of $K$, and $f$ a pliable
function on a subset $S$ of $L^n$, then $f|(S\cap K^n)$ is pliable as a function
on the subset $S\cap K^n$ of $K^n$.
If $K$ is a number field or a $\mathfrak{p}$-adic field, 
$L$ a finite extension of $K$, and $f$ an affinely
 pliable
function on a subset $S$ of $L$, then $f|(S\cap K)$ is affinely
pliable as a function
on the subset $S\cap K$ of $K$.
\end{lem}
\begin{proof}
The intersection of $K$ and an open subgroup of $L^*$ is an open subgroup
of $K^*$.
\end{proof}
\begin{lem}\label{lem:oneoflst}
Let $f$ be a pliable function on a subset $S$ of $\mathbb{Z}^2$. Let
$m$ be a positive integer. Then
\[(x,y)\mapsto f\left(\frac{x}{\gcd(x,y,m)},\frac{y}{\gcd(x,y,m)}\right)\]
is a pliable function on $S'=\{(x,y)\in \mathbb{Z}^2:(x/\gcd(x,y,m),
y/\gcd(x,y,m))\in S\}$.
\end{lem}
\begin{proof}
Suppose $f$ is pliable with respect to $\{(v_j,U_j,\vec{q}_j)\}_{j\in J}$.
Then \[(x,y)\mapsto f\left(\frac{x}{\gcd(x,y,m)},\frac{y}{\gcd(x,y,m)}\right)\]
is pliable with respect to $\{(v_j,U_j,\vec{q}_j)\}_{j\in J}
\cup \{(v_\mathfrak{p},\mathfrak{O}_{K_{\mathfrak{p}}}, (1,0))
\}_{\mathfrak{p}|m} 
\cup \{(v_\mathfrak{p},\mathfrak{O}_{K_{\mathfrak{p}}}, (0,1))
\}_{\mathfrak{p}|m}$.
\end{proof}
\begin{lem}\label{lem:redef}
Let $K$ be a number field or a $\mathfrak{p}$-adic field. Let $f$ be a pliable function
from $X\subset K^n$ to $Y$ (resp. an affinely pliable function from $X\subset K$
to $Y$). Let $x_0\in K^n-X$ (resp. $x_0\in K-X$), $y_0\in Y$. 
Define $f':S\cup \{x_0\}\to Y$
by \[f'(x) = \begin{cases} f(x) &\text{if $x\in S$,}\\y_0 &\text{if $x=x_0$.}\end{cases}\]
Then $f$ is pliable (resp. affinely pliable).
\end{lem}
\begin{proof}
If $f$ is pliable with respect to $\{(v_j,U_j,\vec{q}_j)\}_j$ (resp. aff. pliable
with respect to $\{(v_j,U_j,t_j)\}_j$) then $f'$ is pliable with respect to
$\{(v_j,U_j,\vec{q}_j)\}_j \cup \{(v,K,\vec{v})\}$, where $v$ is an
arbitrary place of $K$ and $\vec{v}$ is any vector orthogonal to $x_0$ 
(resp. aff. pliable
with respect to $\{(v_j,U_j,t_j)\}_j \cup \{(v,K,x_0)\}$, where $v$ is
an arbitrary place of $K$).
\end{proof}
\begin{lem}\label{lem:strau1}
Let $f$ be an affinely pliable function on $\mathbb{Z}$.
Then there are integers $a$, $m$ and $t_0$, $m>0$,
such that $f$ is constant on the set $\{t\in \mathbb{Z} :
t\equiv a \mo m,\:t>t_0\}$.
\end{lem}
\begin{proof}
Immediate from Definition \ref{def:rpli}.
\end{proof}
\begin{lem}\label{lem:strau2}
Let $f$ be a pliable function on $\mathbb{Z}^2$. Then there
are a lattice $L\subset \mathbb{Z}^2$ and a sector $S\subset \mathbb{R}^2$
such that $f$ is constant on $L\cap S$. 
\end{lem}
\begin{proof}
Immediate from Definition \ref{def:pli}.
\end{proof}

\subsection{Pliability of local root numbers}

Let $E$ be an elliptic curve over a field $K$. Given an extension $L/K$,
we write $E(L)$ for the set of $L$-rational points of $E$. We define
$E\lbrack m\rbrack \subset E(\overline{K})$ to be the set of points of
order $m$ on $E$. We write $K(E\lbrack m\rbrack)$ for the minimal subextension
of $K$ over which all elements of $E\lbrack m\rbrack$ are rational. The
extension $K(E\lbrack m\rbrack)/K$ is always finite and Galois.

Write $\widetilde{K}$ for the maximal unramified extension of a local field
$K$.

\begin{lem}\label{lem:STgood}
Let $K$ be a $\mathfrak{p}$-adic field. Let $E$ be an elliptic curve over $K$
with potential good reduction. Then there is a minimal
algebraic extension $L$ of $\widetilde{K}$ over which $E$ acquires good 
reduction. Moreover, 
$L=\widetilde{K}(E\lbrack m\rbrack)$ for all $m\geq 3$ prime to the
characteristic of the residue field of $K$.
\end{lem}
\begin{proof}
See \cite{ST}, Section 2, Corollary 3.
\end{proof}

\begin{lem}\label{lem:ugleg}
 Let $K$ be a $\mathfrak{p}$-adic field. Then there is a
finite extension $K'/\widetilde{K}$ such that every elliptic curve over $K$
with potential good reduction acquires good reduction over $K'$.
\end{lem}
\begin{proof}
We can check directly from the explicit formulas for the group law
(see e.g. \cite{Si}, Chap III, 2.3) that $K(E \lbrack 3\rbrack)/K$
is an extension of degree at most $6$ and $K(E \lbrack 4\rbrack)/K$
is an extension of degree at most $12$. Since $K$ is a $\mathfrak{p}$-adic
field, it has only finitely many extensions of given degree (see e.g.
\cite{La}, II, \S 5, Prop. 14). Let $K_{12}/K$ be the composition of
all extensions of $K$ of degree at most $12$. Since $K_{12}/K$ is the
composition of finitely many finite extensions, it is itself a finite
extension. By Lemma \ref{lem:STgood}, every elliptic curve over $K$
with potential good reduction acquires good reduction over $L = K_{12}
\cdot \widetilde{K}$. Since $K_{12}/K$ is a finite extension,
$L/\widetilde{K}$ is a finite extension.
\end{proof} 
\begin{lem}\label{lem:firfew}
 Let $K$ be a local field of ramification degree $e$
over $\mathbb{Q}_p$. Let $\pi$ be a prime of $K$. Then the reduction 
$\mo \mathfrak{p}$ of an elliptic curve $E$ over $K$ depends only on
\[\begin{aligned}
c_4\cdot \pi^{-4 \min(\lfloor v(c_4)/4\rfloor,\lfloor v(c_6)/6\rfloor,
\lfloor v(\Delta)/12 \rfloor)} &\mo \mathfrak{p}^{5 e + 1},\\
c_6\cdot \pi^{-6 \min(\lfloor v(c_4)/4\rfloor,\lfloor v(c_6)/6\rfloor,
\lfloor v(\Delta)/12 \rfloor)} &\mo \mathfrak{p}^{5 e + 1},\end{aligned}\]
where $c_4$, $c_6$ and $\Delta$ are any choice of parameters for $E$.
\end{lem}
\begin{proof} 
Let $k$ be the residue field of $K$. Let $E$ be an elliptic curve over $K$
with parameters $c_4, c_6,\Delta \in K$. Let 
\[\begin{aligned}
c_4' &= c_4\cdot \pi^{-4 \min(\lfloor v(c_4)/4\rfloor,\lfloor v(c_6)/6\rfloor,
\lfloor v(\Delta)/12 \rfloor)},\\
c_6' &= c_6\cdot \pi^{-6 \min(\lfloor v(c_4)/4\rfloor,\lfloor v(c_6)/6\rfloor,
\lfloor v(\Delta)/12 \rfloor)}.\end{aligned}\]
Suppose $\charac(k)\ne 2,3$. Then a minimal Weierstrass equation for $E$ is given

by
\[y^2 = x^3 - \frac{c_4'}{48} - \frac{c_6'}{864} .\]
Both $\frac{-c_4'}{48}$ and $\frac{-c_6'}{864}$ are integral. The reduction
$\mo \mathfrak{p}$ is simply
\[y^2 = x^3 - (c_4' \cdot 48^{-1} \mo \mathfrak{p}) - 
(c_6' \cdot 864^{-1} \mo \mathfrak{p}) .\]
This depends only on $c_4',c_6' \mo \mathfrak{p}$.

Consider now $\charac(k)=2,3$.  Let $m$ be the smallest positive integer
such that there are $r,s,t\in K$, $u\in \mathfrak{O}_K^*$, for which the equation
\begin{equation}\label{eq:hrumph}
(u^3 y' + s u^2 x' + t)^2 = (u^2 x' + r)^3 - 
\frac{\pi^{4 m} c_4'}{48} (u^2 x' + r) - \frac{\pi^{6 m} c_6'}{864}\end{equation}
has integral coefficients when expanded on $x'$ and $y'$.  (Clearly
$m\leq e$.)  Then, for $m$ and any choice of $r,s,t\in K$, $u\in \mathfrak{O}_K^*$,
giving integral coefficients, (\ref{eq:hrumph}) is a minimal Weierstrass
equation for $E$, and its reduction $\mo \mathfrak{p}$ gives us the reduction 
$E \mo \mathfrak{p}$. 

By \cite{Si}, III, Table 1.2, $r,s,t\in K$, 
$u\in \mathfrak{O}_K^*$ can give us integral coefficients only if 
$3 r, s, t\in \mathfrak{O}_K$ (if $\charac(k)=3$) or
$2 r, s, 2 t\in \mathfrak{O}_K$ (if $\charac(k)=2$).
Thus, both the existence of (\ref{eq:hrumph}) and its
coefficients $\mo \mathfrak{p}$ depend only on $\frac{c_4'}{2\cdot 3\cdot 48},
\frac{c_6'}{864} \mo \mathfrak{p}$. 
Since $c_4'$ and $c_6'$ are integral,
$\frac{c_4'}{2\cdot 3 \cdot 48}$,
$\frac{c_6'}{864} \mo \mathfrak{p}$ depend only on
$c_4' \mo \mathfrak{p}^{5 e + 1}$ and $c_6' \mo \mathfrak{p}^{5 e + 1}$
(if $\charac(k)=2$) or on $c_4' \mo \mathfrak{p}^{2 e+1}$ and
$c_6' \mo \mathfrak{p}^{3 e + 1}$ (if $\charac(k)=3$). The statement follows.
\end{proof}
\begin{lem}\label{lem:sundh}
Let $K$ be a $\mathfrak{p}$-adic field of ramification degree $e$ over
$\mathbb{Q}_p$. Let $L$ be an extension of $K$ of finite ramification degree
over $K$. Let $\mathfrak{p}_K$ be the prime ideal of $K$, $\mathfrak{p}_L$
the prime ideal of $L$.
Then the reduction mod $\mathfrak{p}_L$ of an elliptic curve $E$ defined over $K$
depends only on $K$, $L$, $v_K(c_4)$, $v_K(c_6)$, $v_K(\Delta)$,
$c_4\cdot (1+\mathfrak{O}_K \mathfrak{p}_K^{5 e+1})$ and 
$c_6\cdot (1+\mathfrak{O}_K \mathfrak{p}_K^{5 e+1})$, where
$c_4$, $c_6$ and $\Delta$ are any choice of parameters for $E$.
\end{lem}
\begin{proof}
Let $e'$ be the ramification degree of $L$ over $K$.
Let $\pi_L$ be a prime of $L$, $\pi_K = \pi_L^{e'}$ a prime of $K$.
 By Lemma \ref{lem:firfew},
the reduction $E\mo \mathfrak{p}_L$ depends only on 
$c_4' \mo \mathfrak{p}_L^{5 e e' +1}$ and
$c_6' \mo \mathfrak{p}_L^{5 e e' + 1}$,
where
\[\begin{aligned}
c_4' &= c_4\cdot \pi_L^{-4 \min(\lfloor v_L(c_4)/4\rfloor,
\lfloor v_L(c_6)/6\rfloor, \lfloor v_L(\Delta)/12\rfloor)}\\
c_6' &= c_6\cdot \pi_L^{-6 \min(\lfloor v_L(c_4)/4\rfloor,
\lfloor v_L(c_6)/6\rfloor, \lfloor v_L(\Delta)/12\rfloor)} .
\end{aligned}\]
Since $4\lbrack v_L(c_4)/4\rbrack\leq v_L(c_4) = e' v_K(c_4)$ and
$6\lbrack v_L(c_4)/6\rbrack\leq v_L(c_6) = e' v_K(c_6)$, we can tell
$c_4' \mo \mathfrak{p}_L^{5 e e' +1}$ and
$c_6' \mo \mathfrak{p}_L^{5 e e' + 1}$ from
$v_L(c_4)$, $v_L(c_6)$, $v_L(\Delta)$,
\[\begin{aligned}
c_4 &\cdot \pi_L^{- e' v_K(c_4)} \mo \mathfrak{p}_L^{5 e e' + 1} \text{ and}\\
c_6 &\cdot \pi_L^{- e' v_K(c_6)} \mo \mathfrak{p}_L^{5 e e' + 1}.\end{aligned}\]
(Either of the last two may not be defined, but we can tell as much from
whether $v_L(c_4)$ and $v_L(c_6)$ are finite.) Since 
$v_L(c_4) = e' v_K(c_4)$, $v_L(c_6) = e' v_K(c_6)$,
$v_L(\Delta) = e' v_K(\Delta)$, $\pi_K = \pi_L^{e'}$ and
$\mathfrak{p}_K = \mathfrak{p}_L^{e'}$, it is enough to know
$v_K(c_4)$, $v_K(c_6)$, $v_K(\Delta)$, 
$c_4 \cdot \pi^{-v_K(c_4)} \mo \mathfrak{p}^{5 e + 1}$
and $c_6 \cdot \pi^{-v_K(c_6)} \mo \mathfrak{p}^{5 e + 1}$. The
statement follows immediately.
\end{proof}

\begin{lem}\label{lem:pardiez}
 Let $K$ be a Henselian local field. Let $k$ be the residue field
of $K$. Let $m\geq 2$ be an integer prime to $\charac(k)$. Let $E$ be
an elliptic curve defined over $K$ with
 good reduction at $\mathfrak{p}_K$; denote its reduction
by $\widehat{E}$. 
Then the natural map
\[E\lbrack m\rbrack \to \widehat{E}\lbrack m\rbrack\]
is bijective.
\end{lem}
\begin{proof}
The map is injective by \cite{Si}, Ch. VII, Prop. 3.1(b). It remains to show that
it is surjective. We have a commutative diagram
\[ \begin{CD} 
       0 @> >> E_1(\overline{K}) @>{f_1} >> E(\overline{K}) @>{f_2} >> \widehat{E}(\overline{k}) @> >> 0 \\ 
       @. @ VV{\cdot m} V @ VV {\cdot m} V @ VV {\cdot m} V \\ 
       0 @> >> E_1(\overline{K}) @>{f_1} >> E(\overline{K}) @>{f_2} >> \widehat{E}(\overline{k}) @> >> {0,} \\ 
       \end{CD} 
       \] 
where $E_1(\overline{K})$ is the set of points on $E(\overline{K})$ reducing
to $0$. Let $x$ be an element of $\widehat{E}\lbrack m\rbrack$. Let
$y\in f_2^{-1}(\{x\})$. Let $z\in f_1^{-1}(\{m\cdot y\})$. By \cite{Si}, Ch. VII,
Prop. 2.2 and Ch. IV, Prop. 2.3(b), 
the map $E_1(\overline{K})\stackrel{\cdot m}{\longrightarrow} E_1(\overline{K})$ is surjective. Choose
$w\in E_1(\overline{K})$ such that $m w = z$. Then $m\cdot f_1(w) = f_1(m\cdot w) =
f_1(z) = m\cdot y$. Hence $m\cdot (y-f_1(w))=0$. Since $f_2\circ f_1 = 0$,
$f_2(y-f_1(w)) = f_2(y) = x$. Thus $(y-f_1(w))$ is an element of 
$E(\overline{K})\lbrack m\rbrack$ mapping to $x$.
\end{proof}
\begin{lem}\label{lem:thud}
 Let $K$ be a $\mathfrak{p}$-adic field. Let $L$ be a finite Galois 
extension of $\widetilde{K}$. Let $E_1$, $E_2$ be elliptic curves over $K$ with good
reduction over $L$. Suppose that $E_1$ and $E_2$ reduce to the same curve over the residue field
of $L$. Then $W_\mathfrak{p}(E_1) = W_\mathfrak{p}(E_2)$.
\end{lem}
\begin{proof}
Let $p$ be the characteristic of the residue field of $K$. 
Let $k$ and $l$ be the residue fields of $K$ and $L$, respectively.
The root number $W_\mathfrak{p}(E)$
of an elliptic curve $E$ over $K$ is determined by the 
canonical representation of the Weil-Deligne
group $\mathcal{W}'(\overline{K}/K)$ on the Tate module $T_\ell(E)$, where 
$\ell$ is any prime different from $p$. If $E$ has potential good reduction, we 
can consider the Weil group $\mathcal{W}(\overline{K}/K)$ together with its natural 
representation on $T_\ell(E)$ instead of the Weil-Deligne group and its representation.

Now let $E$ have good reduction over $L$. Let $\mathfrak{q}$ be the prime ideal of $L$.
The natural map $f$ from $E\lbrack \ell^n\rbrack$, $n\geq 1$, to 
$(E \mo \mathfrak{q})\lbrack \ell^n \rbrack$
commutes with the natural actions of
$\mathcal{W}(\overline{K}/K)$ on $E\lbrack \ell^n\rbrack$ and on 
$(E \mo \mathfrak{q})\lbrack \ell^n \rbrack$. 
By Lemma \ref{lem:pardiez}, $f$ is bijective. Hence the
action of $\mathcal{W}(\overline{K}/K)$ on
$E\lbrack \ell^n\rbrack$ is given
by the action of $\mathcal{W}(\overline{K}/K)$ on
$(E \mo \mathfrak{q})\lbrack \ell^n\rbrack$.
Since $l$ is algebraically closed, $(E \mo \mathfrak{q})\lbrack \ell^n\rbrack$
is a subset of $E\mo \mathfrak{q}$. Therefore, the action of $W(\overline{K}/K)$
on $E\lbrack \ell^n\rbrack$ is given by the action of
$W(\overline{K}/K)$ on $E \mo \mathfrak{q}$. 
The action of $\mathcal{W}(\overline{K}/K)$
on \[T_\ell(E) = \lim_{\leftarrow} E\lbrack \ell^n\rbrack\] is thus given
by its action on 
$E\mo \mathfrak{q}$.

Therefore, if $E_1$ and $E_2$ have the same reduction $\mo \mathfrak{q}$,
they have the same local root number $W_\mathfrak{p}(E_1) = W_\mathfrak{p}(E_1)$.
\end{proof} 

\begin{lem}\label{lem:coronet}
Let $K$ be a $\mathfrak{p}$-adic field. Let $\mathcal{E}$ be an elliptic
curve over $K(t)$. Let $S$ be the set of all $t\in K$ such that 
$\mathcal{E}(t)$ is an elliptic curve over $K$ with potential good reduction.
Then the map \[t\mapsto W_\mathfrak{p}(\mathcal{E}(t))\] on $S$ is affinely pliable.
\end{lem}
\begin{proof}
Let $K'/\widetilde{K}$ be as in Lemma \ref{lem:ugleg}. Let $L/\widetilde{K}$
be the Galois closure of $K'/\widetilde{K}$. Since $L/\widetilde{K}$
is the Galois closure of a finite extension, it is itself a finite
extension. The statement then follows immediately from Lemmas
\ref{lem:sundh} and \ref{lem:thud}.
\end{proof}

\begin{lem}\label{lem:rootofall}
Let $K$ be a $\mathfrak{p}$-adic field. 
Let $E$ be an elliptic curve over $K$ given by $c_4, c_6\in K$. Assume $E$ has
 potentially multiplicative reduction. 
Then
\[\begin{aligned}
W_{\mathfrak{p}}(E) &= \quadrec{-1}{p} &\text{if $E$ has additive reduction over $K$,}\\
W_{\mathfrak{p}}(E) &= -\quadrec{-c_6(E) \pi^{-v(c_6(E))}}{\mathfrak{p}}
 \;\;&\text{if $E$ has multiplicative reduction over $K$,}
\end{aligned}\]
where $\pi$ is any prime element of $K$.
\end{lem}
\begin{proof}
This is a classical result that we will translate from the terms presented in
\cite{Ro2}, Section 19. The statement there is as follows.
If $E$ has additive reduction over $K$, then 
$W_{\mathfrak{p}} = \chi(-1)$, where $\chi$ is the ramified character of $K^*$.
If $E$ has multiplicative reduction over $K$, then
\[W_{\mathfrak{p}} = \begin{cases} -1 &\text{if $E$ has split multiplicative reduction,}\\
1 &\text{if $E$ has non-split multiplicative reduction.}\end{cases}\]

Suppose $E$ has additive reduction over $K$. Since $v_K(-1)=0$, $\chi(-1)$
equals $\quadrec{-1}{\mathfrak{p}}$ and we are done.

Suppose that $E$ has multiplicative reduction over $K$ and $\mathfrak{p}$ does not lie
over $2$. Then the reduced curve $E \mo \mathfrak{p}$ has an equation of the form
\[y^2 = x^3+a x^2,\;\;\;\;a\in (\mathfrak{O}_K/\mathfrak{p})^*\]
 (see, e.g., \cite{Si}, App. A, Prop. 1.1). The tangents of the curve at the node 
$(x,y) = (0,0)$ are $\pm \sqrt{a}$. Thus, the reduction is split if and only if 
$a$ is a square. Since the parameter $\overline{c}_6$ of $E \mo \mathfrak{p}$
equals $-64 a^3$, we have that $a$ is a square if and only if 
$\quadrec{-\overline{c}_6}{\mathfrak{p}}=1$. Now $\overline{c}_6$ is the
reduction $\mo \mathfrak{p}$ of the parameter $c_6'$ of a minimal Weierstrass
equation for $E$. Since $E$ has multiplicative reduction, we can take 
$c_6'=c_6\cdot \pi^{-v_{\mathfrak{p}}(c_6)}$. (Notice that $v_{\mathfrak{p}}$ is
even, and thus the choice of $\pi$ is irrelevant.) The statement follows immediately.

Suppose that $E$ has multiplicative reduction over $K$ and $\mathfrak{p}$ lies over $2$.
Then every element of $\mathfrak{O}_K/\mathfrak{p}$ is a square, and thus 
(a) the reduction must be split, and (b) 
$\quadrec{-c_6(E) \pi^{-v(c_6(E))}}{\mathfrak{p}} = 1$. The statement follows.
\end{proof}
\begin{lem}\label{lem:facilnomas}
Let $K$ be a $\mathfrak{p}$-adic field. Let $\mathcal{E}$ be an elliptic
curve over $K(t)$. Let $S$ be the set of all $t\in K$ such that 
$\mathcal{E}(t)$ is an elliptic curve over $K$ with potential multiplicative
 reduction.
Then the map \[t\mapsto W_\mathfrak{p}(t)\] on $S$ is affinely pliable.
\end{lem}
\begin{proof}
For $t\in S$, the curve $\mathcal{E}(t)$ has multiplicative reduction over $K$
if and only if $v(c_6(\mathcal{E}(t)))$ is divisible by $6$. If $\mathcal{E}(t)$
has multiplicative reduction over $K$, its root number
\[W_\mathfrak{p}(E) = -\quadrec{-c_6(E) \pi^{-v(c_6(E))}}{\mathfrak{p}}\]
depends only on 
the coset $c_6(\mathcal{E}(t))\cdot (1 + \pi \mathfrak{O}_K)$.
If $\mathcal{E}(t)$ has additive reduction over $K$, its root number
equals the constant $\quadrec{-1}{\mathfrak{p}}$. 

Therefore $W_{\mathfrak{p}}(\mathcal{E}(t))$ depends only on the coset of
$c_6(\mathcal{E}(t))\cdot (1 + \pi \mathfrak{O}_K)$ in which $c_6(\mathcal{E}(t))$. 
By Proposition
\ref{prop:simil} it follows that $W_\mathfrak{p}(\mathcal{E}(t))$ is 
affinely pliable.
\end{proof}

\begin{lem}\label{lem:choose}
 Let $K$ be a $\mathfrak{p}$-adic field. Let $\mathcal{E}$
be an elliptic curve over $K(t)$. For $t\in K$, let
\[\begin{aligned}
f_1(t) &= \lbrack \text{$\mathcal{E}(t)$ has potential good reduction}\rbrack,\\
f_2(t) &= \lbrack \text{$\mathcal{E}(t)$ has potential multiplicative reduction}\rbrack,\\
f_3(t) &= \lbrack \text{$\mathcal{E}(t)$ is singular}\rbrack .\end{aligned}\]
Then $f_1,f_2,f_3:K\to \{0,1\}$ are affinely pliable.
\end{lem}
\begin{proof}
Since $\mathcal{E}(t)$ is singular for finitely many $t\in K$, $f_1$ is affinely
pliable. 
If $\mathcal{E}(t)$ is non-singular, then $\mathcal{E}(t)$ has potential multiplicative
reduction if and only if $v(j(\mathcal{E}(t)))>0$. Thus, for all but
finitely many $t$, both $f_2(t)$ and $f_3(t)$ depend only on $v(j(\mathcal{E}(t)))$.
By Proposition \ref{prop:simil}, $f_2$ and $f_3$ are affinely pliable.
\end{proof} 

\begin{prop}\label{prop:corona}
Let $K$ be a $\mathfrak{p}$-adic field. Let $\mathcal{E}$ be an elliptic
curve over $K(t)$. Then the map \[t\mapsto W_\mathfrak{p}(\mathcal{E}(t))\]
on $K$ is affinely pliable.
\end{prop}
\begin{proof}
Immediate from Lemmas \ref{lem:coronet}, \ref{lem:facilnomas} and \ref{lem:choose}.
\end{proof}

\begin{prop}\label{prop:krone}
Let $K$ be a number field. Let $\mathfrak{p}\in I_K$ be a prime ideal. Let $\mathcal{E}$
be an elliptic curve over $K(t)$. Then the map
\[t\mapsto W_{\mathfrak{p}}(\mathcal{E}(t))\]
on $K$ is affinely pliable.
\end{prop}
\begin{proof}
Denote by $\mathcal{E}_{\mathfrak{p}}$ be the elliptic curve over $K_{\mathfrak{p}}(t)$
defined by the same equation as $\mathcal{E}$. For $t\in K$, the elliptic curve
$\mathcal{E}_{\mathfrak{p}}(t)$ is the localization $(\mathcal{E}(t))_{\mathfrak{p}}$
of $\mathcal{E}(t)$ at $\mathfrak{p}$. The local root number $W_{\mathfrak{p}}(E)$
of an elliptic curve over $K$ is by definition equal to the root number 
$W_{\mathfrak{p}}(E_{\mathfrak{p}})$ of the localization $E_{\mathfrak{p}}$ of $E$
at $\mathfrak{p}$. By Proposition \ref{prop:corona},
$t\mapsto W_{\mathfrak{p}}(\mathcal{E}_{\mathfrak{p}}(t))$
is an affinely pliable map on $K_{\mathfrak{p}}$. Therefore, its restriction
\[t\mapsto W_{\mathfrak{p}}(\mathcal{E}_{\mathfrak{p}}(t)) =
 W_{\mathfrak{p}}((\mathcal{E}(t))_{\mathfrak{p}}) = W_{\mathfrak{p}}(\mathcal{E}(t))\]
to $K$ is an affinely pliable map on $K$. 
\end{proof}
\subsection{Pliable functions and reciprocity}\label{subs:quadr}

For the following it will be convenient to work in a slightly more abstract
fashion. Let $K$ be a number field. Let $\mathcal{C}_i$, $i\geq 0$, be a multiplicatively closed set of
functions from $\mathfrak{O}_K^i$ to a multiplicative  abelian group 
$\mathcal{G}$. Let $\mathcal{D}$
be a multiplicatively closed set of functions from $\mathfrak{O}_K^2$ to $\mathcal{G}$
such that $(x,y)\mapsto
f(F_1(x,y),\dotsc , F_n(x,y))$ belongs to $\mathcal{D}$ for any $f\in 
\mathcal{C}_n$ and any homogeneous polynomials $F_1,\dotsc ,F_n\in \mathfrak{O}_K[x,y]$.

We want to define a family of operators $\lbrack,\rbrack$ that we may
manipulate much like reciprocity symbols.
Consider a function $\lbrack ,\rbrack_{\mathfrak{d}} : 
\{(x,y)\in (\mathfrak{O}_K-\{0\})^2 :
\gcd(x,y)|\mathfrak{d}^\infty\} \to \mathcal{G}$ 
for every non-zero ideal $\mathfrak{d}\in I_K$. Assume
that $\lbrack , \rbrack_{\mathfrak{d}}$ satisfies the following conditions:

\begin{enumerate}
\item $\lbrack a b, c\rbrack _{\mathfrak{d}} = \lbrack a,c\rbrack _{\mathfrak{d}} \cdot \lbrack
b,c\rbrack _{\mathfrak{d}},$

\item $\lbrack a, b c\rbrack _{\mathfrak{d}} = \lbrack a,b\rbrack _{\mathfrak{d}} \cdot \lbrack
a,c\rbrack _{\mathfrak{d}},$

\item $\lbrack a , b\rbrack _{\mathfrak{d}} = \lbrack a + b c,b\rbrack _{\mathfrak{d}}$ provided that $a + b c\ne 0$,

\item $\lbrack a, b\rbrack _{\mathfrak{d}} = f_{\mathfrak{d}}(a,b)\cdot \lbrack b,a\rbrack _{\mathfrak{d}}, \text{%
\;where $f_{\mathfrak{d}}$ is a function in $\mathcal{C}_2$},$

\item $\lbrack a, b\rbrack _{\mathfrak{d}} = f_{\mathfrak{d},b}(a),
 \text{\;where $f_{\mathfrak{d},b}$ is a function in 
$\mathcal{C}_1$,}$

\item $\lbrack a, b\rbrack _{\mathfrak{d}_1} = 
 f_{\mathfrak{d}_1,\mathfrak{d}_2}(a,b) \lbrack a, b\rbrack _{\mathfrak{d}_2} \text{%
\:for $\mathfrak{d}_1|\mathfrak{d}_2$, \;where $f$ is a function in $\mathcal{C}_2$.}$
\end{enumerate}

\begin{prop}
\label{prop:wichtig} Let $F, G\in \mathfrak{O}_K\lbrack x,y\rbrack$ be
homogeneous polynomials without common factors. 
Let $\mathfrak{d}$ be a non-zero
ideal of $\mathfrak{O}_K$
 such that $\gcd(F(x,y),G(x,y)) | \mathfrak{d}^\infty$ for all coprime 
$x, y\in \mathfrak{O}_K$. Then there is a function $f$ in $\mathcal{D}$ such that 
\begin{equation*}
\lbrack F(x,y),G(x,y) \rbrack_{\mathfrak{d}} = f(x,y) \lbrack x,y \rbrack_1^{ (\deg F)
(\deg G)}
\end{equation*}
for all but finitely many elements $(x,y)$ of $\{(x,y)\in (\mathfrak{O}_K-\{0\})^2: \gcd(x,y)=1\}$.
\end{prop}

\begin{proof}
If $\deg(G) = 0$ the result follows from condition (5). If $\deg(F)=0$ the
result follows from (4) and (5). If $F$ and $G$ is reducible, the statement
follows by (1) or (2) from cases with lower $\deg(F) + \deg(G)$. If $F$ is
irreducible and $G = c x$, $c$ non-zero, then by (1), (2), (3) and (4), 
\begin{equation*}
\begin{aligned} \njac{F(x,y)}{G(x,y)}_{\mathfrak{d}} &=
\njac{a_0 x^k + a_1 x^{k-1} y + \dotsb + a_k y^k }{c x}_{\mathfrak{d}} \\&=
\njac{F(x,y)}{c}_{\mathfrak{d}} \cdot \njac{a_k y^k}{x}_{\mathfrak{d}} \\ &=
\njac{F(x,y)}{c}_{\mathfrak{d}} \cdot \njac{a_k}{x}_{\mathfrak{d}}\cdot \njac{y}{x}_{\mathfrak{d}}^k \\ &=
\njac{F(x,y)}{c}_{\mathfrak{d}} \cdot \njac{a_k}{x}_{\mathfrak{d}}\cdot f_{\mathfrak{d}}^{-k}(x,y) g_{1,\mathfrak{d}}^k(x,y)\njac{x}{y}_1^k
\end{aligned}
\end{equation*}
for some $f_{\mathfrak{d}},g_{1,\mathfrak{d}}\in \mathcal{C}$, and the result follows from (5), the definition
of $\mathcal{D}$ and the already treated case of $\lbrack {\text{constant}},{%
x}\rbrack_{\mathfrak{d}}$. The same works for $F$ irreducible, $G = c y$. The case of $G$
irreducible, $F = \text{$c x$ or $c y$}$ follows from (4) and the foregoing.
For $F$, $G$ irreducible, $\deg(F)<\deg(G)$, we apply (4). We are left with
the case of $F$, $G$ irreducible, $F,G \ne cx, cy$, $\deg(F)\geq \deg(G)$.
Write $F = a_0 x^k + \dotsb + a_k y^k$, $G = b_0 x^l + b_1 x^{l-1} y +
\dotsb + b_l y^l$. Then 
\begin{equation*}
\begin{aligned}
\njac{F(x,y)}{G(x,y)}_{\mathfrak{d}} &= f_{\mathfrak{d},b_0 \mathfrak{d}}(x,y) \njac{F(x,y)}{G(x,y)}_{\mathfrak{d} b_0} \\&=
f_{\mathfrak{d},b_0 \mathfrak{d}}(x,y) \njac{b_0}{G(x,y)}_{b_0 \mathfrak{d}} \njac{b_0 F(x,y)}{G(x,y)}_{b_0 \mathfrak{d}} \\ &=
f_{\mathfrak{d},b_0 \mathfrak{d}}(x,y) \njac{b_0}{G(x,y)}_{b_0 \mathfrak{d}} \njac{b_0 F(x,y) - a_0 G(x,y)}{G(x,y)}_{b_0 \mathfrak{d}}
\end{aligned}
\end{equation*}
for all coprime $x$, $y$ such that $b_0 F(x,y) - a_0 G(x,y)\ne 0$.
(Since $b_0 F(x,y) - a_0 G(x,y)$ is a non-constant homogeneous polynomial, there are only
finitely many such pairs $(x,y)$.)
The coefficient of $x^k$ in $b_0 F(x,y)
- a_0 G(x,y)$ is zero. Hence $b_0 F(x,y) - a_0 G(x,y)$ is a multiple of $y$.
Either it is reducible or it is a constant times $y$. Both cases have
already been considered.
\end{proof}

Now let $\mathcal{G}$ be the group $\{-1,1\}$,
$\mathcal{C}_1$ the set of pliable functions on $\mathfrak{O}_K$,
$\mathcal{C}_2$ the set of pliable functions on $\mathfrak{O}_K^2$ with
$\vec{q}_j \in \{(1,0),(0,1)\}$ for every $j$ and $\mathcal{D}$ the set of
pliable functions on $\mathfrak{O}_K^2$. Let
\begin{equation}  \label{eq:brack}
\lbrack {a},{b}\rbrack_{\mathfrak{d}} = \prod_{\mathfrak{p}\nmid 2 \mathfrak{d}}
\left(\frac{a}{\mathfrak{p}}\right)^{v_\mathfrak{p}(b)}, 
\end{equation}
where $\left(\frac{\cdot}{\mathfrak{p}}\right)$ is the quadratic reciprocity symbol.
The defining condition on $\mathcal{D}$ holds by Proposition \ref%
{prop:plinvar}. Properties (1), (2) and (3) are immediate. Property (5) follows immediately
from 
the fact that $\left(\frac{a}{\mathfrak{p}}\right)$ depends on $a$ only as an element of $K^*/(K^*)^2$;
clearly $(K^*)^2$ is an open subgroup of $K^*$.
It remains to prove (4) and (6).

\begin{lem}
\label{lem:chunky} Given a non-zero ideal 
$\mathfrak{d}$ of $\mathfrak{O}_K$, there is a pliable function $f$ on 
$\mathfrak{O}_K^2$ with 
$q_j\in \{(1,0),(0,1)\}$
such that 
\begin{equation*}
\prod_{\mathfrak{p}\nmid 2 \mathfrak{d}} \left(\frac{a}{\mathfrak{p}}\right)^{v_\mathfrak{p}(b)}
= f(a,b) \prod_{\mathfrak{p}\nmid 2 \mathfrak{d}} \left(\frac{b}{\mathfrak{p}}\right)^{v_\mathfrak{p}(a)}
\end{equation*}
for all non-zero $a,b \in \mathfrak{O}_K$ with $\gcd(a,b)|\mathfrak{d}$.
\end{lem}
\begin{proof}
Let $\left(\frac{a,b}{\mathfrak{p}}\right)$ be the quadratic Hilbert symbol. For $a$, $b$ coprime,
\[\prod_{\mathfrak{p}\nmid 2 \mathfrak{d}} \left(\frac{a}{\mathfrak{p}}\right)^{v_\mathfrak{p}(b)}
= \mathop{\prod_{\mathfrak{p}\nmid 2 \mathfrak{d}}}_{\mathfrak{p}|b}
\left(\frac{a}{\mathfrak{p}}\right)^{v_\mathfrak{p}(b)}
= \mathop{\mathop{\prod_{\mathfrak{p}\nmid 2 \mathfrak{d}}}_{\mathfrak{p}|b}}_{\mathfrak{p}\nmid a}
\left(\frac{a}{\mathfrak{p}}\right)^{v_\mathfrak{p}(b)}
= \mathop{\mathop{\prod_{\mathfrak{p}\nmid 2 \mathfrak{d}}}_{\mathfrak{p}|b}}_{\mathfrak{p}\nmid a}
\left(\frac{b,a}{\mathfrak{p}}\right)
= \mathop{\mathop{\prod_{\mathfrak{p}\nmid 2 \mathfrak{d}}}_{\mathfrak{p}|b}}_{\mathfrak{p}\nmid a}
\left(\frac{a,b}{\mathfrak{p}}\right) .\]
Similarly
\[\prod_{\mathfrak{p}\nmid 2 \mathfrak{d}} \left(\frac{b}{\mathfrak{p}}\right)^{v_{\mathfrak{p}}(a)}
= \mathop{\mathop{\prod_{\mathfrak{p}\nmid 2 \mathfrak{d}}}_{\mathfrak{p}|a}}_{\mathfrak{p}\nmid b} 
\left(\frac{a,b}{\mathfrak{p}}\right) .\]
Hence \[\begin{aligned}
\prod_{\mathfrak{p}\nmid 2 \mathfrak{d}} \left(\frac{a}{\mathfrak{p}}\right)^{v_{\mathfrak{p}}(b)}
 \prod_{\mathfrak{p}\nmid 2 \mathfrak{d}} \left(\frac{b}{\mathfrak{p}}\right)^{v_{\mathfrak{p}}(a)} &= 
\mathop{\prod_{\mathfrak{p}\nmid 2 \mathfrak{d}}}_{\mathfrak{p}| a b} \left(\frac{a,b}{\mathfrak{p}}\right) =
\mathop{\prod_{\mathfrak{p}\nmid 2}}_{\text{$\mathfrak{p}\nmid b$ or $\mathfrak{p}\nmid a b$}}
\left(\frac{a,b}{\mathfrak{p}}\right) \\ &=
\left(\frac{a,b}{\infty}\right) \left(\frac{a,b}{2}\right) 
\prod_{\mathfrak{p}|\gcd(\mathfrak{d}, a b)} \left(\frac{a,b}{\mathfrak{p}}\right) .
\end{aligned}\]
Now note that $\left(\frac{a,b}{\mathfrak{p}}\right)$ and
$\left(\frac{a,b}{\infty}\right)$ are pliable 
on $(\mathfrak{O}_K - \{0\})^2$ with \[\{(v_j,U_j,\vec{q}_j)\} = 
\{(v,(K^*)^2,(1,0)),(v,(K^*)^2,(0,1))\}.\] Therefore
 \[\left(\frac{a,b}{\infty}\right) \left(\frac{a,b}{2}\right) 
\prod_{\mathfrak{p}|\gcd(\mathfrak{d}, a b)} \left(\frac{a,b}{\mathfrak{p}}\right) \]
is pliable on $\{\mathfrak{O}_K-\{0\}\}^2$ with $q_j\in \{(1,0),(0,1)\}$. Set
  \[f(a,b) = \left(\frac{a,b}{\infty}\right) \left(\frac{a,b}{2}\right) 
\prod_{\mathfrak{p}|\gcd(\mathfrak{d}, a b)} \left(\frac{a,b}{\mathfrak{p}}\right) \]
\end{proof}

\begin{lem}
Given non-zero $\mathfrak{d}_1$, 
$\mathfrak{d}_2$ with $\mathfrak{d}_1|\mathfrak{d}_2$, 
there is a pliable function $f$
such that 
\begin{equation*}
\prod_{\mathfrak{p}\nmid 2 \mathfrak{d}_1}
\left(\frac{a}{\mathfrak{p}}\right) = f(a,b) 
\prod_{\mathfrak{p}\nmid 2 \mathfrak{d}_2}
\left(\frac{a}{\mathfrak{p}}\right) 
\end{equation*}
for all $a$, $b$ with $\gcd(a,b)|\mathfrak{d}_1$.
\end{lem}

\begin{proof}
We have 
\[
\prod_{\mathfrak{p}\nmid 2 \mathfrak{d}_1}
\left(\frac{a}{\mathfrak{p}}\right) = \mathop{\prod_{\mathfrak{p}|2 \mathfrak{d}_2}}_{\mathfrak{p}\nmid 2 \mathfrak{d}_1}
 \left(\frac{a}{\mathfrak{p}}\right)^{v_{\mathfrak{p}(b)}} \cdot
\prod_{\mathfrak{p}\nmid 2 \mathfrak{d}_2}
\left(\frac{a}{\mathfrak{p}}\right) .\]
Since $a\to \left(\frac{a}{\mathfrak{p}}\right)$ is pliable, we are done.
\end{proof}

Hence we obtain
\begin{cor}[to Proposition \ref{prop:wichtig}]
\label{cor:important} Let $F, G\in \mathfrak{O}_K\lbrack x,y\rbrack$ be
homogeneous polynomials without common factors. Let $\mathfrak{d}$ be a 
non-zero ideal of $\mathfrak{O}_K$
 such that \[\gcd(F(x,y),G(x,y)) | \mathfrak{d}^\infty\] 
for all coprime integers $x$, $y$. 
Let $\lbrack , \rbrack$ be as in (\ref{eq:brack}). Then there is a
pliable function $f$ on $\mathfrak{O}_K^2$ such that 
\begin{equation*}
\lbrack F(x,y),G(x,y) \rbrack_{\mathfrak{d}} = f(x,y) \text{\;\;(if $\deg F$ or $\deg G$
is even)}
\end{equation*}
\begin{equation*}
\lbrack F(x,y),G(x,y) \rbrack_{\mathfrak{d}} = f(x,y) \lbrack x,y \rbrack_1 \text{\;\;(if 
$\deg F$ and $\deg G$ are odd)}
\end{equation*}
for all coprime $x, y\in \mathfrak{O}_K$ (if $\deg F$ or $\deg G$ is even)
or all coprime, non-zero $x,y\in \mathfrak{O}_K$ (if $\deg F$ and $\deg G$ are odd).
\end{cor}
\begin{proof}
By Proposition \ref{prop:wichtig}, the statement holds for all but finitely many
elements $(x,y)$ of $\{(x,y)\in (\mathfrak{O}_K-\{0\})^2 : \gcd(x,y)=1\}$.
By Lemma \ref{lem:redef}, $f$ can be redefined for finitely many elements of the
domain and still be pliable.
\end{proof}
\subsection{Averages and pliable functions}

What we will now show is essentially that, given a pliable function $f$ and
a function $g$ whose average over lattices of small index is well-known, we
can tell the average of $f\cdot g$ over $\mathbb{Z}^2$. By Corollary \ref%
{cor:important} this will imply, for example, that $\sum \lbrack x^2 + 3 x y
- 2 y^2, 4 x ^3 - x y^2 + 7 y^3 \rbrack_{\mathfrak{d}}\, g(x,y) = o(N^2)$ provided that $%
\sum_{(x,y)\in L} g(x,y) = o(N^2)$ for $L$ small. 

We may start with the parallel statements for affinely pliable functions.

\begin{lem}\label{lem:marks}
 Let $U$ be an open subgroup of $\mathbb{R}^*$. Let $t_1<t_2<\dotsb < t_n$
be real numbers. If $t$, $t'$ are real numbers with $t<t_1$, $t'<t_1$ or $t>t_n$, $t'>t_n$, then
$t-t_i$ and $t'-t_i$ lie in the same coset of $U$ for every $1\leq i\leq n$.
\end{lem}
\begin{proof}
If $U = \mathbb{R}^*$, the statement is trivially true. If $U = \mathbb{R}^+$, note
that $t-t_i$ and $t'-t_i$ lie in the same coset of $U$ if and only if
$\sgn(t-t_i) = \sgn(t'-t_i)\ne 0$. The statement is then obvious.
\end{proof}
\begin{lem}\label{lem:affparta}
Let $p$ be a prime. Let $U$ be an open subgroup of $\mathbb{Q}_p^*$. Let 
$t_1,\dotsc,t_n\in \mathbb{Q}_p$.
Then there is a partition
\[\mathbb{Z} = A_{\infty} \cup \bigcup_{i\geq 0} \bigcup_{k\in K} A_{i,k}\]
such that
\begin{enumerate}
\item $K$ is a finite set,
\item $A_{\infty}$ is a finite subset of $\mathbb{Z}$,
\item $A_{i,k}$ is a disjoint union of at most $c_1$ arithmetic progressions of modulus 
$p^{i+c_2}$,
\item for every $i_0\geq 0$, $A_{\infty} \cup \bigcup_{i\geq i_0} \bigcup_{k\in K} A_{i,k}$ is a disjoint
union of at most $c_1$ arithmetic progressions of modulus $p^{i_0}$,
\item for any choice of $i\geq 0$, $j=1,\dotsc,n$, $k\in K$ and all $t,t'\in A_{i,k}$,
$t-t_j$ and $t'-t_j$ lie in the same coset of $U$.
\end{enumerate}
The positive integers $c_1$, $c_2$ depend only on $p$, $U$ and $t_1,\dotsc, t_n$.
\end{lem}
\begin{proof}
We can assume that $U = 1 + p^l \mathbb{Z}_p$, $l\geq 1$. If $t$, $t'$ lie in the same coset
of $U$, then $t-t_j$ and $t'-t_j$ lie in the same coset of $U$ for all $t_j\in \mathbb{Q}_p -
\mathbb{Z}_p$. Hence we can assume $t_j\in \mathbb{Z}_p$ for all $1\leq j\leq n$.

Let $d = 1 + \max_{j_1\ne j_2} v_p(t_{j_1} - t_{j_2})$. Define
\[\begin{aligned}
K &= ((\mathbb{Z}_p/U)^* \times \{0,1,\dotsc,d\})^n,\\
A_i &= \{t\in \mathbb{Z} : \max_j v_p(t-t_j) = i\},\\
A_{\infty} &= \{t_1,\dotsc,t_n\} \cap \mathbb{Z},\\
A_{i,((k_{1 1},k_{1 2}),\dotsc,(k_{n 1},k_{n 2}))} &=
\{t\in A_i : \frac{t-t_j}{p^{v_p(t-t_j)}} \equiv k_{j 1} \mo p^l,
\min(v_p(t-t_j),d) = k_{j 2}\}.\end{aligned}\]
Statements (1) and (2) hold by definition. We can write $A_i$ in the form
\[A_i = \bigcup_{1\leq j\leq n} (t_j + p^i \mathbb{Z})\]
Since any two arithmetic progressions $t_j + p^i \mathbb{Z}$, 
$t_{j'} + p^i \mathbb{Z}$ of
the same modulus are either
disjoint or identical, it follows that $A_i$ is the union of at most $n$ disjoint arithmetic
progressions of modulus $p^i$. Clearly 
$A_{i_0} = A_{\infty} \cup \bigcup_{i\geq i_0} \bigcup_{k\in K} A_{i,k}$. Hence (4) holds.

For $i<d$,
\[
A_{i,((k_{1 1},k_{1 2}),\dotsc,(k_{n 1},k_{n 2}))} =
\{t\in A_i : \frac{t-t_j}{p^{v_p(t-t_j)}} \equiv k_{j 1} \mo p^l,\:
 v_p(t-t_j) = k_{j 2}\}.\]
If $\max_j k_{j 2} \ne i$, then $A_{i,((k_{1 1},k_{1 2}),\dotsc,(k_{n 1},k_{n 2}))} = \emptyset$.
Otherwise, 
\[\begin{aligned}
A_{i,((k_{1 1},k_{1 2}),\dotsc,(k_{n 1},k_{n 2}))} &=
\bigcap_{1\leq j\leq n} \{t\in \mathbb{Z} : t\equiv p^{k_{j 2}} k_{j 1} + t_j \mo 
 p^{l + k_{j 2}}\}\\
&= \bigcap_{1\leq j\leq n} \{t\in \mathbb{Z} : t - t_j \in k_{j 1} p^{k_{j 2}} U\}.\end{aligned}\]
Both (3) and (5) follow immediately.

For $i\geq d$,
\[
A_{i,((k_{1 1},k_{1 2}),\dotsc,(k_{n 1},k_{n 2}))} =
\{t\in A_i : \frac{t-t_j}{p^{v_p(t-t_j)}} \equiv k_{j 1} \mo p^l,\:
 v_p(t-t_j) = k_{j 2}'\},\]
where \[k_{j 2}' = \begin{cases} k_{j 2} & \text{if $k_{j 2}<d$,}\\
i &\text{if $k_{j 2}\geq d$.}\end{cases}\]
Then \[\begin{aligned}
A_{i,((k_{1 1},k_{1 2}),\dotsc,(k_{n 1},k_{n 2}))} &=
\bigcap_{1\leq j\leq n} \{t\in \mathbb{Z} : t\equiv p^{k_{j 2}'} k_{j 1} + t_j \mo
 p^{l + k_{j 2}'}\}\\
&= \bigcap_{1\leq j\leq n} \{t\in \mathbb{Z} : t - t_j \in k_{j 1} p^{k_{j 2}'} U\}.\end{aligned}\]

Again, (3) and (5) follow.
\end{proof}
\begin{lem}\label{lem:affpartb}
Let $p$ be a prime. Let $U$ be an open subgroup of $\mathbb{Q}_p^*$. Let 
$t_1,\dotsc,t_n\in \mathbb{Q}_p$. Let $a$ be an integer, $m$ a non-negative integer.
Then there is a partition
\[\{t\in \mathbb{Z} : t\equiv a \mo p^m\} = B_{\infty} \cup \bigcup_{i\geq m} \bigcup_{k\in K} B_{i,k}\]
such that
\begin{enumerate}
\item $K$ is a finite set,
\item $B_{\infty}$ is a finite subset of $\mathbb{Z}$,
\item $B_{i,k}$ is a disjoint union of at most $c_1$ arithmetic progressions of modulus 
$p^{i+c_2}$,
\item for every $i_0\geq m$, 
$B_{\infty} \cup \bigcup_{i\geq i_0} \bigcup_{k\in K} B_{i,k}$ is a disjoint
union of at most $c_1$ arithmetic progressions of modulus $p^{i_0}$,
\item for any choice of $i\geq m$, $j=1,\dotsc,n$, $k\in K$ and all $t,t'\in B_{i,k}$,
$t-t_j$ and $t'-t_j$ lie in the same coset of $U$.
\end{enumerate}
The positive integers $c_1$, $c_2$ depend only on $p$, $U$ and $t_1,\dotsc, t_n$.
\end{lem}
\begin{proof}
Let $A_{\infty}$, $A_{i,k}$ be as in Lemma \ref{lem:affparta}.
By Lemma \ref{lem:affparta}, (4), 
\[A_{\infty} \cup \bigcup_{i\geq i_0} \bigcup_{k\in K} A_{i,k}\] is a union of
arithmetic progressions of modulus $p^{i_0}$. Hence, for
$i_0\leq m$, either
\[\{t\in \mathbb{Z}: t\equiv a \mo p^m\} \cap
(A_{\infty} \cup \bigcup_{i\geq i_0} \bigcup_{k\in K} A_{i,k}) = \emptyset\]
or
\[\{t\in \mathbb{Z}: t\equiv a \mo p^m\} \subset
A_{\infty} \cup \bigcup_{i\geq i_0} \bigcup_{k\in K} A_{i,k} .\]

Suppose
\[\{t\in \mathbb{Z}: t\equiv a \mo p^m\} \cap
(A_{\infty} \cup \bigcup_{i\geq m} \bigcup_{k\in K} A_{i,k}) = \emptyset .\]
Let $i_0\geq 0$ be the largest integer such that
\[\{t\in \mathbb{Z}: t\equiv a \mo p^m\} \subset
A_{\infty} \cup \bigcup_{i\geq i_0} \bigcup_{k\in K} A_{i,k} .\]
Then \[\{t\in \mathbb{Z}: t\equiv a \mo p^m\} = \bigcup_{k\in K} (A_{i_0,k} \cap 
\{t\in \mathbb{Z}: t\equiv a \mo p^m\}) .\]
Set $B_{m,k} = A_{i_0,k} \cap
\{t\in \mathbb{Z}: t\equiv a \mo p^m\}$, $B_{i,k}=\emptyset$ for $i\ne m$,
$B_{\infty} = \{t\in A_{\infty}:t\equiv a \mo p^m\}$.

Suppose now
\[\{t\in \mathbb{Z}: t\equiv a \mo p^m\} \subset
A_{\infty} \cup \bigcup_{i\geq m} \bigcup_{k\in K} A_{i,k} .\]
For every $i\geq m$, $k\in K$, $A_{i,k} \cap \{t\in \mathbb{Z}: t\equiv a \mo p^m\}$
is equal to either the empty set or to $A_{i,k}$.
Set $B_{i,k}=\emptyset$ for $i<m$, 
$B_{i,k} = A_{i,k} \cap \{t\in \mathbb{Z}: t\equiv a \mo p^m\}$ for $i\geq m$,
$B_{\infty} = \{t\in A_{\infty}: t\equiv a \mo p^m\}$.
\end{proof}
\begin{lem}
\label{lem:sivcor} Let $M$, $R$ and $C$ be positive integers. Let $%
\{a_n\}_{n=1}^\infty$ be such that
\begin{enumerate}
\item $a_n=0$ for all $n$ for which $\rad(n)\nmid R$,
\item $s_d = \sum_n |a_{d n}|$ converges for every $d$,
\item $s_d = O(C/d)$ for $M<d\leq p_0 M$, where $p_0$ is the
largest prime factor of $R$.
\end{enumerate}
Then 
\begin{equation*}
\sum_n a_n = \sum_{n\leq M} a_n + O\left(\frac{C (\log p_0 M)^{\omega(R)}}{M}%
\right),
\end{equation*}
where the implied constant is absolute.
\end{lem}

\begin{proof}
Every $n>M$ satisfying $\rad(n)|R$ has a divisor $M<d\leq p_0 M$. Hence 
\begin{equation*}
\begin{aligned}
\sum_n a_n &= \sum_{n\leq M} a_n + O\left(\sum_{n>M} |a_n|\right)
=\sum_{n\leq M} a_n + O\left(\mathop{\sum_{M<d\leq p_0 M}}_{\rad(d)|R}
\mathop{\sum_n}_{d|n} \; |a_n|\right)\\
&= \sum_{n\leq M} a_n + O\left(\mathop{\sum_{M<d\leq p_0 M}}_{\rad(d)|R} C/d\right) .
\end{aligned}
\end{equation*}
There are at most $\prod_{p|R} (1 + \log_p p_0 M)$ terms in $\sum_{M<d\leq
p_0 M,\, \rad(d)|R}$. Hence 
\begin{equation*}
\sum_n a_n = \sum_{n\leq M} a_n + O\left(\frac{C (\log p_0 M)^{\omega(R)}}{M}%
\right) .
\end{equation*}
\end{proof}

\begin{lem}\label{lem:weneed}
Let $f, g:\mathbb{Z} \to \mathbb{C}$ be given with $\max |f(x)|\leq 1$, $\max |g(x)|\leq 1$.
Let $f$ be affinely pliable with respect to $\{(v_j, U_j, t_j)\}$. Assume that there are
$\eta_N\leq N$, $\epsilon_N\geq 0$ such that for any $a,m\in \mathbb{Z}$, $0<m\leq \eta_N$,
\begin{equation}\label{eq:knot}
\mathop{\sum_{1\leq x\leq N}}_{x\equiv a \mo m} g(x) \ll \frac{\epsilon_N N}{m}.
\end{equation}
Then, for any $a, m\in \mathbb{Z}$, $0<m\leq \eta_N$,
\[\mathop{\sum_{1\leq x\leq N}}_{x\equiv a \mo m} f(x) g(x) \ll 
 \left( \frac{\epsilon_N}{m} + \frac{(\log \eta_N)^c}{\eta_N}\right) N,\]
where $c$ is the number of distinct finite places among $\{v_j\}$ and the implied
constant depends only on the implied constant in (\ref{eq:knot}) and on $\{v_j,U_j,t_j\}$.
\end{lem}
\begin{proof}
Let $\{p_l\}$ be the set of all finite places among $\{v_j\}$. 
Let $\{t_{l,1},\dotsb ,t_{l,n_l}\}$ be
the set of all $t_j$ such that $v_j$ is induced by $p_l$. 
For every $p_l$, Lemma \ref{lem:affpartb} yields a partition
\[\{x\in \mathbb{Z} : x\equiv a \mo p_l^{v_{p_l}(m)}\} = 
 B_{l,\infty} \cup \bigcup_{i\geq v_{p_l}(m)} \bigcup_{k\in K_l} B_{l,i,k}\]
such that $t-t_{l,j}$ and $t'-t_{l,j}$ lie in the same coset of $U_l$ for any
$t,t'\in B_{l,i,k}$ and any $i$, $j$, $k$. Let 
\[
m_0 = \frac{m}{\prod_l p_l^{v_{p_l}(m)}} .\]
Clearly
\begin{equation}\label{eq:onestar}\begin{aligned}
\{x\in \mathbb{Z} : x\equiv a \mo m\} &= \bigcap_l \left(
B_{l,\infty} \cup \bigcup_{i\geq v_{p_l}(m)} \bigcup_{k\in K_l} B_{l,i,k}\right)
\cap (a + m_0 \mathbb{Z})\\
 &= \left(\bigcap_l B_{l,\infty} \right) \cup \mathop{\bigcup_{n\geq 1}}_{\rad(n)|R}
\bigcup_{\{k_l\} \in \prod_l K_l} \bigcap_{l} B_{l,v_p(m n), k_l} \cap
(a + m_0 \mathbb{Z}),
\end{aligned}\end{equation}
where $R = \prod_l p_l$. Let $t_0$ be the largest of all $t_j$ such that $v_j$ is 
an infinite place; see Lemma \ref{lem:marks}. Since $f$ is affinely pliable with 
respect to $\{v_j,U_j,t_j\}$,
it is constant on 
\begin{equation}\label{eq:squig}
\{x\in \mathbb{Z}: x>t_0\} \cap \bigcap_l B_{l,v_p(m n),k_l}\end{equation} 
for any $n\geq 1$ and any
$\{k_l\}\in \prod_l K_l$. Denote the value of $f$ on 
(\ref{eq:squig})
by $f_{n,\{k_l\}}$. Thanks to (\ref{eq:onestar}), we
can write
\[
\begin{aligned}
\mathop{\sum_{1\leq x\leq N}}_{x\equiv a \mo m} f(x) g(x) &= 
 \sum_{1\leq x\leq t_0} f(x) g(x) +
 \mathop{\sum_{t_0< x\leq N}}_{x\in \cap_l B_{l,\infty}} f(x) g(x) \\ &+
 \mathop{\sum_{n\geq 1}}_{\rad(n)|R} \sum_{\{k_l\}\in \prod_l K_l}
 \mathop{\mathop{\sum_{t_0<x\leq N}}_{x\in \cap_l B_{l,v_p(m n),k_l}}}_{
 x\in a + m_0 \mathbb{Z}} f(x) g(x)\\
&= O(1) + \sum_{\{k_l\}\in \prod_l K_l} \mathop{\sum_{n\geq 1}}_{\rad(n) | R}
 f_{n,\{k_l\}} 
 \mathop{\mathop{\sum_{1\leq x\leq N}}_{x\in \cap_l B_{l,v_p(m n),k_l}}}_{
 x\in a + m_0 \mathbb{Z}} g(x) .\end{aligned}
\]
Fix $\{k_l\}\in \prod_l K_l$. Set
\[a_n = f_{n,\{k_l\}} 
\mathop{\mathop{\sum_{1\leq x\leq N}}_{x\in \cap_l B_{l,v_p(m n),k_l}}}_{
 x\in a + m_0 \mathbb{Z}} g(x)\]
if $\rad(n)|R$, $a_n=0$ otherwise.  
Then
\[
\mathop{\sum_{1\leq x\leq N}}_{x\equiv a \mo m} f(x) g(x) =
\sum_n a_n .\]
Let $s_d = \sum_n |a_{d n}|$. From Lemma \ref{lem:affpartb}, (4),
and the fact that $\max |g(x)|\leq 1$, we get that $s_d\ll \frac{N}{m n}$.
Set $C=N/m$. By Lemma \ref{lem:affpartb}, (3),
$\bigcap_l B_{l,v_p(m n),k_l} \cap (a + m_0 \mathbb{Z})$ is the union of at most $c_3 = c_1^{\#\{p_l\}}$ arithmetic
progressions of modulus $c_4 m n$, where $c_4 = \prod_l p_l^{c_2}$. 
Set
$M=\min\left(\frac{\eta_N}{c_4 m}, \frac{N}{p_0 m}\right)$, where 
$p_0 = \max_l p_l$.
We can now
apply Lemma \ref{lem:sivcor}, obtaining
\begin{equation}\label{eq:aligtw}\begin{aligned}
\sum_n a_n &= \sum_{n\leq M} a_n + O\left(\frac{C (\log p_0 M)^{\omega(R)}}{M}\right)\\
 &= \sum_{n\leq M} a_n + \max\left(\frac{N}{\eta_N} (\log \eta_N/m)^{w(R)},
 (\log N/m)^{w(R)}\right)\\
 &= \sum_{n\leq M} a_n + O\left(\frac{N}{\eta_N} (\log \eta_N)^{\omega(R)}\right)
.\end{aligned}\end{equation}
By (\ref{eq:knot}),
\begin{equation}\label{eq:aligth}\begin{aligned}
 \sum_{n\leq M} a_n &\ll \mathop{\sum_{n\leq M}}_{\rad(n)|R} \frac{\epsilon_N N}{m n}
\leq \sum_{\rad(n)|R} \frac{\epsilon_N N}{m n} \\
&= \frac{\epsilon N}{m} \cdot \prod_{p|R} \left(1 + \frac{1}{p} + \frac{1}{p^2} +\dotsb\right)
\ll \frac{\epsilon N}{m} .\end{aligned}\end{equation}
We conclude that
\[\mathop{\sum_{1\leq x\leq N}}_{x\equiv a \mo m} f(x) g(x)\ll
 \frac{\epsilon N}{m} + \frac{N (\log \eta_N)^{\omega(R)}}{\eta_N} .\]
\end{proof}
\begin{lem}
\label{lem:rsub} Let $U$ be an open subgroup of $\mathbb{R}^*$. Let $\{\vec{q%
}_j\}$ be a finite subset of $\mathbb{R}^n$. Then there is a partition 
\begin{equation*}
\mathbb{R}^n = T_1 \cup \dotsb \cup T_k \cup S_1 \cup \dotsb \cup S_l
\end{equation*}
such that

\begin{enumerate}
\item $T_j$ is a hyperplane,

\item $S_i$ is a sector,

\item $\vec{q}_j \cdot \vec{v}_1$ and $\vec{q}_j \cdot \vec{v}_2$ lie in the
same coset $r U$ of $U$ for any $\vec{v}_1, \vec{v}_2 \in S_i$ and all $j$.
\end{enumerate}
\end{lem}
\begin{proof}
We can assume $U=\mathbb{R}^+$. Set $T_i=\{(x,y)\in \mathbb{R}^n :
(x,y)\cdot \vec{q}_j = 0\}$. Let $S_1,\dotsc ,S_l$ be the connected
components of $\mathbb{R}^n - (T_1\cap T_2\cap \dotsb \cap T_k)$.
\end{proof}

We define $\mathbb{A}_p = \{(x,y)\in \mathbb{Z}^2 : p\nmid \gcd(x,y)\}$.

\begin{lem}\label{lem:oques}
Let $p$ be a prime. Let $n$ be a non-negative integer. For any two distinct lattices
$L, L'\subset \mathbb{Z}^2$ of index 
$\lbrack \mathbb{Z}^2 : L\rbrack = \lbrack \mathbb{Z}^2 : L'\rbrack = p^n$,
the two sets $\mathbb{A}_p\cap L$, $\mathbb{A}_p\cap L'$ are disjoint.
\end{lem}
\begin{proof}
Both L and L' contain $(p^n,0)$ and $(0,p^n)$. Suppose $(x,y)\in L \cap L'$,
$p\nmid \gcd(x,y)$. Then the lattice $L''$ generated by $(p^n,0)$, $(0,p^n)$ and
$(x,y)$ is contained in $L\cap L'$. Since the index $\lbrack \mathbb{Z}^2 : L''\rbrack$
of $L''$ is $p^n$, it follows that $L=L'$. Contradiction.
\end{proof}

\begin{lem}\label{lem:pancakestack}
 Let $p$ be a prime. Let $U$ be an open subgroup of $\mathbb{%
Q}_p^*$. Let $\{\vec{q}_j\}_{j\in J}$ be a finite subset of $\mathbb{Q}_p^2$. Then
there is a partition 
\begin{equation*}
\mathbb{A}_p  = A_{\infty} \cup \bigcup_{i\geq 0}
\bigcup_{k\in K} A_{i,k}
\end{equation*}
such that

\begin{enumerate}
\item $K$ is a finite set,

\item $A_{\infty}$ is the union of finitely many sets of the form $A_{x,y} = \{(n x,
n y) : n\in \mathbb{Z}, p\nmid n\}$,

\item $A_{i,k}$ is a disjoint union of at most $c_1$ lattice cosets of index 
$p^{i+c_2}$, 

\item for every $i_0\geq 0$, the set $A_{\infty}\cup \bigcup_{i\geq i_0}
\bigcup_{k\in K} A_{i,k}$ is a disjoint union of at most $c_1$ sets of the form 
$R\cap \mathbb{A}_p$, where $R$ is a lattice of index $p^{i_0}$; any given
$A_{i,k}$, $i\geq i_0$, lies entirely within one such set $R\cap \mathbb{A}_p$;

\item for any choice of $i\geq 0$, $j\in J$, $k\in K$ and all 
$(x_1,y_1), (x_2,y_2)\in A_{i,k}$,
the inner products
$\vec{q}_j \cdot (x_1,y_1)$ and $\vec{q}_j \cdot (x_2,y_2)$ lie in the
same coset of $U$.
\end{enumerate}
\end{lem}
\begin{proof}
We can assume that $U\subset \mathbb{Z}_p^*$ and $\vec{q}_j \in \mathbb{Z}%
_p^2 - (p \mathbb{Z}_p)^2$. Furthermore we can suppose that for every
pair of indices $j_1$, $j_2$, $j_1\ne j_2$, there is no rational number
$c$ such that $\vec{q}_{j_1} = c \vec{q}_{j_2}$. Hence the determinant
\[D_{j_1,j_2} = \left|\begin{matrix} q_{j_1,1} &q_{j_1,2}\\q_{j_2,1} &q_{j_2,2}\end{matrix}
\right|\] is non-zero. Take $(x,y)\in \mathbb{Z}^2$ with $p\nmid x$.
Then 
\[\begin{aligned}
\min(v_p(\vec{q}_{j_1}\cdot (x,y)), v_p(\vec{q}_{j_2}\cdot (x,y)))&\leq
v_p\left(\left|\begin{matrix} \vec{q}_{j_1}\cdot (x,y) & q_{j_1,2}\\
\vec{q}_{j_2}\cdot (x,y) & q_{j_2,2}\end{matrix}\right|\right) \\ &=
v_p\left(
\left|\begin{matrix} q_{j_1,1} & q_{j_1,2}\\q_{j_2,1}& q_{j_2,2}\end{matrix}
\right|  \cdot 
\left|\begin{matrix} x & 0\\y & 1\end{matrix}\right|\right) =
v_p(D_{j_1,j_2}) .\end{aligned}\]
In the same way \[
\min(v_p(\vec{q}_{j_1}\cdot (x,y)), v_p(\vec{q}_{j_2}\cdot (x,y))) \leq 
v_p(D_{j_1,j_2}) \] for $(x,y)\in \mathbb{Z}^2$ with $p\nmid y$.
Setting $d = \max_{j_1\ne j_2} v_p(D_{j_1,j_2})$ we obtain that
for any given pair $(x,y)\in \mathbb{Z}^2$ with $p\nmid \gcd(x,y)$
there can be at most one index $j$ for which
$v_p(\vec{q}_j\cdot (x,y)) > d$.

Let the cosets of $U$ in $\mathbb{Z}_p^*$ be $U_1, U_2,\dotsc U_m$. 
Let $r$
be the least positive integer such that $p^r \mathbb{Z}_p + 1 \subset U$.
Define  
\begin{equation*}
\begin{aligned}
K &= \{(x_0,y_0,a)\in (\mathbb{Z}/p^{d+r})^2 \times \{1,2,\dotsc , m\} :
p\nmid x_0 \vee p\nmid y_0 \},\\
A_{\infty} &= \{(x,y) \in \mathbb{Z}^2 : \exists j \:\text{\st}\,
(x,y)\cdot \vec{q}_j = 0\} \cap 
\{(x,y)\in \mathbb{Z}^2 : p\nmid \gcd(x,y)\}
.\end{aligned}
\end{equation*}
 For $i>d$, let $A_{i,(x_0,y_0,a)}$ be the set of all $(x,y)\in \mathbb{Z}^2$ such that 
$x\equiv x_0 \mo p^{d+r}$, $y\equiv y_0 \mo p^{d+r}$, $\max_j v_p((x,y)\cdot 
\vec{q}_j) = i$ and $p^{-i} (\vec{q}_{j_0}\cdot (x,y)) \in U_a$, where $j_0$
is the only $j$ for which the maximum $\max_j v_p((x,y)\cdot \vec{q}_j) = i$
is attained. For $i\leq d$ and $a>1$, let $A_{i,(x_0,y_0,a)}$ be the empty
set. For $i\leq d$ and $a=1$, let $A_{i,(x_0,y_0,a)}$ be the set of all 
$(x,y) \in \mathbb{Z}^2$ such that $x\equiv x_0 \mo p^{d+r}$, 
$y\equiv y_0 \mo p^{d+r}$ and $\max_j v_p(\vec{q}_j\cdot (x,y)) = i$.
These definitions for $A_{i,k}$, $k\in K$, give us that
\begin{equation}\label{eq:pancake}
A_{\infty}\cup \bigcup_{i\geq i_0} \bigcup_{k\in K} A_{i,k}
= \{(x,y)\in \mathbb{Z}^2 : p\nmid \gcd(x,y),
\max_j v_p((x,y)\cdot \vec{q}_j)\geq i_0\}.\end{equation}

Properties (1) and (2) follow immediately from our definitions of $K$,
$A$ and $A_{i,(x_0,y_0,a)}$. Let us verify properties (3) and (4).
For $i_0\geq 0$,
\begin{equation}\label{eq:ystar}
A\cup \bigcup_{i\geq i_0} \bigcup_{k\in K} A_{i,k} =
 \bigcup_{j\in J} \left(\{(x,y)\in \mathbb{Z}^2 : v_p((x,y)\cdot \vec{q}_j)\geq i_0\}
\cap \mathbb{A}_p\right) .\end{equation}
By Lemma \ref{lem:oques}, any two distinct sets in the union on the right hand side of
(\ref{eq:ystar}) are disjoint. Since 
 $\{(x,y)\in \mathbb{Z}^2 : v_p((x,y)\cdot \vec{q}_j)\geq i_0\}$ is a lattice of
index $p^{i_0}$, we have proven the first half of (4).
Let $(x,y)\in A_{i,(x_0,y_0,a)}$, $i\geq i_0$, 
$j\in J$. To prove the second half of (4), we must show
that we can tell whether $v_p((x,y)\cdot \vec{q}_k)\geq i_0$ from
$i$, $i_0$, $x_0$, $y_0$, $a$ and $j$ alone. If $i_0\leq d$, this is clear: 
$x_0, y_0 \mo p^d$ give us $x, y \mo p^d$.
If $i_0> d$, then $v_p((x,y)\cdot \vec{q}_j)\geq i_0$ if and only if
$v_p((x,y)\cdot \vec{q}_j)>d$. We can tell whether
$v_p((x,y)\cdot \vec{q}_j)>d$ from $x_0, y_0 \mo p^{d+1}$. Hence (4) holds.

For $i\leq d$, each set $A_{i,(x_0,y_0,a)}$ is either empty
or a lattice coset of index $p^{2(d+r)}$. 
Then $A_{i,(x_0,y_0,a)}$ can be written as a disjoint
union $A_{i,(x_0,y_0,a)} = \bigcup_{j\in J} A'_{i,j,(x_0,y_0,a)}$,
where $A'_{i,j,(x_0,y_0,a)}$ is the set of all $(x,y)\in \mathbb{Z}^2$
such that
\[x\equiv x_0
\mo p^{d+r},\, y\equiv y_0 \mo p^{d+r},\, v_p(\vec{q}_j\cdot (x,y)) = i,\,
p^{-i} (\vec{q_j}\cdot (x,y))\in U_a.\]
The union is disjoint because $v_p((x,y)\cdot \vec{q}_j)=i$ cannot
hold for two different $j$ when $i>d$. Since $1 + p^r \mathbb{Z}_p
\subset U$, we can write $A'_{i,j,(x_0,y_0,a)}$ as a disjoint union
of at most $p^r$ sets of the form
\[\begin{aligned}
L_{i,j,(x_0,y_0,b)} &= \{(x,y)\in \mathbb{Z}^2 : x\equiv x_0
\mo p^{d+r}, y\equiv y_0 \mo p^{d+r}\} \\ &\cap
\{(x,y)\in \mathbb{Z}^2 : \vec{q_j}\cdot (x,y) \equiv b \mo p^{i+r}\} .
\end{aligned}\]
Since this is the intersection of a lattice coset of index $p^{2 d + 2 r}$
and a lattice coset of index $p^{i+r}$, 
$L_{i,j,(x_0,y_0,b)}$ must be a lattice coset of
index $n_i$ satisfying $p^{i+r} | n_i | p^{i + 2 d + 3 r}$.
Hence (3) is satisfied for any $i\geq 0$.

It remains to prove (5). For $i>d$, this is immediate from the
definition of $A_{i,(x_0,y_0,a)}$. Let $i\leq d$. Any two elements 
$(x_1,y_1)$, $(x_2,y_2)$ of
$A_{i,(x_0,y_0,a)}$ must satisfy
$x_1\equiv x_2 \mo p^{d+r}$, $y_1\equiv y_2 \mo p^{d+r}$. Hence
$\vec{q}_j\cdot (x_1,y_1) \cong \vec{q}_j\cdot (x_2,y_2) \mo p^{d+r}$
for every $j$. Since $\max_j v_p(\vec{q}_j\cdot (x_1,y_1))
= \max_j v_p(\vec{q}_j\cdot (x_2,y_2)) = i \leq d$,
we can conclude that $\vec{q}_j\cdot (x_1,y_1)$ and
$\vec{q}_j\cdot (x_1,y_1)$ lie in the same coset of $1 + p^r \mathbb{Z}_p$. 
Hence $\vec{q}_j\cdot (x_1,y_1)$ and
$\vec{q}_j\cdot (x_1,y_1)$ lie in the same coset of $U$.
\end{proof}
\begin{lem}\label{lem:nestle}
Let $L\subset \mathbb{Z}^2$ be a lattice. Let $L', L''\subset L$ be lattice cosets
contained in $L$. Then the intersection $L'\cap L''$ is either the empty set or
a lattice coset of index $\lbrack \mathbb{Z}^2 : L'\cap L''\rbrack$ dividing
$\frac{\lbrack \mathbb{Z}^2 : L'\rbrack \cdot \lbrack \mathbb{Z}^2 : L''\rbrack}{
 \lbrack \mathbb{Z}^2 : L\rbrack}$.
\end{lem}
\begin{proof}
Since $L$ and $\mathbb{Z}^2$ are isomorphic, it is enough to prove the statement for 
$L=\mathbb{Z}^2$. It holds in general that, given two subgroup cosets $L'$,
$L''$ of an abelian
group $Z$, the intersection $L'\cap L''$ is either the empty set or a subgroup coset
of index dividing $\lbrack Z: L'\rbrack \cdot \lbrack Z : L\rbrack$.
\end{proof}
\begin{lem}\label{lem:morepancakes}
Let $p$ be a prime. Let $U$ be an open subgroup of $\mathbb{Q}_p^*$. Let 
$\{\vec{q}_j\}_{j\in J}$ be a finite subset of $\mathbb{Q}_p^2$. Let $L$ be a lattice
of index $\lbrack \mathbb{Z}^2 : L\rbrack = p^m$. Then there is a partition 
\[L\cap \mathbb{A}_p = B_{\infty} \cup \bigcup_{i\geq m} \bigcup_{k\in K} B_{i,k}\]
such that
\begin{enumerate}
\item $K$ is a finite set,

\item $B_{\infty}$ is the union of finitely many sets of the form $A_{x,y} = \{(n x,
n y) : n\in \mathbb{Z}, p\nmid n\}$,

\item $B_{i,k}$ is a disjoint union of at most $c_1$ lattice cosets of index 
$p^{i+c_2}$, 

\item for every $i_0\geq 0$, the set $B_{\infty}\cup \bigcup_{i\geq i_0}
\bigcup_{k\in K} B_{i,k}$ is a disjoint union of at most $c_1$ sets of the form 
$R\cap \mathbb{A}_p$, where $R$ is a lattice of index $p^{i_0}$,

\item for any choice of $i\geq 0$, $j\in J$, $k\in K$ and all 
$(x_1,y_1), (x_2,y_2)\in A_{i,k}$,
the inner products
$\vec{q}_j \cdot (x_1,y_1)$ and $\vec{q}_j \cdot (x_2,y_2)$ lie in the
same coset of $U$.
\end{enumerate}
\end{lem}
\begin{proof}
Let $A_{\infty}$, $A_{i,k}$
be as in Lemma \ref{lem:pancakestack}. By Lemma \ref{lem:pancakestack}, (4),
\[A_{\infty} \cup \bigcup_{i\geq i_0} \bigcup_{k\in K} A_{i,k}\] is a disjoint union
of at most $c_1$ lattices of index $p^{i_0}$. Hence, for $i_0\leq m$, it follows
from Lemma \ref{lem:oques} that either
\[(L\cap \mathbb{A}_p)\cap (A_{\infty} \cup 
 \bigcup_{i\geq i_0} \bigcup_{k\in K} A_{i,k}) = \emptyset\]
or 
\[L\cap \mathbb{A}_p\subset (A_{\infty} \cup
 \bigcup_{i\geq i_0} \bigcup_{k\in K} A_{i,k}) \]
must hold.
Suppose $(L\cap \mathbb{A}_p)\cap (A_{\infty} \cup 
 \bigcup_{i\geq m} \bigcup_{k\in K} A_{i,k}) = \emptyset$. Let $i_0\geq 0$ be
the largest integer such that
\[(L\cap \mathbb{A}_p)\subset (A_{\infty} \cup 
 \bigcup_{i\geq i_0} \bigcup_{k\in K} A_{i,k})  .\]
Then
\[L\cap \mathbb{A}_p = \bigcup_{k\in K} (A_{i_0,k} \cap L) .\]
Set $B_{m,k} = A_{i_0,k}\cap L$, $B_{i,k} = \emptyset $ for $i\ne m$,
$B_{\infty} = A_{\infty} \cap L$. Conditions (1), (2), (4) and (5) follow trivially from
the definitions of $A_{i_0,k}$ and $m$. By Lemma \ref{lem:pancakestack}, (3), $A_{i_0,k}$
is the disjoint union of at most $c_1$ lattice cosets of index $p^{i_0+c_2}$. Take
one such lattice coset and call it $R_0$. By Lemma \ref{lem:pancakestack}, (4),
$R_0$ is contained in a set of the form $R\cap \mathbb{A}_p$, where $R$ is a lattice
of index $p^{i_0}$. Since $p^{i_0} | p^m$, $L$ is contained in a lattice $R'$ of index
$p^{i_0}$. By Lemma \ref{lem:oques}, either $R\cap R'\cap \mathbb{A}_p = \emptyset$ or
$R = R'$. In the former case, $R_0\cap (L \cap \mathbb{A}_p) = \emptyset$. In the latter
case, Lemma \ref{lem:nestle} yields that $R_0 \cap L$ is a lattice coset of index dividing
$p^{(i_0+c_2)+m-i_0} = p^{m+c_2}$ and divided by $\lbrack \mathbb{Z}^2 : R\cap L \rbrack
= \lbrack \mathbb{Z}^2 : L\rbrack = p^m$. Condition (4) follows.

Now suppose 
\[
L\cap \mathbb{A}_p\subset (A_{\infty} \cup
 \bigcup_{i\geq m} \bigcup_{k\in K} A_{i,k}).\]
By Lemma \ref{lem:pancakestack}, 
$A_{\infty} \cup \bigcup_{i\geq i_0} \bigcup_{k\in K} A_{i,k}$
is a disjoint union of sets of the form $R\cap \mathbb{A}_p$, $R$ a lattice of
index $p^m$. By Lemma \ref{lem:oques}, one such $R$ is equal to $L$.
For $i\geq m$, set $B_{i,k} = A_{i,k}$ if $A_{i,k}\subset L$, $B_{i,k}=0$ otherwise.
Set $B_{\infty} = A_{\infty} \cap L$. Conditions (1) to (5) follow easily.
\end{proof}
\begin{prop}\label{prop:pliblock}
Let $f,g: \mathbb{Z}^2 \to \mathbb{C}$ be given with $\max |f(x,y)|, |g(x,y)| \leq 1$. 
Let $f$ be pliable with respect to $\{(v_j,U_j,\vec{q}_j)\}
$. Assume that there are $\eta_N\leq N$, $\epsilon_N\geq 0$ such that for any sector $S$
and any lattice coset $L$ of index $\lbrack \mathbb{Z}^2 : L\rbrack \leq \eta_N$, 
\begin{equation}\label{eq:pligiv}
\mathop{\sum_{(x,y) \in S\cap \lbrack -N,N\rbrack^2 \cap L}}_{\gcd(x,y)=1} 
g(x,y) \ll \frac{\epsilon_N N^2}{\lbrack \mathbb{Z}^2 : L\rbrack}.
\end{equation}
Then, for any sector $S$ and any lattice $L$, 
\begin{equation*}
\mathop{\sum_{(x,y) \in S\cap \lbrack -N,N\rbrack ^2 \cap L}}_{
\gcd(x,y)=1} f(x,y) g(x,y) \ll 
\left(\frac{\epsilon_N}{\lbrack \mathbb{Z}^2 : L\rbrack} + 
\frac{(\log \eta_N)^c}{\eta_N}\right) N^2,
\end{equation*}
where $c$ is the number of distinct finite places among $\{v_j\}$ and the
implied constant depends only on 
the implied constant in (\ref{eq:pligiv}) and on
$\{(v_j, U_j, \vec{q}_j)\}$.
\end{prop}
\begin{proof}
By Lemma \ref{lem:rsub} we can partition $\mathbb{R}^2$ into 
\begin{equation*}
\mathbb{R}^2 = T_1\cup \dotsb \cup T_k \cup S_1 \cup \dotsb \cup S_l
\end{equation*}
such that $\vec{q}_j \cdot (x_1,y_1)$ and $\vec{q}_j \cdot (x_2,y_2)$ lie in
the same coset of $\bigcup_{j} U_j$ for all $(x_1,y_1)$, $(x_2,y_2)$ in $S_i$ and
all $j$ with $v_j = \infty$. The contribution of $T_1, T_2, \dotsc, T_k$
to the final sum is $O(1)$. As there is a finite number of $S_i$'s, it is enough
to prove the desired bound for every $S_i$ separately. Fix $i$ and 
let $S'=S_i\cap S$.

Let $\{p_l\}$ be the set of all finite places among $\{v_j\}$. Let 
$\{\vec{q}_{l,j}\}$ be the set
of all $\vec{q}_j$ such that $v_j$ is induced by $p_l$. Let 
$m = \lbrack \mathbb{Z}^2 : L\rbrack$. We can write
\[L = \bigcap_l L_{p_l} \cap L_{m_0},\]
where $L_{p_l}$ is a lattice of index $p_l^{v_{p_l}(m)}$ and $L_{m_0}$ is 
a lattice of index $m_0 = \frac{m}{\prod_l p_l^{v_{p_l}(m)}}$.

For every $p_l$, 
Lemma \ref{lem:pancakestack} yields a partition 
\begin{equation*}
L_{p_l} \cap \mathbb{A}_{p_l}  = B_{\infty} \cup \bigcup_{i\geq v_{p_l}(m)}
\bigcup_{k\in K_l} B_{l,i,k}
\end{equation*}
such that $\vec{q}_{l,j} \cdot (x_1,y_1)$ 
and $\vec{q}_{l,j} \cdot (x_2,y_2)$ lie in
the same coset of $U$ for any $(x_1,y_1)$, $(x_2,y_2)$ in 
$B_{l,i,k}$ and any $i$, $j$, $k$.

Let $\mathbb{A}=\{x,y\in \mathbb{Z}^2 : \gcd(x,y)=1\}$. Clearly
\begin{equation}\label{eq:quotab}\begin{aligned}
L\cap \mathbb{A} &= \bigcup_l \left(B_{l,\infty}\cup \bigcup_{i\geq v_{p_l}(m)}
  \bigcup_{k\in K_l} B_{l,i,k}\right) \cap \mathbb{A}\\
 &= \left(\bigcup_l B_{l,\infty} \cap \mathbb{A}\right) \cup
\mathop{\bigcup_{n\geq 1}}_{\rad(n)|R} \bigcup_{\{k_l\}\in \prod_l K_l}
 \left(\bigcap_l B_{l,v_p(m n),k_l} \cap L_{m_0}\right) \cap \mathbb{A} .
\end{aligned}\end{equation}
Note that $(\bigcup_l B_{l,\infty} \cap \mathbb{A})$ is a finite set. Since $f$ is affinely
pliable with respect to $\{v_j, U_j, \vec{q}_j\}$, it is constant on
$S'\cap \bigcap_l B_{l,v_p(m n),k_l}$ for any $n\geq 1$ and any 
 $\{k_l\} \in \prod_l K_l$. Denote the value of $f$ on $\bigcap_l B_{l,v_p(m n),k_l}$
by $f_{n,\{k_l\}}$. Thanks to (\ref{eq:quotab}), we can write
\[\begin{aligned}
\mathop{\sum_{(x,y)\in S' \cap \lbrack -N,N\rbrack^2 \cap L}}_{\gcd(x,y)=1}
 f(x,y) g(x,y) &=
 \sum_{(x,y)\in \cap_l B_{l,\infty}\cap \mathbb{A}} f(x,y) g(x,y) \\ &+
\mathop{\sum_{n\geq 1}}_{\rad(n)|R} 
\sum_{\{k_l\}\in \prod_l K_l}
 \mathop{\sum_{(x,y)\in B_{l,v_p(m n), k_l}\cap L_{m_0}}}_{(x,y)\in \mathbb{A}}
 f(x,y) g(x,y) \\
&= \sum_{\{k_l\}\in \prod_l K_l} \mathop{\sum_{n\geq 1}}_{\rad(n)|R} 
 f_{n,\{k_l\}} 
\mathop{\sum_{(x,y)\in B_{l,v_p(m n), k_l}\cap L_{m_0}}}_{(x,y)\in \mathbb{A}}
g(x,y) \\ &+O(1).\end{aligned}\]
Fix $\{k_l\}\in \prod_l K_l$. Set
\[a_n = f_{n,\{k_l\}} 
 \mathop{\sum_{(x,y)\in B_{l,v_p(m n), k_l} \cap L_{m_0}}}_{(x,y)\in \mathbb{A}}
 g(x,y)\]
if $\rad(n)|R$, $a_n=0$ otherwise. 
Then
\[
\mathop{\sum_{(x,y) \in S\cap \lbrack -N,N\rbrack ^2 \cap L}}_{
\gcd(x,y)=1} f(x,y) g(x,y) =
\sum_n a_n .\]
Let $s_d = \sum_n |a_{d n}|$. From Lemma \ref{lem:morepancakes}, (4),
Lemma \ref{lem:ramsay}
and $|g(x,y)|\leq 1$, we get that $s_d\ll \frac{N^2}{m n}$.
Set $C=N^2/m$. By Lemma \ref{lem:morepancakes}, (3),
$\bigcap_l B_{l,v_p(m n),k_l} \cap L_{m_0}$ 
is the union of at most $c_3 = c_1^{\#\{p_l\}}$ lattice cosets of modulus
$c_4 m n$, where $c_4 = \prod_l p_l^{c_2}$. 
Set
$M=\min\left(\frac{\eta_N}{c_4 m}, \frac{N}{p_0 m}\right)$, where 
$p_0 = \max_l p_l$.
We can now
apply Lemma \ref{lem:sivcor}, obtaining
\[\begin{aligned}
\sum_n a_n = \sum_{n\leq M} a_n + O\left(\frac{N}{\eta_N} (\log \eta_N)^{\omega(R)}\right)
.\end{aligned}\]
By (\ref{eq:pligiv}),
\[\begin{aligned}
 \sum_{n\leq M} a_n &\ll \mathop{\sum_{n\leq M}}_{\rad(n)|R} \frac{\epsilon_N N}{m n}
\leq \sum_{\rad(n)|R} \frac{\epsilon_N N}{m n} \\
&= \frac{\epsilon N}{m} 
\cdot \prod_{p|R} \left(1 + \frac{1}{p} + \frac{1}{p^2} +\dotsb\right)
\ll \frac{\epsilon N}{m} = \frac{\epsilon N}{
\lbrack \mathbb{Z}^2 : L\rbrack} .\end{aligned}\]
We conclude that
\[
\mathop{\sum_{(x,y) \in S'\cap \lbrack -N,N\rbrack ^2 \cap L}}_{
\gcd(x,y)=1} f(x,y) g(x,y) \ll 
\left(\frac{\epsilon_N}{\lbrack \mathbb{Z}^2 :
L\rbrack} + 
\frac{(\log \eta_N)^c}{\eta_N}\right) N^2 .\]
As said in the beginning of the proof, it follows immediately that
\[
\mathop{\sum_{(x,y) \in S\cap \lbrack -N,N\rbrack ^2 \cap L}}_{
\gcd(x,y)=1} f(x,y) g(x,y) \ll 
\left(\frac{\epsilon_N}{\lbrack \mathbb{Z}^2 : L\rbrack} + 
\frac{(\log \eta_N)^c}{\eta_N}\right) N^2 .\]
\end{proof}
\section{Using the square-free sieve}
We will now state the results we need
from Chapter \ref{chap:square}, as well as some simple consequences.
\subsection{Conditional results}\label{subs:conres}
We introduce the following quantitative versions of Conjectures
$\mathfrak{A}_1$ and $\mathfrak{A}_2$. 
\begin{conja1q}
The polynomial $P\in \mathfrak{O}_K\lbrack x\rbrack$ obeys
\[\#\{1\leq x\leq N:
\exists \mathfrak{p} \text{\,\st\,} \rho(\mathfrak{p})>N^{1/2}, \mathfrak{p}^2|P(x)\} \ll \delta(N),\]
where $1\ll \delta(N)\ll N$ and $\rho(\mathfrak{p})$ is the rational prime
lying under $\mathfrak{p}$.
\end{conja1q}
\begin{conja2q}
The homogeneous polynomial $P\in \mathfrak{O}_K\lbrack x,y\rbrack$ obeys
\[\#\{-N\leq x,y\leq N:
\exists \mathfrak{p} \text{\,\st\,} \rho(\mathfrak{p})>N, \mathfrak{p}^2|P(x)\} \ll \delta(N),\]
where $1\ll \delta(N)\ll N$ and $\rho(\mathfrak{p})$ is the rational prime
lying under $\mathfrak{p}$.
\end{conja2q}

We can now restate Propositions \ref{prop:hjuhju} and \ref{prop:hjohjo}
as conditional results.
\begin{prop}[$\mathfrak{A}_1(K,P,\delta(N))$]\label{prop:sqfra}
Let $K$ be a number field. Let $f:I_K\times \mathbb{Z}\to \mathbb{C}$,
$g:\mathbb{Z}\to \mathbb{C}$ be given with $\max |f(\mathfrak{a},x)|\leq 1$,
$\max |g(x)|\leq 1$. Assume that $f(\mathfrak{a},x)$ depends only on
$\mathfrak{a}$ and on $x \mo \mathfrak{a}$. Let 
$P\in \mathfrak{O}_K\lbrack x \rbrack$. 
Suppose there are $\epsilon_{1,N}, \epsilon_{2,N}\geq 0$
such that for any integer $a$ and any positive integer $m$,
\begin{equation}\label{eq:whynot2}
\mathop{\sum_{1\leq x\leq N}}_{x\equiv a \mo m} g(x) \ll 
\left(\frac{\epsilon_{1,N}}{m} + \epsilon_{2,N}\right) N .
\end{equation}
Then, for any integer $a$ and any positive integer $m$,
\[\begin{aligned}
\mathop{\sum_{1\leq x\leq N}}_{x\equiv a \mo m} f(\sq_K(P(x)),x) g(x) &\ll 
\left(\frac{\epsilon_{1,N}}{m}
+ (\log N)^{c_1} \sqrt{\max(\epsilon_{2,N},m/N^{1/2})}
\right) \\ &\cdot \tau_{c_2}(m) N + \delta(N),\end{aligned}\]
where 
$c_1$ and $c_2$ depend only on $P$ and $K$, and
the implied constant depends only on $P$, $K$ and
the implied constant in (\ref{eq:whynot2}).
\end{prop}
\begin{prop}[$\mathfrak{A}_2(K,P,\delta(N))$]\label{prop:volver}
Let $K$ be a number field.
Let $f:I_K\times \{(x,y)\in \mathbb{Z}^2 : \gcd(x,y)=1\} \to \mathbb{C}$, $g:\{(x,y)\in \mathbb{Z}^2 : \gcd(x,y)=1\} \to \mathbb{C}$
be given with $\max |f(\mathfrak{a},x,y)| \leq 1$, $\max |g(x,y)|\leq 1$. Assume
that $f(\mathfrak{a},x,y)$ depends only on $\mathfrak{a}$ and on 
$\{\frac{x \mo \mathfrak{p}}{y \mo \mathfrak{p}}\}_{\mathfrak{p}|\mathfrak{a}}
\in \prod_{\mathfrak{p}|\mathfrak{a}} 
\mathbb{P}^1(\mathfrak{O}_K/\mathfrak{p})$. 
Let $P\in \mathfrak{O}_K\lbrack x,y\rbrack$ be a homogeneous polynomial.
Let $S$ be a convex set.
Suppose there are $\epsilon_{1,N}, \epsilon_{2,N}\geq 0$ such that for
any lattice coset $L\subset \mathbb{Z}^2$,
\begin{equation}\label{eq:pcond2}
\mathop{\sum_{(x,y)\in S\cap \lbrack -N,N\rbrack^2\cap L}}_{\gcd(x,y)=1}
g(x,y) \ll \left(\frac{\epsilon_{1,N}}{\phi(\lbrack \mathbb{Z}^2 :L\rbrack)} +
\epsilon_{2,N}\right) N^2.
\end{equation}
Then, for any lattice coset $L\subset \mathbb{Z}^2$,
\[\begin{aligned}
\mathop{\sum_{(x,y)\in S\cap \lbrack -N,N\rbrack^2 \cap L}}_{\gcd(x,y)=1}
&f(\sq_K(P(x,y)),x,y) g(x,y) \\ &\ll
\left(\frac{\epsilon_{1,N}}{\lbrack \mathbb{Z}^2 : L\rbrack} +
(\log N)^{c_1} \sqrt{\max(\epsilon_{2,N},\lbrack \mathbb{Z}^2 : L\rbrack/N)}\right) \tau_{c_2}(m) N + \delta(N),\end{aligned}\]
where $c_1$ and $c_2$ depend only on $P$ and $K$, and
the implied constant depends only on $P$, $K$ and 
the implied constant in (\ref{eq:pcond2}).
\end{prop}
See Appendices \ref{sec:fapp} and \ref{sec:app} for all proven instances of 
$\mathfrak{A}_i(K,P,\delta(N))$.
\subsection{Miscellanea}
We will need the following simple lemmas.
\begin{lem}\label{lem:glugluglu}
For any positive integer $n$,
\[\prod_{p|n} \left(1 + \frac{1}{p}\right) \ll
\log \log n,\]
where the implied constant is absolute.
\end{lem}
\begin{proof}
Obviously
\[\log \prod_{p|n} \left(1 + \frac{1}{p}\right) \leq
\sum_{p|n} \frac{1}{p} .\]
Define \[S(m,r) = \max_{n\leq r} \mathop{\sum_{p|n}}_{p>m} \frac{1}{p} .\]
Then, for any $r$,
\[S(m,r) \leq \frac{1}{p} + S(p, r/p)\]
for some $p>m$. Clearly
\[\begin{aligned}
S(m_1,n) \geq S(m_2,n) &\text{\; if $m_1\leq m_2$,}\\
S(m,n_1) \geq S(m,n_2) &\text{\; if $n_1\geq n_2$.}\end{aligned}\]
Hence
\[\begin{aligned}
S(1,n) &\leq \frac{1}{2} + S(2,n/2) \leq \frac{1}{2} + \frac{1}{3} + S(3,
n/2\cdot 3) \\
&\leq \frac{1}{2} + \frac{1}{3} + \dotsb + \frac{1}{p} + S\left(m,
\frac{n}{\prod_{p\leq m} p}\right) .\end{aligned}\]
Now
\[\prod_{p\leq m} p = \left(\frac{m}{2}\right)^{O((m/2)/(\log m/2))} =
e^{O(m)} .\]
Thus, the least $m$ such that $\prod_{p\leq m} p > n/2$ is at most
$O(\log n)$. Therefore
\[S(1,n) \leq \sum_{p\leq m} \frac{1}{p} \leq \log \log \log n + o(1) .\]
The statement follows.
\end{proof}
\begin{lem}\label{lem:minmat}
Let $g:\mathbb{Z}^2 \to \mathbb{C}$ be given with $|g(x,y)|
\leq 1$ for all $x,y\in \mathbb{Z}$. Let $\eta(N)\leq N$.
Suppose that, for every sector $S$ and every lattice
$L$ of index $\lbrack \mathbb{Z}^2 : L
\rbrack \leq \eta(N)$,
\begin{equation}\label{eq:whenwi}
\sum_{(x,y)\in S \cap \lbrack -N,N\rbrack^2\cap L} g(x,y) \ll
\frac{\epsilon(N) N^2}{\lbrack \mathbb{Z}^2 :L \rbrack} .\end{equation}
Then, for every sector $S$ and every lattice $L$ of index
$\lbrack \mathbb{Z}^2 : L
\rbrack \leq \eta(N)$,
\[\mathop{\sum_{(x,y)\in S \cap \lbrack -N,N\rbrack^2\cap L}}_{\gcd(x,y)=1}
 g(x,y) \ll
\max\left(\epsilon(N) \log \log N, 
\frac{1}{(\eta(N))^{1/2 - \epsilon}}\right)
\frac{N^2}{\lbrack \mathbb{Z}^2 :L \rbrack} .\]
\end{lem}
\begin{proof}
For every positive integer $a$, let
\[\begin{aligned}
S_a &= \{0\},\\
\gamma(a) &= \lbrack \mathbb{Z}^2 : L\cap a \mathbb{Z}^2\rbrack,\\
f_a(0) &= \begin{cases} 1 &\text{if $a=1$,}\\0 &\text{otherwise,}
\end{cases}\\
g_a(0) &= \mathop{\sum_{(x,y)\in S\cap \lbrack -N,N\rbrack^2 \cap L}}_{
 \gcd(x,y)=a} \lambda_K(P(x,y)).\end{aligned}\]
Clearly
\[\mathop{\sum_{(x,y)\in S\cap \lbrack - N,N\rbrack^2 \cap L}}_{
\gcd(x,y)=1} g(x,y) .\]
By Lemma \ref{lem:wicked},
\[\begin{aligned}
\sum_{a=1}^{\infty} f_a(0) g_a(0) &=
 \sum_{\gamma(d)\leq \eta(N)} \left(\sum_{d'|d} \mu(d')
\lbrack d/d'=1\rbrack \right) \mathop{\sum_a}_{d|a} g_a(0)\\
&+ 2 \sum_{\eta(N)<\gamma(d)\leq \eta(N)^2} \tau_3(a)
 \mathop{\sum_{a}}_{d|a} |g_a(0)|
 + 2 \mathop{\sum_{\text{$p$ prime}}}_{\gamma(p)>\eta(N)} 
 \mathop{\sum_{a}}_{d|a} |g_a(0)| \\
&= \sum_{\gamma(d)\leq M} \mu(d) 
 \mathop{\sum_{(x,y)\in S\cap \lbrack - N, N\rbrack^2 \cap L}}_{
a|x,\, a|y} g(x,y) \\
&+ 2 \sum_{\eta(N)<\gamma(d)\leq \eta(N)^2} \tau_3(d)
 \mathop{\sum_a}_{d|a} \left|
 \mathop{\sum_{(x,y)\in S\cap \lbrack - N, N\rbrack^2 \cap L}}_{
\gcd(x,y)=a} g(x,y) \right| \\ &+
2 \mathop{\sum_{\text{$p$ prime}}}_{\gamma(p)>\eta(N)} \tau_3(d)
 \mathop{\sum_a}_{d|a} \left|
 \mathop{\sum_{(x,y)\in S\cap \lbrack - N, N\rbrack^2 \cap L}}_{
\gcd(x,y)=a} g(x,y) \right| .
\end{aligned}\]
Then, by (\ref{eq:whenwi}),
\[\sum_{a=1}^{\infty} f_a(0) g_a(0) =
\sum_{\gamma(d)\leq \eta(N)} \frac{\epsilon(N) N^2}{\gamma(d)} +
2 \sum_{\eta(N)<\gamma(d)\leq \eta(N)^2} \tau_3(d) 
 \frac{N^2}{\gamma(d)} +
2 \mathop{\sum_{\text{$p$ prime}}}_{\gamma(p)>\eta(N)} 
 \frac{N^2}{\gamma(p)} .\]
We can assume that $L$ is not contained in any set of the form
$a \mathbb{Z}^2$, $a>1$, as otherwise the statement is trivial. Thus
$\gamma(d)=d \cdot \lcm(d,\lbrack \mathbb{Z}^2 : L \rbrack)$. Hence
\[\sum_{d=1}^{\infty} \frac{1}{\gamma(d)} \leq
 \sum_{d'|\lbrack \mathbb{Z}^2 : L\rbrack} 
 \frac{1}{d' \lbrack \mathbb{Z}^2 : L\rbrack}
 \sum_d \frac{1}{d^2} \ll \sum_{d'|\lbrack \mathbb{Z}^2 : L\rbrack}
 \frac{1}{d' \lbrack \mathbb{Z}^2 : L\rbrack},\]
\[\begin{aligned}\sum_{\gamma(d)>\eta(N)} \frac{\tau_3(d)}{\gamma(d)} &=
 \sum_{d'|\lbrack \mathbb{Z}^2 : L\rbrack} 
  \frac{\tau_3(d')}{d' \lbrack \mathbb{Z}^2 : L\rbrack}
\sum_{d>(\eta(N)/d')^{1/2}} \frac{\tau_3(d)}{d^2} \ll
\sum_{d'|\lbrack \mathbb{Z}^2 : L\rbrack}
\frac{\tau_3(d')}{\lbrack \mathbb{Z}^2 : L\rbrack \sqrt{d' \eta(N)}}
.\end{aligned}\]
By \ref{lem:glugluglu},
\[\sum_{d'|\lbrack \mathbb{Z}^2 : L\rbrack} \frac{1}{d'} \ll \log \log N .\]
Clearly
\[\sum_{d'|\lbrack \mathbb{Z}^2 : L\rbrack} \frac{\tau_3(d')}{\sqrt{d'}}
\ll \tau_4(\lbrack \mathbb{Z}^2 : L\rbrack) \ll \lbrack \mathbb{Z}^2 : L
\rbrack^{\epsilon} .\]
The statement follows.
\end{proof}
\begin{lem}\label{lem:repla}
Let $K$ be a number field. Let $F\in \mathfrak{O}_K\lbrack x\rbrack$
be a square-free polynomial. Let $a$ be an integer, $m$ a positive integer.
If $\mathfrak{A}_1(K,F,\delta(N))$ holds, then
$\mathfrak{A}_1(K,F(m x + a),\delta(m N))$ holds.
\end{lem}
\begin{proof}
Immediate from the statement of Conjecture $\mathfrak{A}_1$.
\end{proof}
\begin{lem}\label{lem:replb}
Let $K$ be a number field. Let $F\in \mathfrak{O}_K\lbrack x,y\rbrack$
be a square-free homogeneous polynomial. Let $A\in SL_2(\mathbb{Z})$,
$m_A = \max(|a_{1 1}| + |a_{1 2}|,|a_{2 1}|+|a_{2 2}|)$.
If $\mathfrak{A}_2(K,F,\delta(N))$ holds, then
$\mathfrak{A}_2(K,F(a_{1 1} x + a_{1 2} y,a_{2 1} x + a_{2 2} y),\delta(m_A N))$ 
holds.
\end{lem}
\begin{proof}
Immediate from the statement of Conjecture $\mathfrak{A}_2$.
\end{proof}
\begin{lem}\label{lem:multconj1}
Let $K$ be a number field. Let $F, G\in \mathfrak{O}_K\lbrack x\rbrack$
be square-free polynomials without common factors. Then
$\mathfrak{A}_1(K,F \cdot G,\delta(N))$ holds if and only if
$\mathfrak{A}_1(K,F,\delta(N))$ and $\mathfrak{A}_2(K,G,\delta(N))$ both
hold.
\end{lem}
\begin{proof}
We can assume $N^{1/2}$ to be larger than 
$\max_{\mathfrak{p}|\Disc(F,G)} \rho(\mathfrak{p})$. 
Then, for any $\mathfrak{p}$ such that $\rho(\mathfrak{p})>N$,
we have that $\mathfrak{p}$ cannot divide both $F(x)$ and $G(x)$. Hence 
\[\{1\leq x\leq N : \exists \mathfrak{p} \text{\;\st\;} \rho(\mathfrak{p})>N^{1/2},
 \mathfrak{p}^2|F(x)\}\]
equals
\[\begin{aligned}
\{1\leq x\leq N : \exists \mathfrak{p} \text{\;\st\;} \rho(\mathfrak{p})>N^{1/2}, \mathfrak{p}^2|G(x)\} \;\;\cup\\
 \{1\leq x\leq N : \exists \mathfrak{p} \text{\;\st\;} \rho(\mathfrak{p})>N^{1/2}, \mathfrak{p}^2|F(x) \cdot G(x)\}.
\end{aligned}\]
\end{proof}
\begin{lem}\label{lem:multconj2}
Let $K$ be a number field. Let $F, G\in \mathfrak{O}_K\lbrack x,y\rbrack$
be square-free homogenous polynomials without common factors. Then
$\mathfrak{A}_2(K,F \cdot G,\delta(N))$ holds if and only if
$\mathfrak{A}_2(K,F,\delta(N))$ and $\mathfrak{A}_2(K,G,\delta(N))$ both
hold.
\end{lem}
\begin{proof}
Same as that of Lemma \ref{lem:multconj1}.
\end{proof}
\begin{lem}\label{lem:multsqua}
Let $K$ be a number field. Let $F, G, H\in \mathfrak{O}_K\lbrack x\rbrack$
be square-free polynomials. Assume that $F$, $G$ and $H$
are coprime as elements of $K\lbrack x\rbrack$. Then
there is an ideal $\mathfrak{m}$ such that, for any $\mathfrak{M}\in I_K$,
$\mathfrak{m} | \mathfrak{M}$, we can tell
\[\begin{aligned}
\sq_K(F(x) H(x))&/\gcd(\sq_K(F(x) H(x)),\mathfrak{M}^{\infty})\;\;\;
\text{and}\\
\sq_K(G(x) H(x))&/\gcd(\sq_K(G(x) H(x)),\mathfrak{M}^{\infty})
\end{aligned}\]
from
\[
\sq_K(F(x) G(x) H(x))/\gcd(\sq_K(F(x) G(x) H(x)),\mathfrak{M}^{\infty})\]
and 
$x \mo \mathfrak{p}$ for $\mathfrak{p}|\sq_K(F(x) G(x) H(x))$,
$\mathfrak{p}\nmid \mathfrak{M}$.
\end{lem}
\begin{proof}
Let $\mathfrak{m} = \Disc(F,G) \cdot \Disc(F,H) \cdot \Disc(G,H)$.
Take a prime ideal $\mathfrak{p}\nmid \mathfrak{M}$.
Suppose 
\[\mathfrak{p}|(\sq_K(F(x) G(x) H(x))/\gcd(\sq_K(F(x) G(x) H(x)),\mathfrak{M}^{\infty})).\]
We can tell which one of $\sq_K(F(x))$, $\sq_K(G(x))$
or $\sq_K(H(x))$ is divided by $\mathfrak{p}$ if we know which one of
$F(x)$, $G(x)$, $H(x)$ is divided by $\mathfrak{p}$. The latter
question can be answered given 
$x\mo \mathfrak{p}$. 
\end{proof}
Given two square-free polynomials 
$A,B\in \mathfrak{O}_K\lbrack x\rbrack$, we can always find square-free
polynomials $F, G, H\in \mathfrak{O}_K\lbrack x\rbrack$ such that
\begin{itemize}\item $F$, $G$ and $H$ are pairwise coprime as elements
of $K\lbrack x\rbrack$,
\item $A = F H$, $B = G H$.
\end{itemize}
Write $\Lcm(A,B)$ for $F\cdot G\cdot H$. Notice that $\Lcm(A,B)$ is defined
only up to multiplication by a unit of $\mathfrak{O}_K$.
\begin{cor}\label{cor:multsquo}
Let $K$ be a number field. Let $A, B\in \mathfrak{O}_K\lbrack x\rbrack$
be square-free polynomials. Then there is an ideal 
$\mathfrak{m}_{A,B}$ such that, for any $\mathfrak{M}\in I_K$,
$\mathfrak{m}_{A,B} | \mathfrak{M}$, we can tell
\[\begin{aligned}
\sq_K(A(x))&/\gcd(\sq_K(A(x)),\mathfrak{m}^{\infty})\;\;\;
\text{and}\\
\sq_K(B(x))&/\gcd(\sq_K(B(x)),\mathfrak{m}^{\infty})
\end{aligned}\]
from
\[
\sq_K(\Lcm(A,B)(x))/\gcd(\sq_K(\Lcm(A,B)(x)),\mathfrak{M}^{\infty})\]
and 
$x \mo \mathfrak{p}$ for $\mathfrak{p}|\sq_K(\Lcm(A,B))$,
$\mathfrak{p}\nmid \mathfrak{M}$.
\end{cor}
\begin{proof}
Immediate from Lemma \ref{lem:multsqua}.
\end{proof}
We can define $\Lcm$ for homogeneous polynomials in two variables
in the same way we defined it for polynomials in one variable.
\begin{lem}\label{lem:multsquid}
Let $K$ be a number field. Let $A, B\in \mathfrak{O}_K\lbrack x,y\rbrack$
be homogeneous square-free polynomials. Then there is an ideal 
 $\mathfrak{m}_{A,B}$ such that, for any $\mathfrak{M}\in I_K$,
$\mathfrak{m}_{A,B} | \mathfrak{M}$, we can tell, for $x$, $y$ coprime,
\[\begin{aligned}
\sq_K(A(x,y))&/\gcd(\sq_K(A(x,y)),\mathfrak{M}^{\infty})\;\;\;
\text{and}\\
\sq_K(B(x,y))&/\gcd(\sq_K(B(x,y)),\mathfrak{M}^{\infty})
\end{aligned}\]
from
\[
\sq_K(\Lcm(A,B)(x,y))/\gcd(\sq_K(\Lcm(A,B)(x,y)),
                                    \mathfrak{M}^{\infty})\]
and 
$\frac{x \mo \mathfrak{p}}{y\mo \mathfrak{p}} \in \mathbb{P}^1(\mathfrak{O}_K/
\mathfrak{p})$ for $\mathfrak{p}|\sq_K(\Lcm(A,B))$,
$\mathfrak{p}\nmid \mathfrak{M}$.
\end{lem}
\begin{proof}
Same as for Lemma \ref{lem:multsqua} and Corollary \ref{cor:multsquo}.
\end{proof}
\section{The global root number and its distribution}
\subsection{Background and definitions}
We may as well start by reviewing the valuative criteria for the
reduction type of an elliptic curve. Let $K_v$ be a Henselian
field
 of characteristic neither $2$ nor $3$. Let $E$ be an elliptic curve over 
$K_v$. Let $c_4, c_6, \Delta\in K_v$
be a set of parameters corresponding to $E$. Then the reduction of $E$
at $v$ is
\begin{itemize}
\item {\em good} if $v(c_4) = 4 k$, $v(c_6) = 6 k$, $v(\Delta) = 12 k$ for some integer $k$;
\item {\em multiplicative} if $v(c_4) = 4 k$, $v(c_6) = 6 k$, 
$v(\Delta) > 12 k$ for some integer $k$;
\item {\em additive} and {\em potentially multiplicative} if 
$v(c_4) = 4 k + 2$, $v(c_6) = 6 k + 3$ and $v(\Delta) > 12 k + 6$ for
some integer $k$;
\item {\em additive} and {\em potentially good} in all remaining cases.
\end{itemize}

From now on, $K$ will be a number field.
Let $\mathcal{E}$ be an elliptic curve over $K(t)$ given by $c_4, c_6\in K(t)$. Let
$q_0\in K(t)$ be a generator of the fractional ideal of $K(t)$ consisting of all
$q\in K(t)$ such that $q^4 c_4$ and $q^6 c_6$ are both in $K\lbrack t\rbrack$. Choose
$q_1\in \mathfrak{O}_K-\{0\}$ such that $(q_1 q_0)^4 c_4$, $(q_1 q_0)^6 c_6$ and
$(q_1 q_0)^{12} \Delta = (q_1 q_0)^{12} \frac{c_4^3 - c_6^2}{1728}$ are all in
$\mathfrak{O}_K\lbrack t\rbrack$. Let 
 $Q(x,y) = q_1 q_0(y/x) x^{\max(\lceil \deg(q_0^4 c_4)/4\rceil,\lceil \deg(q_0^6 c_6)/6\rceil)}$.
Then
\[\begin{aligned}
C_4(x,y) &= Q^4(x,y) c_4(y/x),\\
C_6(x,y) &= Q^6(x,y) c_6(y/x),\\
D(x,y) &= Q^{12}(x,y) \Delta(y/x)\end{aligned}\]
are homogeneous polynomials in $\mathfrak{O}_K\lbrack x,y\rbrack$. 
Note that
$\deg C_6(x,y) = 6 \deg Q$, and thus $\deg C_6$ is even.

We define $P_v$ as in the introduction: for
$v$ a place of $K(t)$, let
$P_{v}\in \mathfrak{O}_K\lbrack t_0,t_1\rbrack$ to be $P_{v} = t_0$ if $v$ is the
place $\deg(\den)-\deg(\num)$, $P_{v} = t_0^{\deg Q} Q\left(\frac{t_1}{t_0}\right)$ if
$v$ is given by a primitive irreducible polynomial $Q_v\in \mathfrak{O}_K\lbrack t\rbrack$.
(We now note that, for any $v$, there are several possible choices for $Q_v$,
all the same up to multiplication by elements of $\mathfrak{O}_K^*$; we choose one $Q_v$ for each $v$ arbitrarily and
fix it once and for all.)
We can write
\begin{equation}\label{eq:defd}
\begin{aligned}
C_4(x,y) &= C_{4,0} \prod_{v} (P_{v}(x,y))^{e_{v,4}},\\
C_6(x,y) &= C_{6,0} \prod_{v} (P_{v}(x,y))^{e_{v,6}},\\ 
D(x,y) &= D_0 \prod_{v} (P_{v}(x,y))^{e_{v,D}},
\end{aligned}
\end{equation}
where $C_{4,0},C_{6,0},D_0 \in \mathfrak{O}_K\lbrack x,y\rbrack$, 
$e_{v,4}, e_{v,6}, e_{v,D}\geq 0$. 
For all but finitely many places $v$ of $K(t)$,
we have $e_{v,4}=0$, $e_{v,6}=0$, $e_{v,D}=0$.

For any place $v$ of $K(t)$, we can localize $\mathcal{E}$ at $v$, thus making it
an elliptic curve over
the Henselian field $(K(t))_{v}$, and then reduce it modulo $v$. We can restate the
the standard valuative criteria for the reduction type
in terms of $e_{v,4}$, $e_{v,6}$, $e_{v,D}$. The reduction of $\mathcal{E}$ at $v$
is
\begin{itemize}
\item good if $e_{v,D}=0$,
\item multiplicative if $e_{v,4}=0$, $e_{v,6}=0$, $e_{v,D}>0$,
\item additive and potentially multiplicative if 
$e_{v,4}=2$, $e_{v,6}=3$, $e_{v,D}>6$,
\item additive and potentially good in all remaining cases.
\end{itemize}

As before, let $\mathbb{A} = \{(x,y) \in \mathfrak{O}_K : \text{$x$, $y$
coprime}\}$. Let
\begin{equation}\label{eq:adef}
\mathbb{A}_{\mathcal{E}} = \{(x,y)\in \mathbb{A}: 
x\ne 0, c_4(y/x)\ne \infty, c_6(y/x)\ne \infty, 
\Delta(y/x) \ne 0,\infty,\, q_0(y/x)\ne 0\} .\end{equation}
Let $(x,y)\in \mathbb{A}_{\mathcal{E}}$. 
Then $c_4(y/x)$ (resp. $c_6(y/x)$, $\Delta(y/x)$) differs from
$C_4(x,y)$ (resp. $C_6(x,y)$, $\Delta(x,y)$) by a non-zero fourth power $Q^4(x,y)$ (resp.
a non-zero sixth power $Q^6(x,y)$, 
a non-zero twelfth power $Q^{12}(x,y)$). Hence, for every prime ideal
$\mathfrak{p}\in I_K$, the reduction of $\mathcal{E}(y/x)$ at $\mathfrak{p}$ is
\begin{itemize}
\item good if $v_{\mathfrak{p}}(C_4(x,y)) = 4 k$, 
$v_{\mathfrak{p}}(C_6(x,y)) = 6 k$, $v_{\mathfrak{p}}(D(x,y)) = 12 k$ 
for some integer $k$;
\item multiplicative if $v_{\mathfrak{p}}(C_4(x,y)) = 4 k$, 
$v_{\mathfrak{p}}(C_6(x,y)) = 6 k$, 
$v_{\mathfrak{p}}(D(x,y)) > 12 k$ for some integer $k$;
\item additive and potentially multiplicative if 
$v_{\mathfrak{p}}(C_4(x,y)) = 4 k + 2$, 
$v_{\mathfrak{p}}(C_6(x,y)) = 6 k + 3$ and 
$v_{\mathfrak{p}}(D(x,y)) > 12 k + 6$ for
some integer $k$;
\item additive and potentially good in all remaining cases.
\end{itemize}

The {\em root number} of an elliptic curve over a global field $K$
is the product of its local root numbers
\[W(E) = \prod_v W_v(E)\]
over all places $v$ of $K$.
Similarly, given $\mathfrak{d}\in I_K$,
we define the {\em putative root number} 
$V_{\mathfrak{d}}(\mathcal{E})$ 
of an elliptic curve $\mathcal{E}$
over $K(t)$ to be the product of its local putative root numbers
\[V_{\mathfrak{d}}(\mathcal{E}) = \prod_v V_{\mathfrak{d},v}(E)\]
over all places $v$ of $K(t)$. We will define {\em local putative root
numbers} shortly. Note for now that $V^{\mathfrak{d},v}(\mathcal{E}) =1$ for all
but finitely many places $v$ of $K(t)$, just as $W_v(E)=1$ for all
but finitely many places $v$ of $K$.

\begin{prop}\label{prop:rohr}
Let $K$ be a number field. Let $\mathfrak{p}$ be prime ideal of $K$ unramified
over $\mathbb{Q}$. Assume $\mathfrak{p}$ lies
over a rational prime $p$ greater than three. Let $E$ be an elliptic 
curve over $K$ whose reduction at $\mathfrak{p}$ is additive and potentially
good. 
Then
\begin{enumerate}
\item $W_{\mathfrak{p}}(E) = \quadrec{-1}{\mathfrak{p}}$ if 
$v_{\mathfrak{p}}(\Delta(E))$ is even but not divisible by four,
\item $W_{\mathfrak{p}}(E) = \quadrec{-2}{\mathfrak{p}}$ if 
$v_{\mathfrak{p}}(\Delta(E))$ is odd and divisible by three,
\item $W_{\mathfrak{p}}(E) = \quadrec{-3}{\mathfrak{p}}$ if 
$v_{\mathfrak{p}}(\Delta(E))$ is divisible by four but not by three.
\end{enumerate}
\end{prop}
\begin{proof}
Let $a$ be any rational integer not divisible by $\mathfrak{p}$. 
If $\deg(K_{\mathfrak{p}}/\mathbb{Q}_p)$ is even, then
$\quadrec{a}{\mathfrak{p}} = 1$. If $\deg(K_{\mathfrak{p}}/\mathbb{Q}_p)$ is odd, then
$\quadrec{a}{\mathfrak{p}} = \quadrec{a}{p}$. Apply \cite{Rog}, Theorem 2, to the case of the 
trivial one-dimensional representation.
\end{proof}
Define $M_{\mathcal{E}}$, $B_{\mathcal{E}}$, $B'_{\mathcal{E}}$ as in (\ref{eq:M}) and
(\ref{eq:bb'}). Let $\lbrack a,b\rbrack_{\mathfrak{d}}$ be as in (\ref{eq:brack}). 
Let $\mathfrak{d}_0\in I_K$ be the principal ideal generated by
\begin{equation}\label{eq:toolong}
6 \,D_0
\mathop{\prod_{v_1\ne v_2}}_{\text{$\mathcal{E}$ has bad red. at $v_1$, $v_2$}}
\Res(P_{v_1},P_{v_2}) ,\end{equation}
where $D_0$ is as in (\ref{eq:defd}). 
\begin{defn}\label{def:putative}
Let $K$ be a number field. Let $\mathcal{E}$ be an elliptic curve over $K(t)$.
Let $\mathfrak{d}\in I_K$ be an ideal divisible by $\mathfrak{d}_0$.
Let $v$ be a place of $K(t)$.
Define the {\em local putative root number} $V_v(\mathcal{E})$ to be a map
from $\mathbb{A}_{\mathcal{E}}$ to $\{-1,1\}$
whose values are given as follows:
\begin{enumerate}
\item $V_{\mathfrak{d},v}(\mathcal{E}) = 1$ if the reduction $\mathcal{E} \mo v$ is good,
\item $V_{\mathfrak{d},v}(\mathcal{E}) = \lambda_K(P_v(x,y)) \cdot \lbrack - C_6(x,y),P_v(x,y) \rbrack_\mathfrak{d}$
if the reduction is multiplicative,
\item $V_{\mathfrak{d},v}(\mathcal{E}) = \lbrack -1, P_v(x,y) \rbrack_\mathfrak{d}$
if the reduction is additive and potentially multiplicative,
\item $V_{\mathfrak{d},v}(\mathcal{E}) = \lbrack -1, P_v(x,y) \rbrack_\mathfrak{d}$
if the reduction is additive and potentially good, and 
$v(\Delta)$ is even but not divisible by four,
\item $V_{\mathfrak{d},v}(\mathcal{E}) = \lbrack -2, P_v(x,y) \rbrack_\mathfrak{d}$
if the reduction is additive and potentially good, and 
$v(\Delta)$ is odd and divisible by three,
\item $V_{\mathfrak{d},v}(\mathcal{E}) = \lbrack -3, P_v(x,y) \rbrack_\mathfrak{d}$
if the reduction is additive and potentially good, and 
$v(\Delta)$ is divisible by four but not by three.
\end{enumerate}
\end{defn}

We define {\em
half bad} and {\em quite bad} reduction as in section \ref{sec:isdf}. The reduction of
$\mathcal{E}$ at $v$ is
\begin{itemize}
\item {\em half bad} if $e_{v,4}\geq 2$, $e_{v,6}\geq 3$, $e_{v,D}=6$,
\item {\em quite bad} if it is bad but not half bad.
\end{itemize}
The reduction of $\mathcal{E}(y/x)$ at $\mathfrak{p}$ is
\begin{itemize}
\item {\em half bad} if $v_{\mathfrak{p}}(C_4(x,y)) \geq 4 k + 2$,
$v_{\mathfrak{p}}(C_6(x,y)) \geq 6 k + 3$ and 
$v_{\mathfrak{p}}(D(x,y)) = 12 k + 6$ for some integer $k$,
\item {\em quite bad} if it is bad but not half bad.
\end{itemize}
It should be clear that half-bad reduction is a special case of additive, potentially good
reduction.

As in subsection \ref{ssec:famdef}, 
we set $W(\mathcal{E}(y/x))=1$ when $\mathcal{E}(y/x)$ is undefined or singular. 
Note that the set
$\{x,y\in \mathfrak{O}_K:\gcd(x,y)=1, \text{$\mathcal{E}(y/x)$ undefined or singular}\}$
is finite, as is its superset $\{x,y\in \mathfrak{O}_K:\gcd(x,y)=1\}-\mathbb{A}_{\mathcal{E}}$.

\subsection{From the root number to Liouville's function}
\begin{lem}\label{lem:bigone}
Let $K$ be a number field. Let $\mathcal{E}$ be an elliptic curve over 
$K(t)$. Let $\mathfrak{d}_0$ be as in \ref{eq:toolong}.
Let $\mathfrak{d}\in I_K$ be an ideal divisible by $\mathfrak{d}_0$.
The putative root number $V_{\mathfrak{d}}(\mathcal{E})$ is of the form
\[V_{\mathfrak{d}}(\mathcal{E}) = f(x,y) \cdot \lambda_K(M_\mathcal{E}(x,y)),\]
where $f$ is a pliable function on $\{(x,y)\in \mathfrak{O}_K^2 : \text{$x$, $y$ coprime}\}$.
\end{lem}
\begin{proof} Let $v$ be a place of $\mathcal{E}$. If the reduction
of $\mathcal{E}$ at $v$ is good, then $V_{\mathfrak{d},v}(\mathcal{E})$ is equal
to the constant $1$
and hence is pliable. If the reduction of $\mathcal{E}$ at $v$ is additive,
$V_{\mathfrak{d},v}(\mathcal{E})$ is pliable by properties (4) and (5) of $\lbrack,\rbrack_{\mathfrak{d}}$
(see subsection \ref{subs:quadr}). 
If the reduction of $\mathcal{E}$ at $v$ is multiplicative,
then $V_{\mathfrak{d},v}(\mathcal{E})$ is equal to the product of 
$\lambda_K(P_v(x,y))$
and a pliable function by Corollary \ref{cor:important} and by the fact
that $\deg(C_6(x,y))$ is even.

The reduction of $\mathcal{E}$ at $v$ is bad for only a finite number
of places $v$. Since the product of finitely many pliable functions is
pliable, we obtain
\[V_{\mathfrak{d}}(\mathcal{E}) = f(x,y) \prod_{\text{$\mathcal{E}$ has mult. red. at 
$v$}}
\lambda_K(P_v(x,y)) = f(x,y) \cdot \lambda_K(M_{\mathcal{E}}(x,y)),\]
where $f(x,y)$ is a pliable function on 
 $\{(x,y)\in \mathfrak{O}_K^2 : \text{$x$, $y$ coprime}\}$.
\end{proof}

\begin{lem}\label{lem:mess} Let $K$ be a number field. Let $\mathcal{E}$ be an elliptic
curve over $K(t)$. Let $v$ be a place of $K(t)$ where $\mathcal{E}$
has bad reduction. Let $\mathbb{A}_{\mathcal{E}}$ be as in (\ref{eq:adef}).
Let $\mathfrak{d}$
be as in (\ref{eq:toolong}). Then, for any $(x,y)\in \mathbb{A}_{\mathcal{E}}$,
\[\begin{aligned}
\mathop{\prod_{\mathfrak{p}\nmid \mathfrak{d}}}_{\mathfrak{p}|P_{v}(x,y)}
 W_{\mathfrak{p}}(\mathcal{E}(y/x)) &= g_{v}(x,y) \cdot 
V_{\mathfrak{d},v}(\mathcal{E})(x,y)\;
&\text{if $v$ is half bad},\\
\mathop{\prod_{\mathfrak{p}\nmid \mathfrak{d}}}_{\mathfrak{p}|P_{v}(x,y)}
 W_{\mathfrak{p}}(\mathcal{E}(y/x)) &= g_{v}(x,y) \cdot
 h(\sq_K(P_{v}(x,y)),x,y)\cdot V_{\mathfrak{d},v}(\mathcal{E})(x,y)
&\text{if $v$ is quite bad},
\end{aligned}\]
where $g_{v}:\mathbb{A}_{\mathcal{E}}\to \{-1,1\}$,
$h:I_K\times \mathbb{A}_{\mathcal{E}}\to \{-1,1\}$ satisfy the following conditions:
\begin{enumerate}
\item $g_{v}$ is pliable, 
\item 
$h(\mathfrak{a},x,y)$ depends only on $\mathfrak{a}$ and on 
$\{\frac{x \mo \mathfrak{p}}{y \mo \mathfrak{p}}\}_{\mathfrak{p}|\mathfrak{a}} \in 
\prod_{\mathfrak{p}|\mathfrak{a}} \mathbb{P}^1(\mathfrak{O}_K/\mathfrak{p})$,
\item $h(\mathfrak{a}_1 \mathfrak{a}_2,x,y) = h(\mathfrak{a}_1,x,y)
h(\mathfrak{a}_2,x,y)$ for any $\mathfrak{a}_1, \mathfrak{a}_2\in I_K$,
\item $h(\mathfrak{a},x,y)=1$ for $\mathfrak{a}|\mathfrak{d}^{\infty}$.
\end{enumerate}
\end{lem}
\begin{proof}
The reduction of $\mathcal{E}$ at $v$ can be multiplicative or additive. If it
is additive, it can be potentially multiplicative or potentially good.
If it is additive and potentially good, it can be half bad or quite bad.
If it is additive, potentially good and quite bad, then 
$\gcd(e_{v,D},12)$ is $2$, $3$ or $4$. We speak of reduction type
$pg_2$, $pg_3$, $pg_4$ accordingly.

We will construct 
 $h_m, h_{mp}, h_{pg_2}, h_{pg_3}, h_{pg_4}:I_K\times \mathbb{A}_{\mathcal{E}}\to \{-1,1\}$,
each of them satisfying the conditions (2)-(4) enunciated for $h$ in the
statement. We will also define a pliable function 
 $g_{v}:\mathbb{A}_{\mathcal{E}}\to \{-1,1\}$ depending on $v$. Our aim is to prove that
$\prod_{\mathfrak{p}\nmid \mathfrak{d}, \mathfrak{p}|P_{v}(x,y)} 
 W_{\mathfrak{p}}(\mathcal{E}(y/x))$ equals
\begin{equation}\label{eq:squiggle}\begin{aligned}
g_{v}(x,y)\cdot V_{\mathfrak{d},v}(\mathcal{E})(x,y)
 \;\;&\text{if $\mathcal{E} \mo v$ is half bad,}\\
g_{v}(x,y)\cdot h_{pg_2}(\sq(P_{v}(x,y)),x,y)\cdot V_{\mathfrak{d},v}(\mathcal{E})(x,y)
 \;\;&\text{if $\mathcal{E} \mo v$ is of type $pg_2$,}\\
g_{v}(x,y)\cdot h_{pg_3}(\sq(P_{v}(x,y)),x,y)\cdot V_{\mathfrak{d},v}(\mathcal{E})(x,y)
 \;\;&\text{if $\mathcal{E} \mo v$ is of type $pg_3$,}\\
g_{v}(x,y)\cdot h_{pg_4}(\sq(P_{v}(x,y)),x,y)\cdot V_{\mathfrak{d},v}(\mathcal{E})(x,y)
 \;\;&\text{if $\mathcal{E} \mo v$ is of type $pg_4$,}\\
g_{v}(x,y)\cdot h_{pm}(\sq(P_{v}(x,y)),x,y)\cdot V_{\mathfrak{d},v}(\mathcal{E})(x,y)
 \;\;&\text{if $\mathcal{E} \mo v$ is additive and pot. mult.,}\\
g_{v}(x,y)\cdot h_m(\sq(P_{v}(x,y)),x,y)\cdot V_{\mathfrak{d},v}(\mathcal{E})(x,y)
 \;\;&\text{if $\mathcal{E} \mo v$ is multiplicative.}
\end{aligned}\end{equation}
Then we can define $h:I_K\times \mathbb{A}_{\mathcal{E}}\to \{-1,1\}$ to be the function such that
$h(\mathfrak{p}^n,x,y)=1$ for $\mathfrak{p}|\mathfrak{d}$,
\begin{equation}\label{eq:litstar}
h(\mathfrak{p}^n,x,y) = \begin{cases}
 h_m(\mathfrak{p}^n,x,y) &\text{if $\mathfrak{p}|\prod_{\text{$v$ mult.}} P_{v}(x,y)$,}\\
 h_{mp}(\mathfrak{p}^n,x,y) 
  &\text{if $\mathfrak{p}|\prod_{\text{$v$ add. and pot. mult.}} P_{v}(x,y)$,}\\
 h_{pg_2}(\mathfrak{p}^n,x,y) 
&\text{if $\mathfrak{p}|\prod_{\text{$v$ is $pg_2$}} P_{v}(x,y)$,}\\
 h_{pg_3}(\mathfrak{p}^n,x,y) 
&\text{if $\mathfrak{p}|\prod_{\text{$v$ is $pg_3$}} P_{v}(x,y)$,}\\
 h_{pg_4}(\mathfrak{p}^n,x,y) 
&\text{if $\mathfrak{p}|\prod_{\text{$v$ is $pg_4$}} P_{v}(x,y)$,}\\
1 &\text{otherwise}
\end{cases}\end{equation}
for $\mathfrak{p}\nmid \mathfrak{d}$, and 
 $h(\mathfrak{a}_1 \mathfrak{a}_2,x,y) = h(\mathfrak{a}_1,x,y) h(\mathfrak{a}_2,x,y)$
for any $\mathfrak{a}_1, \mathfrak{a}_2\in I_K$.

First note that no more than one case can hold in (\ref{eq:litstar}), as 
$\mathfrak{p}\nmid \mathfrak{d}$ implies that $\mathfrak{p}$ cannot divide both
$P_{v}(x,y)$ and $P_{u}(x,y)$ for $v$, $u$ distinct (see
(\ref{eq:toolong})). Notice, too, that
condition (2) in the statement is fulfilled: since $P_{v}$ is homogeneous,
whether or not $\mathfrak{p}|P_{v}(x,y)$ for given $x$, $y$ depends only on
$\frac{x \mo \mathfrak{p}}{y \mo \mathfrak{p}}$. Finally, it is an immediate consequence
of (\ref{eq:litstar}) that
\[
h(\sq_K(P_{v}(x,y)),x,y) = \begin{cases}
h_m(\sq_K(P_{v}(x,y)),x,y)\;\;
&\text{if $\mathcal{E} \mo v$ is multiplicative,}\\
h_{pm}(\sq(P_{v}(x,y)),x,y)\;\;
&\text{if $\mathcal{E} \mo v$ is add. and pot. m.,}\\
h_{pg_2}(\sq(P_{v}(x,y)),x,y)\;\;
&\text{if $\mathcal{E} \mo v$ is $pg_2$}\\
h_{pg_3}(\sq(P_{v}(x,y)),x,y)
 \;\;&\text{if $\mathcal{E} \mo v$ is $pg_3$}\\
h_{pg_4}(\sq(P_{v}(x,y)),x,y)\;\;&\text{if $\mathcal{E} \mo v$ is $pg_4$.}
\end{cases}\]
The statement then follows from (\ref{eq:squiggle}). It remains to construct
$g_{v}$, $h_m$, $h_{pm}$, $h_{pg_2}$, $h_{pg_3}$, $h_{pg_4}$ and to
prove (\ref{eq:squiggle}).

Let $e_{v,4}$ $e_{v,6}$, $e_{v,D}$ be as in (\ref{eq:defd}). 
Suppose $\mathfrak{p}\nmid \mathfrak{d}$, $\mathfrak{p}|P_{v}(x,y)$. Then
$\mathfrak{p}\nmid P_{u}(x,y)$ for every $u \ne v$. Hence
\begin{equation}\label{eq:tripstar}\begin{aligned} 
v_{\mathfrak{p}}(C_4(x,y)) &= e_{v,4}\cdot v_{\mathfrak{p}}(P_{v}(x,y)),\\
v_{\mathfrak{p}}(C_6(x,y)) &= e_{v,6}\cdot v_{\mathfrak{p}}(P_{v}(x,y)),\\
v_{\mathfrak{p}}(D(x,y)) &= e_{v,D}\cdot v_{\mathfrak{p}}(P_{v}(x,y)).
\end{aligned}\end{equation}

{\em Case 1: $\mathcal{E}$ has multiplicative reduction at $v$.}
We are given that $e_{v,4}=0$, $e_{v,6}=0$, $e_{v,D}>0$.
Hence $v_{\mathfrak{p}}(C_4(x,y))=0$, $v_{\mathfrak{p}}(C_6(x,y))=0$,
$v_{\mathfrak{p}}(D(x,y))>0$. Therefore, $\mathcal{E}(y/x)$ has multiplicative
reduction at $\mathfrak{p}$. By Lemma \ref{lem:rootofall},
\[W_{\mathfrak{p}}(\mathcal{E}(y/x)) = - \quadrec{-C_6(x,y)}{\mathfrak{p}} .\]
Thus
\[\begin{aligned}
\mathop{\prod_{\mathfrak{p}\nmid \mathfrak{d}}}_{\mathfrak{p}|P_{v}(x,y)}
W_{\mathfrak{p}}(\mathcal{E}(y/x)) &= 
\mathop{\prod_{\mathfrak{p}\nmid \mathfrak{d}}}_{\mathfrak{p}|P_{v}(x,y)}
\left(- \quadrec{-C_6(x,y)}{\mathfrak{p}}\right) \\ 
&= \mathop{\prod_{\mathfrak{p}|\mathfrak{d}}}_{\mathfrak{p}|P_{v}(x,y)} 
(-1)^{v_{\mathfrak{p}}(P_{v}(x,y))} 
\prod_{\mathfrak{p}\nmid \mathfrak{d}, \mathfrak{p}^2|P_{v}(x,y)} 
(-1)^{v_{\mathfrak{p}}(P_{v}(x,y))-1}\\ &\cdot
\mathop{\prod_{\mathfrak{p}\nmid \mathfrak{d}}}_{\mathfrak{p}^2 |P_{v}(x,y)}
\quadrec{- C_6(x,y)}{\mathfrak{p}}^{v_{\mathfrak{p}}(P_{v}(x,y))-1} 
\\
&\cdot \prod_{\mathfrak{p}|P_{v}(x,y)} (-1)^{v_{\mathfrak{p}}(P_{v}(x,y))}
\mathop{\prod_{\mathfrak{p}\nmid \mathfrak{d}}}_{\mathfrak{p}|P_{v}(x,y)}
 \quadrec{- C_6(x,y)}{\mathfrak{p}}^{v_{\mathfrak{p}}(P_{v}(x,y))}
 .\end{aligned}\]
Let \[\label{eq:doubstar}\begin{aligned}
g_{v}(x,y) &= \prod_{\mathfrak{p}|\mathfrak{d}, \mathfrak{p}|
 P_{v}(x,y)} (-1)^{v_{\mathfrak{p}}(P_{v}(x,y))},\\ 
h_m(\mathfrak{a},x,y) &= 
 \lambda_K\left(\frac{\mathfrak{a}}{\gcd(\mathfrak{a},\mathfrak{d}^{\infty})}\right)
\cdot \lbrack -C_6(x,y),\mathfrak{a}\rbrack_{\mathfrak{d}} .
\end{aligned}\]
Then
$\prod_{\mathfrak{p}\nmid \mathfrak{d}, \mathfrak{p}|P_{v}(x,y)}
 W_{\mathfrak{p}}(\mathcal{E}(y/x)) $ is
\begin{equation}\label{eq:twinkle}
g_{v}(x,y) \cdot h_m(\sq_K(P_{v}(x,y)),x,y) \cdot
\lambda_K(P_{v}(x,y))
\lbrack - C_6(x,y),P_{v}(x,y)\rbrack_{\mathfrak{d}}
.\end{equation}

The map $t\mapsto (-1)^{v_{\mathfrak{p}}(t)}$ on $K$ is pliable.
Hence, by Proposition \ref{prop:plinvar}, 
$(x,y)\mapsto (-1)^{v_{\mathfrak{p}}(P_{v}(x,y))}$ is a pliable function
on $A$. Since $g_v(x,y)$ equals 
$\prod_{\mathfrak{p}|\mathfrak{d}} (-1)^{v_{\mathfrak{p}}(P_{v}(x,y))}$,
which is a product of finitely many pliable functions, $g_{v}(x,y)$ is pliable.

It remains to show that $h_m(\mathfrak{a},x,y)$ depends only on $\mathfrak{a}$
and $\{\frac{x \mo \mathfrak{p}}{y \mo \mathfrak{p}}\}_{\mathfrak{p}|\mathfrak{a}}$.
For fixed $\mathfrak{a}$, the first factor
$\lambda_K\left(\frac{\mathfrak{a}}{\gcd(\mathfrak{a},\mathfrak{d}^{\infty})}\right)$
 is a constant. Since
\[\lbrack - C_6(x,y),\mathfrak{a}\rbrack_{\mathfrak{d}} =
 \mathop{\prod_{\mathfrak{p}\nmid \mathfrak{d}}}_{\mathfrak{p}|\mathfrak{a}}
 \quadrec{-C_6(x,y)}{\mathfrak{p}}^{v_{\mathfrak{p}}(\mathfrak{a})},\]
it is enough to show that $\quadrec{-C_6(x,y)}{\mathfrak{p}}$ depends only on $\frac{x \mo \mathfrak{p}}{y \mo \mathfrak{p}}$
for every prime $\mathfrak{p}$ with
$\mathfrak{p}|\mathfrak{a}$, $\mathfrak{p}\nmid \mathfrak{d}$. For every $t\in \mathfrak{O}_K^*$,
\[\quadrec{-C_6(r x, r y)}{\mathfrak{p}} = 
 \quadrec{-r^{\deg C_6} C_6(x,y)}{\mathfrak{p}} =
 \quadrec{r}{\mathfrak{p}}^{\deg C_6} \quadrec{- C_6(x,y)}{\mathfrak{p}} .\]
Since $\deg C_6$ is even, it follows that
\[\quadrec{-C_6(r x, r y)}{\mathfrak{p}} = 
 \quadrec{- C_6(x,y)}{\mathfrak{p}} .\]
Hence $\quadrec{-C_6(x,y)}{\mathfrak{p}}$ depends only on 
$\frac{x \mo \mathfrak{p}}{y \mo \mathfrak{p}}$. Therefore 
$h_m(\mathfrak{a},x,y)$ depends only on $\mathfrak{a}$ and
$\{\frac{x \mo \mathfrak{p}}{y \mo \mathfrak{p}}\}_{\mathfrak{p}|\mathfrak{a}}$.

We have shown that $g_{v}$ and $h_{v}$ in (\ref{eq:twinkle}) satisfy properties
(1) and (2) in the statement. Properties (3) and (4) are immediate from
(\ref{eq:doubstar}). 
Since $V_{\mathfrak{d},v}(\mathcal{E})(x,y) = \lambda_K(P_{v}(x,y)) \lbrack - C_6(x,y),P_{v}(x,y)
\rbrack_{\mathfrak{d}}$, we are done.

{\em Case 2: $\mathcal{E}$ has additive, potentially multiplicative reduction at $v$.}
We are given $e_{v,4}=2$, $e_{v,6}=3$, $e_{v,D}>6$. Let $\mathfrak{p}$ be a prime
ideal dividing $P_{v}(x,y)$ but not $\mathfrak{d}$. 
Then
$v_{\mathfrak{p}}(C_4(x,y))= 4 k$,
$v_{\mathfrak{p}}(C_6(x,y)) = 6 k$,
$v_{\mathfrak{p}}(D(x,y)) > 12 k$ if $v_{\mathfrak{p}}(P_{v}(x,y))=2 k$, $k>0$,
and
$v_{\mathfrak{p}}(C_4(x,y))= 4 k + 2$,
$v_{\mathfrak{p}}(C_6(x,y)) = 6 k + 3$,
$v_{\mathfrak{p}}(D(x,y)) > 12 k + 6$ if $v_{\mathfrak{p}}(P_{v}(x,y))=2 k +1$, 
$k\geq 0$. Thus, $\mathcal{E}(y/x)$ has multiplicative reduction at $\mathfrak{p}$
if $v_{\mathfrak{p}}(P_{v}(x,y))$ is even and positive, but has 
additive, potentially multiplicative
reduction if $v_{\mathfrak{p}}(P_{v}(x,y))$ is odd.

Hence, by Lemma \ref{lem:rootofall},
\[\begin{aligned}
\mathop{\prod_{\mathfrak{p}\nmid \mathfrak{d}}}_{\mathfrak{p}|P_{v}(x,y)}
 W_{\mathfrak{p}}(\mathcal{E}(y/x)) &=
\mathop{\mathop{\prod_{\mathfrak{p}\nmid \mathfrak{d}}}_{\mathfrak{p}|P_{v}(x,y)}
 }_{\text{$v_{\mathfrak{p}}(P_{v}(x,y))$ even}}
- \quadrec{-C_6(x,y)\mathfrak{p}^{-v_{\mathfrak{p}}(C_6(x,y))}
}{\mathfrak{p}}
\mathop{\mathop{\prod_{\mathfrak{p}\nmid \mathfrak{d}}}_{\mathfrak{p}|P_{v}(x,y)}}_{\text{$v_{\mathfrak{p}}(P_{v}(x,y))$ odd}}
 \quadrec{-1}{\mathfrak{p}} \\
&= h_{pm}(\sq_K(P_{v}(x,y)),x,y) \cdot
 \lbrack -1, P_{v}(x,y)\rbrack_{\mathfrak{d}} ,\end{aligned}\]
where $h_{pm}(\mathfrak{a},x,y) = 
 \prod_{\mathfrak{p}|\mathfrak{d}, \mathfrak{p}\nmid \mathfrak{a}}
 \left(-\quadrec{-C_6(x,y)}{\mathfrak{p}}\right)^{v_{\mathfrak{p}}(a)}$.
It is clear that $h_{pm}(\mathfrak{a},x,y)$ is multiplicative on $\mathfrak{a}$
and trivial for $\mathfrak{a}|\mathfrak{d}^{\infty}$. As shown above,
$\quadrec{-C_6(x,y)}{\mathfrak{p}}$ depends only on $\mathfrak{p}$
and $\frac{x \mo \mathfrak{p}}{y \mo \mathfrak{p}}$. Hence 
 $h_{pm}(\mathfrak{a},x,y)$ depends only on $\mathfrak{a}$ and
$\{\frac{x \mo \mathfrak{p}}{y \mo \mathfrak{p}}\}_{\mathfrak{p}|\mathfrak{a}}$.
Set $g_{v}(x,y)=1$, 
Since $V_{\mathfrak{d},v}(\mathcal{E})(x,y)=\lbrack -1,P_{v}(x,y)\rbrack_{\mathfrak{d}}$,
\[\mathop{\prod_{\mathfrak{p}\nmid \mathfrak{d}}}_{\mathfrak{p}|P_{v}(x,y)}
 W_{\mathfrak{p}}(\mathcal{E}(y/x)) =
g_{v}(x,y) \cdot h_{pm}(\sq_K(P_v(x,y)),x,y) \cdot V_{\mathfrak{d},v}(\mathcal{E})(x,y) .\]
{\em Case 3: $\mathcal{E}$ has half-bad reduction at $v$.}
We are given $e_{v,4}\geq 2$, $e_{v,6}\geq 3$, $e_{v,D}=6$. Let $\mathfrak{p}$ be a prime
ideal dividing $P_{v}(x,y)$ but not $\mathfrak{d}$. 
Then
$v_{\mathfrak{p}}(C_4(x,y))\geq 4 k$,
$v_{\mathfrak{p}}(C_6(x,y)) \geq 6 k$,
$v_{\mathfrak{p}}(D(x,y)) = 12 k$ if $v_{\mathfrak{p}}(P_{v}(x,y))=2 k$, $k>0$,
and
$v_{\mathfrak{p}}(C_4(x,y))\geq 4 k + 2$,
$v_{\mathfrak{p}}(C_6(x,y)) \geq 6 k + 3$,
$v_{\mathfrak{p}}(D(x,y)) = 12 k + 6$ if $v_{\mathfrak{p}}(P_{v}(x,y))=2 k +1$, 
$k\geq 0$. Thus, $\mathcal{E}(y/x)$ has half-bad
 reduction at $\mathfrak{p}$
if $v_{\mathfrak{p}}(P_{v}(x,y))$ is odd,
and good reduction if $v_{\mathfrak{p}}(P_{v}(x,y))$ is even.
Hence, by Proposition \ref{prop:rohr},
\[W_{\mathfrak{p}}(\mathcal{E}(y/x)) =\begin{cases} 1
&\text{if $v_{\mathfrak{p}}(P_{v}(x,y))$ is even,}\\
\quadrec{-1}{\mathfrak{p}}
&\text{if $v_{\mathfrak{p}}(P_{v}(x,y))$ is odd.}\end{cases}\]
Thereby
\[\mathop{\prod_{\mathfrak{p}\nmid \mathfrak{d}}}_{\mathfrak{p}|P_{v}(x,y)}
W_{\mathfrak{p}}(\mathcal{E}(y/x)) = \mathop{\prod_{\mathfrak{p}\nmid \mathfrak{d}}}_{\mathfrak{p}|P_{v}(x,y)} 
\quadrec{-1}{\mathfrak{p}}^{v_{\mathfrak{p}}(P_{v}(x,y))}
= \lbrack -1, P_{v}(x,y)\rbrack_{\mathfrak{d}} = V_{\mathfrak{d},v}(\mathcal{E})(x,y) .\]
Set $g_{v}(x,y)=1$.

{\em Case 4: $\mathcal{E}$ has $gp_2$ reduction at $v$.}
We are given that the reduction is additive and $\gcd(e_{v,D},12)=2$. 
Then the reduction of $\mathcal{E}(y/x)$ at $\mathfrak{p}$
is good if $6|v_{\mathfrak{p}}(P_{v}(x,y))$ and additive and potentially good
otherwise if $6\nmid v_{\mathfrak{p}}(P_{v}(x,y))$. Hence
\[\gcd(v_{\mathfrak{p}}(D(x,y)),12) =\begin{cases}
2 &\text{if $v_{\mathfrak{p}}(P_{v}(x,y))\equiv 1, 5 \mo 6$,}\\
4 &\text{if $v_{\mathfrak{p}}(P_{v}(x,y))\equiv 2, 4 \mo 6$,}\\
6 &\text{if $v_{\mathfrak{p}}(P_{v}(x,y))\equiv 3 \mo 6$.}
\end{cases}\]
So, by Proposition \ref{prop:rohr},
\[W_{\mathfrak{p}}(\mathcal{E}(y/x)) = \begin{cases}
1 &\text{if $v_{\mathfrak{p}}(P_{v}(x,y))\equiv 0 \mo 6$,}\\
\quadrec{-2}{\mathfrak{p}} &\text{if $v_{\mathfrak{p}}(P_{v}(x,y))\equiv 3 \mo 6$,}\\
\quadrec{-1}{\mathfrak{p}} &\text{if $v_{\mathfrak{p}}(P_{v}(x,y))\not\equiv 0 \mo 3$.}
\end{cases}\]
Let $H:I_K\to \{-1,1\}$ be the multiplicative function such that $H(\mathfrak{p}^n)=1$
for $\mathfrak{p}|\mathfrak{d}$ and
\[H(\mathfrak{p}^n) = \begin{cases}
1 &\text{if $n\equiv 0,4,5 \mo 6$}\\
\quadrec{-1}{\mathfrak{p}} &\text{if $n\equiv 1,3 \mo 6$}\\
\quadrec{2}{\mathfrak{p}} &\text{if $n\equiv 2 \mo 6$}
\end{cases}\]
for $\mathfrak{p}\nmid \mathfrak{d}$. Then
\[W_{\mathfrak{p}}(\mathcal{E}(y/x)) = \begin{cases}
\quadrec{-1}{\mathfrak{p}}^{v_{\mathfrak{p}}(P_{v}(x,y))}
&\text{if $v_{\mathfrak{p}}(P_{v}(x,y))=1$,}\\
H(\mathfrak{p}^{v_{\mathfrak{p}}(P_{v}(x,y))-1} ) 
\quadrec{-1}{\mathfrak{p}}^{v_{\mathfrak{p}}(P_{v}(x,y))}
&\text{if $v_{\mathfrak{p}}(P_{v}(x,y))>1$.}\end{cases}\]
Hence \[\begin{aligned}
\mathop{\prod_{\mathfrak{p}\nmid \mathfrak{d}}}_{
\mathfrak{p}|P_{v}(x,y)} W_{\mathfrak{p}}(\mathcal{E}(y/x)) &= 
\mathop{\prod_{\mathfrak{p}\nmid \mathfrak{d}}}_{\mathfrak{p}^2 | P_{v}(x,y)}
 H(\mathfrak{p}^{v_{\mathfrak{p}}(P_{v}(x,y))-1})
\mathop{\prod_{\mathfrak{p}\nmid \mathfrak{d}}}_{\mathfrak{p} | P_{v}(x,y)}
 \quadrec{-1}{\mathfrak{p}}^{v_{\mathfrak{p}}(P_{v}(x,y))}
\\ &= H(\sq_K(P_{v}(x,y))) \cdot \lbrack -1,P_{v}(x,y)\rbrack_{\mathfrak{d}} .\end{aligned}\]
Set $g_{v}(x,y)=1$, $h_{gp_2}(\mathfrak{a},x,y)=H(\mathfrak{a})$ and we are done.

{\em Case 5: $\mathcal{E}$ has $gp_3$ reduction at $v$.}
We are given that the reduction is additive and $\gcd(e_{v,D},12)=3$. Then
\[\gcd(v_{\mathfrak{p}}(D(x,y)),12) = \begin{cases}
 3 &\text{if $v_{\mathfrak{p}}(P_{v}(x,y))\equiv 1,3 \mo 4$,}\\
 6 &\text{if $v_{\mathfrak{p}}(P_{v}(x,y))\equiv 2 \mo 4$,}\\
 12 &\text{if $v_{\mathfrak{p}}(P_{v}(x,y))\equiv 0 \mo 4$.}
\end{cases}\]
So, by Proposition \ref{prop:rohr},
\[W_{\mathfrak{p}}(\mathcal{E}(y/x)) = \begin{cases}
1 &\text{if $v_{\mathfrak{p}}(P_{v}(x,y))\equiv 0 \mo 4$,}\\
\quadrec{-1}{\mathfrak{p}} &\text{if $v_{\mathfrak{p}}(P_{v}(x,y))\equiv 2 \mo 4$,}\\
\quadrec{-2}{\mathfrak{p}} &\text{if $v_{\mathfrak{p}}(P_{v}(x,y))\not\equiv 1 \mo 2$.}
\end{cases}\]
Let $H:I_K\to \{-1,1\}$ be the multiplicative function such that $H(\mathfrak{p}^n)=1$
for $\mathfrak{p}|\mathfrak{d}$ and
\[H(\mathfrak{p}^n) = \begin{cases}
\quadrec{2}{\mathfrak{p}} &\text{if $n\equiv 1 \mo 4$}\\
1 &\text{if $n\not\equiv 1 \mo 4$}
\end{cases}\]
for $\mathfrak{p}\nmid \mathfrak{d}$. Then
\[\mathop{\prod_{\mathfrak{p}\nmid \mathfrak{d}}}_{
\mathfrak{p}|P_{v}(x,y)} W_{\mathfrak{p}}(\mathcal{E}(y/x)) =
H(\sq_K(P_{v}(x,y))) \cdot \lbrack -2, P_{v}(x,y)\rbrack_{\mathfrak{d}} .\]
Set $g_{v}(x,y)=1$, $h_{v}(\mathfrak{a},x,y) = H(\mathfrak{a})$ and we are done.

{\em Case 6: $\mathcal{E}$ has $gp_4$ reduction at $v$.}
We are given that the reduction is additive and $\gcd(e_{v,D},12)=4$. Then
\[\gcd(v_{\mathfrak{p}}(D(x,y)),12) = \begin{cases}
4 &\text{if $v_{\mathfrak{p}}(P_{v}(x,y))\equiv 1,2 \mo 3$,}\\
12 &\text{if $v_{\mathfrak{p}}(P_{v}(x,y))\equiv 0 \mo 3$.}\end{cases}\]
So, by Proposition \ref{prop:rohr},
\[W_{\mathfrak{p}}(\mathcal{E}(y/x)) = \begin{cases}
1 &\text{if $v_{\mathfrak{p}}(P_{v}(x,y))\equiv 0 \mo 3$,}\\
\quadrec{-3}{\mathfrak{p}} &\text{if $v_{\mathfrak{p}}(P_{v}(x,y))
 \not\equiv 0 \mo 3$.}
\end{cases}\]
Let $H:I_K\to \{-1,1\}$ be the multiplicative function such that $H(\mathfrak{p}^n)=1$
for $\mathfrak{p}|\mathfrak{d}$ and
\[H(\mathfrak{p}^n) = \begin{cases}
\quadrec{-3}{\mathfrak{p}} &\text{if $n\equiv 1,2,3 \mo 6$}\\
1 &\text{otherwise}
\end{cases}\]
for $\mathfrak{p}\nmid \mathfrak{d}$. Then
\[\mathop{\prod_{\mathfrak{p}\nmid \mathfrak{d}}}_{
\mathfrak{p}|P_{v}(x,y)} W_{\mathfrak{p}}(\mathcal{E}(y/x)) =
H(\sq_K(P_{v}(x,y))) \cdot \lbrack -3, P_{v}(x,y)\rbrack_{\mathfrak{d}} .\]
Set $g_{v}(x,y)=1$, $h_{v}(\mathfrak{a},x,y) = H(\mathfrak{a})$ and we are done.
\end{proof}
\begin{prop}\label{prop:ahem}
Let $K$ be a number field. Let $\mathcal{E}$ be an elliptic curve over $K(t)$.
Let $\mathfrak{M}\in I_K$.
Then
there are $g:\mathbb{A}_{\mathcal{E}}\to \{-1,1\}$,
$h:I_K\times \mathbb{A}_{\mathcal{E}}\to \{-1,1\}$ such that 
, for all $(x,y)\in \mathbb{A}_{\mathcal{E}}$,
\[W(\mathcal{E}(y/x)) = g(x,y)\cdot h(\sq_K(B_\mathcal{E}'(x,y)),x,y)\cdot 
\lambda_K(M_\mathcal{E}(x,y)),\]
and, furthermore,
\begin{enumerate}
\item $g$ is pliable, 
\item 
$h(\mathfrak{a},x,y)$ depends only on $\mathfrak{a}$ and on 
$\{\frac{x \mo \mathfrak{p}}{y \mo \mathfrak{p}}\}_{\mathfrak{p}|\mathfrak{a}} \in 
\prod_{\mathfrak{p}|\mathfrak{a}} \mathbb{P}^1(\mathfrak{O}_K/\mathfrak{p})$.
\item $h(\mathfrak{a}_1 \mathfrak{a}_2,x,y) = h(\mathfrak{a}_1,x,y)
h(\mathfrak{a}_2,x,y)$ for any $\mathfrak{a}_1, \mathfrak{a}_2\in I_K$,
\item $h(\mathfrak{a},x,y)=1$ for $\mathfrak{a}|\mathfrak{M}^{\infty}$.
\end{enumerate}
\end{prop}
\begin{proof}
For all $(x,y)\in \mathbb{A}_{\mathcal{E}}$, we can write
\[W(\mathcal{E}(y/x)) = W_{\infty}(\mathcal{E}(y/x)) \prod_{\mathfrak{p}}
 W_{\mathfrak{p}}(\mathcal{E}(y/x)) .\]
It follows from the definition of local root numbers that
$W_{\mathfrak{p}}(\mathcal{E}(y/x))=1$ when $\mathcal{E}(y/x)$ has
good reduction at $\mathfrak{p}$ (see, e.g., \cite{Ro2}, Sec. 19, Prop (i)).
We also know that $W_{\infty}=-1$ (see, e.g., \cite{Ro2}, Sec. 20). 
Let $\mathfrak{d} = \mathfrak{M} \cdot \mathfrak{d}_0$. Then
\[W(\mathcal{E}(y/x)) = - \prod_{\mathfrak{p}}
 W_{\mathfrak{p}}(\mathcal{E}(y/x))
= \prod_{\mathfrak{p}|\mathfrak{d}} W_{\mathfrak{p}}(\mathcal{E}(y/x))\; \cdot\mathop{\prod_{\mathfrak{p}|\mathfrak{d}}}_{\text{$\mathcal{E}(y/x)$ has bad red.
at $\mathfrak{p}$}} 
W_{\mathfrak{p}}(\mathcal{E}(y/x)) .\]

Let $\mathfrak{p}\nmid \mathfrak{d}$ be a prime at which $\mathcal{E}(y/x)$ has
bad reduction.
Since \[D(x,y) = D_0 \prod_{v} (P_{v}(x,y))^{e_{v,D}}\]
and $D_0|\mathfrak{d}$, we must have $\mathfrak{p}|P_{v}(x,y)$ for
some place $v$ with $e_{v,D}>0$. By the definition (\ref{eq:toolong})
of $\mathfrak{d}$, it follows that $\mathfrak{p}\nmid P_{u}(x,y)$
for every place $u\ne v$ of $K(t)$. Thus
\[\begin{aligned}
W(\mathcal{E}(y/x)) &= -\prod_{\mathfrak{p}|\mathfrak{d}} W_{\mathfrak{p}}(\mathcal{E}(y/x))
\mathop{\prod_{v}}_{e_{v,D}>0} 
 \mathop{\mathop{\prod_{\mathfrak{p}\nmid \mathfrak{d}}}_{\mathfrak{p}|P_{v}(x,y)}}_{
\text{$\mathcal{E}(y/x)$ has bad. red. at $\mathfrak{p}$}}
W_{\mathfrak{p}}(\mathcal{E}(y/x))\\
&= -\prod_{\mathfrak{p}|\mathfrak{d}} W_{\mathfrak{p}}(\mathcal{E}(y/x))
\mathop{\prod_{v}}_{e_{v,D}>0} 
 \mathop{\prod_{\mathfrak{p}\nmid \mathfrak{d}}}_{\mathfrak{p}|P_{v}(x,y)}
W_{\mathfrak{p}}(\mathcal{E}(y/x)) .\end{aligned}\]
By Lemma \ref{lem:mess},
\[\begin{aligned}
\mathop{\prod_{v}}_{e_{v,D}>0} W_{\mathfrak{p}}(\mathcal{E}(y/x)) &=
 \mathop{\prod_{\text{$v$ half-bad}}}_{e_{v,D}>0} g_{v}(x,y) 
V_{\mathfrak{d},v}(\mathcal{E})(x,y) \\
&\cdot \mathop{\prod_{\text{$v$ quite bad}}}_{e_{v,\mathfrak{d}}>0}
 g_{v}(x,y) h(\sq_K(P_{v}(x,y)),x,y) V_{\mathfrak{d},v}(\mathcal{E})(x,y) \\
&= \mathop{\prod_{v}}_{e_{v,D}>0} g_{nu}(x,y)
 \mathop{\prod_{\text{$v$ quite bad}}}_{e_{v,D}>0}
 h(\sq_K(P_{v}(x,y)),x,y) \mathop{\prod_{v}}_{e_{v,D}>0} 
 V_{\mathfrak{d},v}(\mathcal{E})(x,y).
\end{aligned}\]
For every two distinct places $v$, $u$ of $K(t)$ with $e_{v,D}>0$,
$e_{u,D}>0$, we know that 
 \[\gcd(P_{v}(x,y),P_{u}(x,y))|\mathfrak{d}^{\infty},\] and thus
$\gcd(\sq_K(P_{v}(x,y)),\sq_K(P_{v}(x,y)))|\mathfrak{d}^{\infty}$.
By
properties (3) and (4) in the statement of Lemma \ref{lem:mess},
\[\mathop{\prod_{\text{$v$ quite bad}}}_{e_{v,D}>0}
    h(\sq_K(P_{v}(x,y)),x,y) = h(\sq_K(B'(x,y)),x,y) .\]
Since $V_{\mathfrak{d},v}(\mathcal{E})=1$ for $v$ with $e_{v,D}=0$,
\[\mathop{\prod_{v}}_{e_{v,D}>0} V_{\mathfrak{d},v}(\mathcal{E})(x,y) =
\prod_{v} V_{\mathfrak{d},v}(\mathcal{E})(x,y) =
V(\mathcal{E})(x,y) .\]
Hence
\[\mathop{\prod_{v}}_{e_{v,D}>0} W_{\mathfrak{p}}(\mathcal{E}(y/x)) =
\left(\mathop{\prod_{v}}_{e_{v,D}>0} g_{v}(x,y)\right)\cdot
h(B'(x,y),x,y)\cdot V(\mathcal{E}(x,y))\]
and thus
\[W(\mathcal{E}(y/x)) = \left(- \prod_{\mathfrak{p}|\mathfrak{d}}
 W_{\mathfrak{p}}(\mathcal{E}(y/x))
\mathop{\prod_{v}}_{e_{v,D}>0} g_{v}(x,y)\right)\cdot
h(B'(x,y),x,y)\cdot V(\mathcal{E}(x,y))\]
By Lemma \ref{lem:bigone},
\[V(\mathcal{E})(x,y) = f(x,y)\cdot \lambda_K(M_{\mathcal{E}}(x,y)),\]
where $f$ is a pliable function. Therefore,
\[W(\mathcal{E}(y/x)) = 
- f(x,y) \prod_{\mathfrak{p}|\mathfrak{d}}
 W_{\mathfrak{p}}(\mathcal{E}(y/x))
\mathop{\prod_{v}}_{e_{v,D}>0} g_{v}(x,y) \cdot
h(B'(x,y),x,y) \lambda_K(M_{\mathcal{E}}(x,y)) .\]
By Proposition \ref{prop:krone} and Lemma \ref{lem:pingpong}, the map
\[(x,y)\mapsto W_{\mathfrak{p}}(\mathcal{E}(y/x))\]
 is pliable. Hence the map
\[g:(x,y)\mapsto \left(- f(x,y) \cdot \prod_{\mathfrak{p}|\mathfrak{d}}
 W_{\mathfrak{p}}(\mathcal{E}(y/x))
\mathop{\prod_{v}}_{e_{v,D}>0} g_{v}(x,y)\right)
\]
on $\mathbb{A}_{\mathcal{E}}$ is the product of finitely many pliable maps. Therefore,
$g$ is itself pliable. We have obtained
\[W(\mathcal{E}(y/x)) = 
g(x,y) \cdot h(B'(x,y),x,y)\cdot \lambda_K(M_{\mathcal{E}}(x,y)),\]
where $g$ is pliable and $h$ depends only on
$\mathfrak{a}$ and on 
$\{\frac{x \mo \mathfrak{p}}{y \mo \mathfrak{p}}\}_{\mathfrak{p}|\mathfrak{a}} \in 
\prod_{\mathfrak{p}|\mathfrak{a}} \mathbb{P}^1(\mathfrak{O}_K/\mathfrak{p})$.
\end{proof}

\subsection{Averages and correlations}\label{subs:aver}
In order to give explicit estimates for the average of 
$W(\mathcal{E}(y/x))$, we need quantitative versions of Hypotheses
$\mathfrak{B}_1$ and $\mathfrak{B}_2$. 

\begin{conjb1e}
Let $\epsilon(N)\geq 0$, $\eta(N)\leq N$.
The polynomial $P\in \mathfrak{O}_K$ obeys
\[\mathop{\sum_{1\leq x\leq N}}_{x\equiv a \mo m} \lambda_K(P(x))
 \ll \frac{\epsilon(N) N}{m}\]
for every $m\leq \eta(N)$.
\end{conjb1e}

\begin{conjb2e}
Let $\epsilon(N)\geq 0$, $\eta(N)\leq N$.
The homogeneous polynomial $P\in \mathfrak{O}_K\lbrack x,y\rbrack$ obeys
\[\sum_{(x,y)\in S \cap \lbrack -N,N\rbrack^2\cap L} \lambda_K(P(x,y)) \ll
\frac{\epsilon(N) N^2}{\lbrack \mathbb{Z}^2 :L \rbrack}\]
for every sector $S$ and every lattice coset $L$ of index 
$\lbrack \mathbb{Z}^2 : L
\rbrack \leq \eta(N)$.
\end{conjb2e}

We can now prove the results stated in the introduction.
\begin{thm}[$\mathfrak{A}_1(K,B_{\mathcal{E}}'(1,t),\delta(N))$,
$\mathfrak{B}_1(K,M_{\mathcal{E}}(1,t),\eta(N),\epsilon(N))$]\label{thm:mart}
Let $K$ be a number field. Let $\mathcal{E}$ be an elliptic curve
over $K(t)$. Suppose $M_{\mathcal{E}}(1,t)$ is non-constant. Then,
for any integers $a$, $m$, $0<m\leq \eta(N)$,
\begin{equation}\label{eq:colabo}\begin{aligned}
\mathop{\sum_{1\leq x\leq N}}_{x\equiv a \mo m}
W(\mathcal{E}(x)) &\ll 
\left( \frac{\epsilon(N)}{m} + \frac{\epsilon'(N)}{\sqrt{m'}}\right)
N + \delta(N), 
\end{aligned}\end{equation}
where 
\begin{equation}\label{eq:siobhan}\begin{aligned}
\epsilon' &= 
  \sqrt{\max((\log \eta(N))^c/\eta(N),N^{-1/2})} 
    \log(-\max((\log \eta(N))^c/\eta(N),N^{-1/2})),\\
m' &= \min(m, \min(N^{1/2},\eta(N)/(\log \eta(N))^c)) ,\end{aligned}
\end{equation}
and both $c$ and the implied constant in (\ref{eq:colabo})
depend only on $\mathcal{E}$ and the implied
constants in hypotheses $\mathfrak{A}_1$ and $\mathfrak{B}_1$.
\end{thm}
\begin{proof}
Let $\mathbb{A}_{\mathcal{E},\mathbb{Z}} = \{t\in \mathbb{Z} : (1,t)\in 
 \mathbb{A}_{\mathcal{E}}\}$. Let $\mathfrak{M}=1$. 
By Proposition \ref{prop:ahem},
\begin{equation}\label{eq:whkn}W(\mathcal{E}(t)) = 
g(1,t) \cdot h(\sq_K(B_{\mathcal{E}}'(1,t)),1,t) \cdot
\lambda_K(M_{\mathcal{E}}(1,t))\end{equation}
for all $t\in \mathbb{A}_{\mathcal{E},\mathbb{Z}}$, where 
$|g(x,y)|=1$, $|h(\mathfrak{a},x,y)|=1$,
$g$ is pliable and $h(\mathfrak{a},1,t)$ depends only on
$\mathfrak{a}$ and $t \mo \rad(\mathfrak{a})$.
Let $g_0(t) = g(1,t)$, $h_0(\mathfrak{a},t)=h(\mathfrak{a},1,t)$.
By Lemma \ref{lem:pongping}, $g_0$ is affinely pliable. 

By $\mathfrak{B}_1(K,M_{\mathcal{E}}(1,t),\eta(N),\epsilon(N))$
and Lemma \ref{lem:weneed},
\[\mathop{\sum_{1\leq x\leq N}}_{x\equiv a \mo m}
 g_0(t) \lambda_K(M_{\mathcal{E}}(1,t)) \ll
 \left(\frac{\epsilon(N)}{m} +
 \frac{(\log \eta(N))^c}{\eta(N)}\right) N\]
for any $a, m\in \mathbb{Z}$, $0<m\leq N$. Then, by 
$\mathfrak{A}_1(K,B_{\mathcal{E}}'(1,t),\delta(N))$ and
Proposition \ref{prop:sqfra}, 
\[\mathop{\sum_{1\leq x\leq N}}_{x\equiv a \mo m}
 h_0(\mathfrak{a},t) g_0(t) \lambda_K(M_{\mathcal{E}}(1,t))\]
is at most a constant times 
\[\left( \frac{\epsilon(N)}{m} + \frac{\epsilon'(N)}{\sqrt{m'}}\right)
N + \delta(N),\]
where $\epsilon'$ and $m'$ are as in (\ref{eq:siobhan}).
By (\ref{eq:whkn}),
\[W(\mathcal{E}(t)) = 
g_0(t) \cdot h(\sq_K(B_{\mathcal{E}}'(1,t)),t) \cdot
\lambda_K(M_{\mathcal{E}}(1,t))\]
for all $t\in \mathbb{A}_{\mathcal{E},\mathbb{Z}}$. Since there are only finitely
many integers not in $\mathbb{A}_{\mathcal{E},\mathbb{Z}}$, the statement follows.
\end{proof}
\begin{thm}[$\mathfrak{A}_1(K,B_\mathcal{E}'(1,t),\delta(N))$,
$\mathfrak{B}_1(K,M_\mathcal{E}(1,t) M_\mathcal{E}(1,t+k),\eta(N),\epsilon(N))$]\label{thm:martp2}
Let $K$ be a number field.
Let $\mathcal{E}$ be an elliptic curve over $K(t)$.
Let $k$ be a non-zero integer.
 Suppose $M_\mathcal{E}(1,t)$ is not constant. Then,
for any integers $a$, $m$, $0<m\leq \eta(N)$, 
\begin{equation}\label{eq:preso}
\mathop{\sum_{1\leq x\leq N}}_{x\equiv a \mo m}
W(\mathcal{E}(x)) W(\mathcal{E}(x+k)) 
\ll \left( \frac{\epsilon(N)}{m} + \frac{\epsilon'(N)}{\sqrt{m'}}\right)
N + \delta(N), 
\end{equation}
where 
\begin{equation}\label{eq:agamn}\begin{aligned}
\epsilon' &= 
  \sqrt{\max((\log \eta(N))^c/\eta(N),N^{-1/2})} 
    \log(-\max((\log \eta(N))^c/\eta(N),N^{-1/2})),\\
m' &= \min(m, \min(N^{1/2},\eta(N)/(\log \eta(N))^c)) ,\end{aligned}
\end{equation}
and both $c$ and the implied constant in (\ref{eq:preso})
depend only on $\mathcal{E}$ and the implied
constants in hypotheses $\mathfrak{A}_1$ and $\mathfrak{B}_1$.
\end{thm}
\begin{proof}
Let $\mathbb{A}_{\mathcal{E},\mathbb{Z}} = \{t\in \mathbb{Z} : (1,t)\in 
 \mathbb{A}_{\mathcal{E}}\}$. 
Let $\mathfrak{M}=\mathfrak{m}_{B_{\mathcal{E}}'(1,t) 
B_{\mathcal{E}}'(1,t+k)}$,
where $\mathfrak{m}$ is as in Corollary \ref{cor:multsquo}.
 By Proposition \ref{prop:ahem}, $W(\mathcal{E}(t))$ equals
\[\begin{aligned}
g(1,t)\,&h(\sq_K(B_{\mathcal{E}}'(1,t)),1,t) 
g(1,t+k) \\ &h(\sq_K(B_{\mathcal{E}}'(1,t+k)),1,t+k) 
\lambda_K(M_{\mathcal{E}}(1,t) M_{\mathcal{E}}(1,t+k))\end{aligned}\]
for all $t\in \mathbb{A}_{\mathcal{E},\mathbb{Z}}$, where 
$|g(x,y)|=1$, $|h(\mathfrak{a},x,y)|=1$,
$g$ is pliable and $h(\mathfrak{a},1,t)$ depends only on
$\mathfrak{a}$ and $t \mo \rad(\mathfrak{a})$.
Let $g_0(t) = g(1,t) g(1,t + k)$, 
\begin{equation}\label{eq:agha}
h_0(t)= h(\sq_K(B_{\mathcal{E}}'(1,t)),1,t) 
h(\sq_K(B_{\mathcal{E}}'(1,t+k)),1,t+k),\end{equation}
By Lemma \ref{lem:pongping}, $g(1,t)$ and $g(1,t+k)$ are affinely pliable,
and hence so is $g_0(t)$. By Lemma \ref{cor:multsquo},
(\ref{eq:agha}) depends only on 
\[\sq_K(\Lcm(B_{\mathcal{E}}'(1,t),B_{\mathcal{E}}'(1,t+k))(x))/
\gcd(\sq_K(\Lcm(B_{\mathcal{E}}'(1,t),B_{\mathcal{E}}'(1,t+k)
)(x)),\mathfrak{M}^{\infty})\]
and on
$x \mo \mathfrak{p}$ for $\mathfrak{p}|\sq_K(\Lcm(
B_{\mathcal{E}}'(1,t),B_{\mathcal{E}}'(1,t+k))(x))$,
$\mathfrak{p}\nmid \mathfrak{M}$.

The remainder of the proof is as for Theorem \ref{thm:mart}. Notice that,
by Lemma \ref{lem:repla},
$\mathfrak{A}_1(K,B_\mathcal{E}'(1,t),\delta(N))$ implies
 $\mathfrak{A}_1(K,B_\mathcal{E}'(1,t+k),\delta(N))$ and thus, by Lemma
\ref{lem:multconj1}, it implies
 $\mathfrak{A}_1(K,B_\mathcal{E}'(1,t) B_\mathcal{E}'(1,t+k),\delta(N))$
as well.
\end{proof}
\begin{thm}[$\mathfrak{A}_1(K,B_\mathcal{E}'(1,t),\delta(N))$]\label{thm:martp3}
Let $K$ be a number field.
Let $\mathcal{E}$ be an elliptic curve over $K(t)$.
Let $c$ be an integer other than zero.
 Suppose $M_\mathcal{E}(1,t)$ is not constant. If
\[
\mathop{\sum_{1\leq x\leq N}}_{x\equiv a \mo m} W(\mathcal{E}(x)) \ll 
\frac{\epsilon(N) N}{m}
\]
for any integers $a$, $m$, $0<m\leq \eta(N)$, 
then
\begin{equation}\label{eq:hunger}
\mathop{\sum_{1\leq x\leq N}}_{x\equiv a \mo m} \lambda_K(P(x))
\ll \left( \frac{\epsilon(N)}{m} + \frac{\epsilon'(N)}{\sqrt{m'}}\right)
N + \delta(N)
\end{equation}
for any integers $a$, $m$, $0<m\leq \eta(N)$,
where 
\begin{equation}\begin{aligned}
\epsilon' &= 
  \sqrt{\max((\log \eta(N))^c/\eta(N),N^{-1/2})} 
    \log(-\max((\log \eta(N))^c/\eta(N),N^{-1/2})),\\
m' &= \min(m, \min(N^{1/2},\eta(N)/(\log \eta(N))^c)) ,\end{aligned}
\end{equation}
and both $c$ and the implied constant in (\ref{eq:hunger})
depend only on $\mathcal{E}$ and the implied
constant in hypothesis $\mathfrak{A}_1$.
\end{thm}
\begin{proof}
Since $|g(x,y)|=|h(\mathfrak{a},x,y)|=1$ for any $\mathfrak{a}$, $x$, $y$,
we can rewrite (\ref{eq:whkn}) as
\[\lambda_K(M_{\mathcal{E}}(1,t)) = 
g(1,t) \cdot h(\sq_K(B_{\mathcal{E}}'(1,t)),1,t) W(\mathcal{E}(t)) .\]
The rest is as in the proof of Theorem \ref{thm:mart}.
\end{proof}

\begin{thm}[$\mathfrak{A}_2(K,B_{\mathcal{E}}',\delta(N))$,
$\mathfrak{B}_2(K,M_{\mathcal{E}},\eta(N),\epsilon(N))$]\label{thm:mart2}
Let $K$ be a number field. Let $\mathcal{E}$ be an elliptic curve
over $K(t)$. Suppose $M_{\mathcal{E}}$ is non-constant. Then,
for every sector $S$ and every lattice coset $L$
of index $\lbrack \mathbb{Z}^2 : L\rbrack \leq \eta(N)$,
\begin{equation}\label{eq:clb}
\mathop{\sum_{(x,y)\in S\cap \lbrack -N, N\rbrack^2 \cap L}}_{\gcd(x,y)=1}
 W(\mathcal{E}(y/x)) \ll
\left( \frac{\epsilon(N)}{\lbrack \mathbb{Z}^2 : L\rbrack} 
+ \frac{\epsilon'(N)}{\sqrt{m'}}\right)
N^2 + \delta(N),\end{equation}
where 
\begin{equation}\label{eq:gherk}\begin{aligned}
\epsilon' &= 
  \sqrt{\max((\log \eta(N))^c/\eta(N),N^{-1/2})} 
    \log(-\max((\log \eta(N))^c/\eta(N),N^{-1/2})),\\
m' &= \min(\lbrack \mathbb{Z}^2 : L \rbrack, 
\min(N^{1/2},\eta(N)/(\log \eta(N))^c)) ,\end{aligned}
\end{equation}
and both $c$ and the implied constant in (\ref{eq:clb})
depend only on $\mathcal{E}$ and the implied
constants in hypotheses $\mathfrak{A}_2$ and $\mathfrak{B}_2$.
\end{thm}
\begin{proof}
By Proposition \ref{prop:ahem},
\begin{equation}\label{eq:whno}
W(\mathcal{E}(y/x)) = g(x,y)\cdot h(\sq_K(B_\mathcal{E}'(x,y)),x,y)\cdot 
\lambda_K(M_\mathcal{E}(x,y)),\end{equation}
for all $(x,y)\in \mathbb{A}_{\mathcal{E}}$,
where $g:\mathbb{A}_{\mathcal{E}}\to \{-1,1\}$,
$h:I_K\times \mathbb{A}_{\mathcal{E}}\to \{-1,1\}$ are such that
\begin{itemize}
\item $g$ is pliable, 
\item 
$h(\mathfrak{a},x,y)$ depends only on $\mathfrak{a}$ and on 
$\{\frac{x \mo \mathfrak{p}}{y \mo \mathfrak{p}}\}_{\mathfrak{p}|\mathfrak{a}} \in 
\prod_{\mathfrak{p}|\mathfrak{a}} \mathbb{P}^1(\mathfrak{O}_K/\mathfrak{p})$.
\end{itemize}
By $\mathfrak{B}_2(K,M_{\mathcal{E}},\eta(N),\epsilon(N))$ and Lemma 
\ref{lem:minmat},
\[\mathop{\sum_{(x,y)\in S \cap \lbrack -N,N\rbrack^2\cap L}}_{\gcd(x,y)=1}
 \lambda(P_K(x,y)) \ll
\max\left(\epsilon(N) ,
\frac{\sqrt{\log \,\lbrack \mathbb{Z}^2 : L\rbrack}}{\eta(N)}\right)
\frac{N^2}{\lbrack \mathbb{Z}^2 :L \rbrack} 
\]
for every lattice $L$ of index $\lbrack \mathbb{Z}^2 : L\rbrack \leq N$.
We can now apply Proposition \ref{prop:pliblock} with
$\epsilon_N = \max(\epsilon(N),\sqrt{\log \eta(N)}/\eta(N))$,
$\eta_N = \eta(N)$, obtaining
\[
\mathop{\sum_{(x,y)\in S\cap \lbrack -N, N\rbrack^2 \cap L}}_{\gcd(x,y)=1}
g(x,y) \lambda_K(M_\mathcal{E}(x,y)) =
\left(\frac{\epsilon(N)}{\phi(\lbrack \mathbb{Z}^2 : L\rbrack)} +
 \frac{(\log \eta_N)^c}{\eta_N}\right) N^2 \]
for any sector $S$ and any lattice $L$. Then, by 
$\mathfrak{A}_2(K,B_{\mathcal{E}}',\delta(N))$ and
Proposition \ref{prop:volver}, the absolute value of
\[
\mathop{\sum_{(x,y)\in S\cap \lbrack -N, N\rbrack^2 \cap L}}_{\gcd(x,y)=1}
h(\sq_K(B_{\mathcal{E}}'(x,y)),x,y)
g(x,y) \lambda_K(M_\mathcal{E}(x,y))\]
is at most a constant times
\[
\mathop{\sum_{(x,y)\in S\cap \lbrack -N, N\rbrack^2 \cap L}}_{\gcd(x,y)=1}
 W(\mathcal{E}(y/x)) \ll
\left( \frac{\epsilon(N)}{\lbrack \mathbb{Z}^2 : L\rbrack} 
+ \frac{\epsilon'(N)}{\sqrt{m'}}\right)
N^2 + \delta(N),\]
where $\epsilon'$ and $m'$ are as in (\ref{eq:gherk}).
Since the set 
$\{(x,y)\in \mathbb{Z}^2 : \text{$x$, $y$ coprime}\} - \mathbb{A}_{\mathcal{E}}$
is finite, the statement follows by (\ref{eq:whno}).
\end{proof}
\begin{thm}[$\mathfrak{A}_2(K\!,B_{\mathcal{E}}',\delta(N))$,
$\mathfrak{B}_2(K\!,M_{\mathcal{E}}(t_0,t_1)
M_{\mathcal{E}}(k_0 x, k_0 y + k_1 x),\eta(N),\epsilon(N))$]
\label{thm:mart2p2}
Let $K$ be a number field. Let $\mathcal{E}$ be an elliptic curve
over $K(t)$. Suppose $M_{\mathcal{E}}$ is non-constant. 
Let $k = k_1/k_0$ be anon-zero rational number, $\gcd(k_0,k_1)=1$.
Then,
for every sector $S$ and every lattice coset $L$
of index $\lbrack \mathbb{Z}^2 : L\rbrack \leq \eta(N)$,
\begin{equation}\label{eq:caleb}
\mathop{\sum_{(x,y)\in S\cap \lbrack -N, N\rbrack^2 \cap L}}_{\gcd(x,y)=1}
 W(\mathcal{E}(y/x)) W(\mathcal{E}(y/x + k)) \ll
\left( \frac{\epsilon(N)}{\lbrack \mathbb{Z}^2 : L\rbrack} 
+ \frac{\epsilon'(N)}{\sqrt{m'}}\right)
N^2 + \delta(c' N),\end{equation}
where 
\[\begin{aligned}
\epsilon' &= 
  \sqrt{\max((\log \eta(N))^c/\eta(N),N^{-1/2})} 
    \log(-\max((\log \eta(N))^c/\eta(N),N^{-1/2})),\\
m' &= \min(\lbrack \mathbb{Z}^2 : L \rbrack, 
\min(N^{1/2},\eta(N)/(\log \eta(N))^c)) ,\end{aligned}\]
and $c$, $c'$ and the implied constant in (\ref{eq:caleb})
depend only on $\mathcal{E}$ and the implied
constants in hypotheses $\mathfrak{A}_2$ and $\mathfrak{B}_2$.
\end{thm}
\begin{proof}
Let $\mathbb{A}_{\mathcal{E},k} =\{(x,y)\in \mathbb{A}_{\mathcal{E}}
: \left(\frac{k_0 x}{\gcd(k_0 x, k_0 y + k_1 x)},
\frac{k_0 y + k_1 x}{\gcd(k_0 x, k_0 y + k_1 x)}\right)\in A_{\mathcal{E}}\}$. Since
$\frac{k_0 y + k_1 x}{k_0 x} = \frac{y}{x} + k$, we can write
$\mathbb{A}_{\mathcal{E},k}$ in full as the set of all coprime
$x,y\in \mathfrak{O}_K$ such that
\[\begin{aligned}
x\ne 0,\;&c_4(y/x)\ne \infty,\: c_6(y/x)\ne \infty,\: 
\Delta(y/x) \ne 0,\infty,\, q_0(y/x)\ne 0,\\
&c_4(y/x+k)\ne \infty, c_6(y/x+k)\ne \infty, 
\Delta(y/x+k) \ne 0,\infty,\, q_0(y/x+k)\ne 0.\end{aligned}\]
Hence $\mathbb{A}-\mathbb{A}_{\mathcal{E},k}$ is a finite set.

Let $F_1(x,y) = k_0 x$, $F_2(x,y) = k_0 y + k_1 x$. For $x$, $y$
coprime, $\gcd(k_0 x, k_0 y + k_1 x)$ must
divide $k_0^2$. 
Let 
\[\mathfrak{M} = k_0 \mathfrak{m}_{B_{\mathcal{E}}',
B_{\mathcal{E}}'(F_1(x,y),F_2(x,y))},\]
where $\mathfrak{m}_{\cdot}$ is as in Lemma \ref{lem:multsquid}. 
Let $g$, $h$ be as in Proposition
\ref{prop:ahem}. Then
\[
W(\mathcal{E}(y/x)) W(\mathcal{E}(y/x + k))\] equals
\[g_1(x,y) \cdot h_1(\sq_K(B_{\mathcal{E}}'(x,y)),x,y) \cdot
\lambda_K(M_{\mathcal{E}}(x,y) M_{\mathcal{E}}(F_1(x,y),F_2(x,y))),\]
for $(x,y)\in \mathbb{A}_{\mathcal{E},k}$, where 
\[\begin{aligned}
g_0(x,y) &= 
g\left(\frac{F_1(x,y)}{\gcd(F_1(x,y),F_2(x,y))},
\frac{F_2(x,y)}{\gcd(F_1(x,y),F_2(x,y))}\right),\\
g_1(x,y) &= g(x,y)\cdot \lambda_K(\gcd(F_1(x,y),F_2(x,y)))^{\deg M_{\mathcal{E}}} g_0(x,y),\\
h_1(x,y) &= h(\sq_K(B_{\mathcal{E}}'(x,y)),x,y) \cdot h_0(x,y) ,\end{aligned}\]
and $h_0(x,y)$ equals
\[h\left(\sq_K\left(
B_\mathcal{E}'\left(\frac{F_1(x,y)}{\gcd(F_1(x,y),F_2(x,y))},
\frac{F_2(x,y)}{\gcd(F_1(x,y),F_2(x,y))}\right)\right),F_1(x,y),
F_2(x,y)\right).\]

By Lemma \ref{lem:oneoflst},
\[(x,y)\mapsto g\left(\frac{x}{\gcd(x,y,k_0^2)},
  \frac{y}{\gcd(x,y,k_0^2)}\right)\]
is a pliable function on $S' = \{(x,y)\in \mathbb{Z}^2 :
 (x/\gcd(x,y,k_0^2),y/\gcd(x,y,k_0^2))\in \mathbb{A}_{\mathcal{E}}\}$.
Then, by Proposition \ref{prop:plisam}, $g_0$
is a pliable function on
\[\{(x,y)\in \mathbb{Z}^2 : \text{$x$, $y$ coprime,}
\left(\frac{F_1(x,y)}{\gcd(F_1(x,y),F_2(x,y))},
\frac{F_2(x,y)}{\gcd(F_1(x,y),F_2(x,y))}\right)
\in \mathbb{A}_{\mathcal{E}}\},\] which is a subset of
$\mathbb{A}_{\mathcal{E},k}$.
Since $\gcd(F_1(x,y),F_2(x,y))|k^{\infty}$ for $x$, $y$ coprime, 
the map \[(x,y)\to \lambda_K(\gcd(F_1(x,y),F_2(x,y)))\] on 
$\mathbb{A}_{\mathcal{E},k}$ is pliable.
Hence $g_1(x,y) = g(x,y)\cdot \lambda_K(\gcd(F_1(x,y),F_2(x,y)))^{\deg M_{\mathcal{E}}} g_0(x,y)$ is pliable.

By Proposition \ref{prop:ahem}, (2), (3) and (4),
$h(\sq_K(B_\mathcal{E}'(x,y)),x,y)$ depends only on
\[\sq_K(B_\mathcal{E}'(x,y))/\gcd(\sq_K(B_\mathcal{E}'(x,y)),\mathfrak{M}^{\infty})\]
and on
$\frac{x \mo \mathfrak{p}}{y \mo \mathfrak{p}}$
for $\mathfrak{p}|\sq_K(B_\mathcal{E}'(x,y))$,
$\mathfrak{p}\nmid \mathfrak{M}$. Hence $h_0(x,y)$ depends only
on \begin{equation}\label{eq:whtocit}
\sq_K(B_\mathcal{E}'(F_1(x,y),F_2(x,y)))/
\gcd(\sq_K(B_\mathcal{E}'(F_1(x,y),F_2(x,y))),\mathfrak{M}^{\infty})
\end{equation} and on
 \[\frac{F_1(x,y)/\gcd(F_1(x,y),F_2(x,y)) \mo \mathfrak{p}}{
F_2(x,y)/\gcd(F_1(x,y),F_2(x,y)) \mo \mathfrak{p}}\]
for \[\mathfrak{p}|\sq_K
(B_\mathcal{E}'(F_1(x,y)/\gcd(F_1(x,y),F_2(x,y)),
F_2(x,y)/\gcd(F_1(x,y),F_2(x,y)))).\;\;\mathfrak{p}\nmid \mathfrak{M} .\] 
Since $\gcd(F_1(x,y),F_2(x,y))|k_0^2$ and $k_0|\mathfrak{M}$,
\[\frac{F_1(x,y)/\gcd(F_1(x,y),F_2(x,y)) \mo \mathfrak{p}}{
F_2(x,y)/\gcd(F_1(x,y),F_2(x,y)) \mo \mathfrak{p}} = 
\frac{F_1(x,y) \mo \mathfrak{p}}{F_2(x,y) \mo \mathfrak{p}}\]
for all $x$, $y$ coprime, $\mathfrak{p}\nmid \mathfrak{M}$.
In turn, since $F_2(x,y)/F_1(x,y) = y/x+k_0/k_1 =y/x+k$, 
\[\frac{F_1(x,y) \mo \mathfrak{p}}{F_2(x,y) \mo \mathfrak{p}} =
\left(\frac{y}{x} + k\right)^{-1} \mo \mathfrak{p} \]
for all $x$, $y$ coprime, $\mathfrak{p}\nmid \mathfrak{M}$. Since $k$ is fixed,
$\left(\frac{y}{x} + k\right)^{-1} \mo \mathfrak{p}$ depends
only on $\frac{y \mo \mathfrak{p}}{x \mo \mathfrak{p}}$. Thus
\[h\left(\sq_K\left(B_\mathcal{E}'\left(
\frac{F_1(x,y)}{\gcd(F_1(x,y),F_2(x,y))},
\frac{F_2(x,y)}{\gcd(F_1(x,y),F_2(x,y))}\right)\right),F_1(x,y),F_2(x,y)\right)\] depends only on
(\ref{eq:whtocit}) and
on $\frac{x \mo \mathfrak{p}}{y \mo \mathfrak{p}}$ for
$\mathfrak{p}|\sq_K(B_\mathcal{E}'(F_1(x,y),F_2(x,y)))$,
$\mathfrak{p}\nmid \mathfrak{M}$. By Lemma \ref{lem:multsquid}, 
it follows that $h_1$ depends only on 
\[\frac{\sq_K(P(x,y))}{\gcd(\sq_K(P(x,y)),\mathfrak{M}^{\infty})},\;\;\;
\text{$\frac{x \mo \mathfrak{p}}{y \mo \mathfrak{p}}$\, for\,
$\mathfrak{p}|\sq_K(P(x,y))$,}\] where
 $P = \Lcm(B_{\mathcal{E}}'(x,y),B_{\mathcal{E}}'(F_1(x,y),F_2(x,y))))$.
It remains to show the fact that $\mathfrak{A}_2(K,P,\delta(c' N))$ holds
for some $c'$ depending only on the implied constant in $\mathfrak{A}_2$.
This follows immediately from $\mathfrak{A}_2(K,B_{\mathcal{E}}',\delta(N))$ and Lemmas 
\ref{lem:replb} and \ref{lem:multconj2}.
\end{proof}

\begin{thm}[$\mathfrak{A}_2(K,B_{\mathcal{E}}',\delta(N))$]
\label{thm:mart2p3}
Let $K$ be a number field. Let $\mathcal{E}$ be an elliptic curve
over $K(t)$. Suppose $M_{\mathcal{E}}$ is non-constant. If
for every sector $S$ and every lattice coset $L$ of index
$\lbrack \mathbb{Z}^2 : L\rbrack \leq \eta(N)$,
\[
\sum_{(x,y)\in S\cap \lbrack -N, N\rbrack^2 \cap L}
 W(\mathcal{E}(y/x)) \ll \frac{\epsilon(N) N^2}{\lbrack \mathbb{Z}^2 : L\rbrack}
\]
then, for every sector $S$ and every lattice coset $L$ of index
$\lbrack \mathbb{Z}^2 : L \rbrack \leq \eta(N)$,
\begin{equation}\label{eq:gherkin}
\mathop{\sum_{(x,y)\in S\cap \lbrack -N, N\rbrack^2 \cap L}
}_{\gcd(x,y)=1} \lambda_K(P(x,y))
\ll
\left( \frac{\epsilon(N)}{\lbrack \mathbb{Z}^2 : L\rbrack} 
+ \frac{\epsilon'(N)}{\sqrt{m'}}\right)
N^2 + \delta(N),\end{equation}
where 
\[\begin{aligned}
\epsilon' &= 
  \sqrt{\max((\log \eta(N))^c/\eta(N),N^{-1/2})} 
    \log(-\max((\log \eta(N))^c/\eta(N),N^{-1/2})),\\
m' &= \min(\lbrack \mathbb{Z}^2 : L \rbrack, 
\min(N^{1/2},\eta(N)/(\log \eta(N))^c)) ,\end{aligned}\]
and both $c$ and the implied constant in (\ref{eq:gherkin})
depend only on $\mathcal{E}$ and the implied
constants in hypotheses $\mathfrak{A}_2$ and $\mathfrak{B}_2$.
\end{thm}
\begin{proof}
By Proposition \ref{prop:ahem},
\[\lambda_K(M_\mathcal{E}(x,y)) = 
g(x,y)\cdot h(\sq_K(B_\mathcal{E}'(x,y)),x,y)\cdot 
W(\mathcal{E}(y/x)),\]
for all $(x,y)\in \mathbb{A}_{\mathcal{E}}$,
where $g:\mathbb{A}_{\mathcal{E}}\to \{-1,1\}$,
$h:I_K\times \mathbb{A}_{\mathcal{E}}\to \{-1,1\}$ are such that
\begin{itemize}
\item $g$ is pliable, 
\item 
$h(\mathfrak{a},x,y)$ depends only on $\mathfrak{a}$ and on 
$\{\frac{x \mo \mathfrak{p}}{y \mo \mathfrak{p}}\}_{\mathfrak{p}|\mathfrak{a}} \in 
\prod_{\mathfrak{p}|\mathfrak{a}} \mathbb{P}^1(\mathfrak{O}_K/\mathfrak{p})$.
\end{itemize}
Proceed as in the proof of Theorem \ref{thm:mart2}.
\end{proof}
Theorems 1.1', 1.3' and 1.4' follow immediately from 
Theorems
\ref{thm:mart}, \ref{thm:mart2} and \ref{thm:mart2p2}, respectively,
and from
the known cases of $\mathfrak{A}_i$ and $\mathfrak{B}_i$ listed
in Appendix \ref{sec:fapp}.
In order to obtain Theorems 1.1--1.4 and Propositions \ref{prop:forgot1},
\ref{prop:forgot2} from Theorems \ref{thm:mart}--\ref{thm:mart2p3}, it
is enough to show that Conjecture $\mathfrak{A}_i(K,P)$ and
Hypothesis $\mathfrak{B}_i(K,P)$, 
as stated in subsection \ref{subs:whenneedit}, imply
$\mathfrak{A}_i(K,P,\delta(N))$ and 
$\mathfrak{B}_i(K,P,\eta(N),\epsilon(N))$, respectively, for
some $\delta(N)$, $\eta(N)$, $\epsilon(N)$ satisfying $\delta(N)=o(N)$,
$\lim_{N\to \infty} \eta(N) = N$, $\epsilon(N) = o(N)$.

The case of $\mathfrak{A}_i$ is clear: since $\mathfrak{A}_1(K,P)$
states that
\[\lim_{N\to \infty} \frac{1}{N} \,\#\{1\leq x\leq N:
\exists \mathfrak{p} \text{\,\st\,} \rho(\mathfrak{p})>N^{1/2}, \mathfrak{p}^2|P(x)\} = 0,\]
we can take
\[\delta(N) = 
\#\{1\leq x\leq N:
\exists \mathfrak{p} \text{\,\st\,} \rho(\mathfrak{p})>N^{1/2}, \mathfrak{p}^2|P(x)\}\]
and thus obtain $\mathfrak{A}_1(K,P,\delta(N))$;
the same works for $\mathfrak{A}_2$. Now assume
that $\mathfrak{B}_1(K,P)$ holds, i.e.,
\[\lim_{N\to \infty} \frac{1}{N} \mathop{\sum_{1\leq n\leq N}}_{n\equiv
a \mo m} \lambda_K(P(n)) = 0\]
for any $a\geq 0$, $m>0$. For every $n\geq 1$, let $A(n)$ be the
smallest positive integer such that, for every $1\leq m\leq n$,
$0\leq a<n$,
\[\frac{1}{N} \mathop{\sum_{1\leq n\leq N}}_{n\equiv a \mo m}
 \lambda_K(P(n)) < \frac{1}{m\cdot n}\]
for all $N\geq A(n)$. Set $A(0)=0$.
For $x\geq 1$, let $B(x)$ be the largest non-negative integer
$n$ such that  $A(n)\leq x$. For every $n>1$, $B(x)>n$ for
all $x\geq A(n)$. Hence $\lim_{x\to \infty} B(x) = \infty$.
Set $\eta(N) = B(N)$,
$\epsilon(N) = \frac{1}{B(N)}$. Then
$\mathfrak{B}_1(K,P,\eta(N),\epsilon(N))$ holds. The same
argument
is valid for $\mathfrak{B}_2$.
\section{Examples}
\subsection{Specimens and how to find them}
Let $K$ be a number field. For any $j\in K(t)$ other than
$j=0$, $j=1728$, the curve given by the equation
\[y^2 = x^3 - \frac{c_4}{48} x - \frac{c_6}{864},\]
\[c_4 := j (j-1728),\;\;\;c_6:=j (j-1728)^2\]
is an elliptic curve over $K(t)$ with $j$-invariant equal to $j$.
Any two elliptic curves $\mathcal{E}$, 
$\mathcal{E}'$ over $K(t)$ with the same
$j$-invariant $j(\mathcal{E}) = j(\mathcal{E}')\ne 0,1728$ 
must be quadratic twists of
each other. Therefore, every elliptic curve $\mathcal{E}$ over $K(t)$ with 
$j$-invariant $j\ne 0,1728$ is given by
\begin{equation}\label{eq:stern}
c_4 = d^2 j (j-1728),\;\;\;c_6=d^3 j (j-1728)^2\end{equation}
for some $d\in (K(t))^*$. Write $t=y/x$. Then the places of potentially
multiplicative reduction of $\mathcal{E}$ 
are given by the factors in the denominator
of $j(y/x)$, where $j(y/x)$ is written as a fraction whose numerator
and denominator have no common factors. The set of places of multiplicative
reduction of $\mathcal{E}$ is, of course, a subset of the set of places of potentially multiplicative reduction. We can choose which subset it is
by adjusting $d$ accordingly.

Thus we can easily find infinitely many elliptic curves
$\mathcal{E}$ over $K(t)$ having $M_{\mathcal{E}}(x,y)$ equal to a given
square-free homogeneous polynomial. (See (\ref{eq:M}) for the definition of
$M_{\mathcal{E}}(x,y)$.) 
Say, for example, that you wish $M_{\mathcal{E}}(x,y)$ to be $y$.
The set of potentially multiplicative places will have to include
the place of $K(t)$ given by $y$. For simplicity's sake, let us
require the set to have that place as its only element. Then $j$ will
have to be a non-constant polynomial on $t^{-1}$. In order for $y$
to give a place of multiplicative reduction over $K(t)$, and not one
of merely potential multiplicative reduction, $v_t(d)$ must be even 
if the degree of $j$ as a polynomial on $t^{-1}$ is even, and odd 
if the degree of $j$ is odd.
These conditions on $d$ and $j$ are sufficient. Thus, e.g., the families
given by
\begin{equation}\label{eq:grat1}\begin{aligned}
j &= t^{-1},\; d = t,\\ c_4 &= t^2 \cdot t^{-1} (t^{-1} - 1728) 
= 1 - 1728 t,\; c_6 = t^3\cdot t^{-1} (t^{-1} - 1728)^2 = (1 - 1728 t)^2,\\
j &= t^{-2},\; d = 1,\\ c_4 &= t^{-2} (t^{-2} - 1728),\;
c_6 = t^{-2} (t^{-2} - 1728)^2,\\
j &= t^{-4} - 3,\; d=(t+1),\\ 
c_4 &= (t+1)^2
 ( t^{-4} - 3) ( t^{-4} - 1731),\;
 c_6 = (t+1)^3 ( t^{-4} - 3) ( t^{-4} - 1731)^2,
\end{aligned} \end{equation}
all have $M_{\mathcal{E}}(x,y) = y$. Note that
$\deg_{\irr} B_{\mathcal{E}}'(x,y)\leq 3$ for all three families in
(\ref{eq:grat1}). Hence Theorems 1.1', 1.3' and 
1.4' can be applied: for any of the families 
in (\ref{eq:grat1}), 
$W(\mathcal{E}(t))$ averages to zero
over the integers and over the rationals; furthermore, $W(\mathcal{E}(t))$
is white noise over the rationals.

In detail, the general procedure for finding all curves $\mathcal{E}$
with $M_{\mathcal{E}}(x,y) = P(x,y)$, $P$ square-free, is as follows. Let 
$P = P_1 \cdots P_2 \dotsb P_n$, $P_i$ irreducible, $P_i\ne P_j$. Suppose
$P_i\ne x$ for all $i$. Let $Q_i(t)$ be the polynomial on $t$ such that
$P_i(y/x) = Q_i(y/x) \cdot x^{\deg Q_i}$. Choose any positive integers
$k_1,\dotsb,k_n$ and four polynomials $R_1(t)$, $R_2(t)$, $R_3(t)$,
$R_4(t)$ coprime to $Q_1(t),\dotsb ,Q_n(t)$; assume that $R_1$ 
is square-free, that
$R_1$, $R_2$, $R_3$ are pairwise coprime, that
$R_4$ is prime to $R_1$ and $R_2$, and that
$\deg R_3 \leq \sum_i k_i \deg Q_i$. Let $R_5$ be the
product of the irreducible factors of $R_2$. Then
\begin{equation}\label{eq:anoteq}
j = \frac{R_3(t)}{R_1(t) R_2(t)^2 \prod_i Q_i(t)^{k_i}},\;\;
 d = R_4(t) R_5(t) \prod_i Q_i(t)^{k_i}\end{equation}
give us an elliptic curve with $\mathcal{E}$ with 
$M_{\mathcal{E}}(x,y) = P(x,y)$; furthermore, any such
curve can be expressed as in (\ref{eq:anoteq}). If 
$P = x\cdot P_1 \cdot P_2 \dotsb P_n$, proceed as above, but require
$\deg R > \sum_i k_i \deg Q_i$. 

The degree $\deg_{\irr} B_{\mathcal{E}}'(x,y)$ of the largest
irr. factor of the polynomial $B_{\mathcal{E}}'(x,y)$ coming from
(\ref{eq:anoteq}) is equal to the largest of
\begin{equation}\label{eq:superf}
\deg_{\irr} P, \deg_{\irr} R_1, \deg_{\irr} R_2,
\deg_{\irr} R_3, 
\deg_{\irr} (R_3 - 1728\cdot R_1 R_2^2 \prod_i Q_i^{k_i})
\end{equation}
or to 1, should all the expressions in (\ref{eq:anoteq}) be zero.
The degree $\deg_{\irr} B_{\mathcal{E}}'(1,t)$ is equal to
(\ref{eq:superf}). Since we need only know the degrees of $M_{\mathcal{E}}$
and $B_{\mathcal{E}}$ to know whether our results hold conditionally
or unconditionally, we see that we have an explicit description of all
families for which our results hold unconditionally. It only remains
to see a few more examples that may not be quite trivial to find.

Take, for instance, the issue of semisimplicity.
Constructing families with \[\deg(M_{\mathcal{E}}(x,y))\leq 3\] and
$c_4$, $c_6$ coprime
is a cumbersome but feasible matter. The
following are a few characteristic specimina:
\[\begin{array}{lll}
c_4 = 1 + \frac{8}{3} t + t^2, &c_6 = 1 + \frac{25}{6} t + 4 t^2 + t^3,
&M_{\mathcal{E}} = (12 x + 5 y) (3 x + 8 y) y,\\
c_4 = 2 + 4 t + t^2,  &c_6 = 1 + 9 t + 6 t^2 + t^3,
&M_{\mathcal{E}} = (7 x + 2 y) ( x^2 + 4 x y + y^2),\\
c_4 = 2 - 4 t + t^2, &c_6 = 3 + 9 t - 6 t^2 + t^3,
&M_{\mathcal{E}} = x^3 + 102 x^2 y - 63 x y^2 + 10 y^3,\\
c_4=4, &c_6 = 11 + t, &M_{\mathcal{E}} = x (3 x + y) (19 x + y),\\
c_4=3, &c_6 = 2 + 7 t, &M_{\mathcal{E}} = x (-23 x^2 + 28 x y + 49 y^2),\\
c_4 = 1+t, &c_6=-1 + 3t, &M_{\mathcal{E}} = x y (-3 x + y),\\
c_4=-2 + 6 t + t^2,\! &c_6 = -\frac{45}{2} + \frac{21}{2} t + 9 t^2 + t^3,\!
&M_{\mathcal{E}} = -\frac{2057}{4} x^3 + \frac{1089}{2} x^2 y + 
\frac{363}{4} x y^2,\\
c_4=(t+1)(t+3),\! &c_6=(13 x^2+12 x y +3 y^2),\! &M_{\mathcal{E}} = x (13 x^2 + 12 x y + 3 y^2) .\end{array}\]
Note that none of these families is strictly speaking semistable, since
they all have additive reduction at the place $\den - \num$ corresponding
to $x$. 

Thanks to (\ref{eq:anoteq}), it is a simple matter to construct
a family $\mathcal{E}$ such
that $M_{\mathcal{E}}(x,y)$ equals the homogeneous polynomial of degree
three for which the parity problem was first treated \cite{HB}:
\[c_4 = 1 - 1728 (t^3+1),\;
c_6 = (1 - 1728 (t^3+1))^2,\; M_{\mathcal{E}}(x,y)= x^3 + 2 y^3.\]

We may conclude by seeing two families $\mathcal{E}$ over $K(t)$,
$K$ a number field other than $\mathbb{Q}$, for which our results
are unconditional. (See Appendix 
\ref{sec:app}.) 
\[K = \mathbb{Q}(\sqrt{5}),\; 
c_4 = (1 - 1728 (t + \sqrt{5})),\;
c_6 = (1 - 1728 (t + \sqrt{5}))^2,\;
M_{\mathcal{E}} = \sqrt{5} x + y,\] 
\[K = \mathbb{Q}(2^{1/3},\omega),\; 
c_4 = t^2 (t^2 - 1728 (t + \omega)),\;\
c_6 = t^2 (t^2 - 1728 (t + \omega))^2,
\; M_{\mathcal{E}} = x (\omega x + y), 
\]
where $\omega$ is a third root of unity.
\subsection{Pathologies}\label{subs:fallada}
There are three kinds of families to which our results do not apply:
(a) constant families, (b) non-constant families with $M_{\mathcal{E}}=1$,
and (c) families over $K(t)$, $K\ne \mathbb{Q}$, such that
$\mathfrak{B}_i(K,M_{\mathcal{E}})$ fails to hold.
The first kind is well understood; if $K$ is Galois, the third kind behaves
essentially like the second kind. (See Appendix \ref{sec:app}.)
Consider, then, $\mathcal{E}$ over $\mathbb{Q}(t)$ with 
$M_{\mathcal{E}}=1$.  Choosing $\mathfrak{M}$ large enough in Proposition
\ref{prop:ahem}, applying Lemmas \ref{lem:strau1} and \ref{lem:strau2}
and assuming $\mathfrak{A}_i(\mathbb{Q},B_{\mathcal{E}}')$, we can
see that there are intersections $S\cap L$ and arithmetic
progressions $a + m \mathbb{Z}$ 
over which $W(\mathcal{E}(t))$ in fact does not average
 to $0$. We may still have 
$\av_{\mathbb{Z}} W(\mathcal{E}(t)) = 0$ or
$\av_{\mathbb{Q},\mathbb{Z}^2} W(\mathcal{E}(t)) = 0$
by cancellation of some sort. The following is an example
where such cancellation does not occur.

Let
\[f(t) = \frac{t^5-1}{t-1},\;
g(t) = \frac{6 (t^7-1)}{t-1}.\]
Define $\mathcal{E}$ to be the elliptic curve over $\mathbb{Q}(t)$ given
by the equation
\[y^2 = x^3 - 3 f (f^3-g^2)^2 x - 2 g (f^3-g^2)^3 .\]
Bounding $\av_\mathbb{Q} \mathcal{E}(t)$ from below by a positive number is
simply a matter of consulting Halberstadt's tables \cite{Ha}. A short 
computer program yields that
\[\begin{aligned}
\frac{1}{25^2} \mathop{\sum_{x=1}^{100} \sum_{y=1}^{100}}_{\gcd(x,y)=1}
 W(\mathcal{E}(y/x))
&= 0.395,\\
\frac{1}{25^2} \mathop{\sum_{x=1}^{100} \sum_{y=1}^{100}}_{\gcd(x,y)=1}
 W(\mathcal{E}(-y/x))
&= 0.35,\\
\frac{1}{100^2} \mathop{\sum_{x=1}^{100} \sum_{y=1}^{100}}_{\gcd(x,y)=1}
 W(\mathcal{E}(y/x))
&= 0.351.\\
\end{aligned}\]

Finally, there is the curious matter of families $\mathcal{E}$
with $M_{\mathcal{E}}(x,y)=x$: the average of $W(\mathcal{E}(t))$
over the rationals is zero, but $M_{\mathcal{E}}(1,t)=1$, and thus
Theorem 1.1 does not apply. This is indeed the case for any 
$\mathcal{E}$ with $j$ a polynomial, $v_t(d)\not\equiv \deg(j)
\mo 2$, $c_4$ and $c_6$ given by $j$ and $d$ as in (\ref{eq:stern}).
If $j$ is a polynomial and $v_t(d) \equiv \deg(j) \mo 2$, then 
$M_{\mathcal{E}}(x,y)=1$. Thus, for any family $\mathcal{E}$
with polynomial $j$, 
there is an arithmetic progression
$a + m \mathbb{Z}$ 
such that 
$\av_{a + m \mathbb{Z}} W(\mathcal{E}(t))$ is non-zero.

%% file: lambda.tex
\chapter{The parity problem}\label{chap:par}
\section{Outline}
Let $f\in \mathbb{Z}\lbrack x,y\rbrack$ be a non-constant homogeneous
polynomial of degree at most $3$. Let $\alpha$ be the Liouville function
($\alpha=\lambda$) or the Moebius function ($\alpha=\mu$). We show that
$\alpha(f(x,y))$ averages to zero. (If $\alpha=\lambda$, we assume, of course,
that $f$ is not of the form $C\cdot g^2$, $C\in \mathbb{Z}$, $g\in
\mathbb{Z}\lbrack x,y\rbrack$.)

The case $\deg f = 1$ is well-known.
Our solution for the case $\deg f = 2$ can hardly be said to be novel, as the main ideas go
back to de la Vall\'ee-Poussin (\cite{DVP1}, \cite{DVP2}) and Hecke (\cite{Hec}). 
Nevertheless, there seems to be no treatment in the literature displaying both
full generality and a strong bound in accordance with the current state of knowledge on
 zero-free regions. We will treat a completely general quadratic form,
without assuming that the form is positive-definite or that its discriminant is
a field discriminant. Our bounds will reflect the broadest known zero-free regions
of Hecke $L$-functions. We will allow the variables to be confined to given lattice
cosets or to sectors in the plane.

The case $\deg f = 3$ appeared to be completely out of reach until rather recently.
We will succeed in breaking parity by an array of methods; in so far as there is an
overall common method, it may be said to consist in the varied usage of traditional 
sieve-methods in non-traditional ways. The strategy used for reducible polynomials is
clearly different from that for irreducible polynomials. (The latter case has a parallel
in the problem of capturing primes.) Nevertheless, there may be some deep similarities
that have come only indirectly and partially to the fore. Note how there seems to be a uniform
barrier for the error bound at $1/(\log N)$. 
Bilinear conditions lurk everywhere.
\section{Preliminaries}
\subsection{The Liouville function}

The Liouville function $\lambda(n)$ is defined on
the set of non-zero rational integers as follows:
\begin{equation}\label{eq:deflam}
\lambda(n) = \prod_{p|n} (-1)^{v_p(n)} .\end{equation}
The following identities are elementary:
\[\lambda(n) = \mu(n) \;\text{ for $n$ square-free},\]
\[
\sum_{d|n} |\mu(d)| \lambda(n/d) = \begin{cases}1 &\text{if $n=1$}\\
0 &\text{if $n>1$,}\end{cases}\]
\[
\sum_n \lambda(n) n^{-s} = \prod_p \frac{1}{1+p^{-s}} = 
\frac{\zeta(2 s)}{\zeta(s)} .\]

We will find it convenient to choose a value for $\lambda(0)$; we adopt the
convention that $\lambda(0) = 0$. We can easily extend the domain of
$\lambda$ further. We define $\lambda$ on $\mathbb{Q}$ by
\begin{equation}\label{eq:extrat}
\lambda\left(\frac{n_0}{n_1}\right) = \frac{\lambda(n_0)}{\lambda(n_1)}
\end{equation}
and on
ideals in a Galois extension $K/\mathbb{Q}$ of degree $n$ by
\begin{equation}\label{eq:exide}
\lambda(\mathfrak{p}_1^{e_1} \mathfrak{p}_2^{e_2} \dotsb \mathfrak{p}_k^{e_k}) =
\prod_{i} \omega^{f(\mathfrak{p_i})\cdot e_i},\end{equation}
where $\omega$ is a fixed $(2n)$th root of unity and
$f(\mathfrak{p_i})$ is the degree of inertia of $\mathfrak{p}_i$ over 
$\mathfrak{p}_i \cap \mathbb{Q}$. Notice that (\ref{eq:exide}) restricts
to (\ref{eq:extrat}), which, in turn, restricts to (\ref{eq:deflam}).
Notice also that the above {\em extension} is different from the natural
{\em generalization} $\lambda_K$:
\begin{equation}\label{eq:exgen}
\lambda_K(\mathfrak{p}_1^{e_1} \mathfrak{p}_2^{e_2} \dotsb \mathfrak{p}_k^{e_k}) =
\prod_{i} (-1)^{e_i}.\end{equation}
We define, as usual,
\begin{equation}\label{eq:exmob}
\mu_K(\mathfrak{p}_1^{e_1} \mathfrak{p}_2^{e_2} \dotsb \mathfrak{p}_k^{e_k}) =
\begin{cases}\prod_{i} (-1)^{e_i} & \text{if $e_i\leq 1$ for all 
$i=1,2,\dotsb,k$}\\ 0 & \text{otherwise.}\end{cases}\end{equation}
\subsection{Ideal numbers and Gr\"ossencharakters}\label{subs:idgross}

Let $K$ be a number field. Write $\mathfrak{O}_K$ for its ring of integers.
Let $I_K$ be the semigroup
of non-zero ideals of $\mathfrak{O}_K$; let $J_K$ be the group of non-zero
fractional ideals of $\mathfrak{O}_K$.  For every $\mathfrak{d}\in I_K$,
define $\mathfrak{O}_{K,\mathfrak{d}}$ to be the set of
elements of $\mathfrak{O}_K$ prime to $\mathfrak{d}$. Define
$I_{K,\mathfrak{d}}$ to be the semigroup of ideals
of $\mathfrak{O}_K$ prime to $\mathfrak{d}$. 

Since the class group of $\mathfrak{O}_K$ is finite, there are ideals $\mathfrak{a}_1, \mathfrak{a}_2,\dotsc,\mathfrak{a}_{i_0}\in I_K$ and positive integers
$h_1, h_2,\dotsb, h_{i_0}$ such that every $\mathfrak{d}\in \mathfrak{O}_K$
can be expressed in a unique way in the form
\begin{equation}\label{eq:decomp}
\mathfrak{d} = \mathfrak{d}_p \mathfrak{a}_1^{d_1} 
\mathfrak{a}_2^{d_2} \dotsb \mathfrak{a}_{i_0}^{d_{i_0}}, 
\text{\;\;\;\;$\mathfrak{d}_p$ principal, $0\leq d_i<h_i$.}\end{equation}
Fix $\alpha_1,\dotsc,\alpha_{i_0}\in \mathfrak{O}_K$ such that
$(\alpha_i) = \mathfrak{a}_i^{h_i}$. 
Choose $\beta_1,\dotsc,\beta_{i_0}$ in the algebraic completion of
$K$ such that  $\beta_i^{h_i} = \alpha_i$ for every $i=1,\dotsb,i_0$.
Define $L = K(\beta_1,\dotsc,\beta_{i_0})$. Let $\mathcal{I}(K)^{\times}$
be the subgroup of $L^*$ generated by $K^*$ and $\beta_1,\dotsb,\beta_{i_0}$.
We say that $\mathcal{I}(K) = \mathcal{I}(K)^{\times} \cap \{0\}$ is the
set of {\em ideal numbers}. For 
$a = \alpha \beta_1^{d_1} \dotsb \beta_{i_0}^{d_{i_0}} \in 
 \mathcal{I}(K)^{\times}$, $\alpha \in K$, let
$\mathfrak{I}(a) = (\alpha) \mathfrak{a}_1^{d_1} \dotsb \mathfrak{a}_{i_0}^{d_{i_0}}$.
Then $\mathfrak{I}: \mathcal{I}(K)^{\times} \to J_K$ is a 
surjective homomorphism
with kernel $\mathfrak{O}_K^*$. We define $\mathcal{I}(\mathfrak{O}_K)^{\times}$ to be the preimage
$\mathfrak{I}^{-1}(I_K)$.

For $a,b\in \mathcal{I}(\mathfrak{O}_K)^{\times}$, we say that 
$a|b$ 
({\em $a$ divides $b$}) 
if $b = a c$ for some $c\in \mathcal{I}(\mathfrak{O}_K)^{\times}$;
we say that $\gcd(a,b)=1$ ({\em $a$ is prime to $b$}) if
there is no non-unit $c\in \mathcal{I}(\mathfrak{O}_K)^{\times}$ such that
$c|a$, $c|b$.

Let $\mathfrak{d}\in I_K$. Let $\hat{\mathfrak{d}}$ be an arbitrary element
of $\mathfrak{I}^{-1}(\mathfrak{d})$.
Define $\mathcal{I}(\mathfrak{O}_K)_{\mathfrak{d}}$ to be the semigroup
of all $a\in \mathcal{I}(\mathfrak{O}_K)^{\times}$ prime to $\hat{\mathfrak{d}}$.
For $a,b\in \mathcal{I}(\mathfrak{O}_K)_{\mathfrak{d}}$, $a = \alpha \beta_1^{a_1} \dotsb \beta_{i_0}^{a_{i_0}}$, $b = \beta \beta_1^{b_1} \dotsb \beta_{i_0}^{b_{i_0}}$,
we say that $a\sim b$ if 
$a_i = b_i$ for every $i=1,\dotsc,i_0$ and
$\hat{\mathfrak{d}}|(\alpha-\beta) \beta_1^{a_1}\dotsb \beta_{i_0}^{a_{i_0}}$.
Define $\mathcal{C}_{\mathfrak{d}}(K)$ to be the set of equivalence
classes of $\mathcal{I}(\mathfrak{O}_K)_{\mathfrak{d}}$ under $\sim$.
 
For every embedding of $K$ into $\mathbb{C}$, 
choose an embedding of $L$
extending it; since $\mathcal{I}(K)\subset L$, we obtain an embedding
of $\mathcal{I}(K)$ into $\mathbb{C}$. Let $\iota_1,\dotsc,\iota_{\deg_K}$
be the embeddings of $\mathcal{I}(K)$ thus obtained; order them
so that $\iota_1,\dotsc,\iota_{r_1}$ come from the real embeddings of $K$
and $\iota_{r_1+1},\dotsc,\iota_{r_1+2 r_2}$ come from the complex
embeddings of $K$. We can assume $\iota_{r_1+r_2+1} = \overline{
\iota_{r_1 + 1}},\dotsc,\iota_{r_1 + 2 r_2} = \overline{\iota_{r_1 + r_2}}$.

For $a,b\in \mathcal{I}(\mathfrak{O}_K)_{\mathfrak{d}}$, we say
that $a\sim_n b$ if $a\sim b$ and $\sgn \iota_i(a) = \sgn \iota_i(b)$
for every $i=1,\dotsc,\deg_K$.
Define $\mathcal{C}_{\mathfrak{d}}^n(K)$ to be the set of equivalence
classes of $\mathcal{I}(\mathfrak{O}_K)_{\mathfrak{d}}$ under $\sim_n$.

We denote the set of all characters $\chi$ of a finite group $G$ by
$\Xi(G)$. Let $\chi \in \Xi(\mathcal{C}_{\mathfrak{d}}^n(K))$. For
$s_1,\dotsc, s_{r_1+r_2}\in \mathbb{R}$, $n_1,\dotsc, n_{r_2}\in \mathbb{Z}$,
define $\gamma_{s,n}: \mathcal{I}(\mathfrak{O}_K)^{\times} \to S^1$ 
as follows:
\begin{equation}\label{eq:pityofwar}
\gamma_{s,n}(a) = \prod_{j=1}^{r_1+r_2} |\iota_j(a)|^{i s_j}
\prod_{j=1}^{r_2} 
 \left(\frac{\iota_{r_1+j}(a)}{|\iota_{r_1+j}(a)|}\right)^{n_j} .\end{equation}
Assume $\gamma_{s,n}(u) \chi(u) = 1$ for every unit 
$u\in \mathfrak{O}_K^* \subset \mathcal{I}(\mathfrak{O}_K)^{\times}$.
 Then we can define
the {\em Gr\"ossencharakter} $\psi_{\chi,s,n}:I_{K,\mathfrak{d}}\to S^1$ by
$\psi(\mathfrak{a}) = \chi(a) \gamma_{s,n}(a)$,
where $a$ is any element of $\mathfrak{I}^{-1}(\mathfrak{a})$.

Consider now $K/\mathbb{Q}$ quadratic. We can 
describe the Gr\"ossencharakters of $K$ as follows. Let $K/\mathbb{Q}$
be imaginary. Write $\iota$ for the embedding $\iota_1$ of 
$\mathcal{I}(\mathfrak{O}_K)$ in $\mathbb{C}$.
Let $\chi\in \Xi(\mathcal{C}_{\mathfrak{d}}(K))$. If
$n$ is an integer such that $\chi(u) (\iota(u))^n = 1$ for every 
$u\in \mathfrak{O}_K^*$, then there is a Gr\"ossencharakter
\begin{equation}\label{eq:disti}
\psi_n(\mathfrak{a}) = \chi(a) \left(\frac{\iota(s)}{|\iota(s)|}\right)^n .
\end{equation}

Let $K/\mathbb{Q}$ now be real. In the definition of $\mathcal{I}(K)^{\times}$, we can choose
$\alpha_1,\dotsc,\alpha_{i_0}$ positive and $\beta_1,\dotsc,\beta_{i_0}$ real.
Thus we can assume that $\iota_1(a),\iota_2(a)\in \mathbb{R}$ for all 
$a\in \mathcal{I}(K)$.
Let $u_1$ be the primitive unit
of $\mathfrak{O}_K$ such that $\iota_1(u_1) > 1$.
For every $\mathfrak{d}\in I_K$, let $k_{\mathfrak{d}}$ be the smallest
positive integer such that $u_1^{k_{\mathfrak{d}}} \equiv 1 \mo \mathfrak{d}$.
Let \[r_{\mathfrak{d}} = \begin{cases} 
1 &\text{if $\frac{\iota_1(u_1)}{\iota_2(u_1)} > 0$,}\\
2 &\text{if $\frac{\iota_1(u_1)}{\iota_2(u_1)} < 0$.}\end{cases}\]
Let $l_{\mathfrak{d}}$ be the positive real number
$\left(\frac{\iota_1(u_1)}{\iota_2(u_1)}\right)^{r_{\mathfrak{d}} k_{\mathfrak{d}}}$.
Let $\chi\in \Xi(\mathcal{C}_{\mathfrak{d}}(K))$. If $n\in \mathbb{Z}$,
$n_0\in \{0,1\}$ are such that
\[\chi(u_1) 
 \left(\sgn\left(\frac{\iota_1(u_1)}{\iota_2(u_1)}\right)\right)^{n_0}
\left|\frac{\iota_1(u_1)}{\iota_2(u_1)}\right|^{2\pi i n/\log l_{\mathfrak{d}}}
= 1,\] then there is a Gr\"ossencharakter
\begin{equation}\label{eq:tilled}\psi_n(\mathfrak{a}) = \chi(a) 
\left(\sgn\left(\frac{\iota_1(a)}{\iota_2(a)}\right)\right)^{n_0}
\left|\frac{\iota_1(a)}{\iota_2(a)}\right|^{2\pi i n/\log l_{\mathfrak{d}}}
.\end{equation}

We define the {\em size} $\mathcal{S}(\psi)$ 
of a Gr\"ossencharakter $\psi$ to be
\begin{equation}\label{eq:siz}
\sqrt{\sum_{j=1}^{r_1+r_2} s_j^2 + \sum_{j=r_1+r_2+1}^{r_1+2 r_2} n_j^2},
\end{equation}
where $s_j$ and $n_j$ are as in (\ref{eq:pityofwar}). For $K/\mathbb{Q}$
quadratic and imaginary,
\[\mathcal{S}(\psi) = n,\]
where $n$ is as in (\ref{eq:disti}). For $K/\mathbb{Q}$ quadratic and real,
\[\mathcal{S}(\psi) = 2^{3/2} \pi n/\log l_{\mathfrak{d}} ,\]
where $n$ is as in (\ref{eq:tilled}). Thus, if we take $K/\mathbb{Q}$ to
be fixed,
\[\mathcal{S}(\psi) \ll N \mathfrak{d} \cdot n .\]
\subsection{Quadratic forms}
We will consider only quadratic forms $a x^2 + b x y + c y^2$ with 
integer coefficients $a,b,c\in \mathbb{Z}$. A quadratic form $a x^2 + b x y
+ c y^2$ is {\em primitive} if $\gcd(a,b,c)=1$.  

Let $n$ be a rational integer. We denote by $\sq(n)$ the largest 
positive integer whose square divides $n$. Define
\[d_n = \begin{cases} \sq(n) &\text{if $4\nmid n$}\\
\sq(n)/2 &\text{if $4|n$.}\end{cases}\]

\begin{lem}\label{lem:quad}
Let $Q(x,y) = a x^2 + b x y + c y^2$ be a primitive, irreducible quadratic form. Let $K = \mathbb{Q}(\sqrt{b^2 - 4 a c})$. Then there are 
algebraic integers $\alpha_1, \alpha_2\in \mathfrak{O}_K$ linearly
independent over $\mathbb{Q}$
such that 
\[ Q(x,y) = \frac{N(x \alpha_1 + y \alpha_2)}{a}\]
for all $x,y\in \mathbb{Z}$. The subgroup $\mathbb{Z} \alpha_1 + \mathbb{Z} \alpha_2$
of $\mathfrak{O}_K$ has index $\lbrack \mathfrak{O}_K : \mathbb{Z} \alpha_1 +
\mathbb{Z} \alpha_2\rbrack = d_{b^2 - 4 a c}$.
\end{lem}
\begin{proof}
Set $\alpha_1 = a$, $\alpha_2 = \frac{b + \sqrt{b^2 - 4 a c}}{2}$. 
\end{proof}

\subsection{Truth and convention} 
Following Iverson and Knuth \cite{Kn}, we define 
$\lbrack \text{true}\rbrack$ to be $1$ and $\lbrack \text{false}\rbrack$ 
to be zero. Thus, for example, $x\to \lbrack x\in S\rbrack$ is the 
characteristic function of a set $S$.

\subsection{Approximation of intervals}
We denote by $S^1$ the unit circle in $\mathbb{R}^2$. An {\em interval}
$I\subset S^1$ is a connected subset of $S^1$. 

\begin{lem}\label{lem:vin}
Let $I\subset S^1$ be an interval with endpoints $x_0$, $x_1$. 
Let $d(x,y)\in \lbrack 0,\pi\rbrack$ denote the angle between two given 
points $x, y\in S^1$.
Then, for any positive $\epsilon$ and any 
positive integer $k$, there are complex numbers $\{a_n\}_{n=-\infty}^\infty$
such that
\[\begin{aligned}
0\leq \sum_{n=-\infty}^{\infty} a_n x^n &\leq 1
\;\;\;\;\;\;\;\;\;\text{for every $x\in S^1$,}\\
\sum_{n=-\infty}^{\infty} a_n x^n &= \lbrack x\in I\rbrack
\;\;\;\;\text{if $d(x,x_0), d(x,x_1) \geq \epsilon/2$,}\end{aligned}\]
\[|a_n|\ll \left(\frac{k}{\epsilon}\right)^k |n|^{-(k+1)}\;\;\;\;\;
\text{for $n\ne 0$},\]
\[|a_0|\ll 1.\]
The implied constant is absolute.
\end{lem}
\begin{proof}
See \cite{Vi}, Ch. 1, Lemma 12.
\end{proof}

\subsection{Lattices, convex sets and sectors}

A \emph{lattice} is a subgroup of $\mathbb{Z}^n$ of finite index; a \emph{%
lattice coset} 
 is a coset of such a subgroup. By the {\em index} of a lattice
coset we mean the index of the lattice of which it is a coset.
For any lattice cosets $L_1$, $%
L_2$ with $\gcd(\lbrack \mathbb{Z}^n : L_1 \rbrack , \lbrack \mathbb{Z}^n :
L_2 \rbrack) = 1$, the intersection $L_1\cap L_2$ is a lattice coset with 
\begin{equation}
\lbrack \mathbb{Z}^n : L_1\cap L_2\rbrack = \lbrack \mathbb{Z}^n :
L_1\rbrack \lbrack \mathbb{Z}^n : L_2\rbrack .
\end{equation}
In general, if $L_1$, $L_2$ are lattice cosets, then $L_1\cap L_2$ is either
the empty set or a lattice coset such that 
\begin{equation}  \label{eq:intersl} \begin{aligned}
\lcm(\lbrack \mathbb{Z}^n : L_1\rbrack,\lbrack \mathbb{Z}^n : L_2\rbrack) 
&\mid
\lbrack \mathbb{Z}^n : L_1\cap L_2\rbrack,\\ \lbrack \mathbb{Z}^n :
L_1\cap L_2\rbrack &\mid \lbrack \mathbb{Z}^n : L_1\rbrack \lbrack \mathbb{Z}^n
: L_2\rbrack .
\end{aligned} \end{equation}

For $S\subset \lbrack -N,N\rbrack^n$ a convex set and $L\subset \mathbb{Z}^n$ a lattice
coset, 
\begin{equation}  \label{eq:wbrl}
\#(S\cap L) = \frac{\Area(S)}{\lbrack \mathbb{Z}^n :L\rbrack} + O(N^{n-1}),
\end{equation}
where the implied constant depends only on $n$.

By a \emph{sector} we will mean a connected component of a set of the form $%
\mathbb{R}^{n}-(T_{1}\cap T_{2}\cap \dotsb \cap T_{n})$, where $T_{i}$ is a
hyperplane going through the origin. Every sector $S$ is convex. Given
a sector $S\subset \mathbb{R}^2$, we may speak of the {\em angle} 
$\alpha\in (0,2 \pi\rbrack$ spanned by $S$, or, for short, the {\em angle 
$\alpha$ of $S$}.

Call a sector $S$ of $\mathbb{R}^2$ a {\em subquadrant} if its closure 
intersects the $x$- and $y$-axes only at the origin. 
By the {\em hyperbolic angle} $\theta \in (0,\infty\rbrack$ of a subquadrant
$S\subset \mathbb{R}^2$ we mean
\[\sup_{(x,y)\in S} \log |x/y| - \inf_{(x,y)\in S} \log |x/y| .\]
Notice that the area of the region
\[\{(x,y)\in S : x^2 + y^2 \leq R\}\]
equals $\frac{1}{2} \alpha R$, where $\alpha$ is the angle of $S$,
whereas
the area of the region
\[\{(x,y)\in S: x y \leq R\}\]
equals $\frac{1}{2} \theta R$, where $\theta$ is the hyperbolic
angle of $S$.
\subsection{Classical bounds and their immediate consequences}

By Siegel, Walfisz and Vinogradov (vd. \cite{Wa}, V \S 5 and V \S 7),
\begin{equation}\label{eq:sw}
\left|\mathop{\sum_{n\leq x}}_{n\equiv a \mo m} \lambda(n)\right| \ll
x e^{-C (\log x)^{2/3}/(\log \log x)^{1/5}},\end{equation}
\begin{equation}\label{eq:swmu}
\left|\mathop{\sum_{n\leq x}}_{n\equiv a \mo m} \mu(n)\right| \ll
x e^{-C (\log x)^{2/3}/(\log \log x)^{1/5}}\end{equation}
for $m\leq (\log x)^A$, with $C$ and the implied constant depending on $A$.

The following lemma is well-known in essence.
\begin{lem}\label{lem:corknut}
Let $K$ be a finite extension of $\mathbb{Q}$. Let 
$\mathfrak{d}$ be an ideal of $\mathfrak{O}_K$. Let
$\psi$ be a
Gr\"ossencharacter on $I_{K,\mathfrak{d}}$.
Assume 
\begin{equation}\label{eq:staclog}
\mathcal{S}(\psi) \ll e^{(\log x)^{3/5} (\log \log x)^{1/5}} .
\end{equation}
If 
\begin{equation} \label{eq:tuxd}
N \mathfrak{d} \ll e^{(\log x)^{2/5} (\log \log x)^{1/5}}
\end{equation}
and $\psi$ is not a real Dedekind character, or
\begin{equation} \label{eq:tuxd2}
N \mathfrak{d} \ll (\log N)^A
\end{equation}
and $\psi$ is a real Dedekind character, then
\begin{equation} \label{eq:fray}
\mathop{\sum_{\mathfrak{m}\in I_{K,\mathfrak{d}}}}_{ N \mathfrak{m%
}\leq x} \psi(\mathfrak{m}) \mu_K(\mathfrak{m}) \ll x
e^{- C \frac{(\log x)^{2/3}}{(\log \log x)^{1/5}}},
\end{equation}
where $C$ and the implied constant in (\ref{eq:fray})
depend only on $K$, $A$, and the implied
constants in (\ref{eq:staclog}), (\ref{eq:tuxd}) and (\ref{eq:tuxd2}).
\end{lem}
\begin{proof}
Clearly
\[
\sum_{\mathfrak{m}\in I_{K,\mathfrak{d}}} \mu_K(\mathfrak{m}) (N \mathfrak{m})^{-s}
 = \prod_{\mathfrak{p}} (1 - (N \mathfrak{p})^{-s}) = 
\frac{1}{L(\psi,s)} .\]
Given the zero-free region in \cite{Col} and the Siegel-type bound in \cite{Fo}
for the 
exceptional zero, the result follows in the
standard fashion (see e.g. \cite{Dav}, Ch. 20--22, or \cite{Col}, \S 6).
\end{proof}

\begin{lem}\label{lem:lamp}
Let $K$ be a quadratic extension of $\mathbb{Q}$. Let 
$\mathfrak{d}$ be an ideal of $\mathfrak{O}_K$. Let
$\psi$ be a
Gr\"ossencharacter on $I_{K,\mathfrak{d}}$.
 Suppose 
\begin{equation}  \label{eq:txd}
\mathcal{S}(\psi) \ll e^{(\log x)^{3/5} (\log \log x)^{1/5}},\;\;\;
N \mathfrak{d} \ll (\log N)^A .
\end{equation}
Then
\begin{equation} \label{eq:frey}
\mathop{\sum_{\mathfrak{m}\in I_{K,\mathfrak{d}}}}_{ N \mathfrak{m%
}\leq x} \psi(\mathfrak{m}) \lambda(N \mathfrak{m}) \ll x
e^{- C \frac{(\log x)^{2/3}}{(\log \log x)^{1/5}}},
\end{equation}
where $C$ and the implied constant in (\ref{eq:frey})
depend only on $K$, $A$ and the implied
constant in (\ref{eq:txd}).
\end{lem}
\begin{proof}
Define 
\begin{equation*}
\phi(\psi,s) = \sum_{\mathfrak{m}\in I_{K,\mathfrak{d}}} \psi(\mathfrak{m})
\lambda(N \mathfrak{m}) (N \mathfrak{m})^{-s}
\end{equation*}
for $\Re s> 1$. We can express $\phi$ as an Euler product: 
\begin{equation*}
\phi(\psi,s) = \mathop{\prod_{\mathfrak{p}\in I_{K,\mathfrak{d}}}}_{ \text{$%
\mathfrak{p}\cap \mathbb{Q}$ splits in $K$}} \frac{1}{1 + \psi(\mathfrak{p})
(N \mathfrak{p})^{-s}} \mathop{\prod_{\mathfrak{p}\in
I_{K,\mathfrak{d}}}}_{\text{ $\mathfrak{p}\cap \mathbb{Q}$ does not split in 
$K$}} \frac{1}{1 - \psi(\mathfrak{p}) (N \mathfrak{p})^{-s}} .
\end{equation*}
Write 
\begin{equation*}
R(\psi,s) =  \mathop{\prod_{\mathfrak{p}\in I_{K,\mathfrak{d}}}}_{\text{ $%
\mathfrak{p}\cap \mathbb{Q}$ ramifies}} \frac{1 + \psi(\mathfrak{p}) 
(N \mathfrak{p})^{-s}}{ 1 - \psi(\mathfrak{p}) 
(N \mathfrak{p})^{-s}} .
\end{equation*}
Then 
\begin{equation*}
\begin{aligned}\phi(\psi,s) 
&=
R(\psi,s)
\mathop{\prod_{\mathfrak{p}\in I_{K,\mathfrak{d}}}}_{
\text{$\mathfrak{p}\cap \mathbb{Q}$ unsplit \& unram.}}
\frac{1 + \psi(\mathfrak{p}) (N \mathfrak{p})^{-s}}{
1 - \psi(\mathfrak{p}) (N \mathfrak{p})^{-s}}
 \prod_{\mathfrak{p}\in I_{K,\mathfrak{d}}} 
\frac{1}{1 + \psi(\mathfrak{p}) (N \mathfrak{p})^{-s}} \\
&=
R(\psi,s)
\mathop{\prod_{p\nmid d}}_{\text{$p$ unsplit \& unram. in $K$}}
 \frac{1 + \chi(p) p^{-2 s}}{1 - \chi(p) p^{-2 s}}
 \prod_{\mathfrak{p}\in I_{K,\mathfrak{d}}}
 \frac{1 - \psi(\mathfrak{p}) (N \mathfrak{p})^{-s}}{
1 - \psi^2(\mathfrak{p}) (N \mathfrak{p})^{-2 s}},\end{aligned}
\end{equation*}
where $d = N \mathfrak{d}$ and $\chi$ is the restriction of $\psi$
to $\mathbb{Z}^+$. We denote 
\begin{equation*}
\begin{aligned}
\chi'(p) &= \begin{cases} 0 &\text{if $p$ ramifies}\\
1 &\text{if $p$ splits}\\
-1 &\text{if $p$ neither splits nor ramifies,}\end{cases}\\
L(\psi,s) &= \prod_{\mathfrak{p}\in I_{K,\mathfrak{d}}}
\frac{1}{1 - \psi(\mathfrak{p}) (N \mathfrak{p})^{-s}}\end{aligned}
\end{equation*}
and obtain 
\begin{equation*}
\phi(\chi,s) = R(\psi,s) \prod_{\text{$p$ ram. in $K$}}
(1 - \chi(p) p^{-2 s})\:
\frac{L(\chi, 2 s)}{L(\chi\cdot \chi^{\prime},2 s)} 
\frac{L(\psi^2,2 s)}{L(\psi,s)} .
\end{equation*}
Proceed as in Lemma \ref{lem:corknut}.
\end{proof}

\subsection{Bilinear bounds}
We shall need bilinear bounds for the Liouville function. For section
\ref{sec:three}, the following lemma will suffice. It is simply a linear
bound in disguise.
\begin{lem}\label{lem:diverto}
Let $S$ be a convex subset of $\lbrack -N, N\rbrack^2$. Let 
$L\subset \mathbb{Z}^2$ be a lattice coset of index 
\begin{equation}\label{eq:lanc}
\lbrack \mathbb{Z}^2 : L\rbrack \ll (\log N)^A.\end{equation}
Let $f:\mathbb{Z}\to \mathbb{C}$ be a function with $\max_y |f(y)|\leq 1$.
 Then, for every $\epsilon>0$,
\begin{equation}\label{eq:gen}
\left| \sum_{(x,y)\in S\cap L} \lambda(x) f(y) \right| \ll
\Area(S) \cdot
 e^{- C (\log N)^{2/3}/(\log \log N)^{1/5}} + N^{1+\epsilon} ,\end{equation}
where $C$ and the implied constant in (\ref{eq:gen})
depend only on $K$, $\epsilon$, $A$ and the implied constant in 
(\ref{eq:lanc}).
\end{lem}
\begin{proof}
For every $y\in \mathbb{Z}\cap \lbrack - N, N\rbrack$, the set 
$\{x : (x,y)\in L\}$ is either the empty set or an arithmetic progression
$m \mathbb{Z} + a_y$, where $m | \lbrack \mathbb{Z}^2 : L \rbrack$.
Let $y_0$ and $y_1$ 
be the least and the greatest $y\in \mathbb{Z}\cap \lbrack -N,N\rbrack$ 
such that $\{x : (x,y)\in S\}$ is non-empty. Let $y\in \mathbb{Z} \cap
\lbrack y_0,y_1\rbrack$. 
Since $S$ is convex and a subset of $\lbrack - N, N\rbrack^2$, 
the set
$\{x : (x,y) \in S\}$ is an interval $\lbrack N_{y,0},N_{y,1}\rbrack$ contained
in $\lbrack - N, N\rbrack$. Hence
\[\begin{aligned}
\left| \sum_{(x,y)\in S\cap L} \lambda(x) f(y) \right| &=
\left| \mathop{\sum_{y_0\leq y\leq y_1}}_{\{x : (x,y)\in L\} \ne \emptyset}
\mathop{\sum_{-N_{y,0}\leq x\leq N_{y,1}}}_{x\equiv a_y \mo m_y}
 \lambda(x) f(y) \right| \\ &\leq
\mathop{\sum_{y_0\leq y\leq y_1}}_{\{x : (x,y)\in L\} \ne \emptyset}
\left| \mathop{\sum_{-N_{y,0}\leq x\leq N_{y,1}}}_{x\equiv a_y \mo m}
 \lambda(x) \right| .\end{aligned}\]
By (\ref{eq:sw}),
\[\begin{aligned}
 \mathop{\sum_{y_0\leq y\leq y_1}}_{\{x : (x,y)\in L\} \ne \emptyset}
\left| \mathop{\sum_{-N_{y,0}\leq x\leq N_{y,1}}}_{x\equiv a_y \mo m}
 \lambda(x) \right| &= 
 \mathop{\mathop{\sum_{y_0\leq y\leq y_1}}_{\{x : (x,y)\in L\} \ne 
\emptyset}}_{N_{y,1}-N_{y,0}> N^{\epsilon}} 
\left| \mathop{\sum_{-N_{y,0}\leq x\leq N_{y,1}}}_{x\equiv a_y \mo m}
 \lambda(x) \right| \\ &+
 \mathop{\mathop{\sum_{y_0\leq y\leq y_1}}_{\{x : (x,y)\in L\} \ne 
\emptyset}}_{N_{y,1}-N_{y,0} \leq N^{\epsilon}} 
\left| \mathop{\sum_{-N_{y,0}\leq x\leq N_{y,1}}}_{x\equiv a_y \mo m}
 \lambda(x) \right| \\
&\ll \sum_{y_0\leq y\leq y_1} (N_{y_1}-N_{y_0}) e^{-C (\log N^{\epsilon})^{2/3}/(\log \log N)^{1/5}} +
N^{1+\epsilon}
.\end{aligned}\]
Clearly \[\Area(S) = \sum_{y=y_0}^{y_1} (N_{y,1} - N_{y,0}) + O(N) .\]
Therefore
\[
\left| \sum_{(x,y)\in S\cap L} \lambda(x) f(y) \right| \ll
\Area(S) \cdot e^{-C (\log N^{\epsilon})^{2/3}/(\log \log N)^{1/5}} +
N^{1+\epsilon} .\]
\end{proof}
As a special case of, say, Theorem 1 in \cite{Le}, we have the following
analogue of Bombieri-Vinogradov:
\begin{equation}\label{eq:bv}
\sum_{m\leq \frac{N^{1/2}}{(\log N)^{2 A+4}}} \mathop{\max_a}_{(a,m)=1} \max_{x\leq N}
\left|\mathop{\sum_{n\leq x}}_{n\equiv a \mo m} \lambda(n) -
\frac{1}{\phi(m)} \mathop{\sum_{n\leq x}}_{\gcd(n,m)=1} \lambda(n)\right| 
\ll \frac{N}{(\log N)^A},\end{equation}
where the implied constant depends only on $A$.

A simpler statement is true.
\begin{lem}\label{bv2}
For any $A>0$,
\[\sum_{m\leq \frac{N^{1/2}}{(\log N)^{2 A+6}}} \max_a 
\max_{x\leq N} \left|\mathop{\sum_{n\leq x}}_{n\equiv a \mo m} \lambda(n) \right| 
\ll \frac{N}{(\log N)^A},\] where the implied constant depends only on $A$.
\end{lem}
\begin{proof}
Write $\rad(m) = \prod_{p|m} p$. Then
\[\sum_{d|\gcd(\rad(m),n)} \lambda(n/d) = \begin{cases}
\lambda(n) & \text{if $\gcd(m,n)=1$}\\
0 &\text{otherwise.}\end{cases}\]
Therefore
\[\begin{aligned}
\sum_{m\leq N^{1/2}} \frac{1}{\phi(m)}
 \max_{x\leq N} \left|\mathop{\sum_{n\leq x}}_{\gcd(n,m)=1} \lambda(n) 
\right| &=
\sum_{m\leq N^{1/2}} \frac{1}{\phi(m)}
 \max_{x\leq N} \left|\sum_{d|\rad(m)} 
\mathop{\sum_{n\leq x}}_{d|n} \lambda(n/d)\right|\\
&\leq \sum_{m\leq N^{1/2}} \frac{1}{\phi(m)}
 \sum_{d|\rad(m)} \max_{x\leq N/d} \left|\sum_{n\leq x} \lambda(n)\right|
\\ &\ll \sum_{m\leq N^{1/2}} \frac{1}{\phi(m)}
 \sum_{d|\rad(m)} N/d \cdot e^{-C \sqrt{\log N/d}}\;\;
\text{by \ref{eq:sw}}\\
&\leq N e^{- C \sqrt{\log N^{1/2}}} \sum_{m\leq N^{1/2}} \frac{1}{\phi(m)}
\sum_{d|\rad(m)} \frac{1}{d} \\
&\ll \frac{N}{(\log N)^A}.\end{aligned}\]

By (\ref{eq:bv}) this implies
\[\sum_{m\leq \frac{N^{1/2}}{(\log N)^{2 A+6}}} \mathop{\max_a}_{\gcd(a,m)=1} 
\max_{x\leq N} \left|\mathop{\sum_{n\leq x}}_{n\equiv a \mo m} \lambda(n) \right| 
\ll \frac{N}{(\log N)^A} .\]

Now
\[\begin{aligned}
\sum_{m\leq \frac{N^{1/2}}{(\log N)^{2 A+6}}} \max_a 
\max_{x\leq N} \left|\mathop{\sum_{n\leq x}}_{n\equiv a \mo m} \lambda(n) \right| &=
\sum_{m\leq \frac{N^{1/2}}{(\log N)^{2 A+6}}} \max_{r|m}
\max_{(a,m)=1} 
\max_{x\leq N} \left|\mathop{\sum_{n\leq x}}_{n\equiv a r\mo m} \lambda(n) \right| \\
&=
\sum_{m\leq \frac{N^{1/2}}{(\log N)^{2 A+6}}} \max_{r|m}
\max_{(a,m)=1} 
\max_{x\leq \frac{N}{r}} \left|\mathop{\sum_{n\leq x}}_{n\equiv a\mo m/r} \lambda(n) \right|\\
&< \sum_{r\leq N^{1/2}}
\sum_{s\leq \frac{(N/r)^{1/2}}{(\log (N/r))^{2 A+6}}}
\max_{(a,s)=1}
\max_{x\leq \frac{N}{r}} \left|\mathop{\sum_{n\leq x}}_{n\equiv a\mo s} \lambda(n) \right|\\
&\ll \sum_{r\leq N^{1/2}} \frac{N/r}{(\log N/r)^{A+1}} \ll
\frac{N}{(\log N)^A} .\end{aligned}\]
\end{proof}

The following lemma is to Lemma \ref{lem:diverto} what Bombieri-Vinogradov is to (\ref{eq:sw}).

\begin{lem}\label{bomb4}
Let $A$, $K$ and $N$ be positive integers such that
$K\leq N^{1/2}/(\log N)^{2 A + 6}$. For $j=1,2,\dotsc,K$,
let $S_j$ be a convex subset of $\lbrack -N,N\rbrack^2$ and let
$L_j\subset \mathbb{Z}^2$ be a lattice coset of index $j$. Let $f:\mathbb{Z}\to \mathbb{C}$ be
a function with $\max_y |f(x,y)| \leq 1$. Then
\[\sum_{j=1}^{K} \left|
\sum_{(x,y)\in S_j\cap L_j} \lambda(x) f(y) \right|
\ll \frac{N^2}{(\log N)^A},\]
where the implicit constant depends only on $A$.
\end{lem}
\begin{proof}
We start with
\[\begin{aligned}
\sum_{j=1}^{K} \left|
\sum_{(x,y)\in S_j \cap L_j} \lambda(x) f(y) \right| &\leq
\sum_{j=1}^{K} 
 \sum_y \left| \mathop{\sum_x}_{(x,y)\in S_j\cap L_j} \lambda(x) \right|\\
&= 
\sum_{j=1}^{K}
 \sum_{k=0}^{\lceil N/j \rceil} \sum_{y = k j}^{(k+1) j -1} \left|
\mathop{\sum_x}_{(x,y)\in S_j\cap L_j} \lambda(x) \right| .\end{aligned}\]

For any $y\in \mathbb{Z}$, the set 
\[\{x : (x,y)\in L_j\}\]
is either the empty set or
an arithmetic progression of modulus $m_j|j$ independent of $y$.
Thus the set
\[A_j=\{(x,y)\in L_j: k j\leq y\leq (k+1) j-1\}\] is the union of $m_j$ sets
of the form
\[B_{y_0,a} = \{(x,y)\in \mathbb{Z}^2 : x\equiv a \mo m_j,\, y = y_0\}\]
with $k j \leq y_0 \leq (k+1) j - 1$. Since an arithmetic progression
of modulus $d$ is the union of $j/d$ arithmetic progressions of modulus $j$, the 
set $A_j$ is the union of $j$ sets of the form
\[C_{x_0,a} =  \{(x,y)\in \mathbb{Z}^2 : x\equiv a \mo j,\, y = y_0\}.\]
Therefore
\[\begin{aligned}
\sum_{j=1}^K
 \sum_{k=0}^{\lceil N/j \rceil} \sum_{y = k j}^{(k+1) j -1} \left|
\mathop{\sum_x}_{(x,y)\in S_j\cap L_j} \lambda(x) \right| &\leq
\sum_{j=1}^K
 \sum_{k=0}^{\lceil N/j \rceil} \sum_{l = 1}^j \left|
\mathop{\sum_x}_{(x,y_0(k,l))\in S_j\cap C_{y_0(k,l),a(k,l)}} 
 \lambda(x) \right| \\ &\leq \sum_{j=1}^K
(N+j)
\max_{y_0} \max_a \left|
\mathop{\sum_x}_{(x,y_0)\in S\cap C_{y_0,a}} \lambda(x)\right|\\ &\leq
\sum_{j=1}^K
(N+j)
\max_{-N\leq b\leq c\leq N} \max_a \left|
\mathop{\sum_{b\leq x\leq c}}_{x\equiv a \mo j} \lambda(x)\right|\\
&\leq \sum_{j=1}^K
 4 (N+j)
\max_{0< c\leq N} \max_a \left|
\mathop{\sum_{0<x\leq c}}_{x\equiv a \mo j} \lambda(x)\right|.\end{aligned}\]

We apply Lemma \ref{bomb5} and are done.
\end{proof}

\begin{cor}\label{bomb5}
Let $A$, $K$, $N$, $d_0$ and $d_1$ be positive integers such that
$K d_1 \leq N^{1/2}/(\log N)^{2 A + 6}$. For $k=1,2,\dotsc ,K$,
let $S_k$ be a convex subset of $\lbrack -N,N\rbrack^2$ and let
$L_k\subset \mathbb{Z}^2$ be a lattice coset of index $\frac{r_k}{d_0} k$ for
some $r_k$ dividing $d_0 d_1$. Then
\[\sum_{k\leq K} \left|
\sum_{(x,y)\in S_k\cap L_k} \lambda(x) \lambda(y) \right|
\ll \tau(d_0 d_1)\cdot \frac{N^2}{(\log N)^A},\]
where the implicit constant depends only on $A$. 
\end{cor}
\begin{proof}
For every $j\leq K d_1$, there are at most $\tau(d_0 d_1)$ lattice cosets
$L_k$ of index $j$. There are no lattice cosets $R_k$ of index greater than
$K d_1$. The statement then follows
from Lemma \ref{bomb4}.
\end{proof}

\subsection{Anti-sieving}
In the next two lemmas we use an upper-bound sieve not to find almost-primes, but to split the integers multiplicatively, with the almost-primes as an error term. A treatment
by means of a cognate of Vaughan's identity would also be possible, but much more cumbersome.
The error term would be the same.
\begin{lem}\label{sieve}
For any given $M_2>M_1>1$,
there are $\sigma_d\in \mathbb{R}$ with $|\sigma_d|\leq 1$
and support on
\[\{M_1\leq d < M_2 : p<M_1 \Rightarrow p\nmid d\}\]
such that for any $a$, $m$, $N_1$ and $N_2$ with $0\leq m<M_1$ and
$0\leq (N_2 - N_1)/m < M_2$,
\[\mathop{\sum_{N_1\leq n< N_2}}_{n\equiv a \mo m} 
\left|1-\sum_{d|n} \sigma_d\right| \ll
\frac{\log M_1}{\log M_2}
\frac{N_2-N_1}{m}
 + M_2^2,\]
where the implied constant is absolute.
\end{lem}
\begin{proof}
Set $\lambda_d$ as in the Rosser-Iwaniec sieve with sieving set
$\mathfrak{P} = \{\text{$p$ prime} : p\geq M_1, p\nmid m\}$
 and upper cut $z=M_2$. 
Set $\sigma_1=0$, $\sigma_d=-\lambda_d$ for $d\ne 1$. Since
\[\mathop{\sum_{N_1\leq n< N_2}}_{n\equiv a \mo m} 
\left|\sum_{d|n} \lambda_d\right| \ll \frac{\log M_1}{\log M_2}
 \frac{N_2-N_1}{m} ,\]  
the statement follows.
\end{proof}
Note that some of the older combinatorial sieves would be enough for
Lemma \ref{sieve},
provided that $M_2$ were kept greater than a given power of $M_1$.
\begin{lem}\label{sieve2}
Let $K/\mathbb{Q}$ be a number field.
Let $M_2>M_1>1$. 
Let $\jmath:K\to \mathbb{R}^{\deg (K/\mathbb{Q})}$ be a bijective
$\mathbb{Q}$-linear map
taking $\mathfrak{O}_K$ to $\mathbb{Z}^{\deg (K/\mathbb{Q})}$.
Then there are $\sigma_{\mathfrak{d}}
\in \mathbb{R}$ with $|\sigma_{\mathfrak{d}}|\leq 1$
and support on
\begin{equation}\label{eq:schi}
\{\mathfrak{d} : M_1\leq N \mathfrak{d} < M_2,\;
\gcd(\mathfrak{d},\lbrack \mathbb{Z}^2 : L\rbrack) = 1,\; 
(N \mathfrak{p}<M_1 \Rightarrow \mathfrak{p}\nmid \mathfrak{d})\}\end{equation}
such that for any positive integer $N>M_2$, any lattice
coset $L\subset \mathbb{Z}^{\deg(K/\mathbb{Q})}$ with index 
$\lbrack \mathbb{Z}^2 : L\rbrack<M_1$ and
any convex set 
  $S\subset \lbrack -N,N\rbrack^{\deg(K/\mathbb{Q})}$,
\[\sum_{\jmath(x)\in S\cap L} 
\left|1-\mathop{\sum_{\mathfrak{d}}}_{x\in \mathfrak{d}} \sigma_{\mathfrak{d}} \right| \ll
\frac{\log M_1}{\log M_2}
\frac{\Area(S)}{\lbrack \mathfrak{O}_K:L\rbrack}
 + N^{\deg(K/\mathbb{Q})-1} M_2^2,\]
where 
the implied constant depends only on $K$.
\end{lem}
\begin{proof}
Set $\lambda_{\mathfrak{d}}$ as in the generalized 
lower--bound Rosser--Iwaniec sieve 
(\cite{Col2}) with sieving set
$\{\text{$\mathfrak{p}$ prime} : N \mathfrak{p} \geq M_1, (N \mathfrak{p},
\lbrack \mathfrak{O}_K:L\rbrack) = 1\}$ and upper cut $z=M_2$. Set $\sigma_{\mathfrak{O}_K} = 0$, $\sigma_\mathfrak{d}
= -\lambda_\mathfrak{d}$ for $\mathfrak{d}\ne \mathfrak{O}_K$.
\end{proof}
\section{The average of $\lambda$ on integers represented by a 
quadratic form}\label{sec:parquad}

We say that a subset $S$ of $\mathbb{C}$ is a {\em sector} if it is a sector
of $\mathbb{R}^2$ under the natural isomorphism $(x+i y)\mapsto (x,y)$ from
$\mathbb{C}$ to $\mathbb{R}^2$.

\begin{lem}\label{lem:imapp}
Let $K$ be an imaginary quadratic extension of $\mathbb{Q}$. Let 
$\mathfrak{d}\in I_K$, $\chi \in \Xi(C_{\mathfrak{d}}(K))$. 
Let $S$ be a sector of $\mathbb{C}$.
Define the function $\sigma_{S,\chi}:I_{K,\mathfrak{d}}\to \mathbb{Z}$
by
\[\sigma_{S,\chi}(\mathfrak{s}) = 
\mathop{\sum_{s\in \mathfrak{I}^{-1}(\mathfrak{s})}}_{\iota(s)\in S} \chi(s) .\]
Then for any positive $\epsilon$ and any positive integer $k$ there are 
Gr\"ossencharakters \[\{\psi_n\}_{-\infty<n<\infty}\] on 
$I_{K,\mathfrak{d}}$,
sectors $S_1$, $S_2$ of angle $\epsilon$, and
complex numbers $\{c_n\}_{-\infty<n<\infty}$ 
such that
\begin{equation}\begin{aligned}
\sigma_{S,\chi}(\mathfrak{s}) = \sum_{n=-\infty}^\infty 
c_n \psi_n(\mathfrak{s})\;\;\;\;\;
&\text{for every $\mathfrak{s}\in I_{K,\mathfrak{d}}$ with 
$\iota(\mathcal{I}^{-1}(\mathfrak{s}))\cap S_i = \emptyset$,}\\
\left|\sum_{n=-\infty}^{\infty} c_n \psi_n(\mathfrak{s}) \right| \ll 1
\;\;\;\;&\text{for every $\mathfrak{s}\in I_{K,\mathfrak{d}}$,}\\
|c_0|\ll 1,\;\;\;|c_n|\ll (k/\epsilon)^k |n|^{-(k+1)} \;\;\;&\text{for $n\ne 0$.}\end{aligned}\end{equation}
The implied constants are absolute.
\end{lem}
\begin{proof}
For every $s\in \mathcal{I}(\mathfrak{O}_K)_{\mathfrak{d}}$,
\[\sigma_{S,\chi}(\mathfrak{I}(s)) = \sum_{u\in \mathfrak{O}_K^*}
 \lbrack \iota(u s)\in S\rbrack \chi(u s) .\]
Since $S$ is a sector, $\iota(u s)\in S$ if and only if
$\iota(u) \frac{\iota(s)}{|\iota(s)|}\in S$.
Now $S\cap S^1$ is an interval. By Lemma \ref{lem:vin} there are
$\{a_n\}_{n=-\infty}^{\infty}$ such that
\[\begin{aligned}
0\leq \sum_{n=-\infty}^{\infty} a_n x^n \leq 1 \;\;\;&\text{for every
$x\in S^1$,}\\
\sum_{n=-\infty}^{\infty} a_n x^n = \lbrack x\in S \cap S^1\rbrack
\;\;\;\;&\text{if $x\in S^1$, $x\notin S_1, S_2$,}\\
|a_n| \ll (k/\epsilon)^k |n|^{-(k+1)} \;\;\;&\text{for $n\ne 0$},\;\;
|a_0|\ll 1 ,\end{aligned}\]
where $S_1$, $S_2$ are sectors of angle $\epsilon$. 
Hence
\[\sum_{u\in \mathfrak{O}_K^*} \lbrack \iota(u s) \in S\rbrack \chi(u s) =
\sum_{u\in \mathfrak{O}_K^*} \sum_{n=-\infty}^{\infty} a_n 
\left(\iota(u) \frac{\iota(s)}{|\iota(s)|}\right)^n \chi(u s) \]
if $s\notin u S_1, u S_2$ for every $u\in \mathfrak{O}_K^*$. 
 Changing the order of summation,
\[\sum_{u\in \mathfrak{O}_K^*} \sum_{n=-\infty}^{\infty}
a_n \left(\iota(u) \frac{\iota(s)}{|\iota(s)|}\right)^n \chi(u s) =
\sum_{n=-\infty}^{\infty} a_n 
 \sum_{u\in \mathfrak{O}_K^*} 
\iota(u)^n \chi(u) \left(\frac{\iota(s)}{|\iota(s)|}\right)^n \chi(s) .\]
We will have $\sum_{u\in \mathfrak{O}_K^*} u^n \chi(u)\ne 0$ only when
$u^n \chi(u) = 1$ for all $u\in \mathfrak{O}_K^*$. Then
there is a Gr\"ossencharakter $\psi_n$ such
that
\[\psi_n(\mathfrak{s}) = \left(\frac{\iota(s)}{|\iota(s)|}\right)^n \chi(s) \]
for every $s\in \mathcal{I}^{-1}(\mathfrak{s})$.
Hence
\[\sum_{n=-\infty}^{\infty} a_n
 \sum_{u\in \mathfrak{O}_K^*} 
\iota(u)^n \chi(u) \left(\frac{\iota(s)}{|\iota(s)|}\right)^n \chi(s) 
 = \mathop{\sum_{-\infty<n<\infty}}_{\iota(u)^n \chi(u)=1}
(\# \mathfrak{O}_K^*) a_n  \psi_n(\mathfrak{s}) .\]
Set \[c_n = \begin{cases} 
(\# \mathfrak{O}_K^*) a_n &\text{if 
$\iota(u)^n \chi(u)=1$,}\\ 0 &\text{otherwise.}\end{cases}\]
\end{proof}
Let the sector $S\subset \mathbb{R}^2$ be a subquadrant.
Define
\[\begin{aligned} \rho(x,y) &= x/y\\
\gamma_-(S) &= \inf_{(x,y)\in S} x/y ,\\
\gamma_+(S) &= \sup_{(x,y)\in S} x/y .\\
\end{aligned}\]
If $S$ is a subquadrant, $\gamma_-(S)$ and $\gamma_+(S)$ are finite non-zero
real numbers of the same sign. Moreover, $(x,y)\in S$ if and only if
$\rho(x,y) \in (\gamma_-(S),\gamma_+(S))$. The sign $\sgn(x)$ is the same for
all $x\in S$. We call it $\sgn(S)$ and define
\[H_S = \{(x,y)\in \mathbb{R}^2 : \sgn(x) = \sgn(S) \} .\]

For $K/\mathbb{Q}$ a real quadratic extension, let 
 $\iota:\mathcal{I}(K)\to \mathbb{R}^2$
be the embedding given by $\iota(a) = (\iota_1(a),\iota_2(a))$. 
\begin{lem}\label{lem:reapp}
Let $K$ be a real quadratic extension of
$\mathbb{Q}$. Let $\mathfrak{d}\in I_K$,
$\chi\in \Xi(C_{\mathfrak{d}}(K))$. Let $S$ be a subquadrant of $\mathbb{R}^2$.
Define the function $\sigma_{S,\chi}:I_{K,\mathfrak{d}}\to \mathbb{Z}$ by
\[\sigma_{S,\chi}(\mathfrak{s}) = 
\mathop{\sum_{s\in \mathfrak{I}^{-1}(\mathfrak{s})}}_{\iota(s)\in S} \chi(s) .\]
Then for any positive $\epsilon$ and any positive integer $k$ there are
Gr\"ossencharakters \[\{\psi_n\}_{-\infty<n<\infty}\] on $I_{K,\mathfrak{d}}$,
sectors $S_1$, $S_2$ of hyperbolic angle at most $\epsilon$, 
and complex numbers $\{c_n\}_{-\infty<n<\infty}$
such that
\[\begin{aligned}
\sigma_{S,\chi}(\mathfrak{s}) &= \sum_{n=-\infty}^\infty 
c_n \psi_n(\mathfrak{s})\;\;\;\;\;
\text{for every $\mathfrak{s}\in I_{K,\mathfrak{d}}$ with 
$\iota(\mathcal{I}^{-1}(\mathfrak{s}))\cap S_i = \emptyset$,}\\
\left|\sum_{n=-\infty}^{\infty} c_n \psi_n(\mathfrak{s}) \right| &\ll 
\frac{|\log(\gamma_+(\iota(S))/\gamma_-(\iota(S)))|}{
|\log(\iota_1(u_1)/\iota_2(u_1))|} + k_{\mathfrak{d}}
\;\;\;\;\text{for every $\mathfrak{s}\in I_{K,\mathfrak{d}}$,}\\
|c_0|, |c_1| &\ll \frac{|\log(\gamma_+(\iota(S))/\gamma_-(\iota(S)))|}{
|\log(\iota_1(u_1)/\iota_2(u_1))|}\\
|c_n|&\ll (k k_{\mathfrak{d}}/\epsilon)^k 
|n|^{-(k+1)} \;\;\;\text{for $n\ne 0, 1$,}\end{aligned}\]
 where $u_1$, $\iota_1$ and $\iota_2$ are as in 
subsection \ref{subs:idgross}.
The implied constants are absolute.
\end{lem}
\begin{proof}
For every $s\in \mathcal{I}(\mathfrak{O}_K)_{\mathfrak{d}}$ with $\iota(s)\in H_S$,
\[
\sigma_{S,\chi}(\mathfrak{I}(s)) = \sum_{u\in \mathfrak{O}_K^*}
 \lbrack \iota(u s) \in S \rbrack \chi(u s) .\]
Since $\iota_1(u_1)$
is positive, $\iota(s)\in H_S$ implies
$\iota(u_1^k s)\in H_S$, $\iota(- u_1^k s)\notin H_S$ 
for every $k\in \mathbb{Z}$.
Hence
\[\begin{aligned}
\sigma_{S,\chi}(\mathfrak{I}(s)) &=
\sum_{k=-\infty}^{\infty} \lbrack \iota(u_1^k s) \in S\rbrack \chi(u^k s) \\
&= \sum_{k=-\infty}^{\infty}
 \lbrack \iota(u_1^k s)\in (S \cap (-S))\rbrack \chi(u_1^k s)\\
&= \sum_{k=-\infty}^{\infty} \lbrack \rho(\iota(u_1^k s)) \in
 (\gamma_-(S),\gamma_+(S))\rbrack \chi(u_1^k s) .\end{aligned}\]
Let $k_{\mathfrak{d}}$, $l_{\mathfrak{d}}$, $r_{\mathfrak{d}}$ be as in
section \ref{subs:idgross}.
Let $C_0$ is the largest integer smaller than
$\frac{|\log (\gamma_+(\iota(S))/\gamma_-(\iota(S)))|}{\log l_{\mathfrak{d}}}$.
Let $\gamma_0 = \gamma_-(S) l_{\mathfrak{d}}^{C_0}$. 
Then 
\[\sigma_{S,\chi}(\mathfrak{I}(s)) = \chi(s) \left(C_0 k_{\mathfrak{d}}
\lbrack \chi(u_1) = 1\rbrack +
\sum_{k=-\infty}^\infty \lbrack \rho(\iota(u_1^k s))\in 
 (\gamma_0,\gamma_+(S))\rbrack \chi(u_1^k)\right) .
\]
Assume $\sgn(\rho(\iota(s))) = \sgn(\gamma_0)$.
Then
there is exactly one integer $n$ such that 
$l_{\mathfrak{d}}^n s \in (\gamma_0, l_{\mathfrak{d}} \gamma_0\rbrack$. 
Let $\phi : \mathbb{R}^*\to S^1$ be given by
\[\phi(r) = e^{2\pi i \frac{\log |r|}{\log l_{\mathfrak{d}}}} .\]
Define $\Phi = \phi \circ \rho \circ \iota:K\mapsto S^1$. Then
\[\sum_{k=-\infty}^\infty \lbrack \rho(\iota(u_1^k s))\in 
 (\gamma_0,\gamma_+(S))\rbrack \chi(u_1^k)=\sum_{k=0}^{k_{\mathfrak{d}}-1} \lbrack 
\Phi(u_1^{r_{\mathfrak{d}} k} s) \in (\phi(\gamma_0),\phi(\gamma_+(S)))
\rbrack \chi(u_1^{r_{\mathfrak{d}} k}) .\]
By Lemma \ref{lem:vin} there are
$\{a_n\}_{n=-\infty}^{\infty}$ such that
\[\begin{aligned}
0\leq \sum_{n=-\infty}^{\infty} a_n x^n \leq 1 \;\;\;&\text{for every
$x\in S^1$,}\\
\sum_{n=-\infty}^{\infty} a_n x^n = \lbrack x\in S \cap S^1\rbrack
\;\;\;\;&\text{if $d(x,\gamma_0), d(x,\gamma_+(S)) \geq \epsilon/2 k_{\mathfrak{d}}$,}\\
|a_n| \ll (k k_{\mathfrak{d}} /\epsilon)^k |n|^{-(k+1)} \;\;\;&\text{for $n\ne 0$},\;\;
|a_0|\ll 1 .\end{aligned}\] 
Hence
\[
\sigma_{S,\chi}(\mathfrak{I}(s)) = \chi(s) \left(C_0 k_{\mathfrak{d}} 
\lbrack \chi(u_1) = 1\rbrack +
\sum_{n=-\infty}^{\infty} a_n \left(
\sum_{k=0}^{k_{\mathfrak{d}}-1}
\Phi(u_1^{r_{\mathfrak{d}} k})^n 
 \chi(u_1^{r_{\mathfrak{d}} k}) \right) 
\Phi(s)^n\right)  ,\]
provided that $d(\Phi(u_1^{r_{\mathfrak{d}} k} s),\gamma_0)\geq \epsilon/2$,
$d(\Phi(u_1^{r_{\mathfrak{d}} k} s),\gamma_+(S))\geq \epsilon/2$ for
every non-negative $k$ less than $k_{\mathfrak{d}}$.
We will have
\begin{equation}\label{eq:asterix}
\sum_{k=0}^{k_{\mathfrak{d}}-1}
\Phi(u_1^{r_{\mathfrak{d}} k})^n 
 \chi(u_1^{r_{\mathfrak{d}} k}) \ne 0\end{equation}
 only when
$\Phi(u_1^{r_{\mathfrak{d}}})^n \chi(u_1^{r_{\mathfrak{d}}}) = 1$. 

Suppose $\frac{\iota_1(u_1)}{\iota_2(u_1)}<0$. Then
there is a Gr\"ossencharakter 
\[\psi_n(\mathfrak{s}) = \chi(s) \sgn(\rho(\iota(s)))^{n_0} \Phi(s)^n,\]
where \[n_0(n) = \begin{cases}1 &\text{if 
$\chi(u_1) \Phi(u_1)^n = 1$}\\
-1 &\text{if $\chi(u_1) \Phi(u_1)^n = -1$.}\end{cases}\] 
Let 
\[c_n = a_n k_{\mathfrak{d}} \lbrack \Phi(u_1^2)^n
 \chi(u_1^2)=1\rbrack \sgn(\gamma_0)^{n_0(n)} + 
C_0 k_{\mathfrak{d}} \lbrack \chi(u_1)=1\rbrack \lbrack n = 0\rbrack .\]
Thus
\[\sigma_{S,\chi}(\mathfrak{I}(s)) = \sum_{n=-\infty}^{\infty}
c_n \psi_n(\mathfrak{I}(s))  \]
for every $s\in \mathcal{I}(\mathfrak{O}_K)_{\mathfrak{d}}$ with
$\iota(s) \in H_S$, $\sgn(\rho(\iota(s))) = \sgn(\gamma_0)$ and
\[d(\Phi(u_1^{r_{\mathfrak{d}} k} s),\gamma_0)\geq \epsilon/2 k_{\mathfrak{d}},\;\;\;
d(\Phi(u_1^{r_{\mathfrak{d}} k} s),\gamma_+(S))\geq \epsilon/2 k_{\mathfrak{d}}\]
for every $0\leq k<k_{\mathfrak{d}}$.
Since $\frac{\iota_1(u_1)}{\iota_2(u_1)}<0$, for every 
$s\in \mathcal{I}(\mathfrak{O}_K)_{\mathfrak{d}}$ there is a 
$u\in \mathfrak{O}_K^*$ such that $\iota(u s)\in H_S$,
 $\sgn(\rho(\iota(u s))) = \sgn(\gamma_0)$.
Hence
\[\sigma_{S,\chi}(\mathfrak{s}) = \sum_{n=-\infty}^{\infty} 
c_n \psi_n(\mathfrak{s})\]
provided that $d(\Phi(s),\gamma_0)\geq \epsilon/2 k_{\mathfrak{d}}$,
$d(\Phi(s),\gamma_+(S))\geq \epsilon/2 k_{\mathfrak{d}}$
for every $s\in \mathcal{I}^{-1}(\mathfrak{s})$.

Suppose now $\frac{\iota_1(u_1)}{\iota_2(u_1)} > 0$. 
We have (\ref{eq:asterix}) only when $\Phi(u_1) \chi(u_1)=1$.
Then there are Gr\"ossencharakters
\[\begin{aligned}
 \psi_{n +}(\mathfrak{s}) &= \chi(s) \Phi(s)^n,\\
 \psi_{n -}(\mathfrak{s}) &= \chi(s) \sgn(\rho(\iota(s)))
\Phi(s)^n \end{aligned} .\]
Let
\[\begin{aligned}
c_{n +} &= a_n k_{\mathfrak{d}} \lbrack \Phi(u_1)^n 
\chi(u_1) = 1\rbrack  + 
C_0 k_{\mathfrak{d}} \lbrack \chi(u_1) = 1\rbrack \lbrack n = 0\rbrack\\
c_{n -} &= (a_n k_{\mathfrak{d}} \lbrack \Phi(u_1)^n 
\chi(u_1) = 1\rbrack  + 
C_0 k_{\mathfrak{d}} \lbrack \chi(u_1) = 1\rbrack \lbrack n = 0\rbrack)
\sgn(\gamma_0) 
.\end{aligned}\]
Then
\begin{equation}\label{eq:manli}
\sigma_{S,\chi}(\mathfrak{I}(s)) = \sum_{n=-\infty}^{\infty}
 \frac{1}{2} (c_{n +} \psi_{n +}(\mathfrak{I}(s)) + c_{n -} \psi_{n -}(\mathfrak{I}(s)))\end{equation}
for every $s\in \mathcal{I}_{K,\mathfrak{d}}$ with $\iota(s)\in H_S$ and
\[d(\Phi(u_1^{r_{\mathfrak{d}} k} s),\gamma_0)\geq \epsilon/2 k_{\mathfrak{d}},\;\;\;
d(\Phi(u_1^{r_{\mathfrak{d}} k} s),\gamma_+(S))\geq \epsilon/2 k_{\mathfrak{d}}\]
for every $0\leq k<k_{\mathfrak{d}}$.
If $\sgn(\rho(\iota(s)))\ne \sgn(\gamma_0)$, both sides of 
(\ref{eq:manli})
 are equal to zero. Hence, for every $\mathfrak{s}\in I_{K,\mathfrak{d}}$,
\[\sigma_{S,\chi}(\mathfrak{s}) = \sum_{n=-\infty}^{\infty}
 \frac{1}{2} (c_{n +} \psi_{n +}(\mathfrak{s}) + c_{n -} \psi_{n -}(\mathfrak{s})) \]
provided that $d(\Phi(s),\gamma_0)\geq \epsilon/2 k_{\mathfrak{d}}$,
$d(\Phi(s),\gamma_+(S))\geq \epsilon/2 k_{\mathfrak{d}}$
for every $s\in \mathcal{I}^{-1}(\mathfrak{s})$.

Now let $s\in \mathcal{I}(\mathfrak{O}_K)_{\mathfrak{d}}$ be given with
\[d(\Phi(s),\gamma_0)<\epsilon/2 k_{\mathfrak{d}} .\]
Then
\[\left|\frac{\log |\rho(s)|}{\log l_{\mathfrak{d}}} - x\right| <
\epsilon/2 k_{\mathfrak{d}}
\] for some $x\in \phi^{-1}(\gamma_0)$. Let us be given $\mathfrak{s} \in
I_{K,\mathfrak{d}}$. Then
\[d(\Phi(s),\gamma_0) < \epsilon/ 2 k_{\mathfrak{d}} \;\;\;
\text{for some $s\in \mathcal{I}^{-1}(s)$}\]
if and only if
\begin{equation}\label{eq:holgra}
\left|\frac{\log|\rho(s)|}{\log l_{\mathfrak{d}}} - x_0 \right| <
 \epsilon/ 2 k_{\mathfrak{d}} \;\;\;\text{for some $s\in \mathcal{I}^{-1}(
\mathfrak{s})$,}
\end{equation}
where $x_0$ is any fixed element of $\phi^{-1}(\gamma_0)$. Clearly
(\ref{eq:holgra}) is equivalent to
\[(x_0 - \epsilon/ 2 k_{\mathfrak{d}}) \log l_{\mathfrak{d}} <
\log |\rho(s)| < (x_0 + \epsilon/2 k_{\mathfrak{d}}) \log l_{\mathfrak{d}},\]
that is, 
\[x_0 \log l_{\mathfrak{d}} - \frac{\epsilon r_{\mathfrak{d}}}{2}
\log \left(\frac{\iota_1(u_1)}{\iota_2(u_1)}\right) 
< \log |\rho(s)| <
 x_0 \log l_{\mathfrak{d}} +
\frac{\epsilon r_{\mathfrak{d}}}{2} \log 
\left(\frac{\iota_1(u_1)}{\iota_2(u_1)}\right) .\]
Thus $S$ is constrained to a section of hyperbolic angle 
$\frac{\epsilon r_{\mathfrak{d}}}{2} \log 
\left(\frac{\iota_1(u_1)}{\iota_2(u_1)}\right)$. The statement follows.
\end{proof}
Let $Q(x,y)$ be a primitive, irreducible quadratic form. Let $K = 
\mathbb{Q}(\sqrt{b^2 - 4 a c})$. We define $\phi_Q:\mathbb{Q}^2\to K$ to be
the map given by
\[\phi_Q(x,y) = \alpha_1 x + \alpha_2 y,\]
where $\alpha_1$, $\alpha_2$ are as in Lemma \ref{lem:quad}. As before, we define
\[\iota(s) = \begin{cases} (\iota_1(s),\iota_2(s))\in \mathbb{R}^2 
                                        &\text{if $K$ is real}\\
 \iota_1(s)\in \mathbb{C} \sim \mathbb{R}^2 &\text{if $K$ is imaginary.}\end{cases}\]
for $s\in \mathcal{I}(K)$.
The map $\iota \circ \phi_Q:\mathbb{Q}^2 \to \mathbb{R}^2$ is linear. For any sector
$S$ of $\mathbb{R}^2$, there is a sector $S_Q$ of $\mathbb{R}^2$ such that 
$(\iota \circ \phi_Q)(S \cap \mathbb{Q}^2) = S_Q \cap \iota(K)$. 

We recall the definition of $\sigma_{S,\chi} : I_{K,\mathfrak{d}}\to \mathbb{Z}$ in 
the statements of Lemmas \ref{lem:imapp} and \ref{lem:reapp}.
\begin{lem}\label{lem:seth}
Let $Q(x,y) = a x^2 + b x y + c y^2\in \mathbb{Z}\lbrack x,y\rbrack$ be a primitive,
irreducible quadratic form. Let $K = \mathbb{Q}(\sqrt{b^2 - 4 a c})$. Let 
$L\subset \mathbb{Z}^2$ be a lattice coset, $S\subset \mathbb{R}^2$ a sector.
If $K$ is real, assume $S_Q$ is a subquadrant. Let 
$d = a\cdot d_{\sq(b^2 - 4 a c)}
\lbrack \mathbb{Z}^2 : L\rbrack$. Then there are sectors 
$\{S_{\mathfrak{r}}\}_{\mathfrak{r}|d^{\infty}}$, $S_{\mathfrak{r}}\subset \mathbb{R}^2$,
and complex numbers 
$\{a_{\mathfrak{r} \chi}\}_{r|d^{\infty}, \chi\in \Xi(C_d(K))}$,
$|a_{\mathfrak{r} \chi}| \leq \frac{d}{\# C_d(K)}$, such that
\[\#\{x,y\in S\cap L:|Q(x,y)|=m\} = 
\mathop{\sum_{\mathfrak{r}}}_{
N \mathfrak{r} = \gcd(a m,d^{\infty})}
\sum_{\chi\in \Xi(C_d(K))} a_{\mathfrak{r} \chi} \mathop{\sum_{\mathfrak{s}}}_{
 N \mathfrak{s} = \frac{|a m|}{\gcd(a m, d^{\infty})}} 
\sigma_{S_{\tau},\chi}(\mathfrak{s})\]
for every positive integer $m$.  If $K$ is real, then,
for every $\mathfrak{r}|d^{\infty}$, $S_{\mathfrak{r}}$ is a subquadrant
satisfying \[
|\log (\gamma_+(S_{\mathfrak{r}})/\gamma_-(S_{\mathfrak{r}}))| =
|\log (\gamma_+(S_Q)/\gamma_-(S_Q))|.\] 
\end{lem}
\begin{proof}
By Lemma \ref{lem:quad},
\[\#\{x,y\in S\cap L: |Q(x,y)| = m\} = 
    \mathop{\mathop{\sum_{s\in \phi_Q(L)}}_{\iota(s)\in S_Q}}_{|N
s| = |a m|} 1.\]
For every $s\in \mathfrak{O}_K$ of norm $N s = \pm a m$, there
is exactly one ideal of norm $\gcd(a m,d^{\infty})$ containing $x$.
Hence
\[\mathop{\mathop{\sum_{s\in \phi_Q(L)}}_{\iota(s)\in S_Q}}_{|N s| = |a m|} 1
= \mathop{\sum_{\mathfrak{r}}}_{N \mathfrak{r} = \gcd(a m,d^{\infty})}
\mathop{\mathop{\sum_{s\in \phi_Q(L)\cap \mathfrak{r}}}_{\iota(s)\in S_Q}}_{
|N s| = |a m|} 1.
\]
Since, by Lemma \ref{lem:quad}, $\phi_Q(L)$ is an additive
 subgroup of $\mathfrak{O}_K$ of index
$d = \lbrack \mathfrak{O}_K : \phi_Q(L)\rbrack = 
a\cdot d_{sq(b^2 - 4 a c)}\lbrack \mathbb{Z}^2 :
L \rbrack$, $\phi_Q(L)\cap \mathfrak{r}$ is an additive
 subgroup of $\mathfrak{r}$ of index dividing
$d$. Therefore, whether or not a given $s\in \mathfrak{r}$ is an element of 
$\phi_Q(L) \cap \mathfrak{r}$ depends only on $s \mo d \mathfrak{r}$. If
$N \mathfrak{r} = \gcd(a m,d^{\infty})$ and 
$N s = a m$, then 
  $N \mathfrak{r} = \gcd(N s, d^{\infty})$,
and so, given that $s\in \mathfrak{r}$, 
$N s / N r$ is prime to $d$. Choose $r\in \mathfrak{I}^{-1}(
\mathfrak{r})$. Then $s/r\in \mathcal{I}(\mathfrak{O}_K)_{d}$. Moreover,
whether or not $s$ is an element of $\phi_Q(L) \cap \mathfrak{r}$ 
depends only on
the equivalence class $\langle s/r\rangle$ 
of $s/r$ in $C_d(K)$. In other words, there is a subset
$C_{\mathfrak{r}}$ of $C_d(K)$ such that $x/r\in \phi_Q(L) \cap \mathfrak{r}$ if and only
if $\langle x/r\rangle \in C_{\mathfrak{r}}$. Then
\[\#\{x,y\in S\cap L : |Q(x,y)| = m \} = \mathop{\sum_{\mathfrak{r}}}_{N \mathfrak{r}
= \gcd(a m,d^{\infty})} \mathop{\mathop{
\mathop{\sum_{s\in \mathcal{I}(\mathfrak{O}_K)_d}}_{\langle s\rangle\in C_{\mathfrak{r}}}
}_{\iota(r s)\in S_Q}}_{|N(\mathcal{I}(s))| = 
a m/\gcd(a m,r^{\infty})} 1 .\]
For $\chi\in \Xi(C_d(K))$, let $a_{\mathfrak{r} \chi} = \frac{1}{\# C_d(K)}
\sum_{c\in C_{\mathfrak{r}}} \overline{\chi(c)} $. Then
\[\lbrack \langle s\rangle \in C_{\mathfrak{r}}\rbrack = \sum_{\chi\in \Xi(C_d(K))}
 a_{\mathfrak{r} \chi} \chi(s) .\]
Hence $\#\{x,y\in S\cap L:|Q(x,y)|=m\}$ equals
\[\begin{aligned}
\mathop{\sum_{\mathfrak{r}}}_{
N \mathfrak{r} = \gcd(a m,d^{\infty})}
&\sum_{\chi\in \Xi(C_d(K))} a_{\mathfrak{r} \chi} 
\mathop{\mathop{\sum_{s\in \mathcal{I}\left(\mathfrak{O}_K\right)_d}
 }_{\iota(r s)\in S}}_{N(\mathcal{I}(s)) =
|a m|/\gcd(|a m|,d^{\infty})} \chi(s)\\
&= \mathop{\sum_{\mathfrak{r}}}_{
N \mathfrak{r} = \gcd(a m,d^{\infty})}
\sum_{\chi\in \Xi(C_d(K))} a_{\mathfrak{r} \chi} 
\mathop{\sum_{\mathfrak{s}}}_{
 N \mathfrak{s} = \frac{|a m|}{\gcd(a m, d^{\infty})}} 
\sigma_{S_{\tau},\chi}(\mathfrak{s}),\end{aligned}\]
where 
\[S_{\mathfrak{r}} = \begin{cases}
 \iota(r)^{-1} S_Q &\text{if $K$ is imaginary,}\\
 \{(x,y)\in \mathbb{R}^2 : (\iota_i(r) x,\iota_2(r) y)\in S_Q\}
 &\text{if $K$ is real.}\end{cases}\]
\end{proof}

\begin{lem}\label{lem:kuhl}
Let $K$ be a quadratic extension of $\mathbb{Q}$. Let $a$ be a non-zero 
rational
integer. Then, for any rational integer $r$ dividing $a$, any ideal
$\mathfrak{d}\in I_{K,r}$ of norm 
\begin{equation}\label{eq:peng} N \mathfrak{d} \ll (\log N)^A\end{equation} 
and any Gr\"ossencharakter $\psi$ on $I_{K,\mathfrak{d}}$
of size $\mathcal{S}(\psi) \ll e^{(\log x)^{3/5} (\log \log x)^{1/5}}$,
we have
\begin{equation}\label{eq:focus}
\mathop{\mathop{\sum_{\mathfrak{s}\in I_{K,\mathfrak{d}}}}_{N \mathfrak{s} \leq x}}_{
r | N \mathfrak{s}} \psi(\mathfrak{s}) \lambda(N \mathfrak{s}) 
\ll x e^{-C \frac{(\log x)^{2/3}}{(\log \log x)^{1/5}}},
\end{equation}
where $C$ and the implied constant in (\ref{eq:focus})
depend only on $K$, $A$, $r$, and the implied constant in (\ref{eq:peng}).
\end{lem}
\begin{proof}
For any $\mathfrak{s}\in I_{K,\mathfrak{d}}$,
\[\begin{aligned}
\lbrack r | N \mathfrak{s}\rbrack &= \lbrack r | N(\gcd(\mathfrak{s},r^\infty))
\rbrack 
= 1 - \mathop{\sum_{\mathfrak{r} | r^{\infty}}}_{r \nmid N \mathfrak{r}}
 \lbrack \mathfrak{r} = \gcd(\mathfrak{s}, r^{\infty})\rbrack\\
&= 1 - \mathop{\sum_{\mathfrak{r} | r^{\infty}}}_{r \nmid N \mathfrak{r}}
 \sum_{\mathfrak{m} | \rad(r)} \mu_K(\mathfrak{m}) \lbrack \mathfrak{r} 
\mathfrak{m} | \gcd(\mathfrak{s}, r^{\infty})\rbrack \\
&= 1 - \mathop{\sum_{\mathfrak{r} | r^{\infty}}}_{r \nmid N \mathfrak{r}}
 \sum_{\mathfrak{m} | \rad(r)} \mu_K(\mathfrak{m}) \lbrack \mathfrak{r} 
\mathfrak{m} | \mathfrak{s} \rbrack .\end{aligned}\]
Hence 
\begin{equation}\label{eq:wrg}
\mathop{\mathop{\sum_{\mathfrak{s}\in I_{K,\mathfrak{d}}}}_{N \mathfrak{s} \leq x}}_{
r | N \mathfrak{s}} \psi(\mathfrak{s}) \lambda(N \mathfrak{s}) 
= \mathop{\sum_{\mathfrak{s}\in I_{K,\mathfrak{d}}}}_{N \mathfrak{s} \leq x}
 \psi(\mathfrak{s}) \lambda(N \mathfrak{s}) 
- \mathop{\sum_{\mathfrak{r} | r^{\infty}}}_{r\nmid N \mathfrak{r}}
 \sum_{\mathfrak{m} | \rad(r)} \mu_K(\mathfrak{m})
 \mathop{\mathop{\sum_{\mathfrak{s}\in I_{K,\mathfrak{d}}}}_{N \mathfrak{s} 
\leq x}}_{\mathfrak{r} \mathfrak{m} | \mathfrak{s}}
 \psi(\mathfrak{s}) \lambda(N \mathfrak{s}) .\end{equation}
We can rewrite the second term on the right side of (\ref{eq:wrg}) as
\[ \mathop{\sum_{\mathfrak{r} | r^{\infty}}}_{r\nmid N \mathfrak{r}}
 \sum_{\mathfrak{m} | \rad(r)} \mu_K(\mathfrak{m})
\psi(\mathfrak{r} \mathfrak{m}) \lambda(N(\mathfrak{r} \mathfrak{m}))
\mathop{\sum_{\mathfrak{s}\in I_{K,\mathfrak{d}}}}_{N \mathfrak{s} 
\leq x/N(\mathfrak{r} \mathfrak{m})}
 \psi(\mathfrak{s}) \lambda(N \mathfrak{s}) .\]
The statement now follows from Lemma \ref{lem:lamp}.
\end{proof}
\begin{lem}\label{lem:chai}
Let $K$ be a finite extension of $\mathbb{Q}$. Let $d$ be a non-zero
rational integer. Then
\[\mathop{\mathop{\sum_{\mathfrak{r}\in I_K}}_{\mathfrak{r}|d^{\infty}}}_{
X^{1/2} < N \mathfrak{r} \leq X} \frac{1}{N \mathfrak{r}} \ll
 \frac{(\log X)^C}{X^{1/2}} ,\]
where $C$ and the implied constant depend only on $K$ and $d$.
\end{lem}
\begin{proof}
Let $\mathfrak{p}$ be the divisor of $d$ of largest norm. Every 
$\mathfrak{r}\in I_K$ with $\mathfrak{r} | d^{\infty}$ and
$N \mathfrak{r} > X^{1/2}$ has a divisor $\mathfrak{d}|\mathfrak{r}$
of norm $X^{1/2} < N \mathfrak{d} \leq X^{1/2} N \mathfrak{p}$. Hence
\[\begin{aligned}
\mathop{\mathop{\sum_{\mathfrak{r}\in I_K}}_{\mathfrak{r} | d^{\infty}}}_{
X^{1/2} < N \mathfrak{r} \leq X} \frac{1}{N \mathfrak{r}} &\leq
 \mathop{\mathop{\sum_{\mathfrak{d}\in I_K}}_{\mathfrak{d} | d^{\infty}}}_{
X^{1/2} < N \mathfrak{d} \leq X^{1/2} N \mathfrak{p}} 
\frac{1}{N \mathfrak{d}} \mathop{\sum_{\mathfrak{a}\in I_K}}_{N \mathfrak{a} \leq X^{1/2}} \frac{1}{N \mathfrak{a}} \\
&\ll (\log X)^{c_1} \frac{1}{X^{1/2}} 
\mathop{\mathop{\sum_{\mathfrak{d}\in I_K}}_{\mathfrak{d} | d^{\infty}}}_{
 N \mathfrak{d} \leq X^{1/2} N \mathfrak{p}} 1 
\ll (\log X)^{c_1} \frac{1}{X^{1/2}} (\log X)^{c_2} .\end{aligned}\]
\end{proof}
\begin{lem}\label{lem:vira}
Let $Q(x,y) = a x^2 + b x y + c y^2\in \mathbb{Z}\lbrack x,y\rbrack$ be a primitive,
irreducible quadratic form. Let $K = \mathbb{Q}(\sqrt{b^2 - 4 a c})$. Let 
$L\subset \mathbb{Z}^2$ be a lattice coset, $S\subset \mathbb{R}^2$ a sector.
Assume 
\begin{equation}\label{eq:cho1}
\lbrack \mathbb{Z}^2 : L\rbrack \ll (\log X)^A .\end{equation}
If $K$ is real, assume $S_Q$ is a subquadrant satisfying
\begin{equation}\label{eq:cho2}|\log(\gamma_+(S_Q)/\gamma_-(S_Q))|\ll (\log X)^A.\end{equation}
Then
\begin{equation}\label{eq:cho3}
\mathop{\sum_{x,y\in S\cap L}}_{|Q(x,y)| \leq X} \lambda(Q(x,y)) \ll
X e^{-C \frac{(\log X)^{2/3}}{(\log \log X)^{1/5}}} ,\end{equation}
where $C$ and the implied constant depend on $a$, $b$, $c$, $A$
and the implied constants in (\ref{eq:cho1}) and (\ref{eq:cho2}).
\end{lem}
\begin{proof}
By Lemma \ref{lem:seth},
\[\mathop{\sum_{x,y\in S\cap L}}_{|Q(x,y)| \leq X} \lambda(Q(x,y)) =
 \sum_{m\leq X} 
\mathop{\sum_{\mathfrak{r}}}_{
N \mathfrak{r} = \gcd(a m,d^{\infty})}
\sum_{\chi\in \Xi(C_d(K))} a_{\mathfrak{r} \chi} \mathop{\sum_{\mathfrak{s}}}_{
 N \mathfrak{s} = \frac{|a m|}{\gcd(a m, d^{\infty})}} 
\sigma_{S_{\mathfrak{r}},\chi}(\mathfrak{s}) \lambda(m) ,\]
where $d = a\cdot d_{\sq(b^2 - 4 a c)} \lbrack \mathbb{Z}^2 : L\rbrack$.
Since $a_{\mathfrak{r} \chi} \leq \frac{d}{\# C_d(K)}$, it will be enough to bound
\begin{equation}\label{eq:canario}\sum_{m\leq X} 
\mathop{\sum_{\mathfrak{r}}}_{
N \mathfrak{r} = \gcd(a m,d^{\infty})}
\mathop{\sum_{\mathfrak{s}}}_{
 N \mathfrak{s} = \frac{|a m|}{\gcd(a m, d^{\infty})}} 
\sigma_{S_{\tau},\chi}(\mathfrak{s}) \lambda(m) .
\end{equation}
We will take $\epsilon$ to be a positive number whose value we shall set
later. By Lemmas \ref{lem:imapp} and \ref{lem:reapp} with $k=1$,
\[\begin{aligned}
\sigma_{S,\chi}(\mathfrak{s}) &= \sum_{n=-\infty}^\infty 
c_n \psi_n(\mathfrak{s})\;\;\;\;\;
\text{for every $\mathfrak{s}\in I_{K,\mathfrak{d}}$ with 
$\iota(\mathcal{I}^{-1}(\mathfrak{s}))\cap S_i = \emptyset$,}\\
\left|\sum_{n=-\infty}^{\infty} c_n \psi_n(\mathfrak{s}) \right| &\ll 
\frac{|\log(\gamma_+(\iota(S))/\gamma_-(\iota(S)))|}{
|\log(\iota_1(u_1)/\iota_2(u_1))|} + k_{\mathfrak{d}}
\;\;\;\;\text{for every $\mathfrak{s}\in I_{K,\mathfrak{d}}$,}\end{aligned}\]
where $S_1$, $S_2$ are sectors of angle at most
$\epsilon$ (if $K/\mathbb{Q}$ is imaginary) or of hyperbolic angle at most
$\epsilon$ (if $K/\mathbb{R}$ is real), and
\[
 |c_n| \ll \frac{|n|^{-2}}{\epsilon} \;\;\;\text{for $K/\mathbb{Q}$                                                     imaginary,}\]
\[|c_n| \ll \frac{k_{\mathfrak{d}}}{\epsilon} |n|^{-2} \;\;\;\text{for $K/\mathbb{Q}$
real, $n\ne 0,1$,}\]
\[|c_0|, |c_1| \ll \max\left(1,
\frac{|\log(\gamma_+(\iota(S))/\gamma_-(\iota(S)))|}{
|\log(\iota_1(u_1)/\iota_2(u_1))|}\right) \;\;\;\text{for $K/\mathbb{Q}$ real.}\]
Let $B$ be a large number whose value will be set later. Since 
$|\psi_n(\mathfrak{s})|=1$, $d\ll (\log N)^A$ and $C_0\ll (\log N)^A$,
 the absolute value of the difference between (\ref{eq:canario}) and
\begin{equation}\label{eq:comi}
\sum_{m\leq X}
\mathop{\sum_{\mathfrak{r}}}_{
N \mathfrak{r} = \gcd(a m,d^{\infty})}
\mathop{\sum_{\mathfrak{s}}}_{
 N \mathfrak{s} = \frac{|a m|}{\gcd(a m, d^{\infty})}} 
\sum_{|n|\leq B} c_n \psi_n(\mathfrak{s})
\lambda(m) \end{equation}
is at most a constant times 
\[\sum_{m\leq X} 
\mathop{\sum_{\mathfrak{r}}}_{
N \mathfrak{r} = \gcd(a m,d^{\infty})}
\mathop{\sum_{\mathfrak{s}}}_{
 N \mathfrak{s} = \frac{|a m|}{\gcd(a m, d^{\infty})}} 
\frac{k_{\mathfrak{d}}}{B\epsilon} \ll \sum_{m\leq X} \frac{k_{\mathfrak{d}} \tau(m)}{B \epsilon} 
\ll \frac{k_{\mathfrak{d}} X \log X}{B \epsilon} .
\]
By (\ref{eq:cho1}), the absolute value of
\[
\sum_{m\leq X}
\mathop{\sum_{\mathfrak{r}}}_{
N \mathfrak{r} = \gcd(a m,d^{\infty})}
\mathop{\sum_{\mathfrak{s}}}_{
 N \mathfrak{s} = \frac{|a m|}{\gcd(a m, d^{\infty})}} 
\sum_{n=-B}^B
c_n \psi_n(\mathfrak{s}) \lambda(m)
\] 
is at most a constant times
\[\max((\log X)^{3 A},(\log X)^{2 A}/\epsilon) 
\max_{-B\leq n\leq B}
\left|\sum_{m\leq X}
\mathop{\sum_{\mathfrak{r}}}_{
N \mathfrak{r} = \gcd(a m,d^{\infty})}
\mathop{\sum_{\mathfrak{s}}}_{
 N \mathfrak{s} = \frac{|a m|}{\gcd(a m, d^{\infty})}} 
\psi_n(N \mathfrak{s}) \right|.\]
Clearly
\[
\sum_{m\leq X}
\mathop{\sum_{\mathfrak{r}}}_{
N \mathfrak{r} = \gcd(a m,d^{\infty})}
\mathop{\sum_{\mathfrak{s}}}_{
 N \mathfrak{s} = \frac{|a m|}{\gcd(a m, d^{\infty})}} 
\psi_n(\mathfrak{s}) \lambda(m) =
\lambda(a) 
\mathop{\sum_{\mathfrak{r}}}_{
\mathfrak{r} | d^{\infty}} \lambda(N \mathfrak{r})
\mathop{\mathop{\sum_{\mathfrak{s}}}_{\frac{a}{\gcd(a,d^\infty)} |
 N \mathfrak{s}}}_{N \mathfrak{s} \leq \frac{a X}{N \mathfrak{r}}}
 \psi_n(\mathfrak{s}) \lambda(N \mathfrak{s}).\]
Now
\[\begin{aligned}
\mathop{\sum_{\mathfrak{r}}}_{
\mathfrak{r} | d^{\infty}} \left|
\mathop{\mathop{\sum_{\mathfrak{s}}}_{\frac{a}{\gcd(a,d^\infty)} |
 N \mathfrak{s}}}_{N \mathfrak{s} \leq \frac{a X}{N \mathfrak{r}}}
 \psi_n(\mathfrak{s}) \lambda(N \mathfrak{s}) \right| &\leq
\mathop{\mathop{\sum_{\mathfrak{r}}}_{
\mathfrak{r} | d^{\infty}}}_{X^{1/2} < N \mathfrak{r} \leq a X}
\frac{a X}{N \mathfrak{r}} \log\left(\frac{a X}{N \mathfrak{r}}\right) \\ &+
\mathop{\mathop{\sum_{\mathfrak{r}}}_{
\mathfrak{r} | d^{\infty}}}_{N \mathfrak{r} \leq X^{1/2}}
\left| \mathop{\mathop{\sum_{\mathfrak{s}}}_{\frac{a}{\gcd(a,d^\infty)} |
 N \mathfrak{s}}}_{N \mathfrak{s} \leq \frac{a X}{N \mathfrak{r}}}
 \psi_n(\mathfrak{s}) \lambda(N \mathfrak{s}) \right|
  .\end{aligned}\]
Set $B = e^{(\log x)^{3/5} (\log \log x)^{1/5}}/(\log x)^A$.
We bound the first term on the right by Lemma \ref{lem:chai} and
the second term by Lemma \ref{lem:kuhl}, obtaining
\[\begin{aligned}
\mathop{\sum_{\mathfrak{r}}}_{
\mathfrak{r} | d^{\infty}} \left|
\mathop{\mathop{\sum_{\mathfrak{s}}}_{\frac{a}{\gcd(a,d^\infty)} |
 N \mathfrak{s}}}_{N \mathfrak{s} \leq \frac{a X}{N \mathfrak{r}}}
 \psi_n(\mathfrak{s}) \lambda(N \mathfrak{s}) \right| &\ll
\frac{a X}{X^{1/2}} (\log X)^C +
\mathop{\mathop{\sum_{\mathfrak{r}}}_{
\mathfrak{r} | d^{\infty}}}_{N \mathfrak{r} \leq X^{1/2}}
\frac{a X}{N \mathfrak{r}} 
 e^{-C \frac{(\log X)^{2/3}}{(\log \log X)^{1/5}}} \\
&\ll X e^{-C' \frac{(\log X)^{2/3}}{(\log \log X)^{1/5}}} .\end{aligned}\]

It remains to estimate
\[\sum_{m\leq X} \mathop{\sum_{\mathfrak{r}}}_{N \mathfrak{r} = 
\gcd(a m,d^{\infty})} \mathop{\mathop{\sum_{\mathfrak{s}}}_{
N \mathfrak{s} = \frac{|a m|}{\gcd(a m,d^{\infty})}}}_{
iota(\mathcal{I}^{-1}(\mathfrak{s}))\cap (S_1 \cup S_2) \ne 0} 1.\] 
It is enough to bound
\[\mathop{\mathop{\sum_{\mathfrak{a}}}_{N \mathfrak{a} \leq X}}_{
 \iota(\mathcal{I}^{-1}(\mathfrak{a})) \cap S_i = \emptyset} 1
= \mathop{\sum_{\iota(s) \in S_i}}_{N S \leq X} 1 \]
for $i=1,2$.
If $K/\mathbb{Q}$ is imaginary, the angle of $S_i$ is at most
$\epsilon$; if $K/\mathbb{Q}$ is real, the hyperbolic angle of $S_i$
is at most $\epsilon$.
Since \[\# \{s\in \iota^{-1}(S) : N s \leq X\}\]
is invariant when $S$ is multiplied by a unit $u\in \mathfrak{O}_K^*$,
we can assume without loss of generality that
 $\log x/y$ is bounded above and below by constants
depending only on $K$. Then the boundary of 
\[\{s\in \iota^{-1}(S) : N s \leq X\}\]
has length equal to at most a constant times $\sqrt{X}$. Hence
\[\mathop{\sum_{\iota(s) \in S_i}}_{N S \leq X} 1 \ll \epsilon X +
\sqrt{X}.\]
Set $\epsilon = \sqrt{B}$.
Then \[
\mathop{\sum_{x,y\in S\cap L}}_{|Q(x,y)| \leq X} \lambda(Q(x,y)) \ll
X e^{-C'' \frac{(\log X)^{2/3}}{(\log \log X)^{1/5}}} ,\]
as was desired.
\end{proof}
\begin{prop}
Let $Q(x,y) =  a x^2 + b x y + c y^2 \in \mathbb{Z}\lbrack x,y\rbrack$
be a quadratic form. Assume $b^2 - 4 a c \ne 0$.
Let $L\subset \mathbb{Z}^2$ be a lattice
coset, $S\subset \mathbb{R}^2$ a sector. Assume 
\begin{equation}\label{eq:midt}\lbrack \mathbb{Z}^2 : L\rbrack \ll (\log X)^A .\end{equation}
Then
\[\sum_{(x,y)\in S\cap L\cap \lbrack - N,N\rbrack^2}
 \lambda(Q(x,y)) \ll N^2 e^{-C \frac{(\log N)^{2/3}}{(\log \log N)^{1/5}}},\]
where $C$ and the implied constant depend only on $a$, $b$, $c$, $A$ and the implied constant
in (\ref{eq:midt}).
\end{prop}
\begin{proof}
If $Q$ is reducible, the statement follows immediately from (\ref{eq:sw}). Assume $Q$ is irreducible.
Let $K = \mathbb{Q}(\sqrt{b^2 - 4 a c})$.

Suppose $K/\mathbb{Q}$ is imaginary. Then $|Q(x,y)|=1$ describes an ellipse in $\mathbb{R}^2$
centered at the origin. Let $S\subset \mathbb{R}^2$ be a subquadrant.
Write the ellipse in polar coordinates:
\[\theta \in \lbrack 0,2 \pi\rbrack,\; r = r_1(\theta),\]
where $r_1:\lbrack 0,2 \pi \rbrack \to \mathbb{R}^+$ is $C^{\infty}$. Let
\[c_{1 0} = \min_{0\leq \theta\leq 2 \pi} r_1(\theta),\;
c_{1 1} = \max_{0\leq \theta\leq 2 \pi} |r_1'(\theta)|.\] Now consider the ellipse
\[\theta \in \lbrack 0,2 \pi\rbrack,\; r = \sqrt{X} r_1(\theta).\]
Any arc
\[\theta \in (\theta_1,\theta_2),\; r = \sqrt{X} r_1(\theta)\]
will lie within the region $R_{\theta_1,\theta_2}(\sqrt{X})$ bounded by the two arcs
\begin{equation}\label{eq:obe}\begin{aligned}
\theta &\in (\theta_1,\theta_2), \; r = \sqrt{X} (r_1(\theta_1) - c_{1 1}
(\theta_2 - \theta_1)),\\
\theta &\in (\theta_1,\theta_2), \; r = \sqrt{X} (r_1(\theta_1) + c_{1 1}
(\theta_2 - \theta_1)),
\end{aligned}\end{equation}
and the two lines $\theta = \theta_1$, $\theta = \theta_2$. It is easy to show
that
\begin{equation}\label{eq:smtar}
\# (R_{\theta_1,\theta_2}(\sqrt{X}) \cap \mathbb{Z}^2) \ll 
 c_{1 1} (\theta_1 - \theta_2)^2 X + (c_{1 1} + 1) (\theta_2-\theta_1) \sqrt{X} .\end{equation}

Write the boundary of the square $\lbrack -1,1\rbrack^2$ in polar coordinates:
\[\theta \in \lbrack 0,2 \pi\rbrack, \; r = r_2(\theta) .\]
Let \[c_{2 0} = \min_{0\leq \theta\leq 2 \pi} r_2(\theta),\:
c_{2 1} = \max_{0\leq \theta\leq 2 \pi} r_2(\theta),\:
c_{2 2} = \max_{0\leq \theta\leq 2 \pi} |r_2'(\theta)|.\]
For any positive real number $N$, path
\[\theta \in (\theta_1,\theta_2),\; r = N \cdot r_2(\theta)\]
lies in the region $R_{\theta_1,\theta_2}'(N)$ bounded by the arcs
\begin{equation}
\label{eq:owe}\begin{aligned}
\theta &\in (\theta_1,\theta_2), \; r = N (r_2(\theta_1) - c_{2 1}
(\theta_2 - \theta_1)) ,\\
\theta &\in (\theta_1,\theta_2), \; r = N (r_2(\theta_1) + c_{2 1}
(\theta_2 - \theta_1)) 
\end{aligned}\end{equation}
and the lines $\theta = \theta_1$, $\theta = \theta_2$. Clearly
\begin{equation}\label{eq:tarsm}
\# (R_{\theta_1,\theta_2}(N) \cap \mathbb{Z}^2) \ll c_{2 1}
(\theta_2 - \theta_1)^2 N^2 + (c_{2 1} + 1) (\theta_2 - \theta_1) N .
\end{equation}
As can be seem from (\ref{eq:obe}) and (\ref{eq:owe}), the region
\[\theta \in (\theta_1,\theta_2), r\leq N r_2(\theta)\]
contains the region
\[\theta \in (\theta_1,\theta_2), r\leq \sqrt{X} r_1(\theta)\]
for 
\[X = \left(\frac{N (r_2(\theta_1) - c_{2 1} (\theta_2 - \theta_1))}{r_1(\theta_1) + c_{1 1} (\theta_2 - \theta_1)}\right)^2 .\]
If $\theta_2-\theta_1 < \frac{r_{2 0}}{2 c_{2 1}}$, we have
$N^2\ll X\ll N^2$.
By (\ref{eq:smtar}) and (\ref{eq:tarsm}), the area between the two regions 
contains
\begin{equation}\label{eq:wro}
O(c_{1 1} (\theta_2 - \theta_1)^2 X + (c_{1 1} + 1) (\theta_2 - \theta_1)
\sqrt{X} + c_{2 1} (\theta_2 - \theta_1)^2 N^2 + (c_{2 1} + 1) 
(\theta_2 - \theta_1) N) \end{equation}
points with integral coordinates. We can rewrite (\ref{eq:wro}) as
\[O((\theta_2 - \theta_1)^2 N^2 + (\theta_2 - \theta_1) N),\]
where the implied constant depends on $r_i$ and $c_{i j}$.
By Lemma \ref{lem:vira},
\[
\mathop{\mathop{\sum_{(x,y)\in L}}_{\theta_1<\theta(x,y)<\theta_2}}_{|Q(x,y)| \leq X} \lambda(Q(x,y)) \ll
X e^{-C \frac{(\log X)^{2/3}}{(\log \log X)^{1/5}}} ,\]
where $\theta(x,y)$ is the angle $0\leq \theta< 2\pi$ between the $x$-axis
and the vector from $(0,0)$ to $(x,y)$.  Hence
\begin{equation}\label{eq:rty}
\mathop{\sum_{(x,y)\in L\cap \lbrack -N,N\rbrack}}_{\theta_1<\theta(x,y)<\theta_2} \lambda(Q(x,y)) \ll N^2 e^{- C \frac{(\log N)^{2/3}}{(\log \log N)^{1/5}}} +
(\theta_2 - \theta_1)^2 N^2 + (\theta_2 - \theta_1) N .\end{equation}

Let $S$ be a sector. We can assume that $S$ is given by
\[\theta<\theta(x,y)<\theta'\] for some $\theta, \theta' \in \lbrack 0,2\pi\rbrack$. Let 
\[\theta_0 = \theta,\, \theta_1 = \frac{\theta'-\theta}{n} + \theta,\,
\theta_2 = \frac{2 (\theta' - \theta)}{n} + \theta,\, \dotsc,\, \theta_n =
\theta' .\]
Then
$\theta_{i+1} - \theta_i = \frac{\theta' - \theta}{n} \leq \frac{2 \pi}{n}$.
Assume $n\geq \frac{4 \pi c_{2 1}}{r_{2 0}}$.
Hence, by (\ref{eq:rty}),
\[
\begin{aligned}
\sum_{(x,y)\in S\cap L \cap \lbrack - N , N\rbrack^2} \lambda(Q(x,y)) &=
\sum_{i=0}^{n-1} \mathop{\sum_{(x,y)\in L\cap \lbrack - N, N\rbrack^2}}_{\theta_i < \theta(x,y) <\theta_{i+1}} \lambda(Q(x,y))\\
&\ll n e^{- C \frac{(\log N)^{2/3}}{(\log \log N)^{1/5}}} N^2 +
\frac{1}{n} N^2 + N .\end{aligned}\]
Choose $n = \min(e^{\frac{C}{2} \frac{(\log N)^{2/3}}{(\log \log N)^{1/5}}},
\frac{4 \pi c_{2 1}}{r_{2 0}})$.
Then
\[\sum_{(x,y)\in S\cap L\cap \lbrack - N, N\rbrack^2} \lambda(Q(x,y)) \ll
e^{- \frac{C}{2} \frac{(\log N)^{2/3}}{(\log \log N)^{1/5}}} N^2.\]

Now suppose that $K/\mathbb{Q}$ is real. Then $|Q(x,y)|=1$ describes
two hyperbolas sharing two axes going through the origin. We can write the 
union of the two hyperbolas in polar coordinates:
\[\theta \in D,\; r = r_1(\theta),\]
where $\theta=\theta_a$, $\theta = \theta_b$ are the axes and
\[D = \lbrack 0, 2\pi\rbrack - \{\theta_a,\theta_b,\theta_a+\pi,\theta_b+\pi
\} .\] 
For $\theta\in \lbrack 0, 2\pi\rbrack$, define
\[d(\theta) = \min(|\theta-\theta_a|,|\theta-\theta_b|,
|\theta-(\theta_a+\pi)|,|\theta - (\theta_b + \pi)|) .\]
The function $r_1:D\to \mathbb{R}^+$ has a positive minimum $c_{1 0}$.
While $r_1(\theta)$ and $r_1'(\theta)$ are unbounded, 
$r_1(\theta) d(\theta)^{1/2}$ and $r_1'(\theta) d(\theta)^{3/2}$ are bounded;
let
\[c_{1 1} = \max_{\theta} |r_1(\theta)| \cdot d(\theta)^{1/2},\;\;
c_{1 2} = |r_1'(\theta)| \cdot d(\theta)^{3/2} .\]
We can define $r_2$, $c_{2 0}$, $c_{2 1}$ and $c_{2 2}$ as before.
Let $(\theta_1,\theta_2)\in D$. The region
\begin{equation}\label{eq:graw1}
\theta\in (\theta_1,\theta_2),\; r\leq N r_2(\theta)\end{equation}
contains the region
\begin{equation}\label{eq:graw2}
\theta\in (\theta_1,\theta_2),\; r\leq \sqrt{X} r_2(\theta)\end{equation}
for
\[X = \left(\frac{N (r_2 (\theta_1) - c_2 (\theta_2 - \theta_1))}{
r_1(\theta_1) + c_{1 1} \frac{\theta_2 - \theta_1}{\min(d(\theta_1),d(\theta_2))^{3/2}}} \right)^2 .\]
Assume 
\[\theta_2 - \theta_1 < \min\left(\frac{r_{2 0}}{2 c_{2 1}} , 
d(\theta_1), d(\theta_2)\right),\;\;\min(d(\theta_1),d(\theta_2))\ll N^{-\epsilon} .\]
Then \begin{equation}\label{eq:tgrw}N^{2 - 3\epsilon} \ll X\ll N^2 .
\end{equation}
It follows that the area between (\ref{eq:graw1}) and (\ref{eq:graw2})
contains
\[O((\theta_2-\theta_1)^2 N^2/\min(d(\theta_1),d(\theta_2))^2 +
(\theta_2-\theta_1) N/\min(d(\theta_1),d(\theta_2))^{3/2} .\]
By Lemma \ref{lem:vira} and (\ref{eq:tgrw}) we get
\[
\mathop{\mathop{\sum_{(x,y)\in L}}_{\theta_1<\theta(x,y)<\theta_2}}_{|Q(x,y)| \leq X} \lambda(Q(x,y)) \ll
N^2 e^{-C \frac{(\log N)^{2/3}}{(\log \log N)^{1/5}}} .\]
As in the case of $K/\mathbb{Q}$ imaginary, 
we can divide any sector $S$ into slices $(\theta_1,\theta_2)$
with \[\theta_2 - \theta_1 \sim e^{-\frac{C}{2} 
 \frac{(\log N)^{2/3}}{(\log \log N)^{1/5}}} .\] We leave out angles
of size \[e^{-\frac{C}{4} 
 \frac{(\log N)^{2/3}}{(\log \log N)^{1/5}}}\] around $\theta_a$,
$\theta_b$, $\theta_a + \pi$ and $\theta_b + \pi$. The statement follows. 
\end{proof}
\section{The average of $\lambda$ on the product of three linear 
factors}\label{sec:three}
\begin{lem}\label{split1}
 For any $M_2>M_1>1$, there
are  $\sigma_d\in \mathbb{R}$ with $|\sigma_d|\leq 1$
and support on
\[\{M_1\leq d<M_2 : p<M_1 \Rightarrow p\nmid d\}\]
such that
\[\begin{aligned}\sum_{(x,y)\in S\cap L} g(x) f(x,y) &=
  \mathop{\sum_a \sum_b \sum_c}_{(a b,c)\in S\cap L} \sigma_a\,
g(a) g(b) f(a b, c) \\ &+ O\left(\frac{\log M_1}{\log M_2} \frac{\Area(S)}{
\lbrack \mathbb{Z}^2:L\rbrack}
  + N M_2\right) \end{aligned}\]
for any positive integer $N>M_2$,
any convex set $S\subset \lbrack -N, N\rbrack^{\deg(K/\mathbb{Q})}$,
any lattice coset $L\subset \mathbb{Z}^{\deg(K/\mathbb{Q})}$ 
with index
$\lbrack \mathbb{Z}^2 : L\rbrack < M_1$,
any function $f:\mathbb{Z}^2\to \mathbb{C}$ 
and any completely multiplicative function
$g:\mathbb{Z}^2 \to \mathbb{C}$ with 
\[\max_{x,y} |f(x,y)|\leq 1,\;\; \max_y |g(y)|\leq 1.\] The implied constant is absolute.
\end{lem}
\begin{proof}
Let $y_1 = \min(\{y\in \mathbb{Z}:\exists x \,\st\, (x,y)\in S\cap L\})$. 
There is an
$l|\lbrack \mathbb{Z}^2:L\rbrack$ such that, for any $y\in \mathbb{Z}$, 
\[(\exists x \,\st\, (x,y)\in L) \Leftrightarrow
(l|y-y_1) .\]
Let
\[\begin{aligned}
N_{j,0} &= \min(\{x: (x, y_1 + j l) \in S\cap L\})\\
N_{j,1} &= \max(\{x: (x, y_1 + j l) \in S\cap L\})+1 .
\end{aligned}\]
Now take $\sigma_d$ as in Lemma \ref{sieve}. If $N_{j,1}-N_{j,0}>M_2$,
then
\[
\sum_{x : (x, y_1 + j l) \in S\cap L}
\left|1 - \sum_{d|x} \sigma_d\right| \ll
\frac{\log M_1}{\log M_2} \frac{N_{j,1}-N_{j,0}}{
\lbrack \mathbb{Z}^2 : L\rbrack/l}\]
Summing this over all $j$ we obtain
\[
\begin{aligned}
\sum_{(x,y)\in S\cap L}
\left|1 - \sum_{d|x} \sigma_d \right|
&\ll
\frac{\log M_1}{\log M_2} \frac{(\Area(S))/ l}{\lbrack \mathbb{Z}^2 : L
\rbrack} M_2 N\\
&\ll 
\frac{\log M_1}{\log M_2} \frac{\Area(S)}{\lbrack \mathbb{Z}^2 : L \rbrack} 
+ M_2 N .
\end{aligned}
\]
Since
\[
\left|\sum_{(x,y)\in S\cap L} g(y) f(x,y) -
\sum_{(x,y)\in S\cap L} \sum_{d|x} \sigma_d g(x) f(x,y) \right|\]
is at most
\[\sum_{(x,y)\in S\cap L} \left| g(y) f(x,y) -
\sum_{d|x} \sigma_d g(y) f(x,y) \right|
\leq \sum_{(x,y)\in S\cap L} \left|1 - \sum_{d|x} \sigma_d\right|
\]
and
\[\mathop{\sum_a \sum_b \sum_c}_{(a b,c)\in S\cap L} \sigma_a\,
g(a) g(b) f(a b, c) = 
\sum_{(x,y)\in S\cap L} \sum_{d|x} \sigma_d g(x) f(x,y) .\]
we are done.  
\end{proof}
\begin{lem}\label{mat3}
Let $c_1$, $c_2$ be integers. 
Let $L\subset \mathbb{Z}^2$ be a lattice. Then the set
$\{(a,b)\in \mathbb{Z}^2 : (a,b c_1),(a,b c_2)\in L\}$ 
is 
either the empty set or a lattice coset $L'\subset \mathbb{Z}^2$
of index dividing $\lbrack \mathbb{Z}^2:L\rbrack^2$.
\end{lem}
\begin{proof}
The set of all elements of $L$ of the form $(a,b c_1)$ is the intersection
of a lattice coset of index $\lbrack \mathbb{Z}^2:L\rbrack$ and a lattice of index $c_1$.
By (\ref{eq:intersl}) it
 is either the empty set or
a lattice coset of index dividing $c_1 \lbrack \mathbb{Z}^2:L\rbrack$. Therefore
the set of all $(a,b)$ such that $(a,b c_1)$ is in $L$ is either
the empty set or a lattice coset
$L_1$
of index dividing $\frac{1}{c_1} c_1 \lbrack \mathbb{Z}^2:L\rbrack = \lbrack \mathbb{Z}^2:L\rbrack$.
Similarly, the set of all $(a,b)$ such that $(a,b c_2)\in L$ is either the
empty set or a lattice coset
$L_2$ of index dividing $\lbrack \mathbb{Z}^2:L\rbrack$. Therefore $L'=L_1\cap L_2$
is either the empty set or 
a lattice coset of index dividing $\lbrack \mathbb{Z}^2:L\rbrack^2$. 
\end{proof}
\begin{defn}
For $A = \left(\begin{matrix} a_{1 1} & a_{1 2}\\ a_{2 1} & a_{2 2}\\a_{3 1} & a_{3 2}\end{matrix}\right)$ we denote
\[A_{1 2} = \left(\begin{matrix} a_{1 1} & a_{1 2}\\ a_{2 1} & a_{2 2}\end{matrix}\right)\;\;\;
A_{1 3} = \left(\begin{matrix} a_{1 1} & a_{1 2}\\ a_{3 1} & a_{3 2}\end{matrix}\right)\;\;\;
A_{2 3} = \left(\begin{matrix} a_{2 1} & a_{2 2}\\ a_{3 1} & a_{3 2}\end{matrix}\right) .\]
\end{defn}
\begin{prop}
Let $S$ be a convex subset of $\lbrack -N,N\rbrack^2$, $N>1$. Let $L\subset \mathbb{Z}^2$ be a lattice coset. Let $a_{1 1}$, $a_{1 2}$, $a_{2 1}$, $a_{2 2}$, $a_{3 1}$, $a_{3 2}$ be rational integers. Then
\[\sum_{(x,y)\in S\cap L} 
\lambda((a_{1 1} x + a_{1 2} y)(a_{2 1} x + a_{2 2} y)(a_{3 1} x + a_{3 2} y))
\ll \frac{\log \log N}{\log N} \frac{\Area(S)}{
\lbrack \mathbb{Z}^2 : L\rbrack}
+ \frac{N^2}{(\log N)^\alpha}\]
for any $\alpha>0$. The implied constant depends only on $(a_{i j})$
and $\alpha$.
\end{prop}
\begin{proof}
We can assume that $A_{12}$ is non-singular, as otherwise the
statement follows immediately from Lemma \ref{lem:diverto}. 
Changing variables we obtain
\[\begin{aligned}
\mathop{\sum_{(x,y)\in S\cap L}}_{\gcd(a_{1 1} x + a_{1 2} y,a_{2 1} x + a_{2 2} y)=1}
&\lambda(a_{1 1} x + a_{1 2} y) \lambda(a_{2 1} x + a_{2 2} y)
\lambda(a_{3 1} x + a_{3 2} y) \\&=
\mathop{\sum_{(x,y)\in A_{1 2} S \cap A_{1 2} L}}_{\gcd(x,y)=1}
\lambda(x) \lambda(y) \lambda\left((a_{3 1}\, a_{3 2}) A_{1 2}^{-1} 
\left(\begin{matrix}x\\y\end{matrix}\right)\right)\\&=
\mathop{\sum_{(x,y)\in A_{1 2} S \cap A_{1 2} L}}_{\gcd(x,y)=1}
\lambda(x) \lambda(y) \lambda(q_1 x + q_2 y) ,\end{aligned}\]
where $q_1 = -\frac{\det(A_{2 3})}{\det(A_{1 2})}$ and
$q_2 = \frac{\det(A_{1 3})}{\det(A_{1 2})}$. Note that
$q_1 x + q_2 y$ is an integer for all $(x,y)$ in $A_{1 2} L$. We can assume that neither $q_1$ nor $q_2$ is zero. Write $S' = A_{1 2} S$, $L' = A_{1 2} L$.
Clearly $S'\subset \lbrack -N',N'\rbrack^2$ for
  $N' = \max(|a_{1 1}| + |a_{1 2}|,|a_{2 1}|+|a_{2 2}|) N$.

Now set 
\[M_1 = (\log N')^{2\alpha +2},\;M_2 = 
\frac{(N')^{1/2}}{(\log N')^{\alpha}} .\] 
Clearly $M_2>M_1$ for $N>N_0$, $N_0$ depending only
on $(a_{i j})$ and $\alpha$.

 By Lemma \ref{split1},
\[\begin{aligned}\sum_{(x,y)\in S'\cap L'} \lambda(x) \lambda(y) 
\lambda(q_1 x + q_2 y) &=
  \mathop{\sum_a \sum_b \sum_c}_{(a b,c)\in S'\cap L'} \sigma_a\,
\lambda(a) \lambda(b) \lambda(c) 
\lambda(q_1 a b + q_2 c)\\ &+ O\left(\frac{\log M_1}{\log M_2} \frac{\Area(S')}{\lbrack \mathbb{Z}^2:L'\rbrack} + N' M_2\right) .\end{aligned}\]

We need to split the domain:
\[  \mathop{\sum_a \sum_b \sum_c}_{(a b,c)\in S'\cap L'} \sigma_a\,
\lambda(a) \lambda(b) \lambda(c) 
\lambda(q_1 a b + q_2 c)
= \sum_{s=1}^{\lceil M_2/M_1 \rceil} T_s
,\]
where
\[T_s = \mathop{\sum_{a = s M_1}^{(s+1) M_1 -1}
        \sum_{|b|\leq N'/s M_1} \sum_c}_{(a b,c)\in S' \cap L'}
 \sigma_a\,
\lambda(a) \lambda(b) \lambda(c) 
\lambda(q_1 a b + q_2 c) .\]
By Cauchy's inequality, 
\[
T_s^2 \leq \frac{(N')^2}{s M_1}
 \sum_c \sum_{|b|\leq N'/s M_1} 
\left(\mathop{\sum_{s M_1\leq a<(s+1) M_1}}_{(a b, c)\in S'\cap L'}
\sigma_a \lambda(a) \lambda(q_1 a b + q_2 c)\right)^2.\]
Expanding the square and changing the order of summation, we get
\[\frac{(N')^2}{s M_1}
\sum_{a_1 = s M_1}^{(s+1) M_1 -1}
\sum_{a_2 = s M_1}^{(s+1) M_1 -1}
\sigma_{a_1} \sigma_{a_2}
\lambda(a_1) \lambda(a_2) \mathop{\sum_c 
\sum_{|b|\leq N'/s M_1}}_{(a_i b,c)\in S'\cap L'} \lambda(q_1 a_1 b + q_2 c)
 \lambda(q_1 a_2 b + q_2 c) .\]
There are at most $M_1\cdot 2 N' \frac{N'}{s M_1}$ terms with
$c_1=c_2$. They contribute at most $\frac{2 (N')^4}{s^2 M_1}$
to $T_s^2$, and thus no more than $((N')^2/\sqrt{M_1}) \log M_2$ to the
sum $\sum_{s=1}^{\lceil M_2/M_1\rceil} T_s$.
It remains to bound
\[\mathop{\sum_{a_1= s M_1}^{(s+1) M_1 -1}
\sum_{a_2= s M_1}^{(s+1) M_1 -1}}_{a_1\ne a_2} \sigma_{a_1} \sigma_{a_2}
\lambda(a_1) \lambda(a_2) \mathop{\sum_c 
\sum_{|b|\leq N'/s M_1}}_{(a_i b,c)\in S'\cap L'} \lambda(q_1 a_1 b + q_2 c)
 \lambda(q_1 a_2 b + q_2 c) .\]

Since $|\sigma_a|\leq 1$ for all $a$,
the absolute value of this is at most
\[ 
\mathop{\sum_{a_1= s M_1}^{(s+1) M_1 -1}
\sum_{a_2= s M_1}^{(s+1) M_1 -1}}_{a_1\ne a_2} 
 \left|
\mathop{\sum_c \sum_b}_{(a_i b,c)\in S'\cap L'} \lambda(q_1 a_1 b + q_2 c)
 \lambda(q_1 a_2 b + q_2 c) \right|.\]

By Lemma \ref{mat3} we can write
$\{(b,c)\in \mathbb{Z}^2: (a_1 b, c),(a_2 b, c)\in S'\cap L'\}$ as
$S''\cap L''$ with $S''$ a convex subset of
$\lbrack -N'/\max(a_1,a_2),N'/\max(a_1,a_2)\rbrack \times
\lbrack -N',N'\rbrack$ and $L''\subset \mathbb{Z}^2$ a lattice coset
of index dividing $\lbrack \mathbb{Z}^2:L'\rbrack^2$. 
Hence we have the sum
\[
\mathop{\sum_{a_1= s M_1}^{(s+1) M_1 -1}
\sum_{a_2= s M_1}^{(s+1) M_1 -1}}_{a_1\ne a_2} 
\left|\sum_{(b,c)\in S''\cap L''} 
\lambda(q_1 a_1 b + q_2 c)
 \lambda(q_1 a_2 b + q_2 c) \right|.\]

Set $S_{a_1,a_2} = 
  \left(\begin{matrix}q_1 a_1 &q_2\\q_1 a_2 &q_2\end{matrix}\right) S''$,
$L_{a_1,a_2} = 
  \left(\begin{matrix}q_1 a_1 &q_2\\q_1 a_2 &q_2\end{matrix}\right)
L''$, $N'' = (|q_1|+|q_2|) N'$.  Clearly $S_{a_1,a_2}$ is a convex subset of $\lbrack -N'',N''\rbrack^2$ with 
\[\Area(S_{a_1,a_2}) = |q_1 q_2 (a_1-a_2)|
\Area(S'') \leq |q_1 q_2| M_1 \frac{4 (N')^2}{s M_1} \ll \frac{N^2}{s},\]
 whereas $L_{a_1,a_2}\subset \mathbb{Z}^2$ is a lattice coset of index $|q_1 q_2 (a_1 - a_2)| \lbrack L'':\mathbb{Z}^2\rbrack$. (That $L_{a_1,a_2}$ is inside $\mathbb{Z}^2$ follows from
our earlier remark that $q_1 x + q_2 y$ is an integer for all
$(x,y)$ in $A_{1 2} L'$.) Now we have
\[
\mathop{\sum_{a_1= s M_1}^{(s+1) M_1 -1}
\sum_{a_2= s M_1}^{(s+1) M_1 -1}}_{a_1\ne a_2} 
 \left| 
\sum_{(v,w)\in S_{a_1,a_2}\cap L_{a_1,a_2}} 
\lambda(v)
 \lambda(w) \right|.\]

This is at most
\[M_1^2 \max_{s M_1\leq a < (s+1) M_1} 
    \mathop{\max_{-M_1\leq d\leq M_1}}_{d\ne 0}
 \left|
\sum_{(v,w)\in S_{a,a+d}\cap L_{a,a+d}} 
\lambda(v)
 \lambda(w) \right|.\]
We can assume that $\lbrack \mathbb{Z}^2 : L\rbrack < (\log N)^{\alpha}$,
as otherwise the bound we are attempting to prove is trivial. Hence
$\lbrack \mathbb{Z}^2 : L ''\rbrack \ll (\log N)^{2 \alpha}$.
By Lemma \ref{lem:diverto},
\[ \left|
\sum_{(v,w)\in S_{a,a+d}\cap L_{a,a+d}} 
\lambda(v)
 \lambda(w) \right|\ll
\frac{N^2}{s} \cdot e^{-C (\log N'')^{3/5}/(\log \log N'')^{1/5}}
+ N^{1+1/3}
.\]
It is time to collect all terms. The total is at most a constant times
\[\begin{aligned}
\frac{\log M_1}{\log M_2} \frac{\Area(S')}{\lbrack \mathbb{Z}^2 :
L'\rbrack} + N' M_2^2 &+
\frac{(N')^2}{\sqrt{M_1}} \log M_2 \\ 
&+ (N')^2 \sqrt{M_1}
\log M_2 \cdot e^{-C (\log N'')^{3/5}/(\log \log N'')^{1/5}} \\
&+ N^{5/3} \sqrt{M_2}
,
\end{aligned}\]
where the constant depends only on $(a_{i j})$ and $\alpha$. Simplifying
we obtain
\[O\left(\frac{\log \log N}{\log N} \frac{\Area(S)}{\lbrack \mathbb{Z}^2 : L\rbrack}
+ \frac{N^2}{(\log N)^\alpha}\right).\]
\end{proof}
\section{The average of $\lambda$ on the 
product of a linear and a quadratic factor}\label{sec:threetwo}
We will be working with quadratic extensions $K/\mathbb{Q}$. It will be 
convenient to use embeddings $\jmath:K\to \mathbb{R}^2$ as in Lemma
\ref{sieve2} instead of embeddings $\iota:K\to \mathbb{R}^2$ of the kind
employed in section \ref{sec:parquad}. (In Lemma \ref{sieve2}, 
$\jmath:K\to \mathbb{R}^2$ takes $\mathfrak{O}_K$ to $\mathbb{Z}^2$,
whereas $\iota:K\to \mathbb{R}^2$ does not.) We define
\[\begin{aligned}
 \jmath(x + y \sqrt{d}) = (x,y) \;\; &\text{if $d\equiv 1 \mo 4$,}\\
 \jmath(x + y \sqrt{d}) = (x-y, 2y) \; &\text{if $d\not\equiv 1 \mo 4$,}
\end{aligned}\]
where $x, y\in \mathbb{Q}$. 

For every $z\in \jmath^{-1}(\lbrack -N,N \rbrack^2)$,
\begin{equation}
|N_{K/\mathbb{Q}} z| \ll N^2 ,
\end{equation}
where the implied constant depends only on $K$. In general there is no 
implication in the opposite sense, as the norm need not be positive definite.
For $K = \mathbb{Q}(\sqrt{d})$, $d<0$,
\begin{equation}\label{eq:boundn1}
\#\{z\in \mathfrak{O}_K : N_{K/\mathbb{Q}}(z) \leq A\} \ll A .
\end{equation}
For $K = \mathbb{Q}(\sqrt{d})$, $d>1$, $A\leq N^2$,
\begin{equation}\label{eq:boundn2}
\#\{z\in \jmath^{-1}(\lbrack -N,N \rbrack^2) : N_{K/\mathbb{Q}}(z) \leq A\}
\ll A \left(1 + \log \frac{N}{\sqrt{A}}\right) + N .
\end{equation}
In either case the implied constant depends only on $d$.

\begin{lem}\label{lem:divbyno}
Let $\mathfrak{a}$ be an ideal in $\mathbb{Q}(\sqrt{d})/\mathbb{Q}$ divisible by no
rational integer $n>1$. Then for any positive $N$,
$y_0\in \lbrack -N,N \rbrack$,
\[\#\{(x,y_0)\in \lbrack -N,N\rbrack^2 : \jmath^{-1}(x,y_0)\in \mathfrak{a}\}
\leq \lceil N/N_{K/\mathbb{Q}}(\mathfrak{a})\rceil .\]
\end{lem}
\begin{proof}
For every rational integer $r\in \mathfrak{a}$, $N \mathfrak{a}
| r$. Hence
\[\{x : \jmath^{-1}(x,y_0) \in \mathfrak{a}\}\]
is an arithmetic progression of modulus $N \mathfrak{a}$.
\end{proof}
\begin{prop}
Let $S$ be a convex subset of $\lbrack -N,N\rbrack^2$, $N>1$. Let $L\subset \mathbb{Z}^2$ be a lattice
coset. Let $a_1$, $a_2$, $a_3$, $a_4$, $a_5$ be rational integers
such that $a_1 x^2 + a_2 x y + a_3 y^2$ is irreducible. Then
\[\sum_{(x,y)\in S\cap L} \lambda((a_1 x^2 + a_2 x y + a_3 y^2) (a_4 x + a_5 y)) \ll
\frac{\log \log N}{\log N} \frac{\Area(S)}{\lbrack \mathbb{Z}^2 : L\rbrack} +
\frac{N^2}{(\log N)^\alpha}\]
for any $\alpha>0$. The implied constant depends only on $(a_{i j})$ and $\alpha$.
\end{prop}
\begin{proof}
Write $d$ for $a_1^2 - 4 a_0 a_2$,
$K/\mathbb{Q}$ for $\mathbb{Q}(\sqrt{d})/\mathbb{Q}$,
$\mathfrak{N} x$ for $N_{K/\mathbb{Q}} x$ and
$\overline{r + s\sqrt{d}}$ for
$r - s\sqrt{d}$.
By Lemma \ref{lem:quad} there are 
$\alpha_1, \alpha_2 \in
\mathfrak{O}_K$ linearly independent over 
$\mathbb{Q}$ and a non-zero rational number $k$ such that
\[a_1 x^2 + a_2 x y + a_3 y^2 = k \mathfrak{N} (x \alpha_1 + y \alpha_2)
= k (x \alpha_1 + y \alpha_2) \overline{(x \alpha_1 + y\alpha_2)} .\]
Hence
\[ \sum_{(x,y)\in S\cap L} \lambda((a_1 x^2 + a_2 x y + a_3 y^2) (a_4 x + a_5 y))\] equals
\[\lambda(k) \sum_{(x,y)\in S\cap L} \lambda((x \alpha_1 + y\alpha_2)
\overline{(x \alpha_1 + y\alpha_2)} (a_4 x + a_5 y)) .\]

By abuse of language we
write $\Re (r+s \sqrt{d})$ for $r$, $\Im (r + s \sqrt{d})$ for
$s$.
Let $C = \left(\begin{matrix}\Re \alpha_1 &\Re \alpha_2\\
\Im \alpha_1&\Im \alpha_2\end{matrix}
\right)^{-1}$. Then $a_4 x + a_5 y = q z + \overline{q z}$
for $z = x \alpha_1 + y \alpha_2$,
\[ q = \frac{1}{2} (a_4 c_{1 1} + a_5 c_{2 1} + \frac{1}{\sqrt{d}} (a_4 c_{1 2} + a_5 c_{2 2})).\] 
Define $\phi_Q:\mathbb{Z}^2\to \mathfrak{O}_K$ to be the mapping
$(x,y) \mapsto (x \alpha_1 + y \alpha_2)$. Let 
$L' = (\iota \circ \phi_Q)(L)$. Let $S'$ be the sector of $\mathbb{R}^2$
such that $(\iota \circ \phi)(S \cap \mathbb{Q}^2) = S' \cap \mathbb{Q}^2$.
Then
\[  \sum_{(x,y)\in S\cap L} \lambda((x \alpha_1 + y\alpha_2)
\overline{(x \alpha_1 + y\alpha_2)} (a_4 x + a_5 y)) =
\sum_{\jmath(z)\in S'\cap L'} 
 \lambda(z \overline{z} (q z + \overline{q z})) .\]
Note that $q z + \overline{q z}$ is an integer for all $z\in L'$.

Let $N'$ be the smallest integer greater than one such that
$j(S')\subset \lbrack -N',N'\rbrack^2$. (Note that $N'\leq c_1 N$, where
$c_1$ is a constant depending only on $\mathbb{Q}$.) Suppose $K/\mathbb{Q}$
is real. Then, by (\ref{eq:boundn2}),
\[\begin{aligned}
\#\{x\in \jmath^{-1}(S') : \mathfrak{N} x 
       \leq \frac{(N')^2}{(\log N)^{\alpha+1}}\}
&\leq \frac{(N')^2}{(\log N')^{\alpha + 1}}
(1 + \log (\log N')^{\alpha +1}) + N \\ &\leq \frac{N^2}{(\log N)^\alpha}
.\end{aligned}\]
The set
\[
\{x\in \lbrack -N,N\rbrack^2 : \mathfrak{N} (\jmath^{-1}(x)) 
  > \frac{(N')^2}{(\log N)^{\alpha+1}}\}
\]
is the region within a square and outside two hyperbolas. As such it is
the disjoint union of at most four convex sets. Hence the set
\[S'' = S \cap \{x\in \lbrack - N, N\rbrack^2 : \mathfrak{N}(\jmath^{-1}(x))
 > (N')^2/(\log N)^{\alpha+1} \}\]
is the disjoint union of at most four convex sets:
\[S'' = S_1 \cup S_2 \cup S_3 \cup S_4 .\]
In the following, $S^*$ will be $S_1$, $S_2$, $S_3$ or $S_4$, and as such
a convex set contained in $S''$.

Suppose now that $K/\mathbb{Q}$ is imaginary. Then the set
\[\{x\in \lbrack - N,N\rbrack^2 : \mathfrak{N}(\jmath^{-1}(x))
 > (N')^2/(\log N)^{\alpha+1} \}\]
is the region within a square and outside the circle given by
\begin{equation}\label{eq:satan}
\{x : \mathfrak{N}(\jmath^{-1}(x)) = (N')^2/(\log N)^{\alpha+1} \}.
\end{equation}
We can circumscribe about (\ref{eq:satan}) a rhombus containing no more
than
\[O((N')^2/(\log N)^{\alpha+1})\]
integer points, where the implied constant depends only on $Q$. We then
quarter the region inside the square $\lbrack -N, N \rbrack$ and outside
the rhombus, obtaining four convex sets $V_1$, $V_2$, $V_3$, $V_4$ inside
$S$. We let $S^*$ be $S\cap V_1$, $S\cap V_2$, $S\cap V_3$ or $S\cap V_4$.

For $K$ either real or imaginary, we now have a convex set 
$S^* \subset \lbrack -N, N\rbrack$ such that, for any 
$\jmath\in \mathfrak{O}_K$, 
\[\jmath(z) \in S^* \Rightarrow \mathfrak{N} z > N^2/(\log N)^{\alpha} .\]
Our task is to bound
\[\mathop{\sum_{z\in \mathfrak{O}_K}}_{\jmath(z) \in S^* \cap L'}
 \lambda(z \bar{z} (q z + \overline{q z})) .
\]

Set \[M_1 = (\log N)^{20(\alpha +1)},\;\; 
M_2 = \frac{N^{1/2}}{4 d \num(\mathfrak{N} q) 
 \lbrack \mathfrak{O}_K : L' \rbrack^2
(\log N)^{16 \alpha +22}} .\]
 By Lemma \ref{sieve2},
\begin{equation}\label{eq:satun}\begin{aligned}
\mathop{\sum_{z\in \mathfrak{O}_K}}_{\jmath(z)\in S^*\cap L'} 
  \lambda(z \overline{z} (q z + \overline{q z})) &=
\sum_{z\in S''\cap L'} \mathop{\sum_{\mathfrak{d}}}_{z\in \mathfrak{d}}
\sigma_{\mathfrak{d}} \lambda(z \overline{z} (q z + \overline{q z})) \\
&+ O\left(\frac{\log M_1}{\log M_2}
\frac{\Area(S')}{\lbrack \mathfrak{O}_K:L'\rbrack}\right)
 + N' M_2 .\end{aligned}\end{equation}
 Let $N'' = (9/4 + |d|) (N')^2 $. Then $\jmath(z)\in \lbrack - N', N'\rbrack$
implies $|\mathfrak{N} z|\leq N''$. Since $\sigma_\mathfrak{d}=0$
when $N \mathfrak{d}<M_1$, the first term on the right of
(\ref{eq:satun}) equals
\[ \mathop{\sum_{\mathfrak{b}}}_{\mathfrak{N} \mathfrak{b}\leq N''/M_1} 
\lambda(\mathfrak{b} \overline{\mathfrak{b}})
   \mathop{\sum_{\mathfrak{a}}}_{\text{$\mathfrak{a} \mathfrak{b}$ principal}}
\sigma_{\mathfrak{a}} \lambda(\mathfrak{a} \overline{\mathfrak{a}})
\mathop{\sum_{(z) = \mathfrak{a} \mathfrak{b}}}_{z\in S'' \cap L'} 
 \lambda(q z + \overline{q z}) .\]
We need to split the domain:
\[\mathop{\sum_{\mathfrak{b}}}_{\mathfrak{N} \mathfrak{b}\leq N''/M_1} 
\lambda(\mathfrak{b} \overline{\mathfrak{b}})
   \mathop{\sum_{\mathfrak{a}}}_{\text{$\mathfrak{a} \mathfrak{b}$ principal}}
\sigma_{\mathfrak{a}} \lambda(\mathfrak{a} \overline{\mathfrak{a}})
\mathop{\sum_{(z) = \mathfrak{a} \mathfrak{b}}}_{\jmath(z)\in S^* \cap L'} 
 \lambda(q z + \overline{q z}) =
\sum_{s=1}^{\lceil \log_2(N''/M_1)\rceil} T_s,\]
where
\[T_s = \mathop{\sum_{\mathfrak{b}}}_{2^{s-1}\leq \mathfrak{N} 
\mathfrak{b}\leq 2^s} 
\lambda(\mathfrak{b} \overline{\mathfrak{b}})
   \mathop{\sum_{\mathfrak{a}}}_{\text{$\mathfrak{a} \mathfrak{b}$ principal}}
\sigma_{\mathfrak{a}} \lambda(\mathfrak{a} \overline{\mathfrak{a}})
\mathop{\sum_{(z) = \mathfrak{a} \mathfrak{b}}}_{\jmath(z)\in S^* \cap L'} 
 \lambda(q z + \overline{q z}) .\]
Notice that $\lambda(\mathfrak{b} \overline{\mathfrak{b}})$,
$\sigma_{\mathfrak{a}}$,  $\lambda(\mathfrak{a} \overline{\mathfrak{a}})$ and
$\lambda(q z + \overline{q z})$ are all real.
By Cauchy's inequality,
\[\begin{aligned}
T_s^2&\leq 2^{s-1} \mathop{\sum_{\mathfrak{b}}}_{2^{s-1}\leq \mathfrak{N}
 \mathfrak{b}\leq 
2^s}
\left(\mathop{\sum_{\mathfrak{a}}}_{\text{$\mathfrak{a} \mathfrak{b}$ principal}}
\sigma_{\mathfrak{a}} \lambda(\mathfrak{a} \overline{\mathfrak{a}})
\mathop{\sum_{(z) = \mathfrak{a} \mathfrak{b}}}_{\jmath(z)\in S^* \cap L'} 
 \lambda(q z + \overline{q z})\right)^2 \\
&\leq
2^{s-1} \sum_{\mathfrak{b}}
\left(
\mathop{\mathop{\sum_{\mathfrak{a}}}_{\text{$\mathfrak{a} \mathfrak{b}$ principal}}}_{
n_{s 0}< \mathfrak{N} \mathfrak{a} \leq n_{s 1}}
\sigma_{\mathfrak{a}} \lambda(\mathfrak{a} \overline{\mathfrak{a}})
\mathop{\sum_{(z) = \mathfrak{a} \mathfrak{b}}}_{\jmath(z)\in S^* \cap L'} 
 \lambda(q z + \overline{q z})\right)^2 ,
\end{aligned}\]
where $n_{s 0} = \frac{(N')^2}{2^s (\log N)^{\alpha+1}}$ and
$n_{s 1} = \min(\frac{N''}{2^{s-1}}, M_2)$.
Expanding the square and changing the order of summation, we get
\[\begin{aligned}
2^{s-1} \mathop{\sum_{\mathfrak{a}_1}}_{n_{s 0} < \mathfrak{N} \mathfrak{a}_1\leq n_{s 1}}
&\mathop{\sum_{\mathfrak{a}_2}}_{n_{s 0} < \mathfrak{N} \mathfrak{a}_2
\leq n_{s 1}}
\sigma_{\mathfrak{a}_1} \sigma_{\mathfrak{a}_2}
\lambda(\mathfrak{a}_1 \overline{\mathfrak{a}_1})
\lambda(\mathfrak{a}_2 \overline{\mathfrak{a}_2})
\\ &\mathop{\sum_{\mathfrak{b}} }_{
\text{$\mathfrak{a}_1 \mathfrak{b}, \mathfrak{a}_2 \mathfrak{b}$ principal}}
\mathop{\sum_{(z_1) = \mathfrak{a}_1 \mathfrak{b}}}_{\jmath(z_1)\in S^* 
\cap L'} 
\mathop{\sum_{(z_2) = \mathfrak{a}_2 \mathfrak{b}}}_{\jmath(z_2)\in S^* 
\cap L'} 
 \lambda(q z_1 + \overline{q z_1}) \lambda(q z_2 + \overline{q z_2}) .
\end{aligned}\]

Write $\mathcal{S}(x + y\sqrt{d})$ for $\max(|x|,|y|)$. Let
$r = (z_2/z_1)\cdot \mathfrak{N} \mathfrak{a}$. We have
$r\in \overline{\mathfrak{a}_1}$ because 
\[(r) = ((z_2)/(z_1))\cdot \mathfrak{N} \mathfrak{a}_1 =
(\mathfrak{a}_2 / \mathfrak{a}_1)\cdot \mathfrak{N} \mathfrak{a}_1 =
\mathfrak{a}_2 \cdot \overline{\mathfrak{a}_1} .\]
Since $\mathfrak{N} z_1 > \frac{(N')^2}{(\log N)^{\alpha + 1}}$
and $\mathcal{S} (z_2 \overline{z_1}) \ll (N')^2$, where the implied
constant depends only on $\mathbb{Q}$,
\begin{equation}\label{eq:S}
\mathcal{S} (r) =
\mathcal{S}\left(\frac{z_2}{z_1} \mathfrak{N} \mathfrak{a}_1\right) =
\mathcal{S}\left(\frac{z_2 \overline{z_1}}{\mathfrak{N} z_1} 
\mathfrak{N} \mathfrak{a}_1\right) =
\mathcal{S}(z_2 \overline{z_1}) \frac{\mathfrak{N} \mathfrak{a}}{
\mathfrak{N} z_1} \ll n_{s 1} (\log N)^{\alpha + 1} .\end{equation}
Set \[R_s = \jmath^{-1}\left(\left[- k n_{s 1} (\log N)^{\alpha + 1},
k n_{s 1} (\log N)^{\alpha + 1}\right]^2\right) ,\]
where $k$ is the implied constant in (\ref{eq:S}) and as such depends only on $K$.
Changing variables we obtain
\[\begin{aligned}
2^{s-1} \mathop{\sum_{\mathfrak{a}}}_{n_{s 0} < \mathfrak{N}
 \mathfrak{a}_1\leq n_{s 1}}
&\mathop{\sum_{r\in \overline{\mathfrak{a}} \cap R_s}}_{
 n_{s 0} < \mathfrak{N} \left(\frac{(r)}{\mathfrak{a}}\right) \leq n_{s 1}}
\sigma_{\mathfrak{a}} \sigma_{(r)/\mathfrak{a}}
\lambda(\mathfrak{a} \overline{\mathfrak{a}})
\lambda\left(\frac{(r)}{\mathfrak{a}}
\overline{\frac{(r)}{\mathfrak{a}}}\right) \\
&\mathop{\mathop{\sum_z}_{\jmath(z)\in \jmath(\mathfrak{a})\cap S^*\cap L'}}_{
\jmath(r z/\mathfrak{N} \mathfrak{a}) \in S^*\cap L'} \lambda(q z + \overline{q z}) \lambda\left(\frac{q r z}{N \mathfrak{a}} + \overline{\frac{q r z}{N \mathfrak{a}}}\right) ,
\end{aligned}\]
that is, $2^{s-1}$ times
\begin{equation}\label{eq:somename}
\mathop{\sum_{\mathfrak{a}}}_{n_{s 0} < \mathfrak{N} \mathfrak{a} \leq n_{s 1}}
\mathop{\sum_{r\in \overline{\mathfrak{a}} \cap R_s}}_{
 n_{s 0} < \mathfrak{N} \left(\frac{(r)}{\mathfrak{a}}\right) \leq n_{s 1}}
\sigma_{\mathfrak{a}} \sigma_{(r)/\mathfrak{a}}
\lambda(r \overline{r})
\mathop{\mathop{\sum_{z}}_{\jmath(z)\in \jmath(\mathfrak{a})\cap S^*\cap L'}}_{
\jmath(r z/\mathfrak{N} \mathfrak{a})\in S^*\cap L'} \lambda(q z + \overline{q z}) \lambda\left(\frac{q r z}{\mathfrak{N} \mathfrak{a}} + 
\overline{\frac{q r z}{\mathfrak{N} \mathfrak{a}}}\right) .
\end{equation}
For any non-zero rational integer $n$,
\[\begin{aligned}
\mathop{\mathop{\sum_{\mathfrak{a}}}_{n_{s 0} <\mathfrak{N} \mathfrak{a}
\leq n_{s 1}}}_{n|\mathfrak{a}}
\sum_{r\in \overline{\mathfrak{a}} \cap R_s}
\mathop{\sum_{\jmath(z)\in \mathfrak{a}\cap S^*}}_{
2^{s-1}\leq \mathfrak{N}((z)/\mathfrak{a})
< 2^s} 1 &\ll
\mathop{\mathop{\sum_{\mathfrak{a}}}_{n_{s 0} <\mathfrak{N} \mathfrak{a}
\leq n_{s 1}}}_{n|\mathfrak{a}}
\frac{(2 k n_{s 1} (\log N)^{\alpha+1})^2}{\mathfrak{N} \mathfrak{a}} 
2^s \log 2^s \\
&\ll \frac{1}{n^2} \frac{N^4 (\log N)^{2\alpha +5}}{2^s}.
\end{aligned}\]
Since the support of $\sigma_{\mathfrak{d}}$ is a subset of
\[\{\mathfrak{d} : M_1\leq \mathfrak{N} \mathfrak{d} < M_2,
\mathfrak{N} \mathfrak{p} < M_1 \Rightarrow \mathfrak{N} \mathfrak{p}
\nmid \mathfrak{d}\},\]
we have that $n|\mathfrak{a}$ and $\sigma_{\mathfrak{a}} \ne 0$
imply $n\geq \sqrt{M_1}$. Therefore (\ref{eq:somename}) equals
\begin{equation}\label{eq:somename2}
\mathop{
 \mathop{\sum_{\mathfrak{a}}}_{
n_{s 0} < \mathfrak{N} \mathfrak{a}\leq n_{s 1}}}_{n>1
\Rightarrow n\nmid \mathfrak{a}}
\mathop{\sum_{r\in \overline{\mathfrak{a}} \cap R_s}}_{
 n_{s 0} < \mathfrak{N} \left(\frac{(r)}{\mathfrak{a}}\right) \leq
 n_{s 1}}
\sigma_{\mathfrak{a}} \sigma_{(r)/\mathfrak{a}}
\lambda(r \overline{r})
\mathop{\mathop{\sum_{z}}_{\jmath(z)\in \jmath(\mathfrak{a})\cap S^*\cap L'}}_{
\jmath(r z/\mathfrak{N} \mathfrak{a})\in S''\cap L'} \lambda(q z + \overline{q z}) \lambda\left(\frac{q r z}{N \mathfrak{a}} + \overline{\frac{q r z}{N \mathfrak{a}}}\right) 
\end{equation}
plus $O(N^4 (\log N)^{2 \alpha +5}/(2^s \sqrt{M_1}))$. The absolute value of
(\ref{eq:somename2}) is at most
\begin{equation}\label{eq:somename3}
\mathop{\mathop{\sum_{\mathfrak{a}}}_{
n_{s 0} < \mathfrak{N} \mathfrak{a} \leq n_{s 1}}}_{n>1\Rightarrow n\nmid \mathfrak{a}}
\sum_{r\in \overline{\mathfrak{a}} \cap R_s}
\left|\mathop{\mathop{\sum_{z}}_{\jmath(z)\in \jmath(\mathfrak{a})\cap S^*\cap L'}}_{
\jmath(r z/\mathfrak{N} \mathfrak{a})\in S''\cap L'}
 \lambda(q z + \overline{q z}) \lambda\left(\frac{q r z}{N \mathfrak{a}} + \overline{\frac{q r z}{N \mathfrak{a}}}\right) \right| .
\end{equation}
By Lemma \ref{lem:divbyno},
\[\begin{aligned}
\mathop{\sum_{\mathfrak{a}}}_{
n_{s 0} < \mathfrak{N} \mathfrak{a}
\leq n_{s 1}}
\sum_{r\in \overline{\mathfrak{a}} \cap R_s\cap \mathbb{Z}}
\sum_{z\in \mathfrak{a}\cap S^*} 1 \;&\ll
\mathop{\sum_{\mathfrak{a}}}_{
n_{s 0} < \mathfrak{N} \mathfrak{a}
\leq n_{s 1}} \left(\frac{N'}{\mathfrak{N} \mathfrak{a}} + 1\right)
 \left(\frac{(N')^2}{\mathfrak{N} \mathfrak{a}} + N'\right) \\ &\ll
\frac{N^3 \log M_1}{n_{s 0}} + N n_{s 1} .\end{aligned}\]

Thus we are left with
\begin{equation}\label{eq:somename4}
\mathop{\mathop{\sum_{\mathfrak{a}}}_{n_{s 0} < \mathfrak{N} \mathfrak{a} \leq n_{s 1}}}_{n>1\Rightarrow n\nmid \mathfrak{a}}
\mathop{\sum_{r\in \overline{\mathfrak{a}} \cap R_s}}_{\Im r \ne 0}
\left|\mathop{\mathop{\sum_{z}}_{\jmath(z)\in \jmath(\mathfrak{a})\cap S^*\cap L'}}_{
\jmath(r z/\mathfrak{N} \mathfrak{a})\in S''\cap L'}
 \lambda(q z + \overline{q z}) \lambda\left(\frac{q r z}{N \mathfrak{a}} + \overline{\frac{q r z}{N \mathfrak{a}}}\right) \right| .
\end{equation}

Notice that $r\in \overline{\mathfrak{a}}$ and
$z\in \mathfrak{a}$ imply $(r z/\mathfrak{N} \mathfrak{a})\in \mathfrak{O}_K$. Hence
$(r/\mathfrak{N} \mathfrak{a})^{-1} \mathfrak{O}_K \supset \mathfrak{a}$. Therefore
$(r/\mathfrak{N} \mathfrak{a})^{-1} \jmath^{-1}(L') \cap \mathfrak{a}$ is either the empty set or a sublattice of 
$\mathfrak{a}$ of index dividing $\lbrack \mathfrak{O}_K : L' \rbrack$. This means that
\[L_{\mathfrak{a},r} = \{z \in \mathfrak{a}\cap \jmath^{-1}(L'): (r z/\mathfrak{N} \mathfrak{a})\in \jmath^{-1}(L')\}\]
is either the empty set or a sublattice of $\mathfrak{a}$ of index 
$\lbrack \mathfrak{a} : L_{\mathfrak{a},r} \rbrack$ 
dividing $\lbrack \mathfrak{O}_K : L' \rbrack^2$,
whereas \[S_{\mathfrak{a},r} = \{z \in S^* : (r z/\mathfrak{N} \mathfrak{a}) \in S''\}\] 
is a convex subset of $\lbrack -N',N' \rbrack^2$.
The map 
\[\kappa:(x,y) \mapsto (q \cdot \phi_Q(x,y) + \overline{q \cdot \phi_Q(x,y)},
\frac{q r \cdot \phi_Q(x,y)}{\mathfrak{N} \mathfrak{a}} + \overline{\frac{q r 
\cdot \phi_Q(x,y)}{\mathfrak{N} \mathfrak{a}}})\]
is given by the matrix
\[\begin{aligned}
\left(\begin{matrix} 2 & 0\\ 2\frac{\Re r}{\mathfrak{N} \mathfrak{a}}& 2 d 
\frac{\Im r}{\mathfrak{N} \mathfrak{a}} \end{matrix}\right)\cdot
\left(\begin{matrix} \Re q & d\, \Im q \\ \Im q & \Re q\end{matrix}\right) 
\;\; &\text{if $d\not\equiv 1 \mo 4$},\\
\left(\begin{matrix} 2 & 0\\ 2\frac{\Re r}{\mathfrak{N} \mathfrak{a}}& 2 d 
\frac{\Im r}{\mathfrak{N} \mathfrak{a}} \end{matrix}\right)\cdot
\left(\begin{matrix} \Re q & d\, \Im q \\ \Im q & \Re q\end{matrix}\right) 
\cdot
\left(\begin{matrix} 1 & \frac{1}{2} \\ 0 & \frac{1}{2}\end{matrix}\right)
\;\; &\text{if $d \equiv 1 \mo 4$}.\end{aligned} \]
Hence $\kappa(L_{\mathfrak{a},r})$ either the empty set or a lattice 
$L'_{\mathfrak{a},r}$ 
of index
\[\lbrack \mathbb{Z}^2 : L'_{\mathfrak{a},r} \rbrack =
\begin{cases}
4 d \,\Im r N q  \lbrack \mathfrak{a} : L_{\mathfrak{a},r} \rbrack
&\text{if $d \not\equiv 1 \mo 4$},\\
2 d \,\Im r N q  \lbrack \mathfrak{a} : L_{\mathfrak{a},r} \rbrack
&\text{if $d \equiv 1 \mo 4$}.\end{cases}\]
and $\kappa(S_{\mathfrak{a},r})$ is a convex set $S'_{\mathfrak{a},r}$
contained in
\[\lbrack - 3 |d| \frac{\mathcal{S}(r)}{\mathfrak{N} \mathfrak{a}} \mathcal{S} (q) N',
 3 |d| \frac{\mathcal{S}(r)}{\mathfrak{N} \mathfrak{a}} \mathcal{S} (q) N' \rbrack^2,\]
which is contained in 
\[\lbrack - 3 |d| S(q) \frac{n_{s 1} (\log N)^{\alpha+1}}{n_{s 0}},
3 |d| S(q) \frac{n_{s 1} (\log N)^{\alpha+1}}{n_{s 0}}\rbrack ,\]
which is in turn contained in
\[
\lbrack - k' (\log N)^{2 \alpha + 2} N, k' (\log N)^{2 \alpha + 2} N \rbrack,\]
where $k'$ depends only on $d$ and $q$. Write (\ref{eq:somename4}) as
\begin{equation}\label{eq:somename5}
\mathop{\mathop{\sum_{\mathfrak{a}}}_{n_{s 0} < N \mathfrak{a} \leq n_{s 1}}}_{n>1\Rightarrow n\nmid \mathfrak{a}}
\mathop{\sum_{r\in \overline{\mathfrak{a}} \cap R_s}}_{\Im r \ne 0}
\left|\sum_{(v,w)\in L'_{\mathfrak{a},r} \cap  S'_{\mathfrak{a},r}}
\lambda(v) \lambda(w)\right|
\end{equation}
\end{proof}

Since $r$ is in $R_s$, $\Im r$ takes values between $-k n_{s 1} (\log N)^{\alpha+1}$
and $k n_{s 1} (\log N)^{\alpha+1}$. By Lemma \ref{lem:divbyno}, $\Im r$ takes each of
these
values  at most
\[\lceil (k n_{s 1} (\log N)^{\alpha+1})/n_{s 0} \rceil \ll
(\log N)^{2\alpha +2} \]
times. Thus (\ref{eq:somename5}) is bounded by a constant times
\[ \frac{N''}{2^{s-1}} (\log N)^{2\alpha + 2}
\sum_{0<y\leq k M_2 (\log n)^{\alpha +1}}
\max_\mathfrak{a} \max_{r : \Im r = y}
\left|\sum_{(v,w)\in L'_{\mathfrak{a},r} \cap  S'_{\mathfrak{a},y}}
\lambda(v) \lambda(w)\right|.\]
By Corollary \ref{bomb5},
\[\sum_{0<y\leq k M_2 (\log n)^{\alpha +1}}
\max_\mathfrak{a} \max_{r : \Im r = y}
\left|\sum_{(v,w)\in L'_{\mathfrak{a},r} \cap  S'_{\mathfrak{a},y}}
\lambda(v) \lambda(w)\right|\] is
\[O\left( \tau(4 d \num(N q) 
\det(N q) [\mathfrak{O}_K : L']^2)
\frac{((\log N)^{2\alpha +2} N)^2}{(\log N)^{8(\alpha+1)}}\right) .\]

It is time to collect all terms. The total is at most 
\[\begin{aligned} \frac{\log M_1}{\log M_2}
\frac{\Area(S')}{\lbrack \mathfrak{O}_K:L'\rbrack}
 &+ N' M_2^2 + \frac{N^2 (\log N)^{\alpha + \frac{7}{2}}}{\sqrt{M_1}}\\ 
&+ \sqrt{N} M_2 (\log N)^{(\alpha + 1)/2} + \sqrt{N} M_2 + N^2 (\log N)^{\alpha}\end{aligned}\]
times a constant depending only on $(a_{i j})$ and $\alpha$.
This simplifies to
\[O\left(\frac{\log \log N}{\log N} \frac{\Area(S)}{
\lbrack \mathbb{Z}^2 : L\rbrack} +
\frac{N^2}{(\log N)^\alpha}\right) .\]

%% file: lambirr.tex
\section{The average of $\lambda$ on irreducible cubics}\label{sec:irrcu}
In the present section we shall prove that $\mu(P(x,y))$ averages to
zero for any irreducible homogeneous
polynomial $P$ of degree $3$. There are two
main stages in the proof: one is the reduction of the problem to a bilinear 
condition, and the other is the demonstration of the bilinear condition.
The second stage resembles its analogue in Heath-Brown's proof that 
$x^3 + 2 y^3$ captures its primes (\cite{HB}); although it is too early
to speak of the general features of a strategy that was first carried out
in \cite{FI1} and is still developing, one may venture that the bilinear
conditions involved in the strategy carry over between related problems
with relative ease. (See Appendix \ref{sec:gebouw}.) The first stage,
namely, the reduction to the bilinear condition, must be attempted with
much closer regard to the specifics of the problem at hand. The reader may
remark that there are few resemblances between subsection \ref{subs:agarag}
and the corresponding sections in \cite{HB}, \cite{HBM}, \cite{HBM2}.
We do follow the example of \cite{HB} in giving a fictively rational
outline before undertaking the actual procedure over a cubic field.
This explanatory device is appropiate in our case because of the inherent
complications of what is essentially an extension of an approach similar
to that in \cite{FI2} to a density below the natural range of the method.
For the sake of familiarity, we will adopt certain notational conventions
used in \cite{FI2}.
\subsection{Sketch}\label{subs:recabarren}
Let $\{a_n\}_{n=1}^{\infty}$ be a bounded sequence of non-negative
real numbers. Write
\[A(x) = \sum_{1\leq n\leq x} a_n, \;\;\;A_d(x) = 
\mathop{\sum_{1\leq n\leq x}}_{d|n} a_n .\]
Our linear axiom will be
\begin{equation}\label{eq:axl}
A_d = \frac{A(x)}{d} + \text{error}\;\;\;\;\;\text{for $d\ll D(x)$,}
\end{equation}
where the error term is small enough to be irrelevant for our purposes.
We also take the bilinear axiom
\begin{equation}\label{eq:bili}
\mathop{\sum_{1\leq r s\leq x}}_{V\leq s\leq 2 V} f(r) g(s) a_{r s}
\ll A(x) (\log x)^{-c_1} ,\end{equation}
valid for any $V$, $f$, $g$ satisfying
\[
x^{1/2} t(x) \leq V \leq x/v(x),\]
\[f(n),g(n) \ll \tau_{c_2}(n),\]
\[
\mathop{\sum_{s\equiv a \mo m}}_{s\leq S} g(s) \ll S e^{- \kappa \sqrt{\log S}}
\;\;\text{for all $m\ll (\log x)^{c_4}$,}\]
where the constants $c_i$ will be as large as needed, and $\kappa$ denotes
an exponent of no importance. We will assume $v(x) > \sqrt{D(x)}$,
as therein lies the origin of certain difficulties that we must learn to
resolve. Set $z(x) = v(x)/\sqrt{D(x)}$. We assume 
\[\begin{aligned}
\log z(x) &\ll (\log x)^{1/c_5},\\
z(x) &\gg (\log x)^{c_6}\\
v(x) D(x) &\gg x \cdot (z(x))^{-\kappa'},\end{aligned}\]
Set 
\[\begin{aligned}
u(x) &= (z(x))^{\kappa'} v(x), \;\; y(x)= (z(x))^{-\kappa'-2} v(x),\\
w(x) &= \frac{u(x)}{z(x)} 2^{\lbrack \log_2 x^{1/2}/(u(x) t(x))\rbrack} .
\end{aligned}\]
We write $t$, $u$, $v$, $w$, $y$, $z$ instead of
$t(x)$, $u(x)$, \ldots, $z(x)$, for the sake of brevity.

We adopt the symbols in \cite{FI2}:
\[\begin{aligned}
f(n\leq a) &= f(n) \cdot \lbrack n\leq a\rbrack,\\
 f(n>a) &= f(n) \cdot \lbrack n > a\rbrack .\end{aligned}\]
For any integer $n$ and any function $f$,
\[\begin{aligned}
f(n) &= f(n\leq a) + f(n>a),\\
f(n>a) &= \sum_{b c |n} \mu(b) f(c>a) .\end{aligned}\]
Hence
\[\begin{aligned}
\mu(n) &= \mu(n\leq u) + \sum_{b c | n} \mu(b) \mu(c>u)\\
&= \mu(n\leq u) + \sum_{b c | n} \mu(b>u) \mu(c>u) +
\sum_{b c | n} \mu(b\leq u) \mu(c>u) .\end{aligned}\]
By M\"obius inversion,
\[\begin{aligned}\sum_{b c|n} \mu(b\leq v) \mu(c>u) &=
 \sum_{b c|n} \mu(b\leq u) \mu(c) -
\sum_{b c|n} \mu(b\leq u) \mu(c\leq u) \\
&= \mu(n\leq u) - \sum_{b c|n} \mu(b\leq u) \mu(c\leq u) .\end{aligned}\]
Therefore
\[\mu(n) = 2 \mu(n\leq u) + \sum_{b c|n} \mu(b>u) \mu(c>u) -
\sum_{b c|n} \mu(b\leq u) \mu(c\leq u) .\]
We can split our ranges of summation:
\[\begin{aligned}
\sum_{b c|n} \mu(b>u) \mu(c>u) &= \sum_{b c|n} \mu(u<b\leq w) \mu(c>u)
\\ &+ \sum_{b c|n} \mu(b>w) \mu(u<c\leq w) + \sum_{b c|n} \mu(b>w)
\mu(c>w),\end{aligned}\]
\[\begin{aligned}
\sum_{b c|n} \mu(b\leq u) \mu(c\leq u) &= 
  \sum_{b c|n} \mu(b\leq u) \mu(c\leq y) +
  \sum_{b c|n} \mu(b\leq y) \mu(y<c\leq u) \\
&+ \sum_{b c|n} \mu(y<c\leq u) \mu(y< c\leq u) .\end{aligned}\]
Thus
\begin{equation}\label{eq:matdi}\begin{aligned}
\mu(n) &= 2 \mu(n\leq u) +
\sum_{b c|n} \mu(u<b\leq w) \mu(c>u) \\
&+ \sum_{b c|n} \mu(b>w) \mu(u<c\leq w) + \sum_{b c|n} \mu(b>w)
\mu(c>w)\\
&-  \sum_{b c|n} \mu(b\leq u) \mu(c\leq y) -
  \sum_{b c|n} \mu(b\leq y) \mu(y<c\leq u) \\
&- \sum_{b c|n} \mu(y<b\leq u) \mu(y< c\leq u) .\end{aligned}
\end{equation}
We denote the terms on the right side of (\ref{eq:matdi}) by
$\beta_1(n), \beta_2(n),\dotsc ,\beta_7(n)$. Set
\[S_j(x) = \sum_{n=1}^x \beta_j(n) a_n .\]
Then
\begin{equation}\label{eq:lins}
\sum_{n=1}^x \mu(n) a_n = S_1(x) + S_2(x) + S_3(x) + S_4(x) 
- S_5(x) - S_6(x) - S_7(x) .\end{equation}

The term $S_1(x)$ can be bounded trivially by $O(u)$. 
We can estimate $S_5(x)$ by means of the linear axiom (\ref{eq:axl}):
\[\begin{aligned}
S_5(x) &= \mathop{\sum_{1\leq n\leq x}}_{b c | n} \mu(b\leq u) 
\mu(c\leq y) a_n\\
&= \sum_{b,c} \mu(b\leq u) \mu(c\leq y) \frac{A(x)}{b c}\\
&= A(x) \cdot \sum_{b\leq u} \mu(b)/b \cdot \sum_{c\leq y} \mu(c)/c\\
&\ll A(x) \cdot e^{-\kappa \sqrt{\log u}} e^{-\kappa \sqrt{\log y}} 
\ll A(x) e^{-\kappa \sqrt{\log x}} .\end{aligned}\]
In the same way,
\[S_6(x) \ll A(x) e^{\kappa \sqrt{\log x}} .\]
We can easily prepare $S_2$ for an application of the bilinear condition
(\ref{eq:bili}):
\[\begin{aligned}
S_2(x) &= \mathop{\sum_{n\leq x}}_{b c | n} \mu(u<b\leq w)
\mu(c>u) a_n,\end{aligned}\]
\[\begin{aligned}
\mathop{\sum_{x/z\leq n\leq x}}_{b c|n} \mu(u<b\leq w) \mu(c>u) a_n
 &= \mathop{\sum_{x/z \leq r s\leq x}}_{x/z w\leq s\leq x/u} f(r)
g(s) a_{r s},\\
 &= \mathop{\sum_{r s\leq x}}_{x/z w\leq s\leq x/u} f(r)
g(s) a_{r s} + O\left(\sum_{n\leq x/z} \tau_3(n) a_n\right),\end{aligned}\]
where
\[\begin{aligned}
f(r) &= \mu(u<b\leq w),\\
g(s) &= \sum_{c|s} \mu(c>u) .\end{aligned}\]
Clearly
\[\label{eq:anbir}\begin{aligned}
\mathop{\sum_{s\equiv a \mo m}}_{s\leq S} g(s) &= 
\mathop{\sum_{s\equiv a \mo m}}_{s\leq S}
 \sum_{c|s} \mu(c>u) = \sum_{d\leq S/u} \sum_{u<c\leq S/d} \mu(c) \\
&= \sum_{d\leq S/u} \left(\sum_{c\leq S/d} \mu(c) - \sum_{c\leq u} \mu(c)\right) 
\ll S e^{-\kappa \sqrt{\log S}} .\end{aligned}\]
Hence, by (\ref{eq:bili}),
\[\mathop{\sum_{r s\leq x}}_{x/z w\leq s\leq x/u} f(r) g(s) a_{r s} \ll
A(x) (\log x)^{- c_1 + 1},\]
and so
\[
S_2(x) \ll A(x) (\log x)^{-c_1 + 1} + A(x/z) (\log x)^{\kappa''} \ll
x (\log x)^{-c_1 + 1} + x (\log x)^{-c_6 + \kappa''}.\]
The sum $S_3(x)$ can be bounded by $x (\log x)^{-c_1+1} + x(\log x)^{-c_6+
\kappa''}$
in exactly the same fashion. Thus, it remains only to bound $S_4$ and $S_7$.
The complications to follow are due to the gap between $v(x)$ and 
$\sqrt{D(x)}$. When there is no such gap, $S_7$ disappears and $S_4$ can
be bounded much more simply; see Appendix \ref{sec:gebouw}.

We will bound $S_7$ first. Let $\{\lambda_d\}$ be a Rosser-Iwaniec
 sieve for the
primes $\{p: u y^{-1} <p \leq w u^{-1}\}$ with upper cut $w u^{-1}$.
By definition,
\[\begin{aligned}
\lambda_1=1,\;\;\lambda_d&=0 \text{  if $d\leq u y^{-1} 
\text{ or } d> w u^{-1}$.}\\
\lambda_d&=0 \text{  if $p|d$ for some $p\leq u y^{-1}$.}\end{aligned}\]
Hence \[ 1 = \sum_{d|n} \lambda_d -\; 
\mathop{\sum_{u y^{-1}<d\leq w u^{-1}}}_{d|n}
\lambda_d \]
for every $d$. Substituting into $\beta_7(n)$, we obtain
\[\label{eq:adam}\begin{aligned}
\beta_7(n) &= \sum_{b c|n} \sum_{d|b} \lambda_d \mu(y<b\leq u) \mu(y<c\leq u)
\\ &- \sum_{b c|n} \mathop{\sum_{u y^{-1}<d\leq w u^{-1}}}_{d|b}
\lambda_d \mu(y<b\leq u) \mu(y<c\leq u) .\end{aligned}\]
We give the names $\beta_8(n)$ and $\beta_9(n)$ to the terms on the right
side of (\ref{eq:adam}). Let
\[S_8(x) = \sum_{n=1}^x \beta_8(n) a_n,\;\;\;
  S_9(x) = \sum_{n=1}^x \beta_9(n) a_n .\]

Let us begin by bounding $S_8$. The main idea should be clear: since 
$\sum_{d|b} \lambda_d$ is small for most $b$, one would think that
$\beta_8(n)$ is small as well. We must proceed with caution, however.
It is only here, and in the corresponding part for $S_4$, that we will
have to incur in error bound greater than $O(A(x) (\log X)^{-B})$.

We will have to resolve two issues.
The domain $y<c\leq u$ of $\mu(y<c\leq u)$ may be wide enough to ruin
a naive bound, and, in addition, $b c$ may be too large for (\ref{eq:axl})

We write
\[S_8(x) = \sum_{y<b\leq u} \sum_{d|b} \lambda_d \mu(y<b\leq u)
 \sum_{h\leq x/b} \sum_{c|h} \mu(y<c\leq u) a_{b h} .\]
We would like to bound $\sum_{h\leq x/b} \sum_{c|h} \mu(y<c\leq u) a_h$.
Now
\[\begin{aligned}
\sum_{c|h} \mu(y<c\leq u) &= \sum_{c|h} \mu(c\leq u) - 
                                \sum_{c|h} \mu(c\leq y) \\
&= \sum_{c|h} \mu(c) - \sum_{c|h} \mu(c>u) - \sum_{c|h} \mu(c\leq y) \\
&= \lbrack h=1\rbrack - \sum_{c|h} \mu(c>u) - \sum_{c|h} \mu(c\leq y) .
\end{aligned}\]

Since $h\geq c\geq y\geq 1$, we may ignore the case $h=1$. We shall bound
\begin{equation}\label{eq:arc1}
\sum_{h\leq x/b} \left| 
\sum_{c|h} \mu(c\leq y) a_{b h} \right| .\end{equation}
Let us first look at the other term, viz.
$\sum_{h\leq x/b} \left| \sum_{c|h} \mu(c\leq y) a_{b h} \right|$.
Clearly
\[
\sum_{c|h} \mu(c>u) a_{b h} = \sum_{c|h} \lbrack c> u \rbrack \mu(c) a_{b h}
= \sum_{c|h} \lbrack h/c > u\rbrack \mu(h/c) a_{b h} .\]
For $h$ square-free,
\[\sum_{c|h} \lbrack h/c>u\rbrack \mu(h/c) a_{b h} = \mu(h) 
\sum_{c|h} \lbrack h/c > u \rbrack \mu(c) a_{b c} = 
\mu(h) \sum_{c|h} \mu(c<h/u) a_{b c} .\]
(The expression for $h$ having a small square factor is in essence the same;
values of $h$ with large square factors can be eliminated.) Hence
\begin{equation}\label{eq:ych}
\sum_{h\leq x/b} \left|\sum_{c|h} \mu(c>u) a_{b h}\right| =
\sum_{h\leq x/b} \left|\sum_{c|h} \mu(c<h/u) a_{b h}\right| .\end{equation}
Since $b h/u\leq x/u\leq D(x)$, the right side of (see (\ref{eq:axl}))
can be bounded like (\ref{eq:arc1}). Let us proceed to bound
(\ref{eq:arc1}).

Suppose $h$ has a prime divisor $p\leq l$, where $l$ is a fixed positive
integer. Then the set of all square-free divisors of $h$ can be partitioned
into pairs $(c,c p)$. Clearly $\mu(c) = - \mu(c p)$. Moreover, we have
either $c\leq y$, $c p\leq y$ or $c>y$, $c p > y$, unless $c$ lies in the
range $y/l<c\leq y$. Thus, all pairs $(c, c p)$ that make a contribution
to $\sum_{c|h} \mu(c\leq y)$ satisfy $y/l <c\leq y$. Hence
\[\left|\sum_{c|h} \mu(c<y) \right| \leq \mathop{\sum_{c|h}}_{y/l<c\leq y} 1 .\] Now define
\[l_0 = 2 = 2^{2^0},\;l_1 = 3 = 2^{2^1},\dotsc,\, h_j = 2^{2^j},\dotsc\]
Note that $x^{1/2} < h_{\lfloor \log_2 \log_2 x\rfloor} \leq x$. Let
\[\begin{aligned}
L_0 &= \{\text{even numbers}\},\\
L_j &= \{h\in \mathbb{Z} : (\exists p\leq l_j \text{ s.t. } p|h) \wedge
(\forall p\leq l_{j-1}, p\nmid h)\} .\end{aligned}\]
Then, by the above,
\[\label{eq:arcs2}\begin{aligned}
\mathop{\sum_{h\leq x/b}}_{h\in L_j} 
\left|\sum_{c|h} \mu(c\leq y)\right|
a_{b h} &\leq \mathop{\sum_{h\leq x/b}}_{h\in L_j} 
\mathop{\sum_{c|h}}_{y/l_j<c\leq y} a_{b h} \\ &\leq
\mathop{\sum_{y/l_j < c\leq y}}_{p|c \Rightarrow p>l_{j-1}} 
\sum_{k\leq x/ b c} a_{b c k} .\end{aligned}\]
By (\ref{eq:axl}) and the fact that $b c\leq y^2 \leq D$,
\[\sum_{k\leq x/ b c} a_{b c k} = A_{b c}(x) \sim \frac{A(x)}{b c} .\]
Hence
\begin{equation}\label{eq:guev}
 \mathop{\sum_{y/l_j < c\leq y}}_{p|c \Rightarrow p>l_{j-1}} 
\sum_{k\leq x/ b c} a_{b c k} \sim \frac{A(x)}{b}
\mathop{\sum_{y/l_j < c\leq y}}_{p|c \Rightarrow p>l_{j-1}} \frac{1}{c} \ll
\frac{A(x)}{b} \frac{\log l_j}{\log l_{j-1}} = \frac{2 A(x)}{b} .
\end{equation}
Considering all sets $L_0$, $L_1$, \ldots, 
$L_{\lfloor \log_2 \log_2 x\rfloor}$, we obtain
\[\sum_{h\leq x/b} \left|\sum_{c|h} \mu(c\leq y)\right| a_{b h} \ll
\frac{A(x)}{b} \log \log x .\]
We conclude that
\begin{equation}\label{eq:boliv}\begin{aligned}
|S_8(x)|&\leq \sum_{y<b\leq u} \sum_{d|b} \lambda_d \sum_{h\leq x/b} \left|
\sum_{c|h} \mu(y<c\leq u)\right| a_{b h} \\
&\ll \sum_{y<b\leq u} \sum_{d|b} \lambda_d \frac{A(x)}{b} \log \log x .
\end{aligned}\end{equation}
(Notice that $\sum_{d|b} \lambda_d$ is always
non-negative.) Since 
\[\sum_{b\leq a} \left(\sum_{d|b} \lambda_d\right) \ll
\frac{a}{(\log w u^{-1})/(\log u y^{-1})} ,\]
we can easily see that 
\[\sum_{y<b\leq u} \left(\sum_{d|b} \lambda_d\right) \frac{1}{b} \ll
\frac{\log u y^{-1}}{(\log w u^{-1})/(\log u y^{-1})} \ll
\frac{(\log z)^2}{\log x} .\]
Therefore
\begin{equation}\label{eq:waveh}
|S_8(x)| \leq \frac{(\log z)^2 \log \log x}{\log x} A(x) .
\end{equation}
It is time to bound $S_9(x)$. We change the order of summation:
\[\begin{aligned}
S_9(x) &= \sum_{y<c\leq u} \mu(c) \sum_{u y^{-1} \leq d\leq w u^{-1}} \lambda_d
 \sum_{y/d <h \leq u/d} \mu(h d) \sum_{k\leq x/ c d h} a_{c d h k} \\
 &= \sum_{u<s\leq w} \mathop{\sum_{d|s}}_{u y^{-1} \leq d \leq w u^{-1}}
 \lambda_d \mu(s/d)
\mu(d) \mathop{\sum_{y/d <h \leq u/d}}_{\gcd(h,d)=1} \mu(h)
\sum_{k\leq x/ c d h} a_{c d h k} .\end{aligned}\]

Since $d$ has no small factors when $\lambda_d\ne 0$, it is a simple matter
to remove the condition $\gcd(h,d)=1$ with an error of at most
$O((\log x)^3/\log z)$. We can make the intervals of summation of $d$ and $h$
independent from each other by slicing $\lbrack u y^{-1}, w u^{-1}\rbrack$
into intervals of the form $\lbrack l,l (1+(\log x))^{-c} )$. There
are at most $O((\log x)^{c+1})$ such intervals. We obtain
\begin{equation}\label{eq:prith}S_9(x) \ll (\log x)^{-c'} A(x) + (\log x)^{c+1} 
\max_{u y^{-1} \leq K \leq w u^{-1}} \left|
\sum_{u < s\leq w} f_K(r) g_K(s) a_{r s} \right| ,\end{equation}
where
\[\begin{aligned}
f_K(r) &= \mathop{\sum_{h|r}}_{y/K < h \leq u/K} \mu(h) ,\\ 
g_K(s) &= \mathop{\sum_{d|s}}_{K\leq d < K (1 + (\log N)^c)} 
\lambda_d \mu(d) \mu(s/d). \end{aligned}\]
We can check that $g_K(s)$ averages to zero over $s\equiv a \mo m$ 
as we did in (\ref{eq:anbir}). 
Hence we can apply the bilinear axiom (\ref{eq:bili}):
\[
\sum_{u<r\leq w} f(r) g(s) \ll A(x) (\log x)^{- c_1 + 1} .\]
Thus
\[S_9(x) \ll (\log x)^{-c'} A(x) + (\log x)^{- c_1 + c + 1} . \]
Remember that we may set $c_1$ to an arbitrarily high value.

It remains to bound $S_4$. We can write
\[\label{eq:vwill}\begin{aligned}
\beta_4(n) &= \sum_{b c|n} \mu(b>w) \mu(c>w)\\
&= \sum_{b c|n} \sum_{d|b} \lambda_d \mu(b>w) \mu(c>w) \\
&- \sum_{b c | n} \mathop{\sum_{u y^{-1} \leq d \leq w u^{-1}}}_{d|b}
\lambda_d \mu(b>w) \mu(c>w) .\end{aligned}\]
We give the names $\beta_{10}$ and $\beta_{11}$ to the terms on the
right side of (\ref{eq:vwill}). Let
\[\begin{aligned} S_{10}(x) &= \sum_{n=1}^x \beta_{10}(n) a_n,\\
S_{11}(x) &= \sum_{n=1}^x \beta_{11}(n) a_n .\end{aligned}\]
We bound $S_{10}(x)$ as we bounded $S_8(x)$. We can obtain an expression 
similar to
(\ref{eq:prith}) for $S_{11}(x)$:
\[S_{11} \ll (\log x)^{-c'} A(x) + (\log x)^{c+1}
 \max_{x w^{-2} \leq K\leq x w^{-1} u^{-1}} 
\left|\sum_{u<s\leq w} f_K(r) g_K(s)\right|,\]
where
\[\begin{aligned}
f_K(s) &= \mathop{\sum_{d|s}}_{K\leq d < K (1 + (\log N)^c)} 
\lambda_d \mu(d) \mu(s/d)\\
g_K(r) &= \mathop{\sum_{h|r}}_{y/K < h \leq u/K} \mu(h) . \end{aligned}\]
Again, we apply (\ref{eq:bili}) and are done:
\[S_{11}(x) \ll (\log x)^{-c'} A(x) + (\log x)^{- c_1 + 1} . \]

We conclude by (\ref{eq:lins}) that
\[\sum_{n=1}^x \mu(n) a_n \ll \frac{(\log z)^2 \log \log x}{\log x} A(x) .\]

\begin{center}
* * *
\end{center}
In the course of the actual procedure we are about to undertake, we will
come across some technical difficulties not present in the above outline.
For example, we will be forced to sieve over ideals and ideal numbers
rather than over rational integers. Our linear sieve axioms will be valid
only on average, unlike, say, (\ref{eq:axl}). Nevertheless, we will 
be able to follow, in the main, the plan we have traced.

As the method we have devised to eliminate a bothersome interval may have
wider applications, it may be worthwhile to review its main idea. We are
given the task of estimating a sum
\[\sum_{a,b \leq X} F_{a b} .\]
We assume we know how to estimate
\begin{equation}\label{eq:jdc}
\mathop{\sum_{a,b \leq X}}_{a\leq x (z(x))^{-1}} F_{a b} \;\;\text{and}
\;\;\mathop{\sum_{a,b \leq X}}_{a\geq x z(x)} F_{a b} ,\end{equation}
where $\log z(x) = o(\sqrt{\log x})$. In order to eliminate the missing
interval, we apply a sieve to the constant function $a\mapsto 1$ with
respect to the primes larger than $z^2(x)$:
\begin{equation}\label{eq:twinkl}
\mathop{\sum_{a,b \leq X}}_{x (z(x))^{-1} \leq a \leq x z(x)} F_{a b}
= \mathop{\sum_{a,b \leq X}}_{x (z(x))^{-1} \leq a \leq x z(x)} 
\sum_{d|a} \lambda_d F_{a b} - 
\mathop{\sum_{a,b \leq X}}_{x (z(x))^{-1} \leq a \leq x z(x)} 
\mathop{\sum_{d|a}}_{d> z^2(x)} \lambda_d F_{a b} .\end{equation}
(Notice the peculiar 
use of a sieve as an identity rather than an approximation.) The first term
on the right can be seen from sieve theory to be at most
\[O\left(\frac{\log z(x)}{\log x} \cdot (\log x z(x) - \log x (z(x))^{-1})\right)
\cdot X .
\]
The second term on the right of (\ref{eq:twinkl}) can be treated 
analogously to the first sum in (\ref{eq:jdc}) with variables $a'=a/d$
and $b'= b d$; clearly $a' b'\leq X$ and $a'\leq x (z(x))^{-1}$ .

\subsection{Axioms}
Let $K/\mathbb{Q}$ be a cubic extension of $\mathbb{Q}$. 
Let $k_0$ be a fixed rational integer. 
Define
\begin{equation}\label{eq:defr}
\mathcal{R} = \{\mathfrak{r} : \mathfrak{r} \in I_K,
\mu_K(\mathfrak{r})^2 = 1,\, \mu_K(N(\mathfrak{r}/\gcd(k_0,\mathfrak{r})))=1\} 
.\end{equation}
We write $\mu_{\mathcal{R}}$ for the M\"obius function with respect
to $\mathcal{R}$:
\[\begin{aligned}
\mu_{\mathcal{R}}(\mathfrak{a}) &=
\mathop{\prod_{\mathfrak{p}|\mathfrak{a}}}_{\mathfrak{p}\in \mathcal{R}} (-1)
&\text{if $\mathfrak{a}$ is square-free,}\\
\mu_{\mathcal{R}}(\mathfrak{a}) &= 0 &\text{  otherwise.}\end{aligned}\]
We are given a bounded sequence 
$\{a_{\mathfrak{r}}\}_{\mathfrak{r}\in \mathcal{R}}$ of non-negative
real numbers, the properties of whose distribution we will now describe.

We abuse notation by writing $\mathfrak{a} < x$, $\mathfrak{a}>x$
when we mean $N \mathfrak{a} < x$, $N \mathfrak{a} > x$;
$N \mathfrak{a} < N \mathfrak{b}$ will, however, still mean
$N \mathfrak{a} < N \mathfrak{b}$. For $\mathfrak{d}\in \mathcal{R}$,
define
\[A_{\mathfrak{d}}(x) = \mathop{\sum_{\mathfrak{d}|\mathfrak{n}}}_{
\mathfrak{n}\leq x} a_{\mathfrak{n}},\;\;\;
A(x) = \sum_{\mathfrak{n} \leq x} a_{\mathfrak{n}} .\]
Write
\begin{equation}\label{eq:lin}
A_{\mathfrak{d}}(x) = \gamma(\mathfrak{d}) A(x) + r_{\mathfrak{d}},
\end{equation}
where $\gamma$ is a bounded multiplicative function supported on 
$\mathcal{R}$ and $r_{\mathfrak{d}}$ is an error term.
We assume our estimates on $\gamma$ to be quite strong for
all primes above $(\log X)^{\kappa}$:
\[\sum_{\mathfrak{p} \leq x} \gamma(\mathfrak{p}) =
 \log \log x + \alpha + O((\log x)^{-B}) ,\]
for any $x> (\log X)^{\kappa}$, some constant $\alpha$ and
 any constant $B>0$, where the
implied constant depends on $B$. Let us be more precise and make
clear that what we are avoiding the divisors of a fixed rational
integer $\delta \leq (\log X)^{\kappa}$:
\begin{equation}\label{eq:bloodyp}
\mathop{\sum_{\mathfrak{p} \leq x}}_{\mathfrak{p}\nmid \delta}
\gamma(\mathfrak{p}) = \log \log x + \alpha + O((\log x)^{-B}) .
\end{equation}
We will also allow ourselves the relative luxury of the following assumption
on the size of $\gamma(\mathfrak{d})$:
\begin{equation}\label{eq:pyram}
\gamma(\mathfrak{d}) \ll 1 / N \mathfrak{d} .\end{equation}
Condition (\ref{eq:pyram}) will be fulfilled for the sequence we are
ultimately interested in. It is possible to replace (\ref{eq:pyram})
with an average condition; see the remark after (\ref{eq:bonnie}).

We have an average bound for the remainder terms $r_{\mathfrak{d}}$:
for any $B_1,B_2>0$, there is a $C>0$ such that
\begin{equation}\label{eq:rbound}
\sum_{\mathfrak{d} \leq x^{2/3} (\log x)^{-C}} \tau^{B_1}(\mathfrak{d})
r_\mathfrak{d}
\ll (\log x)^{-B_2} A(x) .\end{equation}
Typically, $A(x)$ will be about a constant times $x^{2/3}$.
We will assume the consequences
\begin{equation}\label{eq:svulg}
A(x) \gg x^{1/2},
\end{equation}
\begin{equation}\label{eq:size} A(x/z) \ll (\log x)^{-B} A(x)\end{equation}
for any $z$ such that $\log \log x/\log z = o(1)$.

We assume the following axiom.
\begin{bilinc} Let $f, g:I_K\to \mathbb{R}$ satisfy
\begin{equation}\label{eq:comman}
|f(\mathfrak{a})|, |g(\mathfrak{a})| \ll \tau^2(\mathfrak{a}) .\end{equation}
Assume $g$ is a linear combination of the form
\begin{equation}\label{eq:sarand}
g(\mathfrak{a}) = \sum_{\mathfrak{d}|\mathfrak{a}} c_{\mathfrak{d}}
\mu_{\mathcal{R}}(\mathfrak{d}>\ell)\end{equation}
or
\begin{equation}\label{eq:sus}
g(\mathfrak{a}) = \sum_{\mathfrak{d}|\mathfrak{a}} c_{\mathfrak{d}}
\mu_{\mathcal{R}}'(\mathfrak{d}>\ell),\end{equation}
where \[\mu_{\mathcal{R}}' = \mu_{\mathcal{R}} \cdot (
\mathfrak{p} \leq (\log x)^{10} \Rightarrow \mathfrak{p}\nmid \mathfrak{d}),\]
the sequence $c_{\mathfrak{d}}$ is bounded and $\ell>x^{1/\kappa}$
for some constant $\kappa$. We assume furthermore that either $f$ or
$g$ is zero on all numbers with small prime divisors:
\[\mathfrak{p}|\mathfrak{a}, \mathfrak{q}|\mathfrak{b},
\mathfrak{p},\mathfrak{q} \leq (\log x)^{10} \Rightarrow
f(\mathfrak{a}) g(\mathfrak{b}) = 0.\]

Then
\begin{equation}\label{eq:darme}
\mathop{\sum_{\mathfrak{a} \mathfrak{b}\leq x}}_{
x^{1/2} (\log x)^T < N \mathfrak{b} \leq x^{3/2} (\log x)^{-T}}
f(\mathfrak{a}) g(\mathfrak{b}) \ll A(x) (\log x)^{-2},\end{equation}
where $T$ is a constant depending only on  $B$ and on the implied
constant in (\ref{eq:comman}).
\end{bilinc}

Write $P(z)$ for $\prod_{\mathfrak{p}<z} \mathfrak{p}$. Write $P_{10}$ for
$P((\log x)^{10})$. Let 
\[\sum_{*} \dotsb\]
be short for
\[\mathop{\sum_{\mathfrak{b} \mathfrak{c} | \mathfrak{n}}}_{\gcd(\mathfrak{n}/\mathfrak{b},
P_{10}^{\infty})=1} \dotsb\]

We will follow a convention we have already implicitly
used in this subsection: $\kappa$ is a fixed constant given by
the sequence $\{a_{\mathfrak{n}}\}$, and we should be ready for it to
be arbitrarily large, but fixed; $B$ is a parameter that we
can set to be arbitrarily large given our axioms (example: 
``the number of primes in arithmetic progressions of modulus up to
$(\log x)^B$ is \ldots''); finally,
$C$ is a parameter that may have to be
taken to be large if a condition is to be satisfied for a chosen value of
$B$.
\subsection{Technical lemmas}
\begin{lem}\label{lem:ough}
Assume (\ref{eq:bloodyp}). Then, for any $B>0$,
\[\mathop{\sum_{\mathfrak{d}\leq y}}_{(\mathfrak{d},\mathfrak{m})=1} \mu(\mathfrak{d}) g(\mathfrak{d}) \ll (\log y)^{-B} + (\log y)^3
\mathop{\sum_{y^{(\log \log y)^{-2}} \leq \mathfrak{p} \leq y}}_{\mathfrak{p}|\mathfrak{m}} \frac{1}{N \mathfrak{p}} .\]
\end{lem}
\begin{proof}
As in \cite{FI2}, pp. 1048--1049.
\end{proof}
\begin{lem}\label{lem:pemmican}
Assume (\ref{eq:rbound}) and (\ref{eq:bloodyp}). Then
\[\sum_{\mathfrak{n}\leq x} \tau^4(\mathfrak{n}) a_{\mathfrak{n}} \ll (\log x)^{16} A(x).\]
\end{lem}
\begin{proof}
As in \cite{FI2}, p. 1047.
\end{proof}
\begin{lem}\label{lem:strawb}
Assume (\ref{eq:svulg}), (\ref{eq:pyram}) and (\ref{eq:rbound}). Then
\[\sum_{\mathfrak{n}\leq x} \lbrack \mathfrak{n} \leq x^{4/11} \gcd(\mathfrak{n},P_{10}^{\infty})
\rbrack \gcd(\mathfrak{n},P_{10}^{\infty}) a_{\mathfrak{n}} \ll A(x) (\log x)^{-B} .\]
\end{lem}
\begin{proof}
Clearly
\[\begin{aligned}
\sum_{\mathfrak{n} \leq x} \lbrack \mathfrak{n} \leq x^{4/11} \gcd(\mathfrak{n},P_{10}^{\infty})
\rbrack a_{\mathfrak{n}} &\leq \sum_{\mathfrak{b} \leq x^{4/11}} 
 \mathop{\sum_{\mathfrak{c} | P_{10}^{\infty}}}_{\mathfrak{b} \mathfrak{c} \leq x} a_{\mathfrak{b}
\mathfrak{c}} \\ &\leq
\sum_{\mathfrak{b} \leq x^{4/11}} \sum_{\mathfrak{c}\leq x^{1/11}} a_{\mathfrak{b} \mathfrak{c}} \\
&+
 \sum_{\mathfrak{b} \leq x^{1/11}} \mathop{\sum_{\mathfrak{c} | P_{10}^{\infty}}}_{
x^{1/11} < \mathfrak{c} \leq x^{1/11} (\log x)^{10}} \mathop{\sum_{\mathfrak{d}}}_{\mathfrak{b}
\mathfrak{c} \mathfrak{d} \leq x} a_{\mathfrak{b} \mathfrak{c} \mathfrak{d}} \\
&\leq x^{5/11} + A(x) (\log x)^{-B} \\ &+ A(x) \sum_{\mathfrak{b} \leq x^{4/11}}
 \mathop{\sum_{\mathfrak{c} | P_{10}^{\infty}}}_{x^{1/11} \leq \mathfrak{c} \leq
x^{1/11} (\log x)^{10}} \gamma(\mathfrak{b} \mathfrak{c}) .\end{aligned}\]
The cardinality of $\{\mathfrak{c}\leq x^{1/11} (\log x)^{10} : \mathfrak{c} |P_{10}^{\infty} \}$
can be crudely estimated by means of Rankin's trick:
\[\begin{aligned}
\# \{\mathfrak{c} \leq m : \mathfrak{c} |P_{10}^{\infty}\} &\leq
 \sum_{\mathfrak{c} |P_{10}^{\infty}} \frac{m^{9/10}}{(N \mathfrak{c})^{9/10}}
= m^{9/10} \prod_{\mathfrak{p} | P_{10}^{\infty}} \frac{1}{1 - (N \mathfrak{p})^{-9/10}}\\
&\sim m^{9/10} e^{\sum_{\mathfrak{p}|P_{10}^{\infty}} (N \mathfrak{p})^{-9/10}}
\ll m^{9/10} e^{C (\log x)/(\log \log x)} \ll m^{9/10+\epsilon} .\end{aligned}\]
Hence
\[\mathop{\sum_{\mathfrak{c} | P_{10}^{\infty}}}_{x^{1/11}\leq \mathfrak{c} \leq x^{1/11} 
(\log x)^{10}} \frac{1}{N \mathfrak{c}} \ll x^{-1/110 + \epsilon}\]
and thus
\[\sum_{\mathfrak{b} \leq x^{4/11 + \epsilon}} 
\mathop{\sum_{\mathfrak{c} | P_{10}^{\infty}}}_{x^{1/11}\leq \mathfrak{c} \leq x^{1/11} (\log x)^{10}}
\gamma(\mathfrak{b} \mathfrak{c}) \ll (\log x) x^{-1/110 + \epsilon} .\]
\end{proof}
\subsection{Bounds and manipulations}\label{subs:agarag}
Let $z = e^{(\log \log x) (\log \log \log x)^{1/2}}$,
$y = x^{1/3} z^{-2}$, $u = x^{1/3} z$, $w = x^{1/2} z^{-1}$.
As in
(\ref{eq:matdi}) and (\ref{eq:lins}),
\[\begin{aligned}
\mu_{\mathcal{R}}(\mathfrak{n}) &= \beta_1(\mathfrak{n}) + \beta_2(\mathfrak{n}) +
\beta_3(\mathfrak{n}) + \beta_4(\mathfrak{n}) - \beta_5(\mathfrak{n}) -
\beta_6(\mathfrak{n}) - \beta_7(\mathfrak{n}),\\
\sum_{\mathfrak{n} \leq x} \mu_{\mathcal{R}}(\mathfrak{n}) a_{\mathfrak{n}}
 &= S_1(x) + S_2(x) + S_3(x) + S_4(x) - S_5(x) - S_6(x) - S_7(x),
\end{aligned}\]
where
\[\begin{aligned}
\beta_1(\mathfrak{n}) &= \mu_{\mathcal{R}}(\mathfrak{n}\leq u) +
\sum_* \mu_{\mathcal{R}}(\mathfrak{b})
 \mu_{\mathcal{R}}(\mathfrak{c}\leq u),\\
\beta_2(\mathfrak{n}) &= 
\sum_* \mu_{\mathcal{R}}(u<\mathfrak{b}\leq w) \mu_{\mathcal{R}}(\mathfrak{c}>u),\\
\beta_3(\mathfrak{n}) &= \sum_* \mu_{\mathcal{R}}(\mathfrak{b}> w) 
\mu_{\mathcal{R}}(u<\mathfrak{c}\leq w),
\end{aligned}\]
\[\begin{aligned}
\beta_4(\mathfrak{n}) &= \sum_*
\mu_{\mathcal{R}}(\mathfrak{b}> w) \mu_{\mathcal{R}}(\mathfrak{c}>w),\\
\beta_5(\mathfrak{n}) &= \sum_* 
\mu_{\mathcal{R}}(\mathfrak{b} \leq u) \mu_{\mathcal{R}}(\mathfrak{c} \leq y),\\
\beta_6(\mathfrak{n}) &= \sum_* 
\mu_{\mathcal{R}}(\mathfrak{b} \leq y) \mu_{\mathcal{R}}(y<\mathfrak{c} \leq u),\\
\beta_7(\mathfrak{n}) &= \sum_* 
\mu_{\mathcal{R}}(y<\mathfrak{b} \leq u) \mu_{\mathcal{R}}(y<\mathfrak{c} 
\leq u),
\end{aligned}\]
and
\[S_j(x) = \sum_{\mathfrak{n} \leq x} \beta_j(n) a_{\mathfrak{n}} .\]
Clearly
\[\begin{aligned}
S_1(x) &= \sum_{\mathfrak{n} \leq u} \mu_{\mathcal{R}}(\mathfrak{n}) a_{\mathfrak{n}} +
\sum_{\mathfrak{n}\leq x} \sum_* \mu_{\mathcal{R}}(\mathfrak{b})
 \mu_{\mathcal{R}}(\mathfrak{c} \leq u) a_{\mathfrak{n}}\\
&= O(A(u)) + \sum_{\mathfrak{n}\leq x} \mu_{\mathcal{R}}(\gcd(\mathfrak{n},P_{10}^{\infty}))
\mu_{\mathcal{R}}(\mathfrak{n}/\gcd(\mathfrak{n},P_{10}^{\infty})) a_{\mathfrak{n}} \\
&= O(A(u)) + \sum_{\mathfrak{n}\leq x} \lbrack \mathfrak{n} \leq u \gcd(\mathfrak{n}, 
P_{10}^{\infty})\rbrack a_{\mathfrak{n}} .\end{aligned}\]
By (\ref{eq:size}) and Lemma \ref{lem:strawb}, we can conclude that
\[S_1(x) \ll (\log x)^{-B} A(x) .\]
We can rewrite $S_5$ as follows:
\[
S_5(x) = 
\mathop{\sum_{\mathfrak{n} \leq x}} \sum_*
h(\mathfrak{b} \leq u) \mu_{\mathcal{R}}(\mathfrak{c} \leq y)
\mathop{\sum_{\mathfrak{d}\leq x/u y}}_{\mathfrak{p} | \mathfrak{d}
\Rightarrow \mathfrak{p} > (\log x)^{10}} a_{\mathfrak{b} \mathfrak{c}
\mathfrak{d}} .\]
 Since $\frac{\log(x^{2/3}/ /((\log x)^C u y))}{\log \log x^{10}}
\gg (\log \log x) (\log \log \log x)$, we can apply the fundamental
lemma of sieve theory (vd., e.g., \cite{HR}, Ch. 2, or \cite{Iw2},
Lem 2.5) to obtain
\[\begin{aligned}
\mathop{\sum_{\mathfrak{d}\leq x/u y}}_{\mathfrak{p} | \mathfrak{d}
\Rightarrow \mathfrak{p} > (\log x)^{10}} a_{\mathfrak{b} \mathfrak{c}
\mathfrak{d}} &= V_{\mathfrak{b} \mathfrak{c}}
 X (1 + O(e^{-(\log \log x) (\log \log \log x)})) + \text{error}\\
&= V_{\mathfrak{b} \mathfrak{c}} X ( 1 + O(1/(\log x)^{\log \log \log x}))
+ \text{error},\end{aligned}\]
where the error term is collected by (\ref{eq:rbound}), and the leading term in the main term 
is given by \[V_{\mathfrak{b} \mathfrak{c}} = 
\mathop{\prod_{\mathfrak{p} \leq (\log x)^{10}}}_{\mathfrak{p}\nmid \mathfrak{b} \mathfrak{c}}
         (1 - \gamma'(\mathfrak{p})),\]
where $\gamma'(\mathfrak{p}) = \gamma(\mathfrak{p})$ for $\mathfrak{p}\nmid \mathfrak{b} \mathfrak{c}$,
$\gamma'(\mathfrak{p}) = 0$ for $\mathfrak{p}\nmid \mathfrak{b} \mathfrak{c}$,
$\mathfrak{p}\nmid k_0$. We then apply Lemma \ref{lem:ough} and obtain
 \[S_5(x) \ll A(x)/ (\log x)^{B} .\]
In the same way,
 \[S_6(x) \ll A(x)/(\log x)^{B} .\]
As in subsection \ref{subs:recabarren}, we have
\[S_2(x) = \mathop{\sum_{\mathfrak{r} \mathfrak{s} \leq x}}_{
x/z w \leq \mathfrak{s} \leq x/w} f(\mathfrak{r}) g(\mathfrak{s})
a_{\mathfrak{r} \mathfrak{s}} ,\]
where
\[\begin{aligned}
f(\mathfrak{r}) &= \mathop{\sum_{\mathfrak{b}|\mathfrak{r}}}_{\gcd(\mathfrak{r}/\mathfrak{b},P_{10})
= 1} h(u < \mathfrak{b} \leq w),\\
g(\mathfrak{s}) &= \mu_{\mathcal{R}}'(u<\mathfrak{s}<w) .\end{aligned}\]
By the bilinear condition (\ref{eq:darme}),
\[\mathop{\sum_{\mathfrak{r} \mathfrak{s} \leq x}}_{x/z w\leq
\mathfrak{s} \leq x/u} f(\mathfrak{r}) g(\mathfrak{s}) \ll
 A(x) (\log x)^{-B} .\]
By Lemma \ref{lem:pemmican},
\[\sum_{n\leq x/z} 
 \tau_3(n) a_n \ll A(x/z) (\log x)^{\kappa} .\]
Hence
\[\begin{aligned}
S_2(x) &\ll A(x/z) (\log x)^{\kappa} + A(x) (\log x)^{-B} \\
 &\ll A(x) (\log x)^{\kappa}/z^2 + A(x) (\log x)^{-B} .\end{aligned}\]
In the same way,
\[S_3(x) = 
 \mathop{\sum_{\mathfrak{r} \mathfrak{s} \leq x}}_{x/z w \leq
\mathfrak{s} \leq x/w} f(\mathfrak{r}) g(\mathfrak{s}) a_{
\mathfrak{r} \mathfrak{s}} + O\left(\sum_{n\leq x/z} 
 \tau_3(n) a_n + A(x/z)\right) ,\]
where
\[\begin{aligned}
f(\mathfrak{r}) &= \mathop{\sum_{\mathfrak{c}|\mathfrak{s}}}_{\gcd(\mathfrak{s}/\mathfrak{c},
P_{10}) = 1} \mu_{\mathcal{R}}'(\mathfrak{c} > w),\\
g(\mathfrak{s}) &= \mu_{\mathcal{R}}(u<\mathfrak{b}\leq w) \end{aligned}\]
and consequently
\[\begin{aligned}
S_3(x) &\ll A(x/z) (\log x)^{\kappa} + A(x) (\log x)^{-B} \\
 &\ll A(x) (\log x)^{\kappa}/z^2 + A(x) (\log x)^{-B} .\end{aligned}\]

It is time to bound $S_7$. Let $\{\lambda_{\mathfrak{d}}\}$ be a 
generalized Rosser-Iwaniec sieve (see, e.g., \cite{Col2}) for the primes
\begin{equation}\label{eq:bonnie}
\{\mathfrak{p}\in \mathcal{R} : u y^{-1} < \mathfrak{p}\leq w u^{-1} \} ,\end{equation}
upper cut $w u^{-1}$ and sieved set $\mathcal{R}$.

{\bf Remark.} We could sieve only up to a fractional power of $w u^{-1}$,
and change our bounds only by a constant as a result -- a constant that
would not necessarily be greater than $1$. A Selberg sieve (see
the generalization in \cite{Ri1}--\cite{Ri3}) would do just as well; its main 
defect for our purposes, namely, its having coefficients that may grow
as fast as the divisor function, is immaterial in the present context.
Notice also that, if we did not have (\ref{eq:pyram}), it would be best to use
$\gamma(\mathfrak{d})$ as our input, instead of $1/N \mathfrak{d}$, which
we implicitly use by choosing $\mathcal{R}$ to be our sieved set. 
We have made the latter choice here for the sake of simplicity: it
is elements of $\mathcal{R}$, not elements of $\{a_{\mathfrak{n}}\}$,
that are being sieved here.

By definition,
\[\begin{aligned}
\lambda_1 =1, \lambda_{\mathfrak{d}} &= 0 \text{\;if $\mathfrak{d}\leq u y^{-1}$ or
$\mathfrak{d}> w u^{-1}$}\\
\lambda_{\mathfrak{d}} &= 0 \text{\;if $\mathfrak{p}|\mathfrak{d}$ for some
$\mathfrak{p}\leq u y^{-1}$.} \end{aligned}\]
Hence
\begin{equation}\label{eq:andy}
1 = \sum_{\mathfrak{d}|\mathfrak{n}} \lambda_{\mathfrak{d}}\; -
\mathop{\sum_{u y^{-1} < \mathfrak{d} \leq w u^{-1}}}_{\mathfrak{d}|
\mathfrak{n}} \lambda_{\mathfrak{d}}
\end{equation}
for every $\mathfrak{d}\in \mathcal{R}$. We substitute (\ref{eq:andy})
into $S_7$:
\[\begin{aligned}
S_7(x) &= \sum_*
 \sum_{\mathfrak{d} | \mathfrak{c}} \lambda_{\mathfrak{d}}
 h (y<\mathfrak{b}\leq u) \mu_{\mathcal{R}}(y<c\leq u)\\
&- \sum_*
 \mathop{\sum_{u y^{-1} < \mathfrak{d} \leq w u^{-1}}}_{
\mathfrak{d} | \mathfrak{c}} \lambda_{\mathfrak{d}}
h(y\leq b<u) \mu_{\mathcal{R}}(y<c\leq u) \\
&= S_8(x) + S_9(x),\end{aligned}\]
say. 
The argument between (\ref{eq:arc1}) and (\ref{eq:waveh}) is unchanged;
we use the upper bound (\ref{eq:pyram}) to bound $\gamma(\mathfrak{d})$.
As a result,
\[S_8(x) \ll \frac{(\log z)^2 \log \log x}{\log x} .\]
We can express $S_9$ as before:
\[S_9(x) = (\log x)^{-B} A(x) + (\log x)^{C+1} \max_{u y^{-1}\leq R\leq w u^{-1}}
\left| \sum_{u< \mathfrak{s}\leq w} f_R(\mathfrak{r}) g_R(\mathfrak{s}) 
a_{\mathfrak{r} \mathfrak{s}} \right|,\]
where
\[\begin{aligned}
f_R(\mathfrak{r}) &= \mathop{\mathop{\sum_{\mathfrak{h}|\mathfrak{r}}}_{
y/K < \mathfrak{h} \leq u/K}}_{\mathfrak{p}<(\log x)^{10} \Rightarrow \mathfrak{p}\nmid 
\mathfrak{r}/\mathfrak{h}} 
h(\mathfrak{h})\\
g_R(\mathfrak{r}) &= \mathop{\sum_{\mathfrak{d}|\mathfrak{s}}}_{
K \leq \mathfrak{d} \leq K (1 + (\log N)^{-C})}
\lambda_{\mathfrak{d}} h(\mathfrak{d}) \mu_\mathcal{R}'(\mathfrak{s}/\mathfrak{d}) 
 .\end{aligned}\]
Notice that the support of $\lambda_{\mathfrak{d}}$ excludes $\lbrack 2,(\log x)^{10}\rbrack$.
We apply the bilinear axiom (\ref{eq:darme}) and obtain
\[S_9(x) \ll A(x) (\log x)^{-B} .\]

Hence
\[S_7(x) \ll \frac{(\log z)^2 \log \log x}{\log x} .\]`
The same bound can be obtained for $S_4$ by nearly the same argument;
see subsection \ref{subs:recabarren}. We conclude that
\[\sum_{\mathfrak{n}\leq x} h(\mathfrak{n}) a_{\mathfrak{n}} \ll
\frac{(\log z)^2 \log \log x}{\log x} \ll \frac{(\log \log x)^5 (\log \log \log x)}{\log x} .\]
It is easy to check that the factor $\log \log \log x$ above can be replaced by any increasing
function $f(x)$ such that $\lim_{x\to \infty} f(x) = \infty$.

\subsection{Background and references for axioms}\label{subs:gustoso}
Let $f(x,y)\in \mathbb{Z}\lbrack x,y\rbrack$ be an irreducible homoegeneous cubic polynoial.
By \cite{HBM}, Lemma 2.1, we can construct a number field $K/\mathbb{Q}$ of degree
$\deg(K/\mathbb{Q})=3$ and two elements $\omega_1, \omega_2 \in \mathfrak{O}_K$
linearly independent over $\mathbb{Z}$ such that
\[f(x,y) = N_{K/\mathbb{Q}} (x \omega_1 + y \omega_2) N \mathfrak{d}^{-1} ,\]
where $\mathfrak{d}$ is the ideal of $\mathfrak{O}_K$ 
generated by $\omega_1$ and $\omega_2$. By \cite{HBM}, Lemmas 2.2 and 2.3, there is a fixed
rational integer $k_0$ such that $(x \omega_1 + y \omega_2) \mathfrak{d}^{-1}$ is always
an element of $\mathcal{R}$, where $\mathcal{R}$ is as in (\ref{eq:defr}); moreover,
\[\mu_R((x \omega_1 + y \omega_2) \mathfrak{d}^{-1} ) = \mu(f(x,y)) .\]
Given $\eta, \upsilon > 0$ and a lattice $L\subset \mathbb{Z}^2$, we define
\begin{equation}\label{eq:defbox}\begin{aligned}
S &=\lbrack X, (1+\eta) X\rbrack \times \lbrack \upsilon X, \upsilon (1 + \eta) X\rbrack\\
\mathcal{A}_{L,S,\omega_i} &= \{(x \omega_1 + y \omega_2) \mathfrak{d}^{-1} : (x,y) \in L \cap S,\,
\gcd(x,y)=1\} .
\end{aligned}\end{equation}
Then
\[\mathop{\sum_{(x,y)\in L \cap S}}_{\gcd(x,y)=1} \mu(f(x,y)) = 
\sum_{\mathfrak{n}\in \mathcal{A}_{L,S,\omega_i}} \mu_{\mathcal{R}}(\mathfrak{n}) .\]
Hence it is natural to define
\[a_{\mathfrak{n}} = \begin{cases} 1 &\text{if $\mathfrak{n}\in \mathcal{A}_{L,S,\omega_i}$,}\\
0 & \text{otherwise.}\end{cases} \]

Let $x_0 = \max_{\mathfrak{a}\in \mathcal{A}_{L,s,\omega_i}} N \mathfrak{a} = X^3 (1 + O(\eta))$.
For $x\leq x_0$, let $A(x) = \sum_{N \mathfrak{n} \leq x} a_{\mathfrak{n}}$. Clearly
\[A(x) \sim \frac{\nu \eta^2 X^2}{\zeta(2) \lbrack \mathbb{Z}^2 : L\rbrack}
 \mathop{\prod_{p|\lbrack \mathbb{Z}^2 : L\rbrack}}_{L\cap p \mathbb{Z}^2 =
\emptyset} ( 1 - p^{-2})^{-1} \mathop{\prod_{p|\lbrack \mathbb{Z}^2 : L\rbrack}}_{L\cap p \mathbb{Z}^2 \ne \emptyset} ( 1 + p^{-1})^{-1} ,\]
provided that $L$ is not contained in any set of the form $p \mathbb{Z}^2$;
if $L\subset p \mathbb{Z}^2$,then $A(x)=0$ and all of our results are trivial.

Assume
\begin{equation}\label{eq:yadnu}\begin{aligned} 
-\log \log N \ll \log \upsilon &\ll \log \log N ,\\
\log \eta &\gg -\log \log N,\\
\eta/\min(\upsilon,\upsilon^{-1}) = o(1) ,\end{aligned}\end{equation}
where the second restriction 
on $\eta$ is enough for us to avoid associated elements in $\mathfrak{O}_K$.

Axioms (\ref{eq:lin})-(\ref{eq:rbound}) are proven for $L = \mathbb{Z}^2$,
$\upsilon = 1$ in \cite{HBM}, sections 2--3; they are proven
for general $L$ in \cite{HBM2}, in a slightly different formulation. 
Since the bound (\ref{eq:rbound}) can absorb powers of $\log x$, and
the introduction of $\upsilon \ne 1$ does not require any change in the proofs,
and the bounds are uniform for $\lbrack \mathbb{Z} : L \rbrack \ll (\log N)^B$,
$B>0$ arbitrary.
Axiom (\ref{eq:size}) is clear.
The bilinear axiom is proven in subsection \ref{subs:biliflow} under the condition
(\ref{eq:rebilicon}). It remains to be seen that all linear combinations of the form
(\ref{eq:sarand}) satisfy (\ref{eq:rebilicon}). Thanks to the standard zero-free regions
for Hecke $L$-functions (see Lemma \ref{lem:corknut}) we know that $\mu_{\mathcal{R}}$
satisfies (\ref{eq:rebilicon}) for $\lbrack \mathbb{Z}^2 : L\rbrack \ll (\log N)^B$
(and the far stronger bound $\ll x e^{-(\log x)^{3/5}/(\log \log x)^{1/5}}$ as well.) 
It then follows by the fundamental lemma of sieve theory
that the function $\mu_{\mathcal{R}}'$ satisfies
(\ref{eq:rebilicon}) as well. To see (\ref{eq:rebilicon}) for linear combinations, note
simply that
\[
\sum_{n\leq x} \left(\sum_{d|n} c_d \mu(n/d > n^{1/\kappa})\right) =
\sum_{d\leq x^{1-1/\kappa}} c_d \sum_{d^{(1-1/\kappa)^{-1}-1} \leq m\leq x/d} \mu(m)\]
In each inner sum, $x/d>x^{1/\kappa}$, and thus $\log(x/d) \gg \log x$. Hence
we bound the inner sum by $C (x/d) (\log x)^{-B}$, $C$ independent of $d$, and obtain a total
bound of at most
\[C x (\log x)^{-B+1} .\]
\subsection{The bilinear condition}\label{subs:biliflow}
This subsection is a summarized paraphrase of \cite{HB}, pp. 66--83, and
\cite{HBM}, pp. 275--284. This rephrasing is necessary because 
the said references carry their argument for a specific function, whose
special properties they use in ultimately inessential ways.

We recapitulate the framework set out in \cite{HBM}, p. \!258 and p. \!277.
We let $K/\mathbb{Q}$ be a number field of degree $\deg(K/\mathbb{Q})=3$.
We are given $\omega_1, \omega_2 \in \mathfrak{O}_K$ linearly independent over
$\mathbb{Z}$. Let $\mathfrak{d}\in \mathfrak{O}_K \omega_1 + \mathfrak{O}_K
\omega_2$. Let $\delta$ be an arbitrary element of 
$\mathcal{I}^{-1}(\mathfrak{d})$, that is, an ideal number corresponding
to $\mathfrak{d}$.

Every class $A\in \mathcal{C}_1(K)$ is a $\mathbb{Z}$-module and as such
has a basis
$\{w_{A,1},\dotsb , w_{A,3}\}$ consisting of elements 
of $\mathcal{I}(\mathfrak{O}_K)^{\times}$. For $A_0 = \cl \delta^{-1}$,
we can choose $\{w_{A_0,1}, w_{A_0,2},w_{A_0,3}\}$ so that
$\omega_1 \delta^{-1} = w_{A_0,1}$ and
$\omega_2 \delta^{-1} = z w_{A_0,2}$ for some $z\in \mathbb{Z}$.
For other classes $A\in \mathcal{C}_1(K)$ we make the choice of basis
$\{w_{A,1},\dotsb , w_{A,3}\}$ arbitrarily.

Let $\beta\in \mathcal{I}(\mathfrak{O}_K)^{\times}$. Let 
$A_{\beta} = \cl(\beta \delta)^{-1}$. Write
\[\begin{aligned}
\beta w_{A_{\beta},1} &= q_{1 1} w_{A_0,1} + q_{1 2} w_{A_0,2} +
                         q_{1 3} w_{A_0,3}\\
\beta w_{A_{\beta},2} &= q_{2 1} w_{A_0,1} + q_{2 2} w_{A_0,2} +
                         q_{2 3} w_{A_0,3}\\
\beta w_{A_{\beta},3} &= q_{3 1} w_{A_0,1} + q_{3 2} w_{A_0,2} +
                         q_{3 3} w_{A_0,3},\end{aligned}\]
where $q_{i j} \in \mathbb{Z}$. Define $h(\beta)$ to be
$\hat{\beta} = (q_{1 3},q_{2 3},q_{3 3}) \in \mathbb{Z}^3$. 
 
We have thus defined a
map $h:\mathcal{I}(\mathfrak{O}_K)^{\times}\to \mathbb{Z}^3$. For any
ideal class $A\in \mathcal{C}_1(K)$, the restriction
$h_{|A}:A\to \mathbb{Z}^3$ is a 
$\mathbb{Z}$-linear map whose image is of finite index in $\mathbb{Z}^3$.

We say that $\vec{a} = (a_1,a_2,a_3)\in \mathbb{R}^3$ 
is {\em primitive} if $\gcd(a_1,a_2,a_3)=1$.
Let $\vec{a}, \vec{b}\in \mathbb{R}^3$. By $\vec{a}\times \vec{b}$ we
mean the {\em cross product}
\[\vec{a}\times \vec{b} = (a_2 b_3 - a_3 b_2, a_3 b_1 - a_1 b_3,
a_1 b_2 - a_2 b_1) .\]
Note that, if $\vec{a}$ and $\vec{b}$ are primitive and $n$ is
a non-zero integer, we have $\vec{a}\times \vec{b} \in n \mathbb{Z}^3$
if and only if $\vec{b} \equiv \lambda \vec{a} \mo n$ for some
$\lambda \in (\mathbb{Z}/n)^*$.

By a {\em cube} $C\subset \mathbb{R}^3$ of side $\ell$ we mean a set of the
form $(x,x+\ell\rbrack \times (y,y+\ell\rbrack \times 
(z,z+\ell\rbrack$.

For $\vec{a}\in \mathbb{Z}^2$, let $\mathfrak{A}_{\vec{a}} =
(a_1 \omega_1 + a_2 \omega_2) \mathfrak{d}^{-1} \in I_K$.
Given $\eta, \upsilon > 0$ and a lattice 
$L\subset \mathbb{Z}^2$, 
let 
\[\begin{aligned}
\Psi_{L,\eta,\upsilon}(\vec{a}) &= \lbrack \vec{a} \in 
                             L\cap 
(\lbrack X, (1+\eta) X \rbrack \times \lbrack \upsilon X,
\upsilon (1 + \eta) X\rbrack)\rbrack\\
 \mathcal{A}'_{L,\eta,\upsilon} &= \{\mathfrak{A}_{\vec{a}} : 
 \vec{a} \in L\cap (\lbrack X, (1+\eta) X \rbrack \times
\lbrack \upsilon X, \upsilon (1 + \eta) X \rbrack)\}\\
 \mathcal{A}_{L,\eta,\upsilon} &= \{\mathfrak{A}_{\vec{a}} : 
 \vec{a} \in L\cap (\lbrack X, (1+\eta) X \rbrack \times
\lbrack \upsilon X, \upsilon (1 + \eta) X \rbrack),\;
\gcd(a_1,a_2)=1 \}.\end{aligned}\]

 Let $\mathcal{Q}\in I_K$
be the set of all ideals in $I_K$ that are not divisible by any rational prime.
In the following, we use $\alpha$, $\beta$ to denote ideal numbers and
$\mathfrak{a}$, $\mathfrak{b}$ to denote ideals.
\begin{lem}\label{eq:agrabo}
Let $K/\mathbb{Q}$ be a number field of degree 3. Let $\omega_1,
\omega_2 \in \mathfrak{O}_K$
be linearly independent over $\mathbb{Z}$.
Let $f, g:I_K\to \mathbb{R}$ be given with 
\begin{equation}\label{eq:tribou}
|f(\mathfrak{a})|, |g(\mathfrak{a})| \ll \tau^{\kappa}(\mathfrak{a}) .
\end{equation} Assume that, for any $B_1, B_2 > 0$,
\begin{equation}\label{eq:rebilicon}
\mathop{\sum_{\vec{b}\in C}}_{\vec{b}\in L \cap h(A)}
g(\mathcal{I}(h_A^{-1}(\vec{b}))) \ll_{B_1,B_2} \vol(C)
(\log X)^{-B_2}
 \end{equation}
for any class $A\in \mathcal{C}_1(K)$, any cube 
$C\subset \lbrack X, 2 X\rbrack^3$ of side $\ell\geq X (\log X)^{-B_1}$,
and any lattice coset $L$ of index $\lbrack \mathbb{Z}^2 : L\rbrack
\leq (\log X)^{B_1}$. Let
\begin{equation}\label{eq:wiggl}\begin{aligned} 
-\log \log N \ll \log \upsilon &\ll \log \log N ,\\
\log \eta &\gg -\log \log N,\\
\eta/\min(\upsilon,\upsilon^{-1}) = o(1) ,\end{aligned}\end{equation}
 Then, for any $B>0$,
\begin{equation}\label{eq:resultat}
\mathop{\mathop{\sum_{\mathfrak{a} \mathfrak{b} \in \mathcal{A}'_{L,\eta,\upsilon}}}_{
X (\log X)^T < N \mathfrak{b} \leq X^{3/2} (\log X)^{-T}}}_{\mathfrak{a}, \mathfrak{b} \in \mathcal{Q}}
f(\mathfrak{a}) g(\mathfrak{b}) \ll X^2 (\log X)^{-B},\end{equation}
where the constant $T$ and the implied constant in (\ref{eq:resultat})
depend only on $\kappa$, $B$ and the implied constants in 
(\ref{eq:tribou})--(\ref{eq:wiggl}).
\end{lem}
\begin{proof}
The argument is nearly the same as that in \cite{HBM}, pp. 278--283.
Let $X (\log X)^T < V < X^{3/2} (\log X)^{-T}$.
Define \begin{equation}\label{eq:aghst}S_0 =
\mathop{\mathop{\sum_{\mathfrak{a} \mathfrak{b} \in \mathcal{A}'_{
L,\eta,\upsilon,\omega_i}}}_{
V < N \mathfrak{b} \leq 2 V}}_{\mathfrak{a}, \mathfrak{b} \in \mathcal{Q}}
f(\mathfrak{a}) g(\mathfrak{b}) 
 .\end{equation}
(Notice that $S_0$ is not the same as $\sum_9(V)$ in \cite{HBM}, (6.2);
instead, what we have is the first summand on the right hand of 
\cite{HBM}, (6.2). We are avoiding the argument at the beginning of 
\S 11 in \cite{HB}, as it implicitly uses a lacunarity condition that
we do not demand.)
We can rewrite (\ref{eq:aghst}) as
\[S_0 = 
\mathop{\sum_{\phi(\vec{a}) = \delta \alpha \beta,\, 
 \mathcal{I}(\alpha)\in \mathcal{Q}}}_{\vec{a}\in \mathbb{Z}^2,\,
V<N \beta\leq 2 V}
f(\mathcal{I}(\alpha)) G_0(\beta) \Psi_{L,\eta}(\vec{a}),\]
where $\phi(\vec{a}) = a_1 \omega_1 + a_2 \omega_2$,
\[G(\beta) = \begin{cases} g(\mathcal{I}(\beta))
&\text{if $\mathcal{I}(\beta) \in \mathcal{Q}_0$}\\
0 &\text{otherwise,}\end{cases}\;\;\;\;\;\;\;\;
G_0(\beta) = \begin{cases} G(\beta) &\text{if $\beta\in \mathcal{Q}_0$}\\
0 &\text{otherwise,}\end{cases}\]
and $\mathcal{Q}_0$ is defined as in \cite{HBM}, p. 278. (In short, 
$\mathcal{Q_0}$
is the set of all ideal numbers $\beta$ satisfying $\mathcal{I}(\beta)\in
\mathcal{Q}$ and a geometrical condition necessary to exclude multiplication
by units.) In the following we will use $\kappa$ to mean a constant
depending only on the value of $\kappa$ in the statement and
 the implied constants in (\ref{eq:tribou})--(\ref{eq:wiggl}).
We now apply Cauchy's inequality:
\begin{equation}\label{eq:metafun}\begin{aligned}
S_0^2 &\ll
\mathop{\sum_{\alpha}}_{\mathcal{I}(\alpha)\in \mathcal{Q}} 
\left| 
\mathop{\sum_{\phi(\vec{a}) = \delta \alpha \beta}}_{\vec{a}\in \mathbb{Z}^2,\,
V<N \beta\leq 2 V}
G_0(\beta) \Psi_{L,\eta}(\vec{a})\right|^2
\cdot \mathop{\sum_{\mathfrak{a}}}_{N \mathfrak{a} \ll X^3/V} 
 |f(\mathfrak{a})|^2 \\ &\ll X^3 V^{-1} (\log X)^{\kappa}
\mathop{\sum_{\alpha}}_{\mathcal{I}(\alpha)\in \mathcal{Q}} 
\left| 
\mathop{\sum_{\phi(\vec{a}) = \delta \alpha \beta}}_{\vec{a}\in \mathbb{Z}^2,\,
V<N \beta\leq 2 V}
G_0(\beta) \Psi_{L,\eta}(\vec{a})\right|^2 .\end{aligned}\end{equation}
As in \cite{HBM}, p. 279, we expand 
(\ref{eq:metafun}) and remove the diagonal terms:
\[S_0 \ll (X^3 V^{-1} (\log X)^{\kappa} \cdot
(S_1 + O(X^2 (\log X)^{\kappa})))^{1/2} ,\]
where
\[S_1 = \mathop{\sum_{\beta_1\ne \beta_2,\, \vec{a}_i\in \mathbb{Z}^2}}_{
 V<N \beta_i\leq 2 V, i=1,2} G_0(\beta_1) G_0(\beta_2) \Psi_{L,\eta}(\vec{a}_1)
\Psi_{L,\eta}(\vec{a}_2) \psi(\vec{a}_1,\vec{a}_2,\beta_1,\beta_2)\]
with \[\psi(\vec{a}_1,\vec{a}_2,\beta_1,\beta_2) = \# \{\alpha :
 \mathcal{I}(\alpha) \in \mathcal{Q},\;
 \phi(\vec{a}_i) = \delta \alpha \beta_i\; \text{for $i=1,2$}\} .\]
As in \cite{HBM}, Lemma 6.2, we remove a small area and obtain
\[S_0 \ll X^2 Y^{-1/2} (\log X)^{\kappa} + X^{3/2} V^{-1/2} S_2^{1/2}\]
with
\[S_2 = \mathop{\mathop{\sum_{\vec{a_i}\in \mathbb{Z}^2,\, \beta_i\in A}}_{
V < N \beta_i \leq 2 V}}_{d(h(\beta_1) \times 
h(\beta_2))> V X^{-1} Y^{-1}}
 G_0(\beta_1) G_0(\beta_2) \Psi_{L,\eta}(\vec{a}_1) \Psi_{L,\eta}(\vec{a}_2)
 \psi(\vec{a}_1,\vec{a}_2,\beta_1,\beta_2),\]
where $A$ is a class of ideal numbers, $Y$ is a parameter
between $1$ and $(\log X)^{T/3}$ chosen at our pleasure,
 and $d((c_1,c_2,c_3)) = \gcd(c_1,c_2,c_3)$.
(Here we have implicitly used Lemma 6.1 of \cite{HBM}.)

We can now proceed as in \cite{HBM}, pp. 280--282, and obtain the following 
analogue of \cite{HBM}, (6.9):
\[S_0\ll X^2 Y^{-1/2} (\log X)^{\kappa} 
+ X^{3/2} V^{-1/2} Y^7 S_3^{1/2} (\log X)^{\kappa},\]
with
\[S_3 = 
\sum_{d_1\in I} \left|\mathop{\mathop{\sum_{\beta_i \in B}}_{
\hat{\beta}_i \in C_i\cap L_{d_1,i}}}_{d(\hat{\beta}_1\times
\hat{\beta}_2)=d}
G(\beta_1) G(\beta_2)\right|,\]
where $A\in \mathcal{C}_1(K)$ is a class of ideal numbers, $I$
is an interval contained in $\lbrack V X^{-1},\infty\rbrack$,
the lattices $L_{d,i}$ have indices 
$\lbrack \mathbb{Z}^3 : L_{d,i}\rbrack | \lbrack \mathbb{Z}^3 : L\rbrack^3$, and
$C_1, C_2\subset \lbrack V X^{-1}, 2 V X^{-1}\rbrack^3$ are
cubes of side about $V X^{-1} (\log X)^{-2 T/3}$. As in
\cite{HBM}, (6.0)--(6.12), we can conclude that
\[S_0\ll X^2 Y^{-1/2} (\log X)^{\kappa} + X^{3/2} V^{-1/2} Y^7 S_4^{1/2}
(\log X)^{\kappa} ,\]
where
\[S_4 = \mathop{\sum_{d_1\in I}}_{d_1 d < d_0} \left|
\mathop{\mathop{\sum_{\beta_i \in B}}_{
\hat{\beta}_i \in C_i\cap L_{d_1,i}}}_{d_1 d | \hat{\beta}_1\times
\hat{\beta}_2}
G(\beta_1) G(\beta_2)\right|\]
with $d_0 = X^{-1} V Y^{15} + V^{1/6}$. We can bound $S_4$ by means
of a large-sieve argument as in \cite{HB}, p. 78--83, and
\cite{HBM}, p. 283; the contribution from small moduli is estimated by
(\ref{eq:rebilicon}). We obtain
\[\begin{aligned}
S_4 \ll &X V \lbrack \mathbb{Z}^2 : L\rbrack^{\kappa} (\log X)^{\kappa}
\\ &\cdot (Y^{\kappa} (\log X)^{-T/2} + Y (\log X)^{-B_1/2} + 
Y (\log X)^{4 B_1} (\log X)^{-B_2}),\end{aligned}\]
where $B_1$ and $B_2$ are arbitrarily large. (See \cite{HBM2}
for an optimization of the exponent
$\kappa$ in $\lbrack \mathbb{Z}^2 : L\rbrack^{\kappa}$.) Set 
$Y = (\log X)^{2 B + 2 \kappa + 2}$, 
$T = 1000 \kappa^2 (B+\kappa+1)$ (say),
$B_1 = T$, $B_2 = 9 B_1$. Then
\[S_0\ll X^2 (\log X)^{-(B+1)} .\]
The statement follows immediately.
\begin{cor}
Let $K/\mathbb{Q}$ be a number field of degree 3. Let $\omega_1,
\omega_2 \in \mathfrak{O}_K$
be linearly independent over $\mathbb{Z}$.
Let $\eta, \upsilon \in \mathbb{R}^+$,
$f, g:I_K\to \mathbb{R}$ satisfy conditions
(\ref{eq:tribou})--(\ref{eq:wiggl}).
Assume furthermore that
\[\mathfrak{p}|\mathfrak{a}, \mathfrak{q}|\mathfrak{b},
\mathfrak{p},\mathfrak{q} \leq (\log x)^{10} \Rightarrow
f(\mathfrak{a}) g(\mathfrak{b}) = 0.\]
Then
\begin{equation}
\mathop{\mathop{\sum_{\mathfrak{a} \mathfrak{b} \in \mathcal{A}_{L,\eta,\upsilon}}}_{
X (\log X)^T < N \mathfrak{b} \leq X^{3/2} (\log X)^{-T}}}_{\mathfrak{a}, \mathfrak{b} \in \mathcal{Q}}
f(\mathfrak{a}) g(\mathfrak{b}) \ll X^2 (\log X)^{-B},\end{equation}
where the constant $T$ and the implied constant in (\ref{eq:resultat})
depend only on $\kappa$, $B$ and the implied constants in 
(\ref{eq:tribou})--(\ref{eq:wiggl}).
\end{cor}
By Lemma \ref{eq:agrabo} and \cite{HB}, p 67. We are simply removing the coprimality
condition on $\mathfrak{a}$ and $\mathfrak{b}$, given that $\mathfrak{a}$ and $\mathfrak{b}$ are
still kept from having small common factors.
\end{proof}
\section{Final remarks and conclusions}
In section \ref{sec:irrcu}, we used the small-boxes formalism of \cite{HB} and \cite{HBM} rather
than our own convex-subset formalism. It is easy to see that boxes such as $S$ in (\ref{eq:defbox})
satisfying (\ref{eq:yadnu}) can cover convex sets with an error of at most $x (\log x)^{-B}$, where
$B$ is arbitrarily large. 

We saw it fit to work with $\lambda$ in sections \ref{sec:three} and \ref{sec:threetwo},
and with $\mu$ in section \ref{sec:irrcu}. (The first choice was due to complete multiplicativity,
the second one to symmetry.) Thanks to Propositions \ref{prop:hjohjo} and
\ref{prop:diezca} for
$\deg P = 3$, a result on $\lambda$ implies one for $\mu$, and vice versa, without any degradation
in our bounds. 
Notice, lastly, that the condition $\gcd(x,y)=1$ implicit in section \ref{sec:irrcu} (see
$\mathcal{A}_{L,S,\omega_i}$ in (\ref{eq:defbox})) can be removed as in Lemma \ref{lem:minmat}.

We collect all our results on cubic polynomials in the following statement.

\begin{thm}
Let $f(x,y)\in \mathbb{Z}\lbrack x,y\rbrack$ be a homogeneous polynomial of degree $3$. Let
$\alpha$ be the M\"obius function ($\alpha = \mu$) or the Liouville function $(\alpha = \lambda)$.
Let $S$ be a convex subset of $\lbrack - N, N\rbrack^2$. 
Let $L\subset \mathbb{Z}^2$ be a lattice coset of
index $\lbrack \mathbb{Z}^2 : L\rbrack \leq (\log N)^{A}$, where $A$
is an arbitrarily high constant. Then
\[\sum_{(x,y)\in S\cap L} \alpha(f(x,y)) \ll
\begin{cases} \frac{(\log \log N)^5 (\log \log \log N)}{\log N} \frac{\Area(S)}{
\lbrack \mathbb{Z}^2 : L\rbrack} + \frac{N^2}{(\log N)^A}
&\text{ if $f$ is irreducible,}\\
\frac{\log \log N}{\log N} \frac{\Area(S)}{\lbrack \mathbb{Z}^2 : L\rbrack}
+ \frac{N^2}{(\log N)^A} 
&\text{ if $f$ is reducible,}\end{cases}\]
where the implied constant depends only on $f$ and on $A$.
\end{thm}

%% file: square.tex
\chapter{The square-free sieve}\label{chap:square}
\epigraph{They sought it with thimbles, they sought it with care;\\
     They pursued it with forks and hope;\\
They threatened its life with a railway--share;\\
     They charmed it with smiles and soap.}{Lewis Carroll,
{\em The Hunting of the Snark}}

A \emph{square-free
sieve} is a result that gives an upper bound for how often a square-free 
polynomial may adopt values that are not square-free.
More generally, we may wish to approximate the cardinality of the set of
arguments $x_1,\dotsc,x_n$ for which the largest square divisor of the value acquired
by $P(x_1,\dotsc,x_n)$ equals a given $\mathfrak{d}$, or, as in 
Chapter \ref{chap:ell}, we may
wish to control the behavior of a function depending on $\sq(P(x_1,\dotsc,x_n))$.

We may aim at obtaining an asymptotic expression
\begin{equation}\label{eq:bartok}
\text{main term } + O(\text{error term}),
\end{equation}
where the main term will depend on the application; in general, the
error term will depend only on the polynomial $P$ in question,
not on the particular quantity being estimated. We can split the
error term further into one term that can be bounded easily for
any $P$, and a second term, say, $\delta(P)$, which may be
rather hard to estimate, and which is unknown for polynomials $P$
of high enough degree. Given this framework, the strongest
results in the literature may be summarized as follows:

\begin{center} \begin{tabular}{l|ll}
$\deg_{\irr}(P)$ &$\delta(P(x))$ &$\delta(P(x,y))$\\ \hline
$1$ & $\sqrt{N}$ & $1$\\
$2$ & $N^{2/3}$ & $N$\\
$3$ & $N / (\log N)^{1/2}$ & $N^2 / \log N$ \\
$4$ & & $N^2/\log N$ \\
$5$ & & $N^2/\log N$ \\
$6$ & & $N^2/(\log N)^{1/2}$
\end{tabular} \end{center}

Here $\deg_{\irr}(P)$ denotes the degree of the largest irreducible
factor of $P$. The second column gives $\delta(P)$ for polynomials
$P\in \mathbb{Z}\lbrack x\rbrack$ of given $\deg_{\irr}(P)$,
whereas the third column refers to homogeneous polynomials 
$P\in \mathbb{Z}\lbrack x,y\rbrack$. The trivial estimates would
be $\delta(P(x))\leq N$ and $\delta(P(x,y))\leq N^2$. See Appendix
\ref{sec:fapp} for attributions.

Our task can be divided into two halves. The first one, undertaken in 
section \ref{sec:cuacriba}, consists in estimating all terms but
$\delta(N)$. We do as much in full generality for any $P$, over
any number field, for that matter. The second half regards bounding
$\delta(N)$. We improve on all estimates known for $3\leq \deg P\leq 5$:

\begin{center}\begin{tabular}{l|ll}
$\deg_{\irr}(P)$ &$\delta(P(x))$ &$\delta(P(x,y))$\\ \hline
$3$ & $N/(\log N)^{0.5718\dotsb}$ & $N^{3/2}/\log N$\\
$4$ & & $N^{4/3} (\log N)^A$ \\
$5$ & & $N^{(5 + \sqrt{113})/8+\epsilon}$\\
\end{tabular}\end{center}

Most of our improvements hinge on a change from a local to a global
perspective. Such previous work in the field as was purely
sieve-based can be seen as an series of purely local estimates
on the density of points on curves of non-zero genus. Our techniques
involve a mixture of sieves, elliptic curves, 
sphere packings, and
some of the methods described in the epigraph.
\section{Notation}
Let $n$ be a non-zero integer.
We write $\tau (n)$ for the number of positive 
divisors of $n$, $\omega (n)$ for the
number of the prime divisors of $n$, and $\rad(n)$ for the product of the prime divisors
of $n$. For any $k\geq 2$, we write $\tau_k(n)$ for the number of 
$k$-tuples $(n_1,n_2,\dotsc, n_k)\in (\mathbb{Z}^+)^k$ such that
$n_1\cdot n_2\cdot \dotsb n_k = |n|$. Thus $\tau_2(n) = \tau(n)$.
We adopt the convention that $\tau_1(n) = 1$. 
We let
\[\sq(n) = \prod_{p^2|n} p^{v_p(n)-1} .\]
We call a rational integer $n$ {\em square-full} if $p^2|n$ for every prime $p$
dividing $n$. Given any non-zero rational integer $D$, we say that $n$ is
{\em ($D\!$)-square-full} if $p^2|n$ for every prime $p$ that divides $n$ but not $D$.

We denote by $\mathfrak{O}_K$ the ring of integers of a global or local
field $K$. We let $I_K$ be the semigroup of non-zero ideals of $\mathfrak{O}_K$.
Given a non-zero ideal $\mathfrak{a}\in I_K$, 
we write $\tau_K(\mathfrak{a})$ for the number of ideals dividing
 $\mathfrak{a}$,
$\omega_K(\mathfrak{a})$ for the number of prime ideals dividing 
$\mathfrak{a}$,
and $\rad_K(\mathfrak{a})$ for the product of the prime ideals dividing
$\mathfrak{a}$.
Given a positive integer $k$, we write $\tau_{K,k}(\mathfrak{a})$ for the number of 
$k$-tuples $(\mathfrak{a}_1,\mathfrak{a}_2,\dotsc ,\mathfrak{a}_k)$ of
ideals of $\mathfrak{O}_K$ such that
$\mathfrak{a} = \mathfrak{a}_1 \mathfrak{a}_2 \dotsb \mathfrak{a}_k$. 
Thus $\tau_2(\mathfrak{a}) = \tau(\mathfrak{a})$.
We let
\[\begin{aligned}
\sq_K(\mathfrak{a}) &= \begin{cases}
\prod_{\mathfrak{p}^2 | \mathfrak{a}} 
\mathfrak{p}^{v_\mathfrak{p}(\mathfrak{a}) - 1} 
&\text{if $\mathfrak{a}\ne 0$,}\\0 &\text{if $\mathfrak{a}=0$,}\end{cases}\\
\mu_K(\mathfrak{a}) &= \begin{cases}
\prod_{\mathfrak{p}|\mathfrak{a}} (-1) &\text{if $\sq_K(\mathfrak{a})=1$,}\\
0 &\text{otherwise.}\end{cases}\end{aligned}\]
We define $\rho(\mathfrak{a})$ to be the positive integer generating $\mathfrak{a}\cap \mathbb{Z}$.

When we say that a polynomial $f\in \mathfrak{O}_K\lbrack x\rbrack$ or
$f\in K\lbrack x\rbrack$ is {\em square-free}, we always mean that
$f$ is square-free as an element of $K\lbrack x\rbrack$. In other
words, we say that $f\in \mathbb{Z}\lbrack x\rbrack$ is {\em square-free}
if there is no polynomial $g\in \mathbb{Z}\lbrack x\rbrack$
such that $\deg g \geq 1$ and $g|f$. See section \ref{sec:notpre}
for the definitions of the resultant $\Res$ and the discriminant $\Disc$.

Given an elliptic curve $E$ over $\mathbb{Q}$, we write $E(\mathbb{Q})$
for the set of rational (that is, $\mathbb{Q}$-valued) points of $E$.
We denote by $\rnk(E)$ the algebraic rank of $E(\mathbb{Q})$.
\section{Sieving}\label{sec:cuacriba}
\subsection{An abstract square-free sieve}\label{subs:abstsi}
\begin{lem}\label{lem:wicked}
Let $K$ be a number field. 
Let $\{S_\mathfrak{a}\}_{\mathfrak{a}\in I_K}$ be a collection of finite
sets, one for each non-zero ideal $\mathfrak{a}$ of $\mathfrak{O}_K$.
Let a map 
$\phi_{\mathfrak{a}_1,\mathfrak{a}_2}:S_{\mathfrak{a}_2}\to 
S_{\mathfrak{a}_1}$ be given for any non-zero ideals 
$\mathfrak{a}_1$, $\mathfrak{a}_2$ such that $\mathfrak{a}_1|\mathfrak{a}_2$. 
Assume $\phi_{\mathfrak{a}_1,\mathfrak{a}_2}\circ
\phi_{\mathfrak{a}_2,\mathfrak{a}_3} = \phi_{\mathfrak{a}_1,\mathfrak{a}_3}$ 
for all $\mathfrak{a}_1$, $\mathfrak{a}_2$, $\mathfrak{a}_3$ 
such that $\mathfrak{a}_1|\mathfrak{a}_2|\mathfrak{a}_3$.
Let $\{f_\mathfrak{a}\}_{\mathfrak{a}\in I_K}$, 
$f_\mathfrak{a}:S_\mathfrak{a}\to \mathbb{C}$ be given
with $|f_\mathfrak{a}(r)|\leq 1$ for all $\mathfrak{a}\in I_K$ and all 
$r\in S_\mathfrak{a}$. Let
$\{g_\mathfrak{a}\}_{\mathfrak{a}\in I_K}$, 
$g_\mathfrak{a}:S_\mathfrak{a}\to \mathbb{C}$ be such that
\[\sum_\mathfrak{a\in I_K} \sum_{r\in S_\mathfrak{a}} |g_\mathfrak{a}(r)|\] converges.
Write 
\[\begin{aligned}
s_\mathfrak{d} &= \mathop{\sum_\mathfrak{a\in I_K}}_{\mathfrak{d}|\mathfrak{a}} 
\sum_{r\in S_\mathfrak{a}} |g_\mathfrak{a}(r)| ,\\
t_\mathfrak{d}(r) &= \mathop{\sum_\mathfrak{a}}_{\mathfrak{d}|\mathfrak{a}} 
 \mathop{\sum_{r'\in S_\mathfrak{a}}}_{
\phi_{\mathfrak{d},\mathfrak{a}}(r')=r} g_\mathfrak{a}(r').
\end{aligned}\]
Let $\gamma:I_K\to \mathbb{Z}^+$ be a map such that 
$\gamma(\mathfrak{d}_1)\leq \gamma(\mathfrak{d}_1 \mathfrak{d}_2)\leq
\gamma(\mathfrak{d}_1) \gamma(\mathfrak{d}_2)$ for all $\mathfrak{d}_1,
\mathfrak{d}_2\in I_K$.
Then, for any positive integer $M$,
\begin{equation}\label{eq:maravel}\begin{aligned}
\sum_{\mathfrak{a}\in I_K} \sum_{r\in S_\mathfrak{a}} f_\mathfrak{a}(r) 
g_\mathfrak{a}(r) &\leq
\sum_{\gamma(\mathfrak{d})\leq M} \sum_{r\in S_\mathfrak{d}}
 \left(\sum_{\mathfrak{d}'|\mathfrak{d}} \mu_K(\mathfrak{d}') 
f_{\mathfrak{d}/\mathfrak{d}'}(
\phi_{\mathfrak{d}/\mathfrak{d}',\mathfrak{d}}(r))\right) t_\mathfrak{d}(r) 
\\&+ 
2 \mathop{\sum_{\mathfrak{d}\in I_K}}_{M<\gamma(\mathfrak{d})\leq M^2} 
\tau_{K,3}(\mathfrak{d}) s_\mathfrak{d}+ 2 
\mathop{\sum_{\text{$\mathfrak{p}$ prime}}}_{\gamma(\mathfrak{p})>M}
 s_\mathfrak{p} .\end{aligned}\end{equation}
\end{lem}
\begin{proof}
Let $\sigma(\mathfrak{a}) = 
\prod_{\mathfrak{p}|\mathfrak{a},\,\gamma(\mathfrak{p})\leq M}
\mathfrak{p}^{v_\mathfrak{p}(a)}$.
By M\"obius inversion, for any $r\in S_\mathfrak{a}$,
\[\begin{aligned}
\sum_{\mathfrak{d}|\mathfrak{a}} 
\sum_{\mathfrak{d}'|\mathfrak{d}} 
 \mu_K(\mathfrak{d}') f_{\mathfrak{d}/\mathfrak{d}'}(
\phi_{\mathfrak{d}/\mathfrak{d}',\mathfrak{a}}(r)) &= f_\mathfrak{a}(r),\\
\mathop{\sum_{\mathfrak{d}|\mathfrak{a}}}_{\mathfrak{p}|\mathfrak{d}\Rightarrow \gamma(\mathfrak{p})\leq M} 
\sum_{\mathfrak{d}'|\mathfrak{d}} 
 \mu_K(\mathfrak{d}') f_{\mathfrak{d}/\mathfrak{d}'}(
\phi_{\mathfrak{d}/\mathfrak{d}',\mathfrak{a}}(r)) &= 
f_{\sigma(\mathfrak{a})}(\phi_{\sigma(\mathfrak{a}),
\mathfrak{a}}(r)) .\end{aligned}\]
Hence
\[\begin{aligned}
\sum_\mathfrak{a} \sum_{r\in S_\mathfrak{a}} 
 f_\mathfrak{a}(r) g_\mathfrak{a}(r) &=
\sum_r \sum_{r\in S_\mathfrak{a}} (f_\mathfrak{a}(r) - 
f_{\sigma(\mathfrak{a})}(
\phi_{\sigma(\mathfrak{a}),\mathfrak{a}}(r)) )
g_\mathfrak{a}(r) \\
&+ \sum_\mathfrak{a} \sum_{r\in S_\mathfrak{a}} 
w_{\mathfrak{a},r} g_\mathfrak{a}(r)
\\ &+\sum_{\gamma(\mathfrak{d})\leq M} 
\left(\sum_{r\in S_\mathfrak{d}}
\sum_{\mathfrak{d}'|\mathfrak{d}} \mu_K(\mathfrak{d}') 
f_{\mathfrak{d}/\mathfrak{d}'}(
\phi_{\mathfrak{d}/\mathfrak{d}',\mathfrak{a}}(r))\right) 
t_\mathfrak{d}(r),\end{aligned}\]
where we write
\[w_{\mathfrak{a},r} = \mathop{\sum_{
\mathfrak{d}|\mathfrak{a}}}_{\mathfrak{p}|\mathfrak{d}\Rightarrow 
\gamma(\mathfrak{p})\leq M}
\sum_{\mathfrak{d}'|\mathfrak{d}} \mu_K(\mathfrak{d}') 
f_{\mathfrak{d}/\mathfrak{d}'}(\phi_{\mathfrak{d}/\mathfrak{d}',
\mathfrak{a}}(r)) \:-
\mathop{\sum_{\mathfrak{d}|\mathfrak{a}}}_{
\gamma(\mathfrak{d})\leq M} \sum_{\mathfrak{d}'|\mathfrak{d}} 
\mu_K(\mathfrak{d}') 
f_{\mathfrak{d}/\mathfrak{d}'}(\phi_{\mathfrak{d}/\mathfrak{d}',
\mathfrak{a}}(r)).\]

Since $\mathfrak{a}=\sigma(\mathfrak{a})$ 
unless $\mathfrak{a}$ is divisible by a prime $\mathfrak{p}$
with $\gamma(\mathfrak{p})>M$, we know that
\[\sum_r \sum_{r\in S_\mathfrak{a}} (f_\mathfrak{a}(r) - 
f_{\sigma(\mathfrak{a})}(
\phi_{\sigma(\mathfrak{a}),\mathfrak{a}}(r)) ) 
g_\mathfrak{a}(r) \leq
\mathop{\sum_{\text{$\mathfrak{p}$ prime}}}_{\gamma(\mathfrak{p})>M}
 s_\mathfrak{p} .\]
Now take $\mathfrak{a}$, $r$ such that
\begin{equation}\label{eq:gesund}\mathop{\sum_{
\mathfrak{d}|\mathfrak{a}}}_{\mathfrak{p}|\mathfrak{d}\Rightarrow 
\gamma(\mathfrak{p})\leq M}
\sum_{\mathfrak{d}'|\mathfrak{d}} \mu_K(\mathfrak{d}') 
f_{\mathfrak{d}/\mathfrak{d}'}(\phi_{\mathfrak{d}/\mathfrak{d}',
\mathfrak{a}}(r)) \,\ne
\mathop{\sum_{\mathfrak{d}|\mathfrak{a}}}_{
\gamma(\mathfrak{d})\leq M} \sum_{\mathfrak{d}'|\mathfrak{d}} 
\mu_K(\mathfrak{d}') 
f_{\mathfrak{d}/\mathfrak{d}'}(\phi_{\mathfrak{d}/\mathfrak{d}',
\mathfrak{a}}(r)).\end{equation}
This can happen only if $\gamma(\sigma(\mathfrak{a}))>M$.
Let $\mathfrak{d}$ be a divisor of $\mathfrak{a}$ with
$\gamma(\mathfrak{d})\leq M$. We would like to show that there
is a divisor $\mathfrak{d}'$ of $\mathfrak{a}$ such that
$\mathfrak{d}|\mathfrak{d}'$ and $M<\gamma(\mathfrak{d}')\leq M^2$.
Since $\gamma(\mathfrak{d})\leq M$, all prime divisors $\mathfrak{p}$
of $\mathfrak{d}$ obey $\gamma(\mathfrak{p})\leq M$, and thus
$\mathfrak{d} | \sigma(\mathfrak{a})$.
Write $\sigma(\mathfrak{a}) = \mathfrak{d} 
\mathfrak{p}_1 \dotsb 
\mathfrak{p}_k$, where $\mathfrak{p_1},\dotsc , \mathfrak{p_k}$
are not necessarily distinct. 
Let $\mathfrak{a}_0 = \mathfrak{d}$.
For $1\leq i\leq k$, let 
$\mathfrak{a}_i = \mathfrak{d} \mathfrak{p}_1 \dotsb \mathfrak{p}_i$. 
Then
$\gamma(\mathfrak{a}_0) \leq M$, $\gamma(\mathfrak{a}_k) = \gamma(\sigma(\mathfrak{a})) > M$ and
$\gamma(\mathfrak{a}_{i+1})\leq \gamma(\mathfrak{a}_i) \gamma(\mathfrak{p}_i)
\leq \gamma(\mathfrak{a}_i) \cdot M$ for every $1\leq i<k$. Hence
there is an $0\leq i\leq k$ such that $M<\gamma(\mathfrak{a}_i)\leq M^2$.
Since $\mathfrak{d} | \sigma(\mathfrak{a})_i$ and
$\mathfrak{a}_i | \sigma(\mathfrak{a})$, we can set $\mathfrak{d}' =
\mathfrak{a}_i$.

Now bound the right hand side of (\ref{eq:gesund}) trivially:
\[\mathop{\sum_{\mathfrak{d}|\mathfrak{a}}}_{
\gamma(\mathfrak{d})\leq M} \sum_{\mathfrak{d}'|\mathfrak{d}} 
\mu_K(\mathfrak{d}') 
f_{\mathfrak{d}/\mathfrak{d}'}(\phi_{\mathfrak{d}/\mathfrak{d}',
\mathfrak{a}}(r))\leq
\mathop{\sum_{\mathfrak{d}|\mathfrak{a}}}_{\gamma(\mathfrak{d})\leq M}
\tau_K(\rad(\mathfrak{d})) .\]
By the foregoing discussion,
\[\mathop{\sum_{\mathfrak{d}|\mathfrak{a}}}_{\gamma(\mathfrak{d})\leq M}
\tau_K(\rad(\mathfrak{d})) \;\leq
\mathop{\sum_{\mathfrak{d}'|\mathfrak{a}}}_{M<\gamma(\mathfrak{d}')\leq M^2} 
 \sum_{\mathfrak{d}|\mathfrak{d}'} \tau_K(\rad(\mathfrak{d})) \;\;= 
\mathop{\sum_{\mathfrak{d}'|\mathfrak{a}}}_{M<\gamma(\mathfrak{d}')\leq M^2}
\tau_{K,3}(\mathfrak{d}')  .\]

Since
\begin{equation*}
\left|\mathop{\sum_{\mathfrak{d}|\mathfrak{a}}}_{\mathfrak{p}|\mathfrak{d}\Rightarrow \gamma(\mathfrak{p})\leq M} \sum_{\mathfrak{d}'|\mathfrak{d}}
\mu_K(\mathfrak{d}') f_{\mathfrak{d}/\mathfrak{d}'}(
\phi_{\mathfrak{d}/\mathfrak{d}',\mathfrak{a}}(r))\right| = 
|f(\sigma(\mathfrak{a}))|\leq 1\end{equation*}
and since for all terms such that
$\gamma(\sigma(\mathfrak{a}))>M$
we have
\[\mathop{\sum_{\mathfrak{d}'|\mathfrak{a}}}_{
M<\gamma(\mathfrak{d}')\leq M^2}
\tau_{K,3}(\mathfrak{d}') \,\geq 1,\] we can conclude that
\[\left|\sum_\mathfrak{a} \sum_{r\in S_\mathfrak{a}} 
w_{\mathfrak{a},r} g_\mathfrak{a}(r)\right|\]
is less than or equal to twice
\[\sum_{\mathfrak{a}} \sum_{r\in S_\mathfrak{a}} 
\mathop{\sum_{\mathfrak{d}|\mathfrak{a}}}_{M<\gamma(\mathfrak{d})
\leq M^2} \tau_{K,3}(\mathfrak{d}) 
\,|g_\mathfrak{a}(r)| .
\]
Since
\[\sum_{\mathfrak{a}} \sum_{r\in S_\mathfrak{a}} 
\mathop{\sum_{\mathfrak{d}|\mathfrak{a}}}_{M<
\gamma(\mathfrak{d})\leq M^2} \tau_{K,3}(\mathfrak{d}) 
\,|g_\mathfrak{a}(r)| \leq
 \sum_{M< \gamma(\mathfrak{d}) \leq M^2} 
\tau_{K,3}(\mathfrak{d}) s_\mathfrak{d},\]
the result follows.
\end{proof}

\subsection{Solutions and lattices}
\begin{lem}\label{lem:sols}
Let $K$ be a $\mathfrak{p}$-adic field.
Let $P\in \mathfrak{O}_K\lbrack x \rbrack$ be a square-free
polynomial. Then
\[P(x)\equiv 0 \mo \mathfrak{p}^n\]
has at most $\max(|\Disc P|_\mathfrak{p}^{-1} \cdot \deg P, 
|\Disc P|_\mathfrak{p}^{-3})$ roots in $\mathfrak{O}_K/\mathfrak{p}^n$.
\end{lem}
\begin{proof}
Let $\pi$ be a prime element of $K$.
If $P$ is of the form $P = \pi Q$ for some $Q\in \mathfrak{O}_K \lbrack x \rbrack$, the statement follows
from the statement for $Q$. Hence we can assume $P$ is not of the form $P = \pi G$. Write
$P = P_1 \cdot P_2 \cdot \dotsb \cdot P_l$, $P_i\in \mathfrak{O}_K$, $P_i$ irreducible.

If $n\leq 3 v_\mathfrak{p}(\Disc P)$, there are trivially at most 
$\# (\mathfrak{O}_K/\mathfrak{p}^n) = |\mathfrak{p}^n|_\mathfrak{p}^{-1} 
\leq |\Disc P|_\mathfrak{p}^{-3}$ roots.
Assume $n> 3 v_{\mathfrak{p}}(\Disc P)$.
 Let $x$ be a root of 
$P(x)\equiv 0 \mo \mathfrak{p}^n$. Let $P_i$ be a factor for which $v_\mathfrak{p}(P_i(x))$ is maximal. By
\[v_\mathfrak{p}(P'(x)) = v_\mathfrak{p}(\sum_j P_j'(x) \cdot P_1(x) \cdot \dotsb \widehat{P_j(x)}\dotsb \cdot P_n(x))
\geq \min_j(v_\mathfrak{p}(P(x)) - v_\mathfrak{p}(P_j(x))),\] 
$\min(v_\mathfrak{p}(P'(x)),v_\mathfrak{p}(P(x)))\leq v_\mathfrak{p}(\Disc P)$
and $v_\mathfrak{p}(P(x))>v_\mathfrak{p}(\Disc P)$, we have that
\[\min_j(v_\mathfrak{p}(P(x)) - v_\mathfrak{p}(P_j(x))) \leq v_\mathfrak{p}(\Disc P)\] and hence
\[v_\mathfrak{p}(P_i(x)) \geq v_\mathfrak{p}(P(x)) - 
v_\mathfrak{p}(\Disc P)\geq n - v_\mathfrak{p}(\Disc P)\geq 2 v_\mathfrak{p}(\Disc P) + 1.\]
On the other hand $\gcd(P_i(x),P_i'(x))|\Disc P$, and thus $v_\mathfrak{p}(P_i'(x))\leq v_\mathfrak{p}(\Disc P)$.
By Hensel's lemma we can conclude that $P_i$ is linear. Since $v_\mathfrak{p}(P_i(x))\geq n - v_\mathfrak{p}(\Disc P)$,
$x$ is a root of \[P_i(x)\equiv 0 \mo \mathfrak{p}^{n-v_\mathfrak{p}(\Disc P)} .\]
Since $P_i$ is linear and not divisible by $\mathfrak{p}$, it has at most one root in 
$\mathfrak{O}_K/\mathfrak{p}^{n-v_\mathfrak{p}(\Disc P)}$. There are at most $v_\mathfrak{p}(\Disc P)$ elements of 
$\mathfrak{O}_K/\mathfrak{p}^n$ reducing to this root. Summing over all $i$ we obtain that there are at
most $l\cdot v_\mathfrak{p}(\Disc P)$ roots of $P(x)\equiv 0 \mo \mathfrak{p}^n$ in $\mathbb{Z}/\mathfrak{p}^n$. Since
$l\leq \deg P$, the statement follows.
\end{proof}

\begin{lem}\label{lem:dontaskmemyadvice}
 Let $K$ be a number field. Let $\mathfrak{m}$ be a non-zero ideal
of $\mathfrak{O}_K$. Let $P\in \mathfrak{O}_K\lbrack x\rbrack$ be a square-free
polynomial. Then 
\[\{x\in \mathbb{Z} : P(x)\equiv 0 \mo \mathfrak{m}\}\]
is the union of at most $|\Disc P|^3\cdot 
\tau_{\deg P}(\rad(\rho(\mathfrak{m})))$
arithmetic progressions of modulus $\rho(\mathfrak{m})$.
\end{lem}
\begin{proof}
  By Lemma \ref{lem:sols}, for every $\mathfrak{p}|\mathfrak{m}$,
the equation
\[P(x)\equiv 0 \mo \mathfrak{p}^n\]
has at most 
$|\Disc P|_\mathfrak{p}^{-3} \deg P$ roots in $\mathfrak{O}_K/\mathfrak{p}^n$.
For any ideal $\mathfrak{a}$,
the intersection of $\mathbb{Z}$ with a set of the form
\[\{x\in \mathfrak{O}_K : x\equiv x_0 \mo \mathfrak{a}\}\]
is either the empty set or an arithmetic progression of modulus 
$\rho(\mathfrak{a})$. This is in particular true for $\mathfrak{a} = \mathfrak{p}^n$; the set 
\[\{x\in \mathbb{Z} : x\equiv x_0 \mo \mathfrak{p}^n\}\]
is the union of at most 
$|\Disc P|_\mathfrak{p}^{-3} \deg P$ arithmetic progressions of modulus
$\rho(\mathfrak{p}^n)$. 

Now consider a rational prime $p$ at least one of whose prime ideal
divisors divides $m$. Write $m = \mathfrak{p}_1^{n_1} \mathfrak{p}_2^{n_2}
\dotsb \mathfrak{p}_k^{n_k} \mathfrak{m}_0$, where $\mathfrak{p}_1,\dotsc
\mathfrak{p}_k | p$, $n_1\geq n_2\geq \dotsb \geq n_k$ and
$\mathfrak{m}_0$ is prime to $p$. The set
\[\{x\in \mathbb{Z} : x\equiv x_0 \mo \mathfrak{p}_1^{n_1} \dotsb
\mathfrak{p}_k^{n_k}\}\]
is the intersection of the sets
\[\{x\in \mathbb{Z} : x\equiv x_0 \mo \mathfrak{p}_j^{n_j}\},
\text{\;\; $1\leq j\leq k$.}
\]
At the same time, it is a disjoint union of arithmetic progressions
of modulus \[\rho(\mathfrak{p}_1^{n_1}\dotsb\mathfrak{p}_k^{n_k}) =
\rho(\mathfrak{p}_1^{n_1}).\] Since
\[\{x\in \mathbb{Z} : x\equiv x_0 \mo \mathfrak{p}_1^{n_1}\}
\] is the disjoint union of at most
$|\Disc P|_\mathfrak{p}^{-3} \deg P$ arithmetic progressions
of modulus $\rho(\mathfrak{p}_1^{n_1})$,
\[\{x\in \mathbb{Z} : x\equiv x_0 \mo \mathfrak{p}_1^{n_1} \dotsb
\mathfrak{p}_k^{n_k}\}\]
is the disjoint union of at most
$|\Disc P|_\mathfrak{p}^{-3} \deg P$ arithmetic progressions
of modulus $\rho(\mathfrak{p}_1^{n_1})$.

By (\ref{eq:inters}) the statement follows.
\end{proof}

\begin{lem}\label{lem:sollat} 
Let $K$ be a number field. 
Let $\mathfrak{m}$ be a non-zero ideal of $\mathfrak{O}_K$.
Let $P\in \mathfrak{O}_K\lbrack x,y\rbrack $ be a non-constant and square-free
 homogeneous
polynomial. Then the set
\[S = \{(x,y)\in \mathbb{Z}^2 : \gcd(x,y) = 1, \mathfrak{m}|P(x,y)\}\]
is the union of at most 
$|\Disc P|^3\cdot \tau_{2\deg P}(\rad(\rho(\mathfrak{m})))$ 
disjoint sets of the form
$L \cap \{(x,y)\in \mathbb{Z}^2 : \gcd(x,y) = 1\},$
$L$ a lattice of index $\lbrack \mathbb{Z}^2 : L\rbrack = \rho(\mathfrak{m})$.
\end{lem}
\begin{proof}
Let $\mathfrak{p}|\mathfrak{m}$. Let $n = v_\mathfrak{p}(\mathfrak{m})$.
Let $r_1,r_2,\dotsb,r_k\in \mathfrak{O}_K/\mathfrak{p}^n$ be the roots of
$P(r,1)\cong 0 \mo \mathfrak{p}^n$. Let $r'_1,r'_2,\dotsb,r'_{k'}\in
\mathfrak{O}_K/\mathfrak{p}^n$ be such roots of
$P(1,r)\cong 0 \mo \mathfrak{p}^n$ as satisfy $\mathfrak{p}|r$.
 Then the set of solutions to
$P(x,y)\cong 0 \mo \mathfrak{p}^n$ in 
\[\{(x,y)\in \mathbb{Z}^2: \mathfrak{p}\nmid \gcd(x,y)\}\] is the union of the disjoint sets
\[\{(x,y)\in \mathbb{Z}^2: \mathfrak{p}\nmid \gcd(x,y), x\equiv r_i y \mo \mathfrak{p}^n\},\]
\[\{(x,y)\in \mathbb{Z}^2: \mathfrak{p}\nmid \gcd(x,y), y\equiv r_i x \mo \mathfrak{p}^n\} .\] 
Each of these sets is either the empty set or a set of the form $L\cap (\mathbb{Z}^2 - p \mathbb{Z}^2)$,
where $p$ is the rational prime lying under $\mathfrak{p}$ and $L$ is a lattice of index $\rho(\mathfrak{p}^n)$. 
By Lemma \ref{lem:sols},
$k+k'\leq 2 |\Disc P|_\mathfrak{p}^{-3} \deg P$.
The rest of the argument is as in Lemma
\ref{lem:dontaskmemyadvice}.\end{proof}

\subsection{Square-full numbers}
\begin{lem}\label{lem:bleich}
Let $K$ be a number field. Let $D$ be the product of all rational primes 
ramifying in $K/\mathbb{Q}$. Then, for every $\mathfrak{d}\in I_K$,
the rational integer $\rho(\mathfrak{d} \rad_K(\mathfrak{d}))$ is ($D\!$)-square-full. For any
integer $n$, there are at most 
$C \cdot \tau_{\deg(K/\mathbb{Q})+1}(n)$
ideals $\mathfrak{d}\in I_K$ such that $\rho(\mathfrak{d} \rad_K(\mathfrak{d})) = n$,
where $C$ is the product \[\prod_p e_p^{\deg(K/\mathbb{Q})/e_p}\]
taken over all primes $p$ ramifying in $K/\mathbb{Q}$.
\end{lem}
\begin{proof}
The first statement is clear. It is enough to verify the second statement
for $n$ of the form $p^m$. Let $e$ be the ramification degree of $p$ over 
$K/\mathbb{Q}$.
Then the ideals 
$\mathfrak{d}$
such that $\mathfrak{d} \rad \mathfrak{d}$ divides $p^m$ are of the form
$\mathfrak{p}_1^{a_1} \mathfrak{p}_2^{a_2} \dotsb \mathfrak{p}_k^{a_k}$,
where $a_1, a_2,\dotsc ,a_k$ are non-negative integers less than $e m$ and
$\mathfrak{p}_1, \mathfrak{p}_2, \dotsc, \mathfrak{p}_k$ are
the primes lying above $p$. There are 
$(e m)^k\leq (e m)^{\deg(K/\mathbb{Q})/e}$
choices for $a_1, a_2,\dotsc ,a_e$. Hence there are at most 
$(e m)^{\deg(K/\mathbb{Q})/e}$ ideals $\mathfrak{d}$ such that 
$\gamma(\mathfrak{d})=n$. Now $m^l \leq \binom{m+l}{l}$ for all positive
$m$ and $l$. Since $\tau_l(p^m) = \binom{m+l-1}{l-1}$ for $l\geq 2$, the statement follows.
\end{proof}
\begin{lem}\label{lem:bleich2}
Let $K$ be a number field. Let $m$ be a positive integer. Let $D$ be the product of all rational primes ramifying in $K/\mathbb{Q}$. Then, for every $\mathfrak{d}\in I_K$,
$\lcm(m,\rho(\mathfrak{d} \rad_K(\mathfrak{d})))$ is
$(D m)$-square-full. For any integer $n$, there are at most 
$C\cdot \tau_{\deg(K/\mathbb{Q})+2}(n)$ ideals $\mathfrak{d}\in I_K$ such that
$\lcm(m,\rho(\mathfrak{d} \rad_K(\mathfrak{d}))) = n$, 
where $C$ is the product \[\prod_p e_p^{\deg(K/\mathbb{Q})/e_p}\]
taken over all primes $p$ ramifying in $K/\mathbb{Q}$.
\end{lem}
\begin{proof}
Immediate from Lemma \ref{lem:bleich}.
\end{proof}
\begin{lem}\label{lem:changet}
Let $K$ be a number field. Let $k$ be a positive integer. For any $\mathfrak{d}\in I_K$,
\[\tau_{K,k}(\rad_K(\mathfrak{d})) \leq \tau_{k^{\deg K/\mathbb{Q}}}(\rad(\rho(\mathfrak{d}))) .\]
\end{lem}
\begin{proof}
Let $n\in \mathbb{Z}$ be square-free. For every $\mathfrak{d}\in I_K$ such that 
$\rad(\rho(\mathfrak{d}))|n$, we have $\mathfrak{d}|\rho(\mathfrak{d})$ and hence
$\rad_K(\mathfrak{d})|n$. Thus it is enough to prove $\tau_{K,k}(n)\leq 
\tau_{k^{\deg K/\mathbb{Q}}}(n)$.
Since there are at most $\deg K/\mathbb{Q}$ prime ideals
in $I_K$ above a given rational prime, $\tau_{K,k}(n)\leq k^{\deg K/\mathbb{Q}} = 
\tau_{k^{\deg K/\mathbb{Q}}}(n)$ for $n$ prime. The general case follows by
multiplicativity. 
\end{proof}
The following two lemmas will be used frequently enough that their repeated mention
would be irksome.
\begin{lem}\label{lem:supmul}
For any positive integers $k$, $n$, $n'$,
\[\tau_k(n n') \leq \tau_k(n) \tau_k(n') .\]
\end{lem}
\begin{proof}
Let $S_k(n)$ be the set of all $k$-tuples of integers $(n_1,n_2,\dotsc,n_k)$ with
product $\prod_j n_j = n$. 
There is a map $f_k$ from $S_k(n)\times S_k(n')$ to $S_k(n n')$:
\[((n_1,\dotsc,n_k),(n_1',\dotsc,n_k'))\mapsto (n_1 n_1',\dotsc,n_k n_k') .\]
We can show that $f_k$ is surjective as follows. Let 
$(n_1'',\dotsc,n_k'')$ be given with $\prod_j n_j'' = n n'$. Define $n_1=\gcd(n,n_1'')$,
$n_2 = \gcd(n/n_1,n_2'')$, $n_3 = \gcd(n/(n_1 n_2),n_3'')$, \dots;
$n_1' = n_1''/n_1$, $n_2' = n_2''/n_2$, $n_3' = n_3''/n_3$, and so on. Then
$f((n_1,\dotsc,n_k),(n_1',\dotsc,n_k'))=(n_1'',\dotsc,n_k'')$. Hence $f_k$ is
surjective. Since $\tau_k(n) = \# S_k(n)$, $\tau_k(n') = \# S_k(n')$,
$\tau_k(n'') = \# S_k(n'')$, the statement follows.
\end{proof}
\begin{lem}
For any positive integers $k_1$, $k_2$, $n$,
\[\tau_{k_1}(n) \tau_{k_2}(n) \leq \tau_{k_1 k_2}(n).\]
\end{lem}
\begin{proof}
Let $S_k(n)$ be as in the proof of Lemma \ref{lem:supmul}.
There is a map $f_{k_1,k_2}$ from $S_{k_1 k_2}(n)$
to $S_{k_1}(n)\times S_{k_2}(n)$:
\[(n_1,\dotsc,n_{k_1 k_2})\mapsto \left(\left(\prod_{j_2} n_{(j_2-1) k_1 + j_1}\right)_{j_1},
\left(\prod_{j_1} n_{(j_2-1) k_1 + j_1}\right)_{j_2}\right) .\]
We can show that $f_{k_1,k_2}$ is surjective as follows. See $n = p_1^{e_1}\dotsb p_k^{e_j}$
as a box of $e_1 + \dotsb + e_j$ primes of different colours. Every 
$(m_1,\dotsc,m_{k_1})\in S_{k_1}(n)$ (resp. $(m_1',\dotsc,m_{k_2}')\in S_{k_2}(n)$) gives
us a partition of the box into $k_1$ sets $M_1,\dotsc,M_{k_1}$ (resp.
$k_2$ sets $M_1',\dotsc,M_{k_2}'$). Let 
$n_{(j_2 -1) k_1 +j_1}$ be the product of the primes in $M_{j_1}\cap M_{j_2}'$.
Then $f(n_1,\dotsc,n_{k_1 k_2}) = 
((m_1,\dotsc,m_{k_1}),(m_1',\dotsc,m_{k_2}'))$. Hence $f_{k_1,k_2}$ is surjective.
Since $\tau_{k_1}(n) = \# S_{k_1}(n)$, $\tau_{k_2}(n) = \# S_{k_2}(n)$, $\tau_{k_1 k_2}(n) = 
\# S_{k_1 k_2}(n)$, the statement follows.
\end{proof}
\begin{lem}\label{lem:prebrim}
Let $k$ be a positive integer. Then
\[\mathop{\sum_{n\leq N}}_{\text{$n$ square-full}} \tau_k(n) \leq
 (1 + \log N)^{k^3 + k^2 - 2} N^{1/2}.\]
\end{lem}
\begin{proof}
Every square-full number can be written as a product of a square and a cube. Hence
\[\begin{aligned}
\mathop{\sum_{n\leq N}}_{\text{$n$ square-full}} \tau_k(n) &\leq
 \sum_{n=1}^{\sqrt{N}} \sum_{m=1}^{N^{1/3}/n^{2/3}} \tau_k(n^2 m^3) \leq
 \sum_{n=1}^{\sqrt{N}} \tau_k(n)^2 \sum_{m=1}^{N^{1/3}/n^{2/3}} \tau_k(m)^3 \\ &\leq
 \sum_{n=1}^{\sqrt{N}} \tau_k(n)^2 (1 + \log m)^{k^3-1} (N/n^2)^{1/3} \\ &\leq
 (1 + \log N)^{k^3 - 1} N^{1/3} \sum_{n=1}^{\sqrt{N}} \frac{\tau_k(n)^2}{n^{2/3}}\\
  &\leq (1 + \log N)^{k^3 - 1} N^{1/3} (1 + \log \sqrt{N})^{k^2 -1} (\sqrt{N})^{1/3} \\ &\leq
 (1 + \log N)^{k^3 + k^2 - 2} N^{1/2} .\end{aligned}\]
\end{proof}
\begin{lem}\label{lem:oohoo}
Let $k$ be a positive integer. Then
$\sum_{\text{$n$ square-full}} \frac{\tau_k(n)}{n}$
converges.
\end{lem}
\begin{proof}
\[\mathop{\sum_{n=1}^{\infty}}_{\text{$n$ square-full}} \frac{\tau_k(n)}{n} \leq
\sum_{n=1}^{\infty} \sum_{m=1}^{\infty} \frac{\tau_k(n^2 m^3)}{n^2 m^3} \leq
\left(\sum_{n=1}^{\infty} \frac{\tau_k(n)^2}{n^2}\right)
 \left(\sum_{m=1}^{\infty} \frac{\tau_k(m)^3}{m^3}\right).\]
\end{proof}
\begin{lem}\label{lem:compnr}
Let $k$ be a positive integer. Then
\[\mathop{\sum_{n>N}}_{\text{$n$ square-full}} \frac{\tau_k(n)}{n} \ll
\frac{(\log N)^{k^2 + k^3 - 2}}{N^{1/2}},\]
where the implied constant depends only on $k$.
\end{lem}
\begin{proof}
Since $\sum_{n>x} \tau_k(n)^{l_1}/n^{l_2} \ll (\log x)^{k^{l_1}-1}/x^{l_2-1}$,
\[\begin{aligned}
\mathop{\sum_{n>N}}_{\text{$n$ square-full}} \frac{\tau_k(n)}{n} &\leq
 \sum_{n>\sqrt{N}} \sum_{m=1}^{\infty} \frac{\tau_k(n^2 m^3)}{n^2 m^3} +
 \sum_{n=1}^{\sqrt{N}} \sum_{m\geq (N/n^2)^{1/3}} \frac{\tau_k(n^2 m^3)}{n^2 m^3} \\
&\ll \left(\sum_{n>\sqrt{N}} \frac{\tau_k(n)^2}{n^2}\right)
\left(\sum_{m=1}^{\infty} \frac{\tau_k(m)^3}{m^3}\right) +
 \sum_{n=1}^{\sqrt{N}} \frac{\tau_k(n^2)}{n^2} \frac{(\log N)^{k^3-1}}{(N/n^2)^{2/3}}\\
&\ll \frac{(\log N)^{k^2-1}}{\sqrt{N}} + 
 \frac{(\log N)^{k^3 -1}}{N^{2/3}} \sum_{n=1}^{\sqrt{N}} \frac{\tau_k(n^2)}{n^{2/3}}\\
&\ll \frac{(\log N)^{k^2-1}}{\sqrt{N}} + 
 \frac{(\log N)^{k^3 -1}}{N^{2/3}} (\log N)^{k^2 - 1} N^{1/6} .\end{aligned}\]
\end{proof}
\begin{lem}\label{lem:brim} Let $D$ and $k$ be positive integers. Then
\[\mathop{\sum_{n=1}^{\infty}}_{\text{$n$ is $(D)$-square-full}}
\tau_k(n) \ll \tau(\rad(D)) (\log N)^{k^3+k^2-2} N^{1/2},\]
where the implied constant depends only on $k$.
\end{lem}
\begin{proof}
By Lemmas \ref{lem:supmul} and \ref{lem:prebrim},
\[\begin{aligned}
\mathop{\sum_{n\leq N}}_{\text{$n$ is $(D)$-square-full}} \tau_k(n) 
&= \sum_{m|\rad(D)} \mathop{\sum_{n\leq N/m}}_{\text{$n$ square-full}} \tau_k(m n)\\
&\leq \sum_{m|\rad(D)} \tau_k(m) \mathop{\sum_{n\leq N/m}}_{\text{$n$ square-full}} \tau_k(n)\\
&\ll \sum_{m|\rad(D)} \frac{\tau_k(m)}{m^{1/2}} (\log N)^{k^3+k^2-2} N^{1/2}\\
&\ll \tau(\rad(D)) (\log N)^{k^3+k^2-2} N^{1/2}.\end{aligned}\]
\end{proof}
\begin{lem}\label{lem:bram} Let $D$ and $k$ be positive integers. Then
\[\mathop{\sum_{n=1}^\infty}_{\text{$n$ is $(D)$-square-full}}
\frac{\tau_k(n)}{n} \ll \tau(\rad(D)),\]
where the implied constant depends only on $k$.
\end{lem}
\begin{proof}
We have
\[\begin{aligned}
\mathop{\sum_{n=1}^\infty}_{\text{$n$ is $(D)$-square-full}}
\frac{\tau_k(n)}{n} &=
\sum_{m|\rad(D)} \mathop{\mathop{\sum_{n=1}^\infty}_{m|n}}_{\text{$n/m$ is
square-full}} \frac{\tau_k(m (n/m))}{m (n/m)}\\
&\leq
\sum_{m|\rad(D)}\frac{\tau_k(m)}{m}
 \mathop{\sum_{n=1}^\infty}_{\text{$n$ is square-full}}
\frac{\tau_k(n)}{n} \\ 
&\ll \tau(\rad(D)) \mathop{\sum_{n=1}^\infty}_{\text{$n$ is square-full}}
\frac{\tau_k(n)}{n} .\end{aligned}\]
The statement now follows from Lemma \ref{lem:oohoo}.
\end{proof}
\begin{lem}\label{lem:brum} For any positive integers $k$, $N$, $D$,
\[\mathop{\sum_{n>N}}_{\text{$n$ is $(D)$-square-full}} 
\frac{\tau_k(n)}{n} \ll \tau(\rad(D)) 
 \frac{(\log N)^{k^2+k^3-2}}{\sqrt{N}} ,\]
where the implied constant depends only on $k$.
\end{lem}
\begin{proof}
Clearly
\[\begin{aligned}
\mathop{\sum_{n>N}}_{\text{$n$ is $(D)$-square-full}} 
\frac{\tau_k(n)}{n} &\leq
\sum_{m|\rad(D)} \mathop{\sum_{n>N/m}}_{\text{$n$ square-full}}
 \frac{\tau_k(m n)}{m n}\\
&\leq \sum_{m|\rad(D)} \frac{\tau_k(m)}{m}
 \mathop{\sum_{n>N/m}}_{\text{$n$ square-full}}
\frac{\tau_k(n)}{n} .\end{aligned}\]
Hence, by Lemma \ref{lem:compnr},
\[\begin{aligned}
\mathop{\sum_{n>N}}_{\text{$n$ is $(D)$-square-full}} 
\frac{\tau_k(n)}{n} &\leq
\sum_{m|\rad(D)} \frac{\tau_k(m)}{m}
\frac{(\log (N/m))^{k^2 + k^3 - 2}}{(N/m)^{1/2}}
\\ &\leq \frac{(\log N)^{k^2 + k^3 - 2}}{N^{1/2}} 
 \sum_{m|\rad(D)} \frac{\tau_k(m)}{m^{1/2}} \\
&\ll \tau(\rad(D)) \frac{(\log N)^{k^2 + k^3 - 2}}{N^{1/2}} .\end{aligned}\]
\end{proof}
\subsection{A concrete square-free sieve}
\begin{prop}\label{prop:hjuhju}
Let $K$ be a number field. Let $f:I_K\times \mathbb{Z}\to \mathbb{C}$,
$g:\mathbb{Z}\to \mathbb{C}$ be given with $\max |f(\mathfrak{a},x)|\leq 1$,
$\max |g(x)|\leq 1$. Assume that $f(\mathfrak{a},x)$ depends only on
$\mathfrak{a}$ and on $x \mo \mathfrak{a}$. Let 
$P\in \mathfrak{O}_K\lbrack x \rbrack$. 
Suppose there are $\epsilon_{1,N}, \epsilon_{2,N}\geq 0$
such that for any integer $a$ and any positive integer $m$,
\begin{equation}\label{eq:whynot}
\mathop{\sum_{1\leq x\leq N}}_{x\equiv a \mo m} g(x) \ll 
\left(\frac{\epsilon_{1,N}}{m} + \epsilon_{2,N}\right) N .
\end{equation}
Then, for any integer $a$ and any positive integer $m$,
\begin{equation}\label{eq:agaragar}\begin{aligned}
\mathop{\sum_{1\leq x\leq N}}_{x\equiv a \mo m} f(\sq_K(P(x)),x) g(x) &\ll 
\left(\frac{\epsilon_{1,N}}{m} + \frac{\epsilon'}{m'}\right) N \\&+
\#\{1\leq x\leq N: \exists \mathfrak{p} \text{\:\st\:} \rho(\mathfrak{p})>N^{1/2},
\mathfrak{p}^2|P(x,y)\},\end{aligned}\end{equation}
where 
\begin{equation}\begin{aligned}
\epsilon' &= 
  \sqrt{\max(\epsilon_{2,N},N^{-1/2})} 
    \log(-\max(\epsilon_{2,N},N^{-1/2})),\\
m' &= \min(m, \min(N^{1/2},\epsilon_{2,N}^{-1})) ,\end{aligned}
\end{equation}
and both $c$ and the implied constant in (\ref{eq:agaragar})
depend only on $P$, $K$, and
the implied constant in (\ref{eq:whynot}).
\end{prop}
\begin{proof}
Since the statement is immediate for $P$ constant, we may assume that $P$
is non-constant. Define $S_\mathfrak{a} = \mathfrak{O}_K/\mathfrak{a}$.
Let $\phi_{\mathfrak{a}_1,\mathfrak{a}_2}:S_{\mathfrak{a}_2}\to
S_{\mathfrak{a}_1}$, $\mathfrak{a}_1|\mathfrak{a}_2$, be the natural projection
from $S_{\mathfrak{a}_2}$ to $S_{\mathfrak{a}_1}$. 

For any $\mathfrak{a}\in I_K$, $r\in S_\mathfrak{a}$, set
$f_\mathfrak{a}(r) = f(\mathfrak{a},x)$, where $x$ is any integer with
$x\equiv r \mo \mathfrak{a}$. Let
\[g_\mathfrak{a}(r) = 
\mathop{
\mathop{\mathop{\sum_{1\leq x\leq N}}_{x\equiv a \mo m}}_{\sq_K(P(x)) = \mathfrak{a}}}_{
x\equiv r \mo \mathfrak{a}} g(x) .\]
Then
\[\mathop{\sum_{1\leq x\leq N}}_{x\equiv a \mo m} f(\sq_K(P(x,y)),x) g(x) = 
\sum_{\mathfrak{a}\in I_K} \sum_{r\in S_\mathfrak{a}} f_\mathfrak{a}(r)
g_\mathfrak{a}(r) .\]
Our task is thus to estimate
$\sum_{\mathfrak{a}\in I_K} \sum_{r\in S_\mathfrak{a}} f_\mathfrak{a}(r)
g_\mathfrak{a}(r)$.

Let $s_{\mathfrak{d}}$, $t_{\mathfrak{d}}(r)$ 
 be defined as in the statement of Lemma \ref{lem:wicked}. 
Let \[\gamma(\mathfrak{d}) = \lcm(\rho(\mathfrak{d} \rad_K(\mathfrak{d})),m).\] 
Let $M\leq N^{1/2}$; its optimal
value will be chosen later. We can now apply Lemma \ref{lem:wicked}. What remains
to do is estimate the right side of the inequality it gives us.

By Lemma \ref{lem:dontaskmemyadvice},
\begin{equation}\label{eq:eisler}
s_{\mathfrak{d}}\leq \#\{1\leq x \leq N: \mathfrak{d} \rad_K{\mathfrak{d}} | P(x),\:x\equiv a
\mo m\} \ll \frac{\tau_{\deg P}(\rad_K(\rho(\mathfrak{d}))) N}{
\gamma(\mathfrak{d})}\end{equation}
for $\gamma(\mathfrak{d}) \leq N$.
By definition
\begin{equation}\label{eq:gorey}
t_{\mathfrak{d}}(r) = 
\mathop{\mathop{\mathop{\sum_{1\leq x\leq N}}_{x\equiv a \mo m}}_{
\mathfrak{d} | \sq_K(P(x))}}_{x\equiv r \mo \mathfrak{d}} g(x) .\end{equation}
We can bound
\[
\sum_{\gamma(\mathfrak{d})\leq M} \sum_{r\in S_\mathfrak{d}}
 \left(\sum_{\mathfrak{d}'|\mathfrak{d}} \mu_K(\mathfrak{d}') 
f_{\mathfrak{d}/\mathfrak{d}'}(
\phi_{\mathfrak{d}/\mathfrak{d}',\mathfrak{d}}(r))\right) t_\mathfrak{d}(r) 
\]
trivially by
\[\sum_{\gamma(\mathfrak{d})\leq M} \tau_{K,2}(\rad_K(\mathfrak{d}))
\sum_{r\in S_\mathfrak{d}} |t_\mathfrak{d}(r)|.\]
We then write $\sum_{r\in S_\mathfrak{d}} |t_\mathfrak{d}(r)|$ in full as
\[\sum_{r\in S_\mathfrak{d}} \left|
\mathop{\mathop{\mathop{\sum_{1\leq x\leq N}}_{x\equiv a \mo m}}_{
\mathfrak{d} | \sq_K(P(x))}}_{x\equiv r \mo \mathfrak{d}} g(x) \right| .\]
By Lemma \ref{lem:dontaskmemyadvice}, the set 
$\{x\in \mathbb{Z} : \mathfrak{d}
 | \sq_K(P(x))\}$ is the union of at most \[|\Disc P|^3
\tau_{\deg P}(\rad(\rho(\mathfrak{d})))\] disjoint sets of the form
of the form $L_c = \{x\in \mathbb{Z} : x\equiv c \mo \rho(\mathfrak{d}
\rad_K(\mathfrak{d}))\}$. For every $L_c$, there is an $r\in S_{\mathfrak{d}}$ such that
$x\equiv r \mo \mathfrak{d}$ for every $x\in L_c$. Hence
\[\sum_{r\in S_{\mathfrak{d}}} \left|
\mathop{\mathop{\mathop{\sum_{1\leq x\leq N}}_{x\equiv a \mo m}}_{\mathfrak{d} | \sq_K(P(x))}
}_{x\equiv r \mo \mathfrak{d}} g(x) \right| \leq |\Disc P|^3 \tau_{\deg P}(
\rad(\rho(\mathfrak{d}))) \max_c \left|
 \mathop{\mathop{\sum_{1\leq x\leq N}}_{x\equiv a \mo m}}_{x\equiv c \mo \rho(\mathfrak{d}
\rad_K(\mathfrak{d}))} g(x) \right| .\]
 
We can now apply (\ref{eq:whynot}), obtaining
\[
\sum_{r\in S_\mathfrak{d}} |t_\mathfrak{d}(r)| 
\ll \tau_{\deg P}(\rad(\rho(\mathfrak{d})))
\left(\frac{\epsilon_{1,N}}{\gamma(\mathfrak{d})} + \epsilon_{2,N}\right) N .\]
Lemma \ref{lem:wicked} now yields
\[\begin{aligned}
\sum_{\mathfrak{a}\in I_K} \sum_{r\in S_\mathfrak{a}} 
  f_\mathfrak{a}(r) g_\mathfrak{a}(r) &\leq
\sum_{\gamma(\mathfrak{d}) \leq M}
\sum_{r\in S_\mathfrak{d}}
 \left(\sum_{\mathfrak{d}'|\mathfrak{d}} \mu_K(\mathfrak{d}') 
f_{\mathfrak{d}/\mathfrak{d}'}(
\phi_{\mathfrak{d}/\mathfrak{d}',\mathfrak{d}}(r))\right) t_\mathfrak{d}(r) 
\\&+ 
2 \sum_{M<\gamma(\mathfrak{d})\leq M^2} 
\tau_{K,3}(\mathfrak{d}) s_\mathfrak{d}+ 2 
\mathop{\sum_{\text{$\mathfrak{p}$ prime}}}_{\gamma(\mathfrak{p})>M}
 s_\mathfrak{p} \\
&\leq \sum_{\gamma(\mathfrak{d})\leq M} \tau_{K,2}(\rad_K(\mathfrak{d}))
 \tau_{\deg P}(\rad(\rho(\mathfrak{d}))) 
  \left(\frac{\epsilon_{1,N}}{\gamma(\mathfrak{d})} + \epsilon_{2,N}\right) N
 \\ &+
2 \sum_{M<\gamma(\mathfrak{d})\leq M^2}\frac{\tau_{K,3}(\mathfrak{d}) 
\tau_{\deg P}(\rad(\rho(\mathfrak{d}))) N}{\gamma(\mathfrak{d})}
+ 2 
\mathop{\sum_{\text{$\mathfrak{p}$ prime}}}_{\gamma(\mathfrak{p})>M}
 s_\mathfrak{p} .\end{aligned}\]
By Lemma \ref{lem:changet}, we get
\[\begin{aligned}
\sum_{\mathfrak{a}\in I_K} \sum_{r\in S_\mathfrak{a}} 
  f_\mathfrak{a}(r) g_\mathfrak{a}(r) &\leq
\sum_{\gamma(\mathfrak{d})\leq M} 
 \tau_{2^{\deg K} + \deg P}(\rad(\rho(\mathfrak{d})))
\left(\frac{\epsilon_{1,N}}{\gamma(\mathfrak{d})} + \epsilon_{2,N}\right) N \\ &+
2 \sum_{M<\gamma(\mathfrak{d})\leq M^2}
\frac{\tau_{3^{\deg K} + \deg P}(\rad(\rho(\mathfrak{d})))}{\gamma(\mathfrak{d})}
N
+ 2 
\mathop{\sum_{\text{$\mathfrak{p}$ prime}}}_{\gamma(\mathfrak{p})>M}
 s_\mathfrak{p} .\end{aligned}\]
By Lemma \ref{lem:bleich2}, 
\[
\sum_{\gamma(\mathfrak{d})\leq M} 
\tau_{2^{\deg K} + \deg P}(\rad(\rho(\mathfrak{d})))
\] is at most a constant times
\[\mathop{\mathop{\sum_{n\leq M}}_{\text{$n$ is $(D m)$-square-full}}}_{m|n}
\tau_{2^{\deg K} + \deg P}(\rad(n)) \tau_{\deg(K/\mathbb{Q}) + 2}(n) ,
\]
where $D$ is the product of all rational primes ramifying in $K/\mathbb{Q}$.
Similarly,
\[
\sum_{\gamma(\mathfrak{d})\leq M} 
\frac{\tau_{2^{\deg K} + \deg P}(\rad(\rho(\mathfrak{d})))}{\gamma(
\mathfrak{d})} \]
is at most a constant times
\[\mathop{\mathop{\sum_{n\leq M}}_{\text{$n$ is $(D m)$-square-full}}}_{m|n}
\frac{\tau_{2^{\deg K} + \deg P}(\rad(n)) \tau_{\deg(K/\mathbb{Q}) + 2}(n)}{n}
,\] and
\[\sum_{M<\gamma(\mathfrak{d})\leq M^2}
\frac{\tau_{3^{\deg K} + \deg P}(\rad(\rho(\mathfrak{d})))}{\gamma(\mathfrak{d})}\] is at most a constant times
\[
\mathop{\mathop{\sum_{M<n\leq M^2}}_{\text{$n$ is $(D m)$-square-full}}}_{m|n}
 \frac{\tau_{3^{deg K} + \deg P}(\rad(n))
\tau_{\deg(K/\mathbb{Q}) + 2}(n)}{n} .\]
By Lemma \ref{lem:dontaskmemyadvice},
\[\begin{aligned}
\mathop{\sum_{\text{$\mathfrak{p}$ prime}}}_{\gamma(\mathfrak{p})>M}
 s_\mathfrak{p} &=
\mathop{\mathop{\sum_{\text{$\mathfrak{p}$ prime}}}_{M<\gamma(\mathfrak{p})\leq N}}_{
 \mathfrak{p}\nmid m} s_{\mathfrak{p}} +
\mathop{\mathop{\sum_{\text{$\mathfrak{p}$ prime}}}_{M<\gamma(\mathfrak{p})\leq N}}_{
 \mathfrak{p}|m} s_{\mathfrak{p}} +
\mathop{\mathop{\sum_{\text{$\mathfrak{p}$ prime}}}_{N<\gamma(\mathfrak{p})\leq N m}}_{
 \mathfrak{p}\nmid m} s_{\mathfrak{p}} +
\mathop{\mathop{\sum_{\text{$\mathfrak{p}$ prime}}}_{N<\gamma(\mathfrak{p})\leq N m}}_{
 \mathfrak{p}|m} s_{\mathfrak{p}} +
\mathop{\sum_{\text{$\mathfrak{p}$ prime}}}_{\gamma(\mathfrak{p})> N m}
 s_\mathfrak{p} 
\\ &\ll \mathop{\sum_{\text{$p$ prime}}}_{M< m p^2\leq N} \frac{N}{m p^2} +
\mathop{\mathop{\sum_{\text{$p$ prime}}}_{M< m p \leq N}}_{p|m} \frac{N}{m p} +
\mathop{\sum_{\text{$p$ prime}}}_{N< m p^2\leq N m} \frac{N}{p^2} \\ &+
\mathop{\mathop{\sum_{\text{$p$ prime}}}_{N< m p^2\leq N m}}_{p|m} \frac{N}{p} +
\mathop{\sum_{\text{$\mathfrak{p}$ prime}}}_{\gamma(\mathfrak{p})> N m}
 s_\mathfrak{p} \\
&\leq \frac{N}{\sqrt{M m}} + \frac{N \omega(m)}{M} + \frac{N}{\sqrt{N/m}} + 
m \omega(m) +
\mathop{\sum_{\text{$\mathfrak{p}$ prime}}}_{\gamma(\mathfrak{p})> N m}
 s_\mathfrak{p} .\end{aligned}\]
Write $rem(M) = \frac{N}{\sqrt{M m}} + \frac{N \omega(m)}{M} + \frac{N}{\sqrt{N/m}} + 
m \omega(m)$; it will be swallowed by higher-order terms shortly. (We will assume
$m<N^{1/2}$, as the bound would otherwise be trivial.) Now
\[\begin{aligned}
\sum_{\mathfrak{a}\in I_K} \sum_{r\in S_\mathfrak{a}} 
  f_\mathfrak{a}(r) g_\mathfrak{a}(r) &\ll
\mathop{\mathop{\sum_{n\leq M}}_{\text{$n$ is $(D m)$-square-full}}}_{m|n}
\tau_{q_1}(n) \left(\frac{\epsilon_{1,N}}{n} + \epsilon_{2,N}\right) N  \\ &+
\mathop{\mathop{\sum_{M<n\leq M^2}}_{\text{$n$ is $(D m)$-square-full}}}_{m|n} 
\frac{\tau_{q_2}(n)}{n} N + rem(M) + 
\mathop{\sum_{\text{$\mathfrak{p}$ prime}}}_{\gamma(\mathfrak{p})> N}
 s_\mathfrak{p} \\
&\leq \mathop{\sum_{n\leq M/m}}_{\text{$n$ is $(D m)$-square-full}}
 \tau_{q_1}(n) \tau_{2 q_1}(m) \left(\frac{\epsilon_{1,N}}{m n} + \epsilon_{2,N}\right) N 
\\ &+
 \mathop{\sum_{n>M/m}}_{\text{$n$ is $(D m)$-square-full}}
  \frac{\tau_{q_2}(n) \tau_{2 q_2}(m)}{m n} N +
rem(M)+
 \mathop{\sum_{\text{$\mathfrak{p}$ prime}}}_{\gamma(\mathfrak{p})> N m}
 s_\mathfrak{p},
\end{aligned}\]
where $q_1=(2^{\deg K} + \deg P)(\deg(K/\mathbb{Q}) + 1)$, 
$q_2=(3^{\deg K} + \deg P)(\deg(K/\mathbb{Q}) + 1)$.
Now note that
\[\mathop{\sum_{\text{$\mathfrak{p}$ prime}}}_{\gamma(\mathfrak{p})> N}
 s_\mathfrak{p} \ll
\#\{1\leq x\leq N: \exists \mathfrak{p} \text{\:\st\:} \rho(\mathfrak{p})>N^{1/2},
\mathfrak{p}^2|P(x,y)\}. \]
By Lemmas \ref{lem:brim}, \ref{lem:bram} and \ref{lem:brum} we can conclude that
\[\begin{aligned}
\sum_{\mathfrak{a}\in I_K} \sum_{r\in S_\mathfrak{a}} 
  f_\mathfrak{a}(r) g_\mathfrak{a}(r) &\ll
\left(\frac{\epsilon_{1,N}}{m} + \epsilon_{2,N} (\log M)^{q_3} \sqrt{M/m} +
 \frac{(\log M)^{q_4}}{\sqrt{M m}}\right) \\ &\cdot \tau_{q_5}(m) \tau(\rad(D m)) N\\
&+ rem(M) +
\#\{-N\leq x\leq N: \exists \mathfrak{p} \text{\:\st\:} \rho(\mathfrak{p})>N^{1/2},
\mathfrak{p}^2|P(x,y)\}\\
&\ll
\left(\frac{\epsilon_{1,N}}{m} + \epsilon_{2,N} (\log M)^{q_3} \sqrt{M/m} +
 \frac{(\log M)^{q_4}}{\sqrt{M m}}\right) \tau_{q_6}(m) N\\
&+
\#\{-N\leq x\leq N: \exists \mathfrak{p} \text{\:\st\:} \rho(\mathfrak{p})>N^{1/2},
\mathfrak{p}^2|P(x,y)\},\end{aligned}\]
where $q_3 = q_1^3 + q_1^2 -2$, $q_4 = q_2^3+q_2^2-2$, $q_5 = \max(2 q_1, 2 q_2)$,
$q_6 = 2 q_5$. Set $M=\min\left(N^{1/2},\frac{1}{\epsilon_{2,N}}\right)$,
$c_1 = q_6$, $c_2=\max(q_3,q_4)$. The statement follows.
\end{proof}

\begin{prop}\label{prop:hjohjo}
Let $K$ be a number field.
Let $f:I_K\times \{(x,y)\in \mathbb{Z}^2 : \gcd(x,y)=1\} \to \mathbb{C}$, $g:\{(x,y)\in \mathbb{Z}^2 : \gcd(x,y)=1\} \to \mathbb{C}$
be given with $\max |f(\mathfrak{a},x,y)| \leq 1$, $\max |g(x,y)|\leq 1$. Assume
that $f(\mathfrak{a},x,y)$ depends only on $\mathfrak{a}$ and on 
$\{\frac{x \mo \mathfrak{p}}{y \mo \mathfrak{p}}\}_{\mathfrak{p}|\mathfrak{a}}
\in \prod_{\mathfrak{p}|\mathfrak{a}} 
\mathbb{P}^1(\mathfrak{O}_K/\mathfrak{p})$. 
Let $P\in \mathfrak{O}_K\lbrack x,y\rbrack$ be a homogeneous polynomial.
Let $S$ be a subset of $\mathbb{R}^2$.
Suppose there are $\epsilon_{1,N},\, \epsilon_{2,N}\geq 0$ such that for
any lattice coset $L\subset \mathbb{Z}^2$,
\begin{equation}\label{eq:pcond}
\mathop{\sum_{(x,y)\in S\cap \lbrack -N,N\rbrack^2\cap L}}_{\gcd(x,y)=1}
g(x,y) \ll \left(\frac{\epsilon_{1,N}}{
\lbrack \mathbb{Z}^2 :L\rbrack} +
\epsilon_{2,N}\right) N^2.
\end{equation}
Then, for any lattice coset $L\subset \mathbb{Z}^2$,
\begin{equation}\label{eq:primni}\begin{aligned}
\mathop{\sum_{(x,y)\in S\cap \lbrack -N,N\rbrack^2 \cap L}}_{\gcd(x,y)=1}
&f(\sq_K(P(x,y)),x,y) g(x,y) \\ &\ll
\left(\frac{\epsilon_{1,N}}{\lbrack \mathbb{Z}^2 : L\rbrack} +
\frac{\epsilon'}{\sqrt{m'}}\right) N^2 \\ &+
\{-N\leq x,y\leq N : \exists \mathfrak{p} \text{\:\st\:} \rho(\mathfrak{p})>N,
\mathfrak{p}^2|P(x,y)\},\end{aligned}\end{equation}
where
\[\begin{aligned}
 \epsilon' &= \sqrt{\max(\epsilon_{2,N},N^{-1/2})} \log(-\max(\epsilon_{2,N},
N^{-1/2})),\\
 m' &= \min(\lbrack \mathbb{Z}^2 : L\rbrack, \min(N^{1/2},\epsilon_{2,N}^{-1})) ,\end{aligned}\]
the constants $c_1$ and $c_2$ depend only on $P$ and $K$, and
the implied constant in (\ref{eq:primni}) depends only on $P$, $K$ and 
the implied constant in (\ref{eq:pcond}).
\end{prop}
\begin{proof}
Since the statement is immediate for $P$ constant we may assume 
that $P$ is non-constant.
Define $S_\mathfrak{a} = \prod_{\mathfrak{p}|\mathfrak{a}} 
\mathbb{P}^1(\mathfrak{O}_K/\mathfrak{p})$.
Let $\phi_{\mathfrak{a}_1,\mathfrak{a}_2} : K_{\mathfrak{a}_2} 
\to K_{\mathfrak{a}_1}$, $\mathfrak{a}_1|\mathfrak{a}_2$, be
the natural projection from $S_{\mathfrak{a}_2}$ to 
$S_{\mathfrak{a}_1}$. Write 
$\phi_\mathfrak{a}(x,y) =
\{\frac{x \mo \mathfrak{p}}{y \mo \mathfrak{p}}\}_{\mathfrak{p}|\mathfrak{a}}
\in S_\mathfrak{a}$ 
for any coprime $x$, $y$.

For any $\mathfrak{a}\in I_K$, $r\in S_\mathfrak{a}$,
set $f_\mathfrak{a}(r) = f(\mathfrak{a},x,y)$, where $x$, $y$
are any coprime integers with $\phi_\mathfrak{a}(x,y) = r$. Let
\[
g_\mathfrak{a}(r) = \mathop{\mathop{\mathop{\sum_{(x,y)\in S\cap \lbrack -N,N
\rbrack^2\cap L}}_{\gcd(x,y)=1}}_{\sq_K(P(x,y)) = \mathfrak{a}}}_{
\phi_\mathfrak{a}(x,y) = r}
g(x,y).\]
Then
\[\mathop{\sum_{(x,y)\in S\cap \lbrack -N,N\rbrack^2\cap L}}_{\gcd(x,y)=1}
f(\sq_K(P(x,y)),x,y) g(x,y) = \sum_{\mathfrak{a}\in I_K}
\sum_{r\in S_\mathfrak{a}} f_{\mathfrak{a}}(r) g_{\mathfrak{a}}(r) .\]
The question now is how to estimate
$\sum_{\mathfrak{a}\in I_K}
\sum_{r\in S_\mathfrak{a}} f_{\mathfrak{a}}(r) g_{\mathfrak{a}}(r)$.

Let $s_{\mathfrak{d}}$, $t_{\mathfrak{d}}(r)$ be as in the statement of Lemma \ref{lem:wicked}. Let 
\[\gamma(\mathfrak{d}) = \lcm(\rho(\mathfrak{d} \rad_K(\mathfrak{d})),\lbrack \mathbb{Z}^2 : L
\rbrack) .\] Let $M\leq N$. 
By Lemmas \ref{lem:ramsay} and \ref{lem:sollat},
\[
\begin{aligned}
s_\mathfrak{d} &\leq \#\{(x,y)\in S\cap \lbrack -N,N\rbrack^2\cap L : \gcd(x,y)=1,
\mathfrak{d} \rad_K(\mathfrak{d}) |P(x,y) \}\\
&\ll \frac{\tau_{2 \deg P}(\rad_K(\rho(\mathfrak{d}))) N^2}{\gamma(\mathfrak{d})}
\end{aligned}
\]
for $\gamma(\mathfrak{d})\leq N^2$. By definition,
\begin{equation}\label{eq:deft} t_\mathfrak{d}(r) =
\mathop{\mathop{\mathop{\sum_{(x,y)\in S\cap \lbrack -N,N
\rbrack^2 \cap L}}_{\gcd(x,y)=1}}_{
\mathfrak{d}|\sq_K(P(x,y))}}_{\phi_\mathfrak{d}(x,y) = r}
g(x,y).\end{equation}
We can bound
\[
\sum_{\gamma(\mathfrak{d})\leq M} \sum_{r\in S_\mathfrak{d}}
 \left(\sum_{\mathfrak{d}'|\mathfrak{d}} \mu_K(\mathfrak{d}') 
f_{\mathfrak{d}/\mathfrak{d}'}(
\phi_{\mathfrak{d}/\mathfrak{d}',\mathfrak{d}}(r))\right) t_\mathfrak{d}(r) 
\]
trivially by
\[\sum_{\gamma(\mathfrak{d})\leq M} \tau_{K,2}(\rad_K(\mathfrak{d}))
\sum_{r\in S_\mathfrak{d}} |t_\mathfrak{d}(r)|.\]
We write $\sum_{r\in S_\mathfrak{d}} |t_\mathfrak{d}(r)|$ in full as
\[\sum_{r\in S_\mathfrak{d}} \left|
\mathop{\mathop{\mathop{\sum_{(x,y)\in S\cap \lbrack -N,N
\rbrack^2\cap L}}_{\gcd(x,y)=1}}_{
\mathfrak{d}|\sq_K(P(x,y))}}_{\phi_\mathfrak{d}(x,y) = r}
g(x,y)\right| .\]
By Lemma \ref{lem:sollat}, the set $\{(x,y)\in \mathbb{Z}^2 :
\gcd(x,y)=1, \mathfrak{d}|\sq_K(P(x))\}$ 
is the union of at most $|\Disc P|^3 \tau_{2\deg P}(\rad_K(\mathfrak{d}))$
disjoint sets of the form \[R\cap \{(x,y)\in \mathbb{Z}^2 : \gcd(x,y)=1\},\]
where $R$ is a lattice of index $\rho(\mathfrak{d} \rad_K(\mathfrak{d}))$.
For every $R$ of index $\rho(\mathfrak{d} \rad_K(\mathfrak{d}))$, there
is an $r\in S_{\mathfrak{d}}$ such that $\phi_{\mathfrak{d}}(x,y)=r$ for
every $(x,y)\in R$.
Hence
\[\sum_{r\in S_{\mathfrak{d}}} \left|
\mathop{\mathop{\mathop{\sum_{(x,y)\in S\cap \lbrack -N,N
\rbrack^2\cap L}}_{\gcd(x,y)=1}}_{\mathfrak{d}|\sq_K(P(x,y))}}_{
\phi_\mathfrak{d}(x,y) = r}
g(x,y) \right|\] is equal to at most 
$(\Disc P)^3 \tau_{2\deg P}(\rad(m))$ times
\[
\mathop{\max_R}_{\lbrack \mathbb{Z}^2 : R\rbrack = \gamma(\mathfrak{d})}
\left|
\mathop{\mathop{\sum_{(x,y)\in S\cap \lbrack -N,N
\rbrack^2 \cap L}}_{\gcd(x,y)=1}}_{(x,y)\in R}
g(x,y) \right|.\]
We can now apply (\ref{eq:pcond}), obtaining
\[\mathop{\mathop{\sum_{(x,y)\in S\cap \lbrack -N,N
\rbrack^2 \cap L}}_{\gcd(x,y)=1}}_{(x,y)\in R} g(x,y) \ll
\left(\frac{\epsilon_{1,N}}{\gamma(\mathfrak{d})} +
 \epsilon_{2,N}\right) N^2 .\]
Lemma \ref{lem:wicked} now yields
\[\begin{aligned}
\sum_{\mathfrak{a}\in I_K} \sum_{r\in S_\mathfrak{a}} 
  f_\mathfrak{a}(r) g_\mathfrak{a}(r) &\leq
\sum_{\gamma(\mathfrak{d}) \leq M} \sum_{r\in S_\mathfrak{a}}
 \left(\sum_{\mathfrak{d}'|\mathfrak{d}} \mu_K(\mathfrak{d}') 
f_{\mathfrak{d}/\mathfrak{d}'}(\phi_{\mathfrak{d}/\mathfrak{d}',\mathfrak{d}}(r))\right) t_\mathfrak{d}(r)\\&+ 
2 \sum_{M<\gamma(\mathfrak{d})\leq M^2} 
\tau_{K,3}(\mathfrak{d}) s_\mathfrak{d}+ 2 
\mathop{\sum_{\text{$\mathfrak{p}$ prime}}}_{\gamma(\mathfrak{p})>M}
 s_\mathfrak{p} \\
&\leq \sum_{\gamma(\mathfrak{d})\leq M} \tau_{K,2}(\rad_K(\mathfrak{d}))
 \tau_{2 \deg P}(\rad(\rho(\mathfrak{d}))) 
\left(\frac{\epsilon_{1,N}}{\phi(\gamma(\mathfrak{d}))} + 
 \epsilon_{2,N}\right) N^2
\\ &+
2 \sum_{M<\gamma(\mathfrak{d})\leq M^2}\frac{\tau_{K,3}(\mathfrak{d}) 
\tau_{2 \deg P}(\rad(\rho(\mathfrak{d})) N^2}{\gamma(\mathfrak{d})}
+ 2 
\mathop{\sum_{\text{$\mathfrak{p}$ prime}}}_{\gamma(\mathfrak{p})>M}
 s_\mathfrak{p} .\end{aligned}\]
The remainder of the argument is the same as in Proposition \ref{prop:hjuhju}.
\end{proof}
\begin{Rem}
Proposition \ref{prop:hjohjo} still holds if ``lattice coset'' is
replaced by ``lattice'' throughout the statement.
\end{Rem}
\section{A global approach to the square-free sieve}
\subsection{Elliptic curves, heights and lattices}\label{subs:shlime}
As is usual, we write $\hat{h}$ for the canonical height on an elliptic
curve $E$, and $h_x$, $h_y$ for the height on $E$ with respect to $x$, $y$:
\[h_x((x,y)) = \begin{cases} 0 &\text{if $P=O$,}\\
\log H(x) &\text{if $P=(x,y)$,}\end{cases}\]
\[h_y((x,y)) = \begin{cases} 0 &\text{if $P=O$,}\\
\log H(y) &\text{if $P=(x,y)$,}\end{cases}\]
where $O$ is the origin of $E$, taken to be the point at infinity, and
\[\begin{aligned}
H(y) &= (H_K(y))^{1/\lbrack K : \mathbb{Q}\rbrack},\\
H_K(y) &= \prod_v \max(|y|_v^{n_v},1),\end{aligned}\]
where $K$ is any number field containing $y$,  the product
$\prod_v$ is taken over all places $v$ of $K$, and $n_v$ denotes
the degree of $K_v/\mathbb{Q}_v$.

In particular, if $x$ is a rational number $x_0/x_1$, $\gcd(x_0,x_1)=1$, then
\[\begin{aligned}
H(x) &= H_{\mathbb{Q}}(x) = \max(|x_0|,|x_1|),\\
h_x((x,y)) &= \log(\max(|x_0|,|x_1|)).\end{aligned}\]

The differences $|\hat{h} - \frac{1}{2} h_x|$ and 
$|\hat{h} - \frac{1}{3} h_y|$ are bounded
on the set of all points of $E$ (not merely on $E(\mathbb{Q})$). This
basic property of the canonical height will be crucial in our analysis.

\begin{lem}\label{lem:soide}
Let $f\in \mathbb{Z}\lbrack x\rbrack$ be a cubic polynomial of 
non-zero discriminant. For every square-free rational integer $d$, 
let $E_d$ be
the elliptic curve 
\[E_d : d y^2 = f(x) .\]
Let $P = (x,y) \in E_d(\mathbb{Q})$. Consider the point
 $P'=(x,d^{1/2} y)$ on $E_1$.
Then $\hat{h}(P) = \hat{h}(P')$, where the canonical heights are defined on
$E_d$ and $E_1$, respectively,
\end{lem}
\begin{proof} Clearly $h_x(P') = h_x(P)$. Moreover $(P+P)' = P' + P'$.
Hence \[\hat{h}(P) = \frac{1}{2} \lim_{N\to \infty} 4^{-N}
h_x(\lbrack 2^N\rbrack P) = \frac{1}{2} \lim_{N\to \infty} 4^{-N}
h_x(\lbrack 2^N\rbrack P') = \hat{h}(P') .\]
\end{proof}

\begin{lem}\label{lem:pseusi}
Let $f\in \mathbb{Z}\lbrack x\rbrack$ be an irreducible cubic polynomial
of non-zero discriminant. Let $E$ be the elliptic curve given by
$E:y^2 = f(x)$. Let $d\in \mathbb{Z}$ be square-free. Let $x$, $y$ be rational
numbers, $y\ne 0$, such that $P = (x, d^{1/2} y)$ lies on $E$. Then
\[h_y(P) \geq \frac{3}{8} \log |d| + C_f,\]
where $C_f$ is a constant depending only on $f$.
\end{lem}
\begin{proof}
Write $y = y_0/y_1$, where $y_0$ and $y_1$ are coprime integers. Then
\begin{equation}\label{eq:wraw}
H(y) = \max\left(\frac{|y_0| |d|^{1/2}}{\sqrt{\gcd(d,y^2)}},
  \frac{|y_1|}{\sqrt{\gcd(d,y_1^2)}}\right) .\end{equation}

Write $a$ for the leading coefficient of $f$. Let $p|\gcd(d,y_1^2)$, $p\nmid a$.
Since $d$ is square-free, $p^2\nmid \gcd(d,y^2)$. Suppose $p^2\nmid y_1$. Then $\nu_p(d y^2) = -1$.
However, $d y^2 = f(x)$ implies that, if $\nu_p(x)\geq 0$, then $\nu_p(d y^2)\geq 0$, and
if $\nu_p(x)<0$, then $\nu_p(d y^2)\leq -3$. Contradiction. Hence $p|\gcd(d,y_1^2)$,
$p\nmid a$ imply $p^2\nmid \gcd(d,y_1^2)$, $p^2|y_1$. Therefore
 $|y_1|\geq (\gcd(d,y_1^2)/a)^2$. 

By (\ref{eq:wraw}) it follows that
\[\begin{aligned}
H(P)&\geq \max\left(\frac{|d|^{1/2}}{\sqrt{\gcd(d,y_1^2)}},
  \frac{|y_1|}{\sqrt{\gcd(d,y_1^2)}}\right)\\
&\geq \max\left(\frac{|d|^{1/2}}{\sqrt{\gcd(d,y_1^2)}},
  \frac{(\gcd(d,y_1^2))^{3/2}}{a^2}\right) .\end{aligned}\]
Since $\max(|d|^{1/2} z^{-1/2}, z^{3/2}/a_3^2)$ is minimal when
$|d|^{1/2} z^{-1/2} = z^{3/2}/a_3^2$, i.e., when
$z = a_3 |d|^{1/4}$, we obtain
\[H(P)\geq |d|^{3/8} |a_1|^{-1/2} .\]
Hence 
\[h_y(P) = \log H(P) \geq \frac{3}{8} \log |d| - \frac{1}{2} \log |a| .\]
\end{proof}
\begin{cor}\label{cor:nosilv}
Let $f\in \mathbb{Z}\lbrack x\rbrack$ be a cubic polynomial of non-zero discriminant.
For every square-free rational integer $d$, let $E_d$ be the elliptic curve
\[E_d : d y^2 = f(x) .\]
Let $P = (x,y)\in E_d(\mathbb{Q})$. Then
\[\hat{h}(P)\geq \frac{1}{8} \log |d| + C_f ,\]
where $C_f$ is a constant depending only on $f$.
\end{cor}
\begin{proof}
Let $P' = (x, d^{1/2} y)\in E_1$. By Lemma \ref{lem:soide}, 
$\hat{h}(P) = \hat{h}(P')$. The difference $|\hat{h} - h_x|$ is bounded on $E$. The
statement follows from Lemma \ref{lem:pseusi}.
\end{proof}

The following crude estimate will suffice for some of our purposes.
\begin{lem}\label{lem:salom}
Let $Q$ be a positive definite quadratic form on $\mathbb{Z}^r$. Suppose
$Q(\vec{x})\geq c_1$ for all non-zero $\vec{x}\in \mathbb{Z}^r$. Then there are
at most 
\[(1 + 2 \sqrt{c_2/c_1})^{r}\]
values of $\vec{x}$ for which $Q(\vec{x})\leq c_2$.
\end{lem}
\begin{proof}
There is a linear bijection $f:\mathbb{Q}^r\to \mathbb{Q}^r$ taking $Q$
to the square root of the Euclidean norm: $Q(\vec{x}) = |f(\vec{x})|^2$ for
all $\vec{x}\in \mathbb{Q}^r$. Because $Q(\vec{x})> c_1$ for all
non-zero $\vec{x}\in \mathbb{Z}^r$, we have that $f(\mathbb{Z}^r)$ is a lattice
$L\subset \mathbb{Q}^r$ such that $|\vec{x}|\geq c_1^{1/2}$ for all
$\vec{x}\in L$, $\vec{x}\ne 0$. We can draw a sphere $S_{\vec{x}}$ of
radius $\frac{1}{2} c_1^{1/2}$ around each point $\vec{x}$ of $L$. The
spheres do not overlap. If $\vec{x}\in L$, $|\vec{x}|\in c_2^{1/2}$,
then $S_{\vec{x}}$ is contained in the sphere $S'$ of radius 
$c_2^{1/2} + c_1^{1/2}/2$ around the origin. The total volume
of all spheres $S_{\vec{x}}$ within $S'$ is no greater than the volume of
$S'$. Hence
\[\# \{\vec{x}\in L: |\vec{x} | \leq c_2^{1/2}\} \cdot (c_1^{1/2}/2)^r
\leq (c_2^{1/2} + c_1^{1/2}/2)^r .\]
The statement follows.
\end{proof}
\begin{cor}\label{cor:magd}
  Let $E$ be an elliptic curve over $\mathbb{Q}$. Suppose there are no non-torsion points
$P\in E(\mathbb{Q})$ of canonical height $\hat{h}(P)<c_1$. Then there are at most
\[O\left((1 + 2\sqrt{c_2/c_1})^{\rnk(E)}\right)\]
points $P\in E(\mathbb{Q})$ for which $\hat{h}(P) < c_2$. The implied constant is absolute.
\end{cor} \begin{proof}
The canonical height $\hat{h}$ is a positive definite quadratic form on the free part
$\mathbb{Z}^{\rnk(E)}$ of $E(\mathbb{Q}) \sim \mathbb{Z}^{\rnk(E)} \times T$. A classical
theorem of Mazur's \cite{Maz} states that the cardinality of $T$ is at most $16$. 
 Apply Lemma \ref{lem:salom}.
\end{proof}
Note that we could avoid the use of Mazur's theorem, since Lemmas \ref{lem:soide} and 
\ref{lem:pseusi} imply that the torsion group of $E_d$ is either $\mathbb{Z}/2$ or trivial
for large enough $d$.

\subsection{Twists of cubics and quartics}

Let $f(x) = a_4 x^4 + a_3 x^3 + a_2 x^2 + a_1 x + a_0 \in 
\mathbb{Z}\lbrack x\rbrack$ be an irreducible polynomial of degree $4$.
For every square-free $d\in \mathbb{Z}$, consider the curve
\begin{equation}\label{eq:conco1} C_d : d y^2 = f(x) .\end{equation}
If there is a rational point $(r,s)$ on $C_d$, then there is a birational
map from $C_d$ to the elliptic curve
\begin{equation}
E_d : d y^2 = x^3 + a_2 x^2 + (a_1 a_3 - 4 a_0 a_4) x - 
(4 a_0 a_2 a_4 - a_1^2 a_4 - a_0 a_3^2) .\end{equation}
Moreover, we can construct such a birational map in terms of $(r,s)$
as follows. Let $(x,y)$ be a rational point on $C_d$. We can rewrite
(\ref{eq:conco1}) as
\[y^2 = \frac{1}{d} f(x) .\]
We change variables:
\[x_1 = x-r, \;\; y_1 = y\]
satisfy
\[y^2 = \frac{1}{d} \left(\frac{1}{4!} f^{(4)}(r) x_1^4 +
\frac{1}{3!} f^{(3)}(r) x_1^3 + \frac{1}{2!} f''(r) x_1^2 +
\frac{1}{1!} f'(r) x_1 + f(r)\right) .\]
We now apply the standard map for putting quartics in Weierstrass form:
\[\begin{aligned}
x_2 &= (2 s (y_1 + s) + f'(r) x_1/d)/x_1^2,\\
y_2 &= (4 s^2 (y_1 + s) + 2 s (f'(r) x_1/d + f''(r) x_1^2/(2 d)) -
(f'(r)/d)^2 x_1^2/(2 s))/x_1^3 \end{aligned}\]
satisfy
\begin{equation}\label{eq:pfauen1}
y_2^2 + A_1 x_2 y_2 + A_3 y_2 = x_2^3 + A_2 x_2^2 + A_4 x_2 + A_6
\end{equation}
with
\[\begin{aligned}
A_1 &= \frac{1}{d} f'(r)/s,\;\;&A_2 = 
                    \frac{1}{d}(f''(r)/2 - (f'(r))^2/(4 f(r))),\\
A_3 &= \frac{2 s}{d} f^{(3)}(r)/3!,\;\;&A_4 = -\frac{1}{d^2} \cdot
4 f(r) \cdot \frac{1}{4!} f^{(4)}(r) ,\\
A_6 &= A_2 A_4 .\end{aligned}\]
To take (\ref{eq:pfauen1}) to $E_d$, we apply a linear change of variables:
\[\begin{aligned}
x_3 &= d x_2 + r (a_3 + 2 a_4 r),\\
y_2 &= \frac{d}{2} (2 y_2 + a_1 x_2 + a_3) \end{aligned}\]
satisfy
\[d y_3^2 = 
x_3^3 + a_2 x_3^2 + (a_1 a_3 - 4 a_0 a_4) x_3 - 
(4 a_0 a_2 a_4 - a_1^2 a_4 - a_0 a_3^2) .\]
We have constructed a birational map $\phi_{r,s}(x,y)\mapsto (x_3,y_3)$
from $C_d$ to $E_d$.

Now consider the equation
\begin{equation}\label{eq:threest}
d y^2 = a_4 x^4 + a_3 x^3 z + a_2 x^2 z^2 + a_1 x z^3 + a_0 z^4 .
\end{equation}
Suppose there is
a solution $(x_0,y_0,z_0)$ to (\ref{eq:threest}) with 
$x_0, y_0, z_0 \in \mathbb{Z}$, $|x_0|, |z_0|\leq N$, $z_0\ne 0$. Then 
$(x_0/z_0,y_0/z_0^2)$ is a rational point on (\ref{eq:conco1}). We can
set $r = x_0/z_0$, $s = y_0/z_0^2$ and define a map $\phi_{r,s}$ from
$C_d$ to $E_d$ as above. Now let $x, y, z\in \mathbb{Z}$, $|x|, |z|\leq N$,
$z_0\ne 0$,
be another solution to (\ref{eq:threest}). Then
\[P = \phi_{r,s}(x_0/z_0,y_0/z_0^2)\]
is a rational point on $E_d$. Notice that $|y_0|, |y| \ll (N^4/d)^{1/2}$.
Write \[\phi_{r,s}(P) = (u_0/u_1,v),\]
where $u_0,u_1\in \mathbb{Z}$, $v\in \mathbb{Q}$, $\gcd(u_0,u_1)=1$. By a 
simple examination of the construction of $\phi_{r,s}$ we can determine
that $\max(u_0,u_1) \ll N^7$, where the implied constant depends only
on $a_0, a_1, \dotsb, a_4$.
In other words,
\begin{equation}\label{eq:quadsta}
h_x(P) \leq 7 \log N + C,\end{equation}
where $C$ is a constant depending only on $a_j$. Notice that (\ref{eq:quadsta})
holds even for $(x,y,z) = (x_0,y_0,z_0)$, as then $P$ is the origin of $E$.

The value of $h_x(P)$ is independent of whether $P$ is considered as a rational
point of $E_d$ or as a point of $E_1$. Let $\hat{h}_{E_1}(P)$ be the 
canonical height of $P$ as a point of $E_1$. Then
\[|\hat{h}_{E_1} (P) - \frac{1}{2} h_x(P)|\leq C',\]
where $C'$ depends only on $f$. By Lemma \ref{lem:soide}, the canonical height
$\hat{h}_{E_1}(P)$ of $P$  as a point of $E_1$ equals the canonical height
$\hat{h}_{E_d}(P)$ of $P$ as a point of $E_d$. Hence 
\[|\hat{h}_{E_d}(P) - \frac{1}{2} h_x(P)|\leq C' .\]
Then, by (\ref{eq:quadsta}),\[\hat{h}_{E_d}(P)\leq \frac{7}{2} \log N + (C/2 + C') .\]

We have proven
\begin{lem}\label{lem:ari} Let 
$f(x,z) = 
a_4 x^4 + a_3 x^3 z + a_2 x^2 z^2 + a_1 x z^3 + a_0 z^4\in \mathbb{Z}\lbrack x,z\rbrack$ 
be an irreducible homogeneous polynomial. Then there is a constant $C_f$
such that the following holds. Let $N$ be any positive integer. Let $d$
be any square-free integer. Let $S_{d,1}$ be the set of all
solutions $(x,y,z)\in \mathbb{Z}^3$ to
\[d y^2 = f(x,z) \]
satisfying $|x|, |z| \leq N$, $\gcd(x,z)=1$. 
Let $S_{d,2}$ be the set of all rational points $P$
on
\begin{equation}\label{eq:vova} E_d : d y^2 = 
x^3 + a_2 x^2 + (a_1 a_3 - 4 a_0 a_4) x - 
(4 a_0 a_2 a_4 - a_1^2 a_4 - a_0 a_3^2) \end{equation}
with canonical height
\[\hat{h}(P) \leq \frac{7}{2}\log N + C_f .\]
Then there is an injective map from $S_{d,1}$ to $S_{d,2}$.
\end{lem}
We can now apply the results of subsection \ref{subs:shlime}.
\begin{prop}\label{prop:clinton}
Let $f(x,z) =
a_4 x^4 + a_3 x^3 z + a_2 x^2 z^2 + a_1 x z^3 + a_0 z^4\in 
\mathbb{Z}\lbrack x,z\rbrack$ 
be an irreducible homogeneous polynomial. Then there are constants 
$C_{f,1}$, $C_{f,2}$, $C_{f,3}$
such that the following holds. Let $N$ be any positive integer. Let $d$
be any square-free integer. Let $S_{d}$ be the set of all 
solutions $(x,y,z)\in \mathbb{Z}^3$ to
\[d y^2 = f(x,z) \]
satisfying $|x|, |z| \leq N$, $\gcd(x,z)=1$. Then
\[\# S_d \ll \begin{cases} 
\left(1 + 2 \sqrt{(\frac{7}{2} \log N + C_{f,1})/
(\frac{1}{8} \log |d| + C_{f,2})}\right)^{\rnk(E_d)}
&\text{if $|d|\geq C_{f,4}$,}\\
\left(1 + 2 C_{f,3} \sqrt{\frac{7}{2} \log N + C_{f,1}}\right)^{\rnk(E_d)}
&\text{if $|d| < C_{f,4}$,}\end{cases}\]
where $C_{f,4} = e^{9 C_{f,2}}$, $E_d$ is as in (\ref{eq:vova}),
and the implied constant depends only on $f$.
\end{prop}
\begin{proof}
If $|d|\leq C_{f,4}$, apply Corollary \ref{cor:magd} and Lemma
\ref{lem:ari}. If $|d|> C_{f,4}$, apply Corollary \ref{cor:nosilv},
Corollary \ref{cor:magd} and Lemma
\ref{lem:ari}. 
\end{proof}
\subsection{Divisor functions and their averages}
As is usual, we denote by $\omega(d)$ the number of prime divisors of
a positive integer $d$.
Given an extension $K/\mathbb{Q}$,
we define \[\omega_K(d) = 
\mathop{\sum_{\mathfrak{p}\in I_K}}_{\mathfrak{p}|d} 1 .\]
\begin{lem}\label{lem:notse}
Let $f(x)\in \mathbb{Z}\lbrack x\rbrack$ be an irreducible polynomial of
degree $3$ and non-zero discriminant.
Let $K = \mathbb{Q}(\alpha)$, where $\alpha$ is a root of $f(x)=0$.
For every square-free rational integer $d$, let $E_d$ be the elliptic curve given by
\[ d y^2 = f(x) .\]
Then 
\[\rnk(E_d) = C_f + \omega_K(d) - \omega(d),\]
where $C_f$ is a constant depending only on $f$.
\end{lem}
\begin{proof}
Write $f(x) = a_3 x^3 + a_2 x^2 + a_1 x + a_0$.
Let $f_d(x) = a_3 x^3 + d a_2 x^2 + d^2 a_1 x + d^3 a_0$. 
Then $d \alpha$ is a root of $f_d(x)=0$.
Clearly $\mathbb{Q}(d \alpha) = \mathbb{Q}(\alpha)$. If $p$ is
a prime of good reduction for $E_1$, then $E_d$ will have additive reduction
at $p$ if $p|d$, and good reduction at $p$ if $p\nmid d$.
 The statement now follows immediately from the standard bound in,
say, \cite{BK}, Prop. 7.1.
\end{proof}

\begin{lem}\label{lem:patabran}
Let $K/\mathbb{Q}$ be a non-Galois extension of $\mathbb{Q}$ of degree $3$.
Let $L/\mathbb{Q}$ be the normal closure of $K/\mathbb{Q}$.
Let $K'/\mathbb{Q}$ be the quadratic subextension of $K/\mathbb{Q}$.
 Then the following
statements are equivalent:
\begin{itemize} 
\item 
$p$ splits as $p = \mathfrak{p}_1 \mathfrak{p}_2$ in $K/\mathbb{Q}$, where 
$\mathfrak{p}_1$ and $\mathfrak{p}_2$ are prime ideals of $K$,
\item $p$ does not split in $K'/\mathbb{Q}$.
\end{itemize}
\end{lem}
\begin{proof}
Clearly $\Gal(K/\mathbb{Q}) = S_3$. Consider the Frobenius element
$\Frob_p$ as a conjugacy class in $S_3$. There are three conjugacy classes in
$S_3$; we shall call them $C_1$ (the identity), $C_2$ (the transpositions)
and $C_3$ (the $3$-cycles). 
If $\Frob_p=C_1$, then $p$ splits completely in $K$ and in $K'$. It remains
to consider the other two cases, $\Frob_p=C_2$ and $\Frob_p = C_3$.

Suppose $\Frob_p = C_2$. Then $p$ splits as $p = \mathfrak{q}_1 \mathfrak{q}_2 \mathfrak{q}_3$
in $L/\mathbb{Q}$. We have
\[C_2 = \{\Frob_{\mathfrak{q}_1}, \Frob_{\mathfrak{q}_2}, \Frob_{\mathfrak{q}_3}\} .\]
Hence exactly one of $\Frob_{\mathfrak{q}_1}$, 
$\Frob_{\mathfrak{q}_2}$, $\Frob_{\mathfrak{q}_3}$ is the transposition fixing
$K$. Say $\Frob_{\mathfrak{q}_1}$ fixes $K$. Let $\mathfrak{p}_1, 
\mathfrak{p}_2, \mathfrak{p}_3\in I_K$ be the primes (not 
distinct) lying under $\mathfrak{q}_1$, $\mathfrak{q}_2$ and $\mathfrak{q}_3$.
Then 
$\deg(K_{\mathfrak{p}_1}/\mathbb{Q}_p) = 1$, whereas
$\deg(K_{\mathfrak{p}_i}/\mathbb{Q}_p) = 2$ for $i=2,3$. Hence 
$p$ splits as $p = \mathfrak{p}_1 \mathfrak{p}_2$ in $K/\mathbb{Q}$.
Since $\deg(L/K')=3$ is odd and
 $N \mathfrak{q}_i = p^2$ is an even power of $p$, we can see that
$p$ cannot split in $K'/\mathbb{Q}$. 

Finally, consider $\Frob_p = C_3$. Then $p$ splits as
$p = \mathfrak{q}_1 \mathfrak{q}_2$ in $L/\mathbb{Q}$.
Since $\deg(L/K')$ and $\deg(K/\mathbb{Q})$ are both odd, it follows
that $p$ splits in $K'/\mathbb{Q}$ but not in $K/\mathbb{Q}$.
\end{proof}

\begin{lem}\label{lem:taube}
 Let $K/\mathbb{Q}$ be an extension of $\mathbb{Q}$ of degree $3$.
Let $\alpha$ be a positive real number.
Let \[
S_{\alpha}(X) = \sum_{n\leq X} 2^{\alpha \omega_K(n) - \alpha \omega(n)} .\]
Then
\begin{equation}\label{eq:atya}\begin{aligned}
 S_{\alpha}(X) &\sim C_{K,\alpha} 
 X (\log X)^{(2^{2\alpha}-1)/3} \text{ if $K/\mathbb{Q}$ is Galois,}\\
 S_{\alpha}(X) &\sim C_{K,\alpha} X 
(\log X)^{\frac{1}{2} (2^{\alpha} - 1) +
\frac{1}{6} (2^{2 \alpha} - 1)} 
 \text{ if $K/\mathbb{Q}$ is not Galois,}
\end{aligned}\end{equation}
where $C_{K,\alpha}>0$ depends only on $K$ and $\alpha$, and
the dependence on $\alpha$ is continuous.
\end{lem}
\begin{proof}
Suppose $K/\mathbb{Q}$ is Galois. Then, for $\Re s>1$,
\[\sigma_{K/\mathbb{Q}}(s) = \prod_{\mathfrak{p}\in I_K} 
 \frac{1}{1-(N \mathfrak{p})^{-s}} = \prod_{\text{$p$ ramified}}
 \frac{1}{1- p^{-s}} \mathop{\prod_{\text{$p$ unsplit}}}_{
\text{\& unram.}} \frac{1}{1 - p^{-3 s}}
\prod_{\text{$p$ split}} \frac{1}{(1- p^{-s})^3} .\]
Hence
\begin{equation}\label{eq:comc}
 \prod_{\text{$p$ split}} (1 + \beta p^{-s}) = L_1(s) 
(\zeta_{K/\mathbb{Q}}(s))^{\beta/3} ,\end{equation}
where $L_1(s)$ is continuous and bounded on $\{s : \Re s > 1 - 1/4\}$.
Now
\[\begin{aligned}
2^{\alpha \omega_K(n) - \alpha \omega(n)} &=
\mathop{\prod_{p|n}}_{\text{$p$ split in $K/\mathbb{Q}$}} 2^{2 \alpha} 
= \mathop{\prod_{p|n}}_{\text{$p$ split in $K/\mathbb{Q}$}}
(1 + (2^{2 \alpha} - 1)) \\
&= \mathop{\sum_{a b = n}}_{p|a \Rightarrow \text{$p$ split}}
\prod_{p|a} (2^{2\alpha}-1) .\end{aligned}\]
Hence
\[\begin{aligned}
\sum_n 2^{\alpha \omega_K(n) - \alpha \omega(n)} n^{-s} &=
\left(\sum_n n^{-s}\right) \cdot \mathop{\sum_n}_{p|n\Rightarrow
\text{$p$ split}} \prod_{p|n} (2^{2\alpha}-1) n^{-s} \\
&= \zeta(s) \cdot \prod_{\text{$p$ split}} (1 + (2^{2 \alpha}-1) p^{-s}) .
\end{aligned}\]
By (\ref{eq:comc}) it follows that
\[\sum_n 2^{\alpha \omega_K(n) - \alpha \omega(n)} n^{-s} = L_1(s)
 (\zeta_{K/\mathbb{Q}}(s))^{(2^{2 \alpha}-1)/3} \zeta(s) .\]
Both $\zeta(s)$ and $\zeta_{K/\mathbb{Q}}$ have a pole of order $1$ at $s=1$.
By a Tauberian theorem (see, e.g., \cite{PT}, Main Th.) we can conclude that
\[\frac{1}{X} \sum_{n\leq X} 2^{\alpha \omega_K(n) - \alpha \omega(n)} 
\sim C_{K,\alpha} (\log X)^{(2^{2\alpha}-1)/3}\]
for some positive constant $C_{K,\alpha}>0$.

Now suppose that $K/\mathbb{Q}$ is not Galois. Denote the splitting
type of a prime $p$ in $K/\mathbb{Q}$ by $p = \mathfrak{p}_1 \mathfrak{p}_2$,
$p = \mathfrak{p}_1 \mathfrak{p}_2 \mathfrak{p}_3$, 
$p = \mathfrak{p}_1^2 \mathfrak{p}_2$, etc.
Then 
\[\zeta_{K/\mathbb{Q}}(s) = \prod_{\mathfrak{p}\in I_K} \frac{1}{(1 - (N \mathfrak{p})^{-s})} = L_2(s) \prod_{p = \mathfrak{p}_1 \mathfrak{p}_2} 
\frac{1}{(1-p^{-s})} 
\prod_{p = \mathfrak{p}_1 \mathfrak{p}_2 \mathfrak{p}_3}
\frac{1}{(1 - p^{-s})^{3}},\]
where $L_2(s)$ is continuous, non-zero
 and bounded on $\{s: \Re s > 1 - \frac{1}{4}\}$.
Let $L/\mathbb{Q}$ be the Galois closure of $K/\mathbb{Q}$. Let $K'/\mathbb{Q}$
be the quadratic subextension of $L/\mathbb{Q}$. Then we obtain from 
Lemma \ref{lem:patabran} that
\[
\prod_{p = \mathfrak{p}_1 \mathfrak{p}_2} 
\frac{1}{(1-p^{-s})} =
\prod_{\text{$p$ unsplit in $K'/\mathbb{Q}$}} \frac{1}{(1 - p^{-s})}
= L_3(s) \zeta(s) \zeta_{K'/\mathbb{Q}}^{- 1/2}(s), \]
where $L_3(s)$ is continuous
 and bounded on $\{s: \Re s > 1 - \frac{1}{4}\}$.

Now \[\begin{aligned}
2^{\alpha \omega_K(n) - \alpha \omega(n)} &=
 \mathop{\prod_{p|n}}_{p = \mathfrak{p}_1 \mathfrak{p}_2} 2^{\alpha}
 \mathop{\prod_{p|n}}_{p = \mathfrak{p}_1 \mathfrak{p}_2 \mathfrak{p}_3} 
2^{2 \alpha} \\ &=
 \mathop{\prod_{p|n}}_{p = \mathfrak{p}_1 \mathfrak{p}_2} (1+ (2^{\alpha}-1))
 \mathop{\prod_{p|n}}_{p = \mathfrak{p}_1 \mathfrak{p}_2 \mathfrak{p}_3} 
(1 + (2^{2 \alpha}-1)) \\ &=
\mathop{\mathop{\sum_{a b c = n}}_{p|a \Rightarrow p = \mathfrak{p}_1 \mathfrak{p}_2}}_{p|b \Rightarrow p = \mathfrak{p}_1 \mathfrak{p}_2 \mathfrak{p}_3}
\prod_{p|a} (2^{\alpha}-1) \prod_{p|b} (2^{2 \alpha}-1) .
\end{aligned}\]
Hence
\[\begin{aligned}
\sum_n 2^{\alpha \omega_K(n) - \alpha \omega(n)} n^{-s} &=
\left(\sum_n n^{-s}\right) \cdot
\mathop{\sum_n}_{p|n \Rightarrow 
 p = \mathfrak{p}_1 \mathfrak{p}_2} \prod_{p|n} (2^{\alpha}-1) n^{-s} \\
&\cdot \mathop{\sum_n}_{p|n \Rightarrow 
 p = \mathfrak{p}_1 \mathfrak{p}_2 \mathfrak{p}_3} 
\prod_{p|n} (2^{2 \alpha}-1) n^{-s} \\
&= \zeta(s) \prod_{p = \mathfrak{p}_1 \mathfrak{p}_2} ( 1 + (2^{\alpha}-1)
p^{-s})^{-1} \prod_{p = \mathfrak{p}_1 \mathfrak{p}_2 \mathfrak{p}_3}
 (1 + (2^{2 \alpha} - 1) p^{-s})^{-1} \\
&= L_4(s) \zeta(s) (\zeta_{K/\mathbb{Q}}(s))^{(2^{2 \alpha} - 1)/3}
(\zeta(s) \zeta_{K'/\mathbb{Q}}^{-1/2}(s))^{(2^{\alpha} - 1) - (2^{2 \alpha} - 1)/3}.\end{aligned}\]
Since $\zeta(s)$, $\zeta_{K/\mathbb{Q}}$ and $\zeta_{K'/\mathbb{Q}}$ each
have a pole of order $1$ at $s=1$, we can apply a Tauberian theorem
as before, obtaining
\[\frac{1}{X} \sum_{n\leq X} 2^{\alpha \omega_K(n) - \alpha \omega(n)} 
\sim C_{K,\alpha} (\log X)^{\frac{1}{2} (2^{\alpha} - 1) +
\frac{1}{6} (2^{2 \alpha} - 1)} .\]
\end{proof}
\subsection{The square-free sieve for homogeneous quartics}
We need the following simple lemma.
\begin{lem}\label{lem:facil}
 Let $f\in \mathbb{Z}\lbrack x,z\rbrack$ be a homogeneous
polynomial. Then there is a constant $C_f$ such that the following holds.
Let $N$ be a positive integer larger than $C_f$. Let $p$ be a 
prime larger than $N$. Then there are at most $12 \deg(f)$ pairs $(x,y)\in 
\mathbb{Z}^2$, $|x|, |z|\leq N$, $\gcd(x,z) = 1$, such that
\begin{equation}\label{eq:maud} p^2 | f(x,z) .\end{equation}
\end{lem}
\begin{proof}
If $N$ is large enough, then $p$ does not divide the discriminant of $f$.
Hence
\begin{equation}\label{eq:bablon}
f(r,1) \equiv 0 \mo p^2\end{equation}
has at most $\deg(f)$ solutions in $\mathbb{Z}/p^2$. If $N$ is large enough for
$p^2$ not to divide the leading coefficients of $f$, then $(x,z) = (1,0)$
does not satisfy (\ref{eq:maud}). Therefore, any solution (x,z) to 
(\ref{eq:maud}) gives us a solution $r = x/z$ to (\ref{eq:bablon}). 
We can focus on solutions $(x,y)\in \mathbb{Z}^2$ to (\ref{eq:maud})
with $x$, $y$ non-negative, as we need only flip signs to repeat the 
procedure for the other quadrants.

Suppose we have two solutions $(x_0,z_0), (x_1, z_1)\in \mathbb{Z}^2$,
$0\leq |x_0|, |x_1|, |z_0|, |z_1|\leq N$, $\gcd(x_0,z_0) = \gcd(x_1,z_1) =1$,
such that \[x_0/z_0 \equiv r \equiv x_1/z_1 \mo p^2.\] Then
\[x_0 z_1 - x_1 z_0 \equiv 0 \mo p^2 .\]
Since $0\leq x_j, z_j\leq N$ and $p>N$, we have that 
\[-p^2 < x_0 z_1 - x_1 z_0 < p^2,\]
and thus $x_0 z_1 - x_1 z_0$ must be zero. Hence $x_0/z_0 = x_1/z_1$.
Since $\gcd(x_0,z_0) = \gcd(x_1,z_1) = 1$ and $\sgn(x_0) = \sgn(x_1)$,
it follows that $(x_0,z_0) = (x_1,z_1)$.
\end{proof} 
\begin{Rem} It was pointed out by Ramsay \cite{Ra} that
an idea similar to that in Lemma \ref{lem:facil} suffices to improve
Greaves's bound for homogeneous sextics \cite{Gr} from
$\delta(N) = N^2 (\log N)^{-1/3}$ to $\delta(N) = N^2 (\log N)^{1/2}$.
\end{Rem}
\begin{prop}\label{prop:fourcoffees}
Let $f\in \mathbb{Z}\lbrack x,z\rbrack$ be a homogeneous irreducible
polynomial of degree $4$. Let
\[\delta(N) = \{x,z\in \mathbb{Z}^2 : |x|, |z|\leq N, \gcd(x,z)=1,
 \exists p>N \text{ s.t. } p^2 | f(x,y)\} .\]
Then
\[\delta(N)\ll N^{4/3} (\log N)^A,\]
where A and the implied constant depend only on $f$.
\end{prop}
\begin{proof}
Write $A = \max_{|x|, |z|\leq N} f(x,z)$. Clearly $A\ll N^4$. We can
write 
\[\begin{aligned}\delta(N) &\leq \sum_{0 < |d|\leq M} 
\#\{x,y,z\in \mathbb{Z}^3, |x|, |z|\leq N, \gcd(x,z)=1 :
d y^2 = f(x,z) \} \\
&+\sum_{N<p\leq \sqrt{A/M}} \# \{x,z\in \mathbb{Z}^2, |x|, |z|\leq N,
\gcd(x,z) = 1: p^2 | f(x,z) \} .\end{aligned}\]
Let $M\leq N^3$. By Lemma \ref{lem:facil},
\[\sum_{N<p\leq \sqrt{A/M}} \# \{x,z\in \mathbb{Z}^2, |x|, |z|\leq N,
\gcd(x,z) = 1: p^2 | f(x,z) \} \ll \frac{1}{\log N} \sqrt{N^{4 - \beta}},\]
where $\beta = (\log M)/(\log N)$. It remains to estimate
\[\sum_{0<|d|\leq M} S(d),\]
where we write
\[S(d) = \#\{x,y,z\in \mathbb{Z}^3, |x|, |z|\leq N, \gcd(x,z)=1 :
d y^2 = f(x,z) \}. \]

Let $C_{f,1}$, $C_{f,2}$, $C_{f,3}$, $C_{f,4}$ be as in Proposition 
\ref{prop:clinton}. Let $K$, $C_f$, $\omega$ and $\omega_K$ be as in
Lemma \ref{lem:notse}. Write $C_{f,5}$ for $C_f$.

By Proposition \ref{prop:clinton},
\[
\sum_{0<|d|< C_{f,4}} S(d) \ll \left(1 + 2 C_{f,3} \sqrt{\frac{7}{2} \log N +
 C_{f,1}}\right)^{C_1} \ll (\log N)^{C_2},\]
where $C_1 = \max_{0<d<C_{f,4}} \rnk(E_d)$, $C_2$ and the implied
constant depend only on $f$. Let $\epsilon$ be a small positive real
number. By Proposition \ref{prop:clinton} and Lemma \ref{lem:notse},
\[\begin{aligned}
\sum_{C_{f,4}\leq |d|<N^{\epsilon}} S(d) &\ll
\sum_{C_{f,4}\leq |d|<N^{\epsilon}} \left(1 + 2 \sqrt{\frac{7}{2} \log N +
C_{f,1}}\right)^{\rnk(E_d)} \\
&\ll \sum_{C_{f,4}\leq |d|<N^{\epsilon}} \left(1 + 2 \sqrt{\frac{7}{2} \log N +
C_{f,1}}\right)^{C_{f,5} + \omega_K(d) - \omega(d)}.\end{aligned}\]
We have the following crude bounds:
\begin{equation}\label{eq:cru}
\omega(d) \leq \frac{\log |d|}{\log \log |d|},\;\;
\omega_K(d)\leq 3 \omega(d) .\end{equation}
Hence
\[\begin{aligned}
\sum_{C_{f,4}\leq |d|<N^{\epsilon}} 
S(d) &\ll
\sum_{C_{f,4}\leq d<N^{\epsilon}} (\log N)^{C_{f,5} + 2 \log d/\log \log d}
\\ &\leq
N^{\epsilon} (\log N)^{C_1} (\log N)^{2 \epsilon \log N/\log \log N} \\
&\leq (\log N)^{C_1} N^{3 \epsilon} ,\end{aligned}\]
where $C$ depends only on $f$ and $\epsilon$.
For any $d$ with $|d|>N^{\epsilon}$, 
 Proposition \ref{prop:clinton} and Lemma \ref{lem:notse} give us
\[\begin{aligned}
S(d) 
&\ll \left(1 + 2 \sqrt{(\frac{7}{2} \log N + C_{f,1})/(\frac{1}{8} 
 \epsilon \log N + C_{f,2})}\right)^{\rnk(E_d)} \\
&\ll (12 \epsilon^{-1/2})^{C_{f,5} + \omega_K(d) - \omega(d)} 
\leq 2^{C_2 \omega_K(d) - C_2 \omega_K(d)} ,\end{aligned}\]
where $C_2$ depends only on $f$ and $\epsilon$. By Lemma \ref{lem:taube}
we can conclude that 
\[\begin{aligned}
\sum_{N^{\epsilon}<|d|\leq M} S(d) &\ll \sum_{d=1}^M 2^{C_2
\omega_K(d) - C_2 \omega_K(d)} \\
&\ll C_3 M (\log N)^{C_4},\end{aligned}\]
where $C_3$ and $C_4$ depend only on $f$ and $\epsilon$. Set $M = N^{4/3}$,
$\epsilon = 1/4$.
\end{proof}

\subsection{Homogeneous cubics}

\begin{prop}\label{prop:aulait}
Let $f\in \mathbb{Z}\lbrack x,z\rbrack$ be a homogeneous irreducible
polynomial of degree $3$. Let
\[\delta(N) = \{x,z\in \mathbb{Z}^2 : |x|, |z|\leq N, \gcd(x,z)=1,
 \exists p>N \text{ s.t. } p^2 | f(x,y)\} .\]
Then
\[\delta(N)\ll N^{4/3} (\log N)^A,\]
where A and the implied constant depend only on $f$.
\end{prop}
\begin{proof}
Write $A = \max_{|x|, |z|\leq N} f(x,z)$. Clearly $A\ll N^4$. We can
write 
\[\begin{aligned}\delta(N) &\leq 
\sum_{0 < |d|\leq M} 
\#\{x,y,z\in \mathbb{Z}^3, |x|, |z|\leq N, \gcd(x,z)=1 :
d y^2 = f(x,z) \} \\
&+\sum_{N<p\leq \sqrt{A/M}} \# \{x,z\in \mathbb{Z}^2, |x|, |z|\leq N,
\gcd(x,z) = 1: p^2 | f(x,z) \} .\end{aligned}\]
Let $M\leq N^2$. 
By Lemma \ref{lem:facil}, the second term on the right is
$O(N^{2-\beta/2}/\log N)$. Now notice that any point
$(x,y,z)\in \mathbb{Z}^3$ on $d y^2 = f(x,z)$ gives us a rational
point $(x',y') = (x/z,y/z^2)$ on
\begin{equation}\label{eq:tineq}
d' {y'}^2 = f(x',1),\end{equation}
where $d' = d z$. Moreover, a rational point on (\ref{eq:tineq})
can arise from at most one point $(x,y,z)\in \mathbb{Z}^3$, $\gcd(x,z)=1$,
in the given fashion. 

If $d\leq M$, then $|d'|=|d z|\leq M N$. The height $h_x(P)$
of the point $P=(x/z,y/z^2)$ is at most $N$. It follows by Lemma 
\ref{lem:soide} that $\hat{h}(P)\leq N + C_f$, where $C_f$ is a
constant depending only on $f$. By Corollaries \ref{cor:nosilv} and
\ref{cor:magd}, there are at most
\[O((1 + 2 \sqrt{(\log N + C_f')/(\log |d| + C_f)})^{\rnk(E_d)})\]
rational points $P$ of height $\hat{h}(P) \leq N + C_f$.
We proceed as in Proposition \ref{prop:fourcoffees}, and obtain that
\[\sum_{0 < |d|\leq M} 
\#\{x,y,z\in \mathbb{Z}^3, |x|, |z|\leq N, \gcd(x,z)=1 :
d y^2 = f(x,z) \}\]
is at most $O(M N (\log N))^A$. Set $\beta = 1/3$.
\end{proof}
\subsection{Homogeneous quintics}
We extract the following result from \cite{Gr}. 
\begin{lem}\label{lem:tumba}
Let $f\in \mathbb{Z}\lbrack x,y\rbrack$ be a homogeneous irreducible
polynomial of degree at most $5$. For all $M<N^{\deg f}$, $\epsilon>0$,
\begin{equation}\label{eq:sivest}
\sum_{d=1}^M \#\{x,y,z\in \mathbb{Z}^3, |x|, |z|\leq N, \gcd(x,z)=1 :
d y^2 = f(x,z) \} \ll N^{(18-\frac{1}{2} \beta^2)/(10-\beta) + \epsilon},
\end{equation}
where $\beta = (\log M)/(\log N)$. The implied constant depends only on $f$ 
and $\epsilon$.
\end{lem}
\begin{proof}
By \cite{Gr}, Lemmas 5 and 6, where the parameters $d$ and $z$ (in
the notation of \cite{Gr}, not ours) are set to the values $d=1$ and
 $z = N^{(1-\beta/2)/(5/2 - \beta/4)}$.
\end{proof}
\begin{prop}\label{prop:jocha}
Let $f\in \mathbb{Z}\lbrack x,z\rbrack$ be a homogeneous irreducible 
polynomial of degree $5$. Let
\[\delta(N) = \{x,z\in \mathbb{Z}^2 : |x|, |z|\leq N, \gcd(x,z)=1,
 \exists p>N \text{ s.t. } p^2 | P(x,y)\} .\]
Then, for any $\epsilon>0$,
\[\delta(N) \ll 
N^{(5+\sqrt{113})/8 + \epsilon}\]
where the implied constant depends only on $f$ and $\epsilon$.
\end{prop}
\begin{proof}
Let $A = \max_{|x|, |z|\leq N} f(x,z)$. Clearly $A\ll N^{\deg(f)}$.
We can write
\[\begin{aligned}
\delta(N) &\leq \sum_{0<|d|\leq M} 
\#\{x,y,z\in \mathbb{Z}^3, |x|, |z|\leq N, \gcd(x,z)=1 :
d y^2 = f(x,z) \} \\
&+\sum_{N<p\leq \sqrt{A/M}} \# \{x,z\in \mathbb{Z}^2, |x|, |z|\leq N,
\gcd(x,z) = 1: p^2 | f(x,z) \} .\end{aligned}\]
By Lemmas \ref{lem:tumba} and \ref{lem:facil}, 
\[\delta(N) \ll N^{(18 - \frac{1}{2} \beta^2)/(10 - \beta) + \epsilon} +
\frac{1}{\log N} \sqrt{N^{\deg(f) - \beta}},\]
where $\beta = (\log M)/(\log N)$. Set $\beta = (15 - \sqrt{113})/4$.
\end{proof}

\subsection{Quasiorthogonality, kissing numbers and cubics}
\begin{lem}\label{lem:agamen}
Let $f\in \mathbb{Z}\lbrack x\rbrack$ be a cubic polynomial of non-zero
discriminant. Let $d$ be a square-free integer. Then, for any two 
distinct
integer points $P=(x,y)\in \mathbb{Z}^2$, $P'=(x',y')\in \mathbb{Z}^2$
on the elliptic curve
\[E_d : d y^2 = f(x),\]
we have
\[\hat{h}(P+P') \leq 3 \max(\hat{h}(P),\hat{h}(P')) + C_f,\]
where $C_f$ is a constant depending only on $f$.
\end{lem}
\begin{proof}
Write $f(x) = a_3 x^3 + a_2 x^2 + a_1 x + a_0$. Let $P + P' = (x'',y'')$.
By the group law,
\[\begin{aligned}
x'' &= \frac{d (y_2 - y_1)^2}{a_3 (x_2 - x_1)^2} - \frac{a_2}{a_3} - x_1 - x_2
\\
&= \frac{d (y_2 - y_1)^2 - a_2 (x_2 - x_1)^2 - a_3 (x_2 - x_1)^2 (x_1 + x_2)}{
a_3 (x_2 - x_1)^2} .\end{aligned}\]
Clearly $|a_3 (x_2 - x_1)^2|\leq 4 |a_3| \max(|x_1|^2,|x_2|^2)$. Now
\[|d (y_2 - y_1)^2| \leq 4 |d| \max(y_1^2, y_2^2) = 4 \max(|f(x_1)|,|f(x_2)|) .\]
Hence
\[|d (y_2 - y_1)^2 - a_2 (x_2 - x_1)^2 - a_3 (x_2 - x_1)^2 (x_1 + x_2)
 \leq A \max( |x|^3, |x'|^3) ,\]
where $A$ is a constant depending only on $f$. Therefore
\[\begin{aligned}
h_x(P) &= \log(\max(|\num(x'')|,|\den(x'')|)) \\
 &\leq 3 \max(\log |x|, \log |x'|) + \log A \\
 &\leq 3 \max(h_x(P),h_x(P')) + \log A .\end{aligned}\]
By Lemma \ref{lem:soide}, the difference $|\hat{h}-h_x|$ 
is bounded by a constant
independent of $d$. The statement follows immediately.
\end{proof}
Consider the elliptic curve
\[E_d:d y^2 = f(x) .\]
There is a $\mathbb{Z}$-linear map from $E_d(\mathbb{Q})$ to
$\mathbb{R}^{\rnk(E_d)}$ taking the canonical height to the square of
the Euclidean norm. In other words, any
 given integer point $P=(x,y)\in E_d$ will be taken
to a point $L(P)\in \mathbb{R}^{\rnk(E_d)}$ whose Euclidean norm $|L(P)|$
satisfies
\[|L(P)|^2 = \hat{h}(P) = \log x + O(1) ,\]
where the implied constant depends only on $f$. In particular, the set of all
integer points $P=(x,y)\in E_d$ with
\begin{equation}\label{eq:perdie}
N^{1-\epsilon}\leq x\leq N\end{equation}
will be taken to a set of points $L(P)$ in $\mathbb{R}^{\rnk(E_d)}$ with
\[(1-\epsilon) \log N + O(1) \leq |L(P)|^2 \leq \log N + O(1) .\]
Let $P, P'\in E_d$ be integer points satisfying (\ref{eq:perdie}).
Assume $L(P)\ne L(P')$. By Lemma \ref{lem:agamen},
\[|L(P) + L(P')|^2 = |L(P + P')|^2 \leq 3 \max(|L(P)|^2,|L(P')|^2) + O(1).\]
Therefore, the inner product $L(P) \cdot L(P')$ satisfies
\[\begin{aligned} L(P) \cdot L(P') &= \frac{1}{2} (|L(P) + L(P')|^2 -
(|L(P)|^2 + |L(P')|^2)) \\
&\leq \frac{1}{2} (3 \max(|L(P)|^2,|L(P')|^2) + O(1) - (|L(P)|^2 +
|L(P')|^2))\\
&\leq \frac{1}{2} ((1+\epsilon) \log(N) + O(1)) \\
&\leq \frac{1}{2} \frac{(1+\epsilon) + O((\log N)^{-1})}{(1-\epsilon)^2}
 |L(P)| |L(P')| .\end{aligned}\]
We have proven
\begin{lem} 
Let $f\in \mathbb{Z}\lbrack x\rbrack$ be a cubic polynomial
of non-zero discriminant. Let $d$ be a square-free integer. Consider the
elliptic curve
\[E_d : d y^2 = f(x) .\]
Let $S$ be the set
\[\{(x,y)\in \mathbb{Z}^2 : N^{1-\epsilon} \leq |x|\leq N, d y^2 = f(x)\} .\]
Let $L$ be a linear map taking $E(\mathbb{Q})$ to $\mathbb{R}^{\rnk(E_d)}$
and the canonical height $\hat{h}$ to the square of the Euclidean norm.
Then, for any distinct points $P, P' \in L(S) \subset \mathbb{R}^{\rnk(E_d)}$
with the angle $\theta$ between $P$ and $P'$ is at least
\[\arccos\left(\frac{1}{2} \frac{(1 + \epsilon) + O((\log N)^{-1})}{(1-\epsilon)^2} \right) 
= 60^\circ + O(\epsilon + (\log N)^{-1}) ,\]
where the implied constant depends only on $f$.
\end{lem}

Let $A(\theta,n)$ be the maximal number of points that
can be arranged in $\mathbb{R}^n$ with angular separation no smaller
than $\theta$. Kabatiansky and Levenshtein (\cite{KL}; vd. also
\cite{CS}, (9.6)) show that, for $n$ large enough,
\[\frac{1}{n} \log_2 A(n,\theta) \leq
\frac{1 + \sin \theta}{2 \sin \theta} \log_2 
 \frac{1 + \sin \theta}{2 \sin \theta} -
\frac{1 - \sin \theta}{2 \sin \theta} \log_2 
 \frac{1 - \sin \theta}{2 \sin \theta} .\]
Thus we obtain
\begin{cor}\label{cor:coeur}
Let $f\in \mathbb{Z}\lbrack x\rbrack$ be a cubic polynomial
of non-zero discriminant. Let $d$ be a square-free integer. Consider the
elliptic curve
\[E_d : d y^2 = f(x) .\]
Let $S$ be the set
\[\{(x,y)\in \mathbb{Z}^2 : N^{1-\epsilon} \leq |x|\leq N, d y^2 = f(x)\} .\]
Then
\[\# S \ll 2^{(\alpha + O(\epsilon + (\log N)^{-1})) \rnk(E_d)},\]
where 
\[\alpha = \frac{2 + \sqrt{3}}{2 \sqrt{3}} \log_2
\frac{2 + \sqrt{3}}{2 \sqrt{3}} + \frac{2 - \sqrt{3}}{2 \sqrt{3}} \log_2
\frac{2 - \sqrt{3}}{2 \sqrt{3}}\]
and the implied constants depend only on $f$.
\end{cor}
Notice that we are using the fact that the size of the torsion group
is bounded.

\begin{prop}\label{prop:bethoo}
Let $f\in \mathbb{Z}\lbrack x\rbrack$ be an irreducible cubic polynomial.
Let
\[\delta(N) = \{1\leq x \leq N : \exists p > N^{1/2} \text{ s.t. }
p^2 | f(x) \} .\]
Then
\begin{equation}\label{eq:gerg}\delta(N) \ll N (\log N)^{-\beta},\end{equation}
where
\[\beta = 
-((2^{2 \alpha}-1)/9-2/3) = 0.5839\dotsc\]
if the discriminant of $f$ is a square, 
\[\beta = 
-\left(\frac{1}{6} (2^{\alpha} - 1) +
\frac{1}{18} (2^{2 \alpha} - 1) - 2/3\right) = 0.5718\dotsc\]
if the discriminant of $f$ is not a square, and
\[\alpha = \frac{2 + \sqrt{3}}{2 \sqrt{3}} \log_2
\frac{2 + \sqrt{3}}{2 \sqrt{3}} + \frac{2 - \sqrt{3}}{2 \sqrt{3}} \log_2
\frac{2 - \sqrt{3}}{2 \sqrt{3}} = 0.4014\dotsc .\]
The implied constant in (\ref{eq:gerg}) depends only on $f$.
\end{prop}
\begin{proof}
Let $A = \max_{1\leq x\leq N} f(x)$. Clearly $A\ll N^3$. We can write
\[\begin{aligned}
\delta(N) &\leq \sum_{N^{1/2}<p<\sqrt{A/M}} \# \{1\leq x\leq N : p^2 | f(x)\}\\
&+ \{1\leq x\leq N^{1-\epsilon} : \exists p>N^{1/2} \text{ s.t. }
p^2 | f(x) \} \\
&+ \sum_{1\leq |d|\leq M} \#\{x,y\in \mathbb{Z}^2 : N^{1-\epsilon}\leq x\leq N,
d y^2 = f(x) \} .\end{aligned}\]
Let $M\leq N^2$. Then
the first term is at most
\[\sum_{N^{1/2} < p < \sqrt{A/M}} 3\ll \frac{3 \sqrt{A/M}}{\log \sqrt{A/M}} \ll \frac{N^{3/2} M^{-1/2}}{\log N} .\]
The second term is clearly no greater than $N^{1-\epsilon}$.
It remains to bound
\[\sum_{1\leq |d|\leq M} B(d),\]
where \[B(d) = 
\#\{x,y\in \mathbb{Z}^2 : N^{1-\epsilon}\leq x\leq N,
d y^2 = f(x) \} .\]
By Lemma \ref{lem:notse} and Corollary \ref{cor:coeur}
\[B(d) \ll
2^{(\alpha + O(\epsilon + (\log N)^{-1})) (\omega_K(d) - \omega(d))},\]
where $K$ is as in Lemma \ref{lem:notse}
and $\alpha$ is as in Corollary \ref{cor:coeur}. Thanks to (\ref{eq:cru}),
we can omit the term $O((\log N)^{-1})$ from the exponent.
Hence it remains to estimate
\[S(M) =
\sum_{1\leq d\leq M} 2^{(\alpha + O(\epsilon)) (\omega_K(d) - \omega(d))} .\]
By Lemma \ref{lem:taube},
\[\begin{aligned}
 S(M) &\ll
 M (\log M)^{(2^{2(\alpha+\epsilon)}-1)/3} \text{ if $K/\mathbb{Q}$ is 
Galois,}\\
 S(M) &\ll
 M (\log M)^{\frac{1}{2} (2^{\alpha+\epsilon} - 1) +
\frac{1}{6} (2^{2 (\alpha+\epsilon)} - 1)} 
 \text{ if $K/\mathbb{Q}$ is not Galois.}\end{aligned}\]
Let $\epsilon = (\log \log M)^{-1}$. Note that $K/\mathbb{Q}$ is
Galois if and only if the discriminant of $f$ is a square.
Then
\[\begin{aligned}
 S(M) &\ll
 M (\log M)^{(2^{2 \alpha}-1)/3} \text{ if $\Disc(f)$ is a square,}\\
 S(M) &\ll
 M (\log M)^{\frac{1}{2} (2^{\alpha} - 1) +
\frac{1}{6} (2^{2 \alpha} - 1)}
 \text{ if $\Disc(f)$ is not a square.}\end{aligned}\]
Set
\[\begin{aligned}
M &= N (\log N)^{-2 (2^{2 \alpha}-1)/9-2/3} \text{ if $\Disc(f)$ is a square,}\\
M &= N (\log N)^{-\frac{1}{3} (2^{\alpha} - 1) -
\frac{1}{9} (2^{2 \alpha} - 1) - 2/3}
 \text{ if $\Disc(f)$ is not a square.}\end{aligned}\]
Hence
\[\begin{aligned}
S(M) &= N (\log N)^{(2^{2 \alpha}-1)/9-2/3} \text{ if $\Disc(f)$ is a square,}\\
S(M) &= N (\log N)^{\frac{1}{6} (2^{\alpha} - 1) +
\frac{1}{18} (2^{2 \alpha} - 1) - 2/3}
 \text{ if $\Disc(f)$ is not a square.}\end{aligned}\]
The statement follows.
\end{proof}
\section{Square-free integers}\label{sec:cuanumer}
In Chapter \ref{chap:ell}, we had the chance to employ the
framework from section \ref{sec:cuacriba} in its full generality.
We will now give a simpler and more traditional application.

\begin{thm}\label{thm:croiss}
Let $f\in \mathbb{Z}\lbrack x\rbrack$ be an irreducible polynomial
of degree $3$. Then the number of positive integers $x\leq N$ for
which $f(x)$ is square-free is given by
\begin{equation}\label{eq:atrya}
N \prod_p \left(1 - \frac{\ell(p^2)}{p^2}\right) + O(N (\log N)^{-\beta}),
\end{equation}
where
\[\begin{aligned}
\beta &= \begin{cases} 0.5839\dotsc
&\text{ if the discriminant of $f$ is a square,}\\ 
0.5718\dotsc
&\text{ if the discriminant of $f$ is not a square,}\end{cases}\\
\ell(m) &= \#\{x\in \mathbb{Z}/m : f(x) \equiv 0 \mo m\}.\end{aligned}\]
Note that $\epsilon$ is an arbitrarily small positive number, and
that the implied constant depends in (\ref{eq:atrya}) depends only on
$f$ and $\epsilon$.
\end{thm}
\begin{proof}
Define the terms needed for Lemma \ref{lem:wicked} as follows. Let $K = \mathbb{Q}$. Let $\gamma(d) = d \rad(d)$. Let $S_a = \{\emptyset\}$
for every $a\in \mathbb{Z}^+$; let $\phi_{a_1,a_2}:S_{a_2}\to S_{a_1}$
be the map taking $\emptyset$ to $\emptyset$. Define
\[\begin{aligned}
f_a(\emptyset) &= \begin{cases} 1 & \text{if $a=1$,}\\
0 &\text{otherwise,}\end{cases}\\
g_a(\emptyset) &= \mathop{\sum_{1\leq x\leq N}}_{\sq(f(x))=a} 1 .\end{aligned}\]
Then the cardinality of $\{1\leq x\leq N : \text{$f(x)$ square-free}\}$
equals
\[\sum_{a\in \mathbb{Z}^+} \sum_{r\in S_a} f_a(r) g_a(r),\]
which is the expression on the left side of the inequality
(\ref{eq:maravel}). It remains to estimate the right side.

Write $f(a)$, $g(a)$ instead of $f_a(\emptyset)$, $g_a(\emptyset)$
for the sake of brevity. Then
\[\begin{aligned}
\sum_{\gamma(d)\leq M} \sum_{r\in S_d} \left(\sum_{d'|d} \mu(d') f(d/d')
\right) t_d(r) &= \sum_{\gamma(d) \leq M} \mu(d) t_d(r) 
= \sum_{\gamma(d)\leq M} \mu(d) \mathop{\sum_{1\leq x\leq N}}_{
d|\sq(f(x))} 1 \\
&= \sum_{d^2\leq M} \mu(d)
\mathop{\sum_{1\leq x\leq N}}_{d^2|f(x)} 1 .\end{aligned}\]
Assume $M\leq N$. Then
\[\begin{aligned}
\sum_{d^2\leq M} \mu(d)
 \mathop{\sum_{1\leq x\leq N}}_{d^2 | f(x)} 1 &=
\mathop{\sum_{\text{$d$ square-free}}}_{d^2\leq M} \mu(d)
 \frac{N \ell(d^2)}{d^2} + O(M^{1/2}) \\
&= \sum_{d} \mu(d) \frac{N \ell(d^2)}{d^2}
 - \sum_{d^2>M} \mu(d)
\frac{N \ell(d^2)}{d^2} + O(M^{1/2}) \\
&= N \prod_p \left(1 - \frac{\ell(p^2)}{p^2}\right) +
O\left(N \sum_{d^2>M} \frac{\tau_3(d)}{d^2} + M^{1/2}\right)\\
&= N \prod_p \left(1 - \frac{\ell(p^2)}{p^2}\right) +
O(N M^{-1/2} (\log N)^3) .\end{aligned}\]

Assume $M\leq \sqrt{N}$. We may now bound the second term on the
right side of (\ref{eq:maravel}). By Lemmas 
\ref{lem:dontaskmemyadvice} and \ref{lem:brum},
\[\begin{aligned}
\sum_{M<\gamma(d)\leq M^2} \tau_3(d) s_d &= \sum_{M<\gamma(d)\leq M^2}
 \tau_3(d) \mathop{\sum_{1\leq x\leq N}}_{\gamma(d)|f(x)} 1\\
&\ll \sum_{M<\gamma(d)\leq M^2} \tau_3(d) \tau_3(\rad(d)) \frac{N}{\gamma(d)}\\
&\ll M^{-1/2} N (\log M)^{9^2+9^3-2} .\end{aligned}\]
The remaining term of (\ref{eq:maravel}) is
\[2 \mathop{\sum_p}_{p^2 > M} s_p = 2 \sum_{p>\sqrt{M}}
 \mathop{\sum_{1\leq x\leq N}}_{p^2|f(x)} 1 =
2 \sum_{\sqrt{M}<p\leq N^{1/2}} 
  \mathop{\sum_{1\leq x\leq N}}_{p^2|f(x)} 1 +
2 \sum_{p>N^{1/2}} 
  \mathop{\sum_{1\leq x\leq N}}_{p^2|f(x)} 1 .\]
By Lemma \ref{lem:brum},
\[\sum_{\sqrt{M}<p\leq N^{1/2}} 
  \mathop{\sum_{1\leq x\leq N}}_{p^2|f(x)} 1 \ll
\sum_{p\geq \sqrt{M}} \frac{N}{p^2} \ll M^{-1/2} N .\]

Hence we have
\[\begin{aligned}\#\{1\leq x\leq N: \text{$f(x)$ square-free}\} &=
N \prod_p \left(1 - \frac{\ell(p)}{p^2}\right) + 
2 \mathop{\sum_{p>N^{1/2}}}_{p^2|f(x)} 1 \\ &+ O(N M^{-1/2} (\log M)^{9^2 +
9^3 - 2}) .\end{aligned}\]
Set $M = N^{1/2}$. Notice that, for $N$ large enough,
 no more than three squares of primes
$p^2$, $p>N^{1/2}$, may divide $f(x)$ for any $1\leq x\leq N$.
Thus
\[\mathop{\sum_{p>N^{1/2}}}_{p^2|f(x)} 1 \ll
 \{1\leq x\leq N : \exists p>N^{1/2} \text{\,s.t.\,} p^2|f(x)\} .\]
By Proposition \ref{prop:bethoo}, the statement follows.
\end{proof}
\begin{thm}
Let $f\in \mathbb{Z}\lbrack x,y\rbrack$ be a homogeneous polynomial
of degree no greater than $6$. Then the number of integer pairs
$(x,y)\in \mathbb{Z}^2\cap \lbrack -N,N\rbrack^2$
 for which 
$f(x,y)$ is square-free is given by
\[4 N^2 \prod_p \left(1 - \frac{\ell_2(p^2)}{p^4}\right) +
\begin{cases} O(N (\log N)^{A_1}) &\text{if $\deg_{\irr}(f) = 1, 2$,}\\
O(N^{4/3} (\log N)^{A_2}) &\text{if $\deg_{\irr}(f) = 3, 4$,}\\
O(N^{(5+\sqrt{113})/8 + \epsilon}) 
&\text{if $\deg_{\irr}(f) = 5$,}\\
O(N^2 (\log N)^{-1/2}) &\text{if $\deg_{\irr}(f) = 6$,}\end{cases}\]
where $\epsilon$ is an arbitrarily small positive number, $A_1$ is
an absolute constant, $A_2$ depends only on $f$, the implied constant
depends only on $f$ and $\epsilon$, $\deg_{\irr}$ denotes the
degree of the irreducible factor of $f$ of largest degree, and
\[\ell_2(m) = \#\{(x,y)\in 
(\mathbb{Z}/m)^2 : f(x,y) \equiv 0 \mo m\}.\]
\end{thm}
\begin{proof}
Set $K$, $\gamma$, $S_a$, $\phi_{a_1,a_2}$ and $f_a$ as in the proof
of Theorem \ref{thm:croiss}. Let
\[g_a(\emptyset) = \mathop{\sum_{(x,y)\in \mathbb{Z}^2 \cap 
\lbrack - N, N\rbrack^2}}_{\sq(f(x))=a} 1 .\]
We proceed as in Theorem \ref{thm:croiss}. Let $M\leq N$. Then
\[\begin{aligned}
\sum_{d^2\leq M} \mu(d) \mathop{\sum_{(x,y)\in \mathbb{Z}^2 \cap 
\lbrack - N, N\rbrack^2}}_{d^2|f(x)} 1 &=
\sum_{d^2\leq M} \mu(d) \frac{4 N^2 \ell_2(d^2)}{d^4} + O(M^{1/2} N) \\
&= \sum_d \mu(d) \frac{4 N^2 \ell_2(d^2)}{d^4} - 
\sum_{d^2>M} \mu(d) \frac{4 N^2 \ell_2(d^2)}{d^4} + O(M^{1/2} N) \\
&= 4 N^2 \prod_p \left(1 - \frac{\ell_2(p^2)}{p^4}\right) +
O(N^2 M^{-1/2} (\log N)^3).\end{aligned}\]
Notice that the first equality is justified even for $M>N^{1/2}$,
as the solutions to $d^2|f(x)$ fall into lattices of index $d^2$ with
$d \mathbb{Z}^2$ as their pairwise intersection.
By Lemmas \ref{lem:ramsay} and \ref{lem:brum},
\[\begin{aligned}
\sum_{M<\gamma(d)\leq M^2} \tau_3(d) s_d &= 
\sum_{M<\gamma(d)\leq M^2} \tau_3(d) 
 \mathop{\sum_{(x,y)\in \mathbb{Z}^2 \cap \lbrack - N, N\rbrack^2}}_{
\gamma(d) | f(x)} 1\\
&\ll \sum_{M<\gamma(d)\leq M^2} \tau_3(d) \tau_{12}(\rad(d))
\frac{N^2}{\gamma(d)}\\
&\ll M^{-1/2} N^2 (\log M)^{A_1},\end{aligned}\]
where $A_1 = 36^2 + 36^3 - 2$. The remaining term is
\[2 \mathop{\sum_p}_{p^2 > M} s_p = \sum_{p>\sqrt{M}}
\mathop{\sum_{(x,y)\in \mathbb{Z}^2 \cap \lbrack - N, N\rbrack^2}}_{p^2|f(x)}
1 ,\]
which is at most a constant times
\[M^{-1/2} N^2  + 
\{x,z\in \mathbb{Z}^2 : |x|, |z|\leq N, \gcd(x,z)=1,
 \exists p>N \text{ s.t. } p^2 | f(x,y)\} .\]
Use Prop. \ref{prop:aulait} for $\deg_{\irr}(f)=3$,
Prop. \ref{prop:fourcoffees} for $\deg_{\irr}(f)=4$
and Prop. \ref{prop:jocha} for $\deg_{\irr}(f)=5$. Use the trivial bound for
$\deg_{\irr}(f)=1,2$, and the estimate in \cite{Gr}, Lemma 3, for
$\deg_{\irr}(f)=6$.
\end{proof}

%% file: appa.tex
\chapter{Addenda on the root number}
\section{Known instances of conjectures $\mathfrak{A}_i$ and $\mathfrak{B}_i$ over the rationals}\label{sec:fapp}
The quantitative versions of $\mathfrak{A}_i$ and $\mathfrak{B}_i$ were introduced in
subsections \ref{subs:conres} and \ref{subs:aver}.
As before, we denote by $\deg_{\irr} P$ the degree of the irreducible factor
of $P$ of highest degree.
\begin{prop}
Conjecture $\mathfrak{A}_1(\mathbb{Q},P,\delta(N))$ holds for
\begin{enumerate}
\item $\deg_{\irr} P = 1$, $\delta(N) = \sqrt{N}$,
\item $\deg_{\irr}  P = 2$, $\delta(N) = N^{2/3}$,
\item $\deg_{\irr}  P = 3$, $\delta(N) = N (\log N)^{-0.5839\dotsc}$ if the
discriminants of all irreducible factors of degree $3$ of $P$ are square,
\item $\deg_{\irr} P = 3$, $\delta(N) = N (\log N)^{-0.5718\dotsc}$, in
general.
\end{enumerate}
\end{prop}
\begin{proof}
The case $\deg_{\irr} P = 1$ is trivial.
The result for $\deg_{\irr} P = 2$ is due to Estermann (\cite{Es}).
See Chapter \ref{chap:square} for $\deg_{\irr} P = 3$. The best
previous bound for $\deg_{\irr} P = 3$, namely
$\delta(N) = N (\log N)^{-1/2}$, was due to Hooley
(\cite{Ho}, Ch. IV).
\end{proof}
\begin{prop}\label{prop:diezca}
Conjecture $\mathfrak{A}_2(\mathbb{Q},P,\delta(N))$ holds for
\begin{enumerate}
\item $\deg_{\irr}  P = 1$, $\delta(N) = 1$,
\item $\deg_{\irr}  P = 2$, $\delta(N) = N$,
\item $\deg_{\irr}  P = 3$, $\delta(N) = N^{3/2}/(\log N)$,
\item $\deg_{\irr}  P = 4$, $\delta(N) = N^{4/3} (\log N)^A$,
\item $\deg_{\irr}  P = 5$, $\delta(N) = N^{(5+ \sqrt{113})/8 + \epsilon}$,
\item $\deg_{\irr}  P = 6$, $\delta(N) = N^2/(\log N)^{1/2}$,
\end{enumerate}
where $\epsilon$ is an arbitrarily small positive integer, and $A$
and the implied constant depends only on $\epsilon$.
\end{prop}
\begin{proof}
The cases $\deg_{\irr} P = 1$ and $\deg_{\irr} P = 2$ are trivial.
See Chapter \ref{chap:square} for $3\leq \deg_{\irr} P \leq 6$.
The best previous bound for $\deg_{\irr} P = 3, 4, 5$ was
$N^2 (\log N)^{-1}$, due to Greaves \cite{Gr}.
While, in the cited work, Greaves gives the bound $N^2 (\log N)^{-1/3}$,
his methods suffice to obtain $N^2 (\log N)^{-1/2}$, as was remarked
by Ramsay (\cite{Ra}, 1991, unpublished; see reference in \cite{GM}).
\end{proof}
\begin{prop}\label{prop:prib}
Hypothesis $\mathfrak{B}_1(\mathbb{Q},P,\eta(N),\epsilon(N))$ holds for $\deg P=1$,
$\eta(N) = (\log N)^A$, $\epsilon(N) = C_1 e^{- C_2 (\log N)^{3/5}/(\log \log N)^{1/5}}$, where
$A$ is arbitrarily large and $C_1$, $C_2$ depend on $A$ and $P$.
\end{prop}
\begin{proof}
By Siegel-Walfisz (vd. \cite{Wa}, V \S 5 and V \S 7). (For an elementary proof of equivalence with the Prime Number Theorem,
see, e.g., \cite{A}.)
\end{proof}
\begin{prop}
Hypothesis $\mathfrak{B}_2(\mathbb{Q},P,\eta(N),\epsilon(N))$ holds for 
\begin{enumerate}
\item $\deg(P) = 1$, $\eta(N) = (\log N)^A$, 
 $\epsilon(N) = C_1 e^{- C_2 (\log N)^{3/5}/(\log \log N)^{1/5}}$, $A$ arbitrarily large, $C_1$, $C_2$ depending on $A$ and $P$,
\item $\deg(P) = 2$, $\eta(N) = (\log N)^A$,
 $\epsilon(N) = C_1 e^{- C_2 (\log N)^{3/5-\epsilon}}$, $A$ arbitrarily large, $\epsilon$ an
arbitrarily small positive number, $C_1$, $C_2$ depending on $A$, $P$ and $\epsilon$,
\item $\deg(P) = 3$, $P$ reducible, $\eta(N) = (\log N)^A$,
 $\epsilon(N) = C \frac{\log \log N}{\log N}$, $A$ arbitrarily large, $C$ depending on $A$ and $P$,
\item $\deg(P) = 3$, $P$ irreducible, $\eta(N) = (\log N)^A$,
$\epsilon(N) = C \frac{(\log \log N)^5 \log \log \log N}{\log N}$,
$A$ arbitrarily large, $C$ depending on $A$ and $P$.
\end{enumerate}
\end{prop}
\begin{proof}
The case $\deg P = 1$ follows immediately from Proposition \ref{prop:prib}. 
For $\deg P = 2, 3$,
see Chapter \ref{chap:par}. As was said before, the case $\deg P = 2$ is in essence well-known and classical. 
\end{proof}
\section{\!Reducing hypotheses on number fields to their rational analogues}\label{sec:app}
Given a number field $K$ and a polynomial $P(x) = a_n x^n + a_{n-1} x^{n-1}
+\dotsb + a_0\in \mathfrak{O}_K\lbrack x\rbrack$ (or a homogeneous
polynomial $P(x,y) = a_n x ^n + a_{n-1} x^{n-1} y + \dotsb + a_0 \in
\mathfrak{O}_K\lbrack x,y\rbrack$), we define
\[K_P = \mathbb{Q}\left(\frac{a_{n-1}}{a_n},\frac{a_{n-2}}{a_n},\dotsb,
\frac{a_0}{a_n}\right) .\]
\begin{lem}\label{lem:trommel} 
Let $K$ be a number field. Let $P\in \mathfrak{O}_K\lbrack x\rbrack$ be a monic, irreducible polynomial.
Suppose $K=K_P$.
Then there is a finite set $D$ of rational primes such that for every $x\in \mathbb{Z}$
and every rational prime $p$ not in $D$,
\begin{enumerate}
\item at most one prime ideal $\mathfrak{p}\in I_K$ lying over $p$ divides $P(x)$,
\item if some $\mathfrak{p}\in I_K$ lying over $p$ divides $P(x)$, then 
$\mathfrak{N}_{K/\mathbb{Q}} \mathfrak{p} = p$,
\item $\sum_{\mathfrak{p}\in I_K, \mathfrak{p}|p} v_{\mathfrak{p}}(P(x))
= v_p(\mathfrak{N}_{K/\mathbb{Q}} P(x)).$
\end{enumerate}
\end{lem}
\begin{proof}
Let $L/\mathbb{Q}$ be the Galois closure of $K/\mathbb{Q}$. Let $G = \Gal(L/\mathbb{Q})$,
$H = \Gal(L/K)$. Then for any ideal $\mathfrak{a}\in I_K$,
\[\mathfrak{N}_{K/\mathbb{Q}} \mathfrak{a} = \prod_{\sigma H} \sigma \mathfrak{a},\]
where the product is taken over all cosets $\sigma H \subset G$ of $H$. Let 
$\sigma$ be an element of $G$ not in $H$. By definition, $\sigma$ cannot leave $K$ fixed.
Since the ratios among the coefficients of $P$ generate $K_P=K$, 
$\sigma$ would leave $K$ fixed if $P_{\sigma}$ were a multiple of $P$.
Hence $P_{\sigma}$ is not a multiple of $P$.
Since
$P$ is irreducible, it follows that $P$ and $P_{\sigma}$ are coprime. Let
$D$ be the set of all rational primes lying under prime ideals dividing $\Disc(P,P_{\sigma})$ for some $\sigma \in G$ not in $H$.

Suppose there are two distinct prime ideals $\mathfrak{p}_1, 
\mathfrak{p}_2 \in I_K$ such that $\mathfrak{p}_1, \mathfrak{p}_2 |
P(x)$, $\mathfrak{p}_1, \mathfrak{p}_2 | p$, $p\notin D$. Then
$\mathfrak{p}_1' | \mathfrak{p}_1$, $\mathfrak{p}_2' | \mathfrak{p}_2$
for some prime ideals $\mathfrak{p}_1', \mathfrak{p}_2' \in I_L$.
There is a $\sigma\in G$ such that $\sigma \mathfrak{p}_1'=\mathfrak{p}_2'$.
Then $\mathfrak{p}_2'$ divides both $P$ and $P_\sigma$. Since
$\mathfrak{p}_1 \neq \mathfrak{p}_2$, $\sigma$ does not fix $K$. Hence
$\sigma \notin H$. Therefore $\mathfrak{p}_2'|\Disc(P,P_{\sigma})$,
and thus $\mathfrak{p}_2'$ must lie over a prime in $D$. Contradiction.
Hence (1) is proven.

Now take $\mathfrak{p}\in I_K$ lying over $p\notin D$. Assume
$\mathfrak{p} | P(x)$ for some $x\in \mathbb{Z}$. Obviously
\[\mathfrak{N}_{L/\mathbb{Q}} \mathfrak{p} = \prod_{\sigma \in G}
\sigma \mathfrak{p} = 
\left(\prod_{\sigma H} \sigma \mathfrak{p}\right)^{\deg L/K} .\]
Since $p\notin D$ and $\mathfrak{p}|P(x)$, we have 
$\gcd(\sigma \mathfrak{p},\sigma' \mathfrak{p})=1$ for
$\sigma$, $\sigma'$ with $\sigma H \notin \sigma' H$. Therefore
$\prod_{\sigma H} \sigma \mathfrak{p}$ divides $p$. Hence 
$\mathfrak{N}_{L/\mathbb{Q}} \mathfrak{p} | p^{\deg L/K}$. Since
$\mathfrak{N}_{L/\mathbb{Q}} \mathfrak{p} = 
(\mathfrak{N}_{K/\mathbb{Q}} \mathfrak{p})^{\deg L/K}$,
we have $\mathfrak{N}_{K/\mathbb{Q}} \mathfrak{p} | p$. Therefore
$\mathfrak{N}_{K/\mathbb{Q}} \mathfrak{p} = p$; this is (2).

Finally,
\[\begin{aligned}
v_p(\mathfrak{N}_{K/\mathbb{Q}} P(x)) &= 
v_p\left(\mathfrak{N}_{K/\mathbb{Q}} \left(
 \mathop{\prod_{\mathfrak{p}\in I_K}}_{\mathfrak{p}|p}
 \mathfrak{p}^{v_{\mathfrak{p}}(P(x))} \right) \right) =
v_p\left(\mathop{\prod_{\mathfrak{p}\in I_K}}_{\mathfrak{p}|p}
(\mathfrak{N}_{K/\mathbb{Q}} \mathfrak{p})^{v_{\mathfrak{p}}(P(x))}
\right)\\
&= v_p\left(\mathop{\prod_{\mathfrak{p}\in I_K}}_{\mathfrak{p}|p}
 p^{v_{\mathfrak{p}}(P(x))} \right)
= \mathop{\sum_{\mathfrak{p}\in I_K}}_{\mathfrak{p}|p}
v_{\mathfrak{p}}(P(x)) .\end{aligned}\]
\end{proof}
\begin{lem}\label{lem:trommel2} 
Let $K$ be a number field. Let $P\in \mathfrak{O}_K\lbrack x,y\rbrack$ be an 
irreducible polynomial. Suppose $K=K_P$.
Then there is a finite set $D$ of rational primes such that for 
all coprime $x, y\in \mathbb{Z}$
and every rational prime $p$ not in $D$,
\begin{enumerate}
\item at most one prime ideal $\mathfrak{p}\in I_K$ lying over $p$ divides $P(x,y)$,
\item if some $\mathfrak{p}\in I_K$ lying over $p$ divides $P(x,y)$, then 
$\mathfrak{N}_{K/\mathbb{Q}} \mathfrak{p} = p$,
\item $\sum_{\mathfrak{p}\in I_K, \mathfrak{p}|p} v_{\mathfrak{p}}(P(x,y))
= v_p(\mathfrak{N}_{K/\mathbb{Q}} P(x,y)).$
\end{enumerate}
\end{lem}
\begin{proof}
Same as that of Lemma \ref{lem:trommel}.
\end{proof}
\begin{prop}\label{prop:baram}
Let $K$ be a number field. Let $P\in \mathfrak{O}_K\lbrack x \rbrack$ be
a square-free, non-constant polynomial. Let $P = P_1 P_2 \dotsb P_k$, $P_i$ irreducible
in $\mathfrak{O}_K\lbrack x\rbrack$. Then
Conjecture $\mathfrak{A}_1(K,P,\delta(N))$ is equivalent to
Conjecture $\mathfrak{A}_1(\mathbb{Q},Q,\delta(N))$, where the
polynomial $Q(x)\in \mathbb{Z}\lbrack x\rbrack$ is defined as the
product of the irreducible factors of
$\mathfrak{N}_{K_{P_i}/\mathbb{Q}} (c_i P_i(x)) \in \mathbb{Z}\lbrack x\rbrack$, 
$i=1,\dotsb,k$,
where $c_1,\dotsc,c_k$ are constants in $\mathfrak{O}_K$.
\end{prop}
\begin{proof}
Since $\mathfrak{A}_1(K,P_1\cdot P_2,\delta(N))$ is equivalent to
$\mathfrak{A}_1(K,P_1,\delta(N)) \wedge \mathfrak{A}_1(K,P_2,\delta(N))$,
it is enough to prove the statement for $P$ irreducible. Choose
a non-zero $c\in \mathfrak{O}_K$ such that the leading coefficient of $c P$ lies
in $K_P$. Then all coefficients of $c P$ lie in $\mathfrak{O}_{K_P}$.
Since we can take $N^{1/2}$ to be larger than every prime divisor of $c$,
it follows that we can assume that $P$ has all its coefficients in 
$\mathfrak{O}_{K_P}$. 
Since we can also let $N^{1/2}$ be larger than all primes ramifying in
$K/K_P$, we can assume $K=K_P$.

Let
\[\begin{aligned}
S_1(N) &= \{1\leq x\leq N: \exists \mathfrak{p} \text{\;\st\;} \rho(\mathfrak{p})>N^{1/2},
\mathfrak{p}^2 |P(x)\}\\
S_2(N) &= \{1\leq x\leq N: \exists p \text{\;\st\;} p>N^{1/2}, 
p^2|\mathfrak{N}_{K/\mathbb{Q}} P(x)\}.
\end{aligned}\]
We recall that conjecture $\mathfrak{A}_1(K,P,\delta(N))$\; states that $\# S_1(N)\ll \delta(N)$, whereas
conjecture $\mathfrak{A}_1(\mathbb{Q},
\mathfrak{N}_{K/\mathbb{Q}} P, \delta(N))$
states that $\# S_2(N) \ll \delta(N)$. We can assume 
$N^{1/2}\geq \max_{p|D} p$, where $D$ is as in Lemma \ref{lem:trommel}.
Then, for every prime ideal $\mathfrak{p}\in I_K$ such that 
$\rho(\mathfrak{p})>N^{1/2}$, $\mathfrak{p}^2|P(x)$, Lemma \ref{lem:trommel}
implies that $\mathfrak{N}_{K/\mathbb{Q}} \mathfrak{p} = \rho(\mathfrak{p})
> N^{1/2}$. Obviously, if $\mathfrak{p}^2 |P(x)$, then
$(\mathfrak{N}_{K/\mathbb{Q}} \mathfrak{p})^2|\mathfrak{N}_{K/\mathbb{Q}}
P(x)$. Thus $S_1(N)$ is a subset of $S_2(N)$. Conversely, if there is a rational
prime $p$ such that $p^2|P(x)$, $p>N^{1/2}\geq \max_{p|D} p$,
we obtain from Lemma \ref{lem:trommel} that $\mathfrak{p}^2|P(x)$
for some $\mathfrak{p}$ lying over $p$. Hence $S_2(N)\subset S_1(N)$,
and therefore $S_1(N)=S_2(N)$, for sufficiently large $N$. The statement
follows immediately.
\end{proof}
\begin{prop}\label{prop:cuatro}
Let $K$ be a number field. Let $P\in \mathfrak{O}_K\lbrack x,y\rbrack$ be
a non-constant homogeneous polynomial. Let $P = P_1 P_2 \dotsb P_k$, $P_i$
irreducible in $\mathfrak{O}_K\lbrack x,y\rbrack$.
Then Conjecture $\mathfrak{A}_2(K,P,\delta(N))$ is equivalent
to Conjecture $\mathfrak{A}_2(\mathbb{Q},Q,\delta(N))$, where
the polynomial $Q(x,y)\in \mathbb{Z}\lbrack x,y\rbrack$ 
as the product of the irreducible factors of
$\mathfrak{N}_{K_{P_i}/\mathbb{Q}} (c_i P_i(x,y)) \in \mathbb{Z}\lbrack x\rbrack$, 
$i=1,\dotsb,k$,
where $c_1,\dotsc,c_k$ are constants in $\mathfrak{O}_K$.
\end{prop}
\begin{proof}
Same as that of Proposition \ref{prop:baram}.
\end{proof}
As was pointed out in the introduction, Hypothesis 
$\mathfrak{B}_i(K,P,\eta(N),\epsilon(N))$ is false for some choices of
$K$ and $P$. Thus we cannot hope to reduce it to the case
$K=\mathbb{Q}$ without restrictions. We will, however,
analyse
the situation completely, provided that $K/\mathbb{Q}$
is Galois:
we can then show
$\mathfrak{B}_i(K,P,\eta(N),\epsilon(N))$ to be false in some cases
and equivalent to $\mathfrak{B}_i(K,P,\eta(N),\epsilon(N))$ in all
other cases.

\begin{lem}\label{lem:mitzra}
Let $K$ be a number field. Let $L$ be a finite Galois extension of $K$.
Suppose $\deg(L/K)$ is odd. Then the restriction of $\lambda_L$
to $I_K$ equals $\lambda_K$.
\end{lem}
\begin{proof}
Let $\mathfrak{p}\in I_K$ be a prime ideal. Let $e$ and $f$ be the ramification degree and the inertia degree of $\mathfrak{p}$, respectively.
Write
\[\mathfrak{p} = \mathfrak{P}_1^{e} \dotsb \mathfrak{P}_n^e,\]
where $n$ is the number of primes of $I_L$ lying over $\mathfrak{p}$.
Since $\deg(L/K) = e f n$, both $e$ and $n$ must be odd.
Hence
\[\lambda_L(\mathfrak{p}) = \lambda_L(\mathfrak{P}_1^e \dotsb
\mathfrak{P}_n^e) = (-1)^{n e} = -1 = \lambda_K(\mathfrak{p}) .\]
Since $\lambda_L$ is completely multiplicative, we conclude that
$\lambda_L(\mathfrak{a}) = \lambda_K(\mathfrak{a})$ for all
$\mathfrak{a}\in I_K$.
\end{proof}

Given a non-zero ideal $\mathfrak{m}\in I_K$, we define
$I_K^\mathfrak{m}$ to be the semigroup of ideals prime to $\mathfrak{m}$
and $P_K^{\mathfrak{m}}$ to be the semigroup of principal ideals
$(x)$ with $x\equiv 1 \mo \mathfrak{m}$ and $x$ totally positive.

\begin{lem}\label{lem:mitzrb}
Let $K$ be a number field. Let $L$ be a finite extension of $K$.
Suppose $\deg(L/K)$ is even. Then the restriction of $\lambda_L$
to $\mathfrak{O}_K$ is pliable.
\end{lem}
\begin{proof}
The order $\deg(L/K)$ of $\Gal(L/K)$ is even. Hence there is
an element $\sigma\in \Gal(L/K)$ of order $2$. Let $K'$ be the fixed
field of $\sigma$. Once we show that $\lambda_L|_{\mathfrak{O}_K'}$
is pliable, we will have by Lemma \ref{lem:opres} that
$\lambda_L|_{\mathfrak{O}_K} = 
(\lambda_L|_{\mathfrak{O}_K'})|_{\mathfrak{O}_K}$.

Let $\mathfrak{p}\in I_K'$. Then
\[\lambda_L(\mathfrak{p}) = \begin{cases}
1 &\text{if $\mathfrak{p}$ splits or ramifies,}\\
-1 &\text{if $\mathfrak{p}$ is unsplit.}\end{cases}\]
Let $\mathfrak{m}$ be the conductor of $L/K'$.
Let $H^{\mathfrak{m}} =
(\mathfrak{N}_{L/K'} I_L^{\mathfrak{m}}) P_K^{\mathfrak{m}}$.
By class field theory (see, e.g., \cite{Ne}, p. 428), 
\begin{itemize}
\item $H^{\mathfrak{m}}$ is an open subgroup of 
$I_K^{\mathfrak{m}}$ of index $2$,
\item a
prime ideal $\mathfrak{p}\in I_K^{\mathfrak{m}}$
splits if and only if it lies in $H^{\mathfrak{m}}$.
\end{itemize}
Therefore, given an ideal $\mathfrak{a}\in I_K$, we have
$\lambda_K(\mathfrak{a})=1$ if and only if $\mathfrak{a}_0\in
H^{\mathfrak{m}}$, where we write $\mathfrak{a} =
\mathfrak{a}_{\mathfrak{m}} \mathfrak{a}_{\mathfrak{m},0}$, 
$\mathfrak{a}_{\mathfrak{m}}|\mathfrak{m}^{\infty}$,
$\mathfrak{a}_{\mathfrak{m},0} \in I_K^{\mathfrak{m}}$.
Since $H^{\mathfrak{m}}$ contains $I_K^{\mathfrak{m}}$,
we have that $\lambda_K(\mathfrak{a})$ depends only on
$\mathfrak{a}_{\mathfrak{m},0} P_K^{\mathfrak{m}}$.
Since we can tell $\mathfrak{a}_{\mathfrak{m}}$ from the
coset of $P_K^{\mathfrak{m}}\subset I_K$ in
which $\mathfrak{a}$ lies, we can say that $\lambda_K(\mathfrak{a})$
depends only on $\mathfrak{a} P_K^{\mathfrak{m}}$.

For every real infinite place $v$ of $K$, let $U_{v} = \mathbb{R}^+$.
For every $\mathfrak{p}|\mathfrak{m}$, let 
$U_{\mathfrak{p}} = 1 + \mathfrak{p}^{v_{\mathfrak{p}}(\mathfrak{m})}
\mathfrak{O}_{K_{\mathfrak{p}}}$. 
Let $x$ be a non-zero element of $\mathfrak{O}_K$. 
Suppose we are 
given $x U_{\mathfrak{p}}$ for every $\mathfrak{p}|\mathfrak{m}$
and $x U_{v}$ for every real infinite place $v$.
Then, by the Chinese remainder theorem, we know
$x P_K^{\mathfrak{m}}$. By the above paragraph, we can tell
$\lambda_K(\mathfrak{a})$ from $x P_K^{\mathfrak{m}}$. We conclude
that $\lambda_K$ is pliable with respect to
$\{v,U_v,0\}_{\text{$v$ real}} \cup 
 \{\mathfrak{p},U_{\mathfrak{p}},0\}_{\mathfrak{p}|\mathfrak{m}}$.
\end{proof}
\begin{prop}\label{prop:hwyhwy}
Let $K$ be a finite Galois extension of $\mathbb{Q}$.
Let $P\in \mathfrak{O}_K\lbrack x\rbrack$ be a square-free, non-constant
polynomial. Let $P = P_1 P_2 \dotsb P_k$, $P_i$ irreducible in
$\mathfrak{O}_K\lbrack x\rbrack$. Then 
\[\lambda_K(P(x)) = f(x)\cdot \lambda\left(\mathop{\prod_{i}}_{\text{
$\deg(K/K_{P_i})$ odd}} \mathfrak{N}_{K_{P_i}/\mathbb{Q}} 
(c_i P_i(x))\right),\]
where $f:\mathbb{Z}\to \{-1 , 0,1\}$ is affinely pliable
and $c_1,\dotsb ,c_k$ are constants in $\mathfrak{O}_K$.
\end{prop}
\begin{proof}
Since (a) $\lambda_K$ and $\lambda$ are completely multiplicative, and
(b) the product of affinely pliable functions is affinely pliable,
it is enough to prove the statement for the case of $P$ irreducible.
Choose a non-zero $c\in \mathfrak{O}_K$ such that the leading coefficient
of $c P$ lies in $K_P$. Then every coefficient of $c P$
lies in $K_P$. 

If $\deg(K/K_{P_i})$ is even, Lemma
\ref{lem:mitzrb} gives us that the restriction of $\lambda_K$
to $\mathfrak{O}_{K_P}$ is pliable. By Proposition \ref{prop:simil},
it follows that the map $x\mapsto \lambda_K(c P(x))$ is pliable
on $\mathfrak{O}_K$. Since $\lambda_K(P(x)) = \lambda_K(c) \lambda_K(c P(x))$, we are done.

 Suppose
$\deg(K/K_{P_i})$ is odd.
By Lemma \ref{lem:mitzra}, $\lambda_K(c P(x)) = \lambda_{K_P}(c P(x))$.
Let $D$ be as in Lemma \ref{lem:trommel}. Then
\[\lambda_{K_P}\left(\prod_{\rho(\mathfrak{p}) \notin D} \mathfrak{p}^{
v_{\mathfrak{p}}(c P(x))}\right) = 
\lambda\left(\prod_{p \notin D} p^{v_p(\mathfrak{N}_{K_P/\mathbb{Q}}
 (c P(x)))} \right) ,\]
where, as before, we write $\rho(\mathfrak{p})$ for the rational
prime lying under $\mathfrak{p}$.
Clearly
\[\lambda_{K_P}(c P(x)) = 
\prod_{\rho(\mathfrak{p}) \in D} (-1)^{v_{\mathfrak{p}}(c P(x))}
\cdot
\lambda_{K_P}\left(\prod_{\rho(\mathfrak{p}) \notin D} \mathfrak{p}^{
v_{\mathfrak{p}}(c P(x))}\right)  .
\]
Set $f(x) =  \prod_{\rho(\mathfrak{p}) \in D} (-1)^{v_{\mathfrak{p}}(c P(x))}$. Since there are finitely many prime ideals lying over elements
of $D$, we conclude that $f$ is a product of finitely many affinely
pliable functions, and is thus pliable itself.
\end{proof}
\begin{prop}\label{prop:hwohwo}
Let $K$ be a finite Galois extension of $\mathbb{Q}$.
Let $P\in \mathfrak{O}_K\lbrack x,y\rbrack$ be a square-free, non-constant
homogeneous polynomial. Let $P = P_1 P_2 \dotsb P_k$, $P_i$ irreducible in
$\mathfrak{O}_K\lbrack x,y\rbrack$. Then 
\[\lambda_K(P(x)) = f(x,y)\cdot \lambda\left(\mathop{\prod_{i}}_{\text{
$\deg(K/K_{P_i})$ odd}} \mathfrak{N}_{K_{P_i}/\mathbb{Q}} 
(c_i P_i(x,y))\right),\]
where $f:\mathbb{Z}^2\to \{-1 , 0,1\}$ is pliable
and $c_1,\dotsb ,c_k$ are constants in $\mathfrak{O}_K$.
\end{prop}
\begin{proof}
Same as that of Proposition \ref{prop:hwyhwy}.
\end{proof}
\begin{cor}\label{cor:aga1}
Let $K$ be a finite Galois extension of $\mathbb{Q}$.
Let $P\in \mathfrak{O}_K\lbrack x\rbrack$ be a square-free, non-constant
polynomial. Let $P = P_1 P_2 \dotsb P_k$, $P_i$ irreducible in
$\mathfrak{O}_K\lbrack x\rbrack$. Let 
\[Q(x) =  
\mathop{\prod_{i}}_{\text{
$\deg(K/K_{P_i})$ odd}} \mathfrak{N}_{K_{P_i}/\mathbb{Q}} 
(c_i P_i(x)),\]
where $c_1,\dotsb ,c_k\in \mathfrak{O}_K$ are is in Proposition
\ref{prop:hwyhwy}. Then
\begin{itemize}
\item $\mathfrak{B}_1(K,P,\eta(N),\epsilon(N))$ is equivalent
to $\mathfrak{B}_1(\mathfrak{Q},Q,\eta(N),\epsilon(N))$ if
$Q$ is not of the form $c R^2$, $c\in \mathfrak{O}_K$, 
$R\in \mathfrak{O}_K\lbrack x\rbrack$,
\item $\mathfrak{B}_1(K,P,\eta(N),\epsilon(N))$ is false if
$Q$ is of the form $c R^2$, $c\in \mathfrak{O}_K$, 
$R\in \mathfrak{O}_K\lbrack x\rbrack$.
\end{itemize}
\end{cor}
\begin{proof}
Immediate from Proposition \ref{prop:hwyhwy} and Lemma \ref{lem:strau1}.
\end{proof}
\begin{cor}\label{cor:aga2}
Let $K$ be a finite Galois extension of $\mathbb{Q}$.
Let $P\in \mathfrak{O}_K\lbrack x,y\rbrack$ be a square-free, non-constant
homogeneous polynomial. Let $P = P_1 P_2 \dotsb P_k$, $P_i$ irreducible in
$\mathfrak{O}_K\lbrack x,y\rbrack$. 
Let
\[Q(x,y) = \mathop{\prod_{i}}_{\text{
$\deg(K/K_{P_i})$ odd}} \mathfrak{N}_{K_{P_i}/\mathbb{Q}} 
(c_i P_i(x,y)),\]
where $c_1,\dotsb ,c_k\in \mathfrak{O}_K$ are is in Proposition
\ref{prop:hwohwo}. Then
\begin{itemize}
\item 
$\mathfrak{B}_2(K,P,\eta(N),\epsilon(N))$ is equivalent
to $\mathfrak{B}_2(\mathfrak{Q},Q,\eta(N),\epsilon(N))$ if
$Q$ is not of the form $c R^2$, $c\in \mathfrak{O}_K$, 
$R\in \mathfrak{O}_K\lbrack x,y\rbrack$,
\item $\mathfrak{B}_2(K,P,\eta(N),\epsilon(N))$ is false if
$Q$ is of the form $c R^2$ for some $c\in \mathfrak{O}_K$, 
$R\in \mathfrak{O}_K\lbrack x,y\rbrack$.
\end{itemize}
\end{cor}
\begin{proof}
Immediate from Proposition \ref{prop:hwohwo} and Lemma \ref{lem:strau2}.
\end{proof}

\section{Ultrametric analysis, field extensions and pliability}\label{subs:ultra}
In this appendix, we show how pliable functions arise naturally
in the context of extensions of local fields. While the rest of the present
work does not depend on the following results, the reader might
find that the following instantiation of pliability illuminates the said
concept.

Let $K$ be a field of characteristic zero. Consider a polynomial $f(x)$ with coefficients
in $K((t))$:
\begin{equation}\label{eq:fueq}
f(x) = x^n + a_{n-1}(t) x^{n-1} + a_{n-2}(t) x^{n-2} + \dotsb + a_0(t).
\end{equation}
The Newton-Puiseux method yields fractional power series $\eta_i(t)$, 
$i=1,2,\dotsc ,n$,
\begin{equation}\label{eq:formsols}
\eta_i(t) = c_{k,i} t^{k/l} + c_{k+1,i} t^{(k+1)/l} + \dotsb
\end{equation}
with coefficients in a finite extension $L/K$, such that
\[f(x) = \prod_i (x-\eta_i(t))\]
formally. In particular, if $f(x)$ is irreducible over $\overline{K}((t))$, we have
\begin{equation}\label{eq:floppy}\begin{aligned}
\eta_0(t) &= c_k t^{k/n} + c_{k+1} t^{(k+1)/n} + \dotsb\\
\eta_j(t) &= c_k \omega^{k j} t^{k/n} + c_{k+1} \omega^{(k+1) j} t^{(k+1)/n} + \dotsb,\,\,
1< j< n,\end{aligned}\end{equation}
where $\omega$ is a primitive $n$th root of unity.

We may rephrase this as follows: any finite extension $R$ of $K((t))$ may be embedded in 
$L((t^{1/k}))$ for some positive integer $l$ and some finite extension $L$ of $K$. Regard
$K((t))$ as a local field with respect to the valuation 
\begin{equation}\label{eq:defvt}v_t(c_k t^k + c_{k+1} t^{k+1} + \dotsb) = k\,\,\text{if $c_k\ne 0$} .\end{equation}
What (\ref{eq:floppy}) then implies is that any totally ramified finite Galois extension of $K((t))$
of degree $n$ can be identified with $K((t^{1/n}))$. An unramified finite Galois extension of
$K((t))$ can be written as $L((t))$, where $L$ is the residue field of the extension, and as such
a finite Galois extension of $L$. Hence an arbitrary finite Galois extension $R$ of $K((t))$
can be identified with $L((t^{1/l}))$, where $l$ is a positive integer and $L$ is a finite
Galois extension of $L$.

Assume from now on that $K$ is a $\mathfrak{p}$-adic field. 
Let $\mathcal{C}_g^\infty(K,t)$ be
the ring of power series $\eta(t) \in K\lbrack \lbrack t\rbrack \rbrack$ that converge
in a neighbourhood of $0$. (In other words, $\mathcal{C}_g^\infty(K,t)$ is the ring of germs of analytic
functions around $0$.)
Let $\mathcal{M}_g^\infty(K,t)$ be the field of fractions of 
$\mathcal{C}_g^\infty(t)$. It is a local field with respect to the
 valuation $v_t$ defined in (\ref{eq:defvt}).

Consider $\eta\in K((t))$. By the {\em radius of convergence} $r(\eta)$
of $\eta$ we mean the largest $r\geq 0$ such that 
$t^{-v_t(\eta)} \eta$ converges inside the open ball $B_0(r)$ of
radius $r$ about zero. We can see $\eta$ as an element of 
$\mathcal{M}_g^\infty(K,t)$ if and only if $r(\eta)>0$. Write
\[\eta = c_{-k} t^{-k} + c_{-k+1} t^{-k+2} + \dotsb .\] Then
$r(\eta)$ is positive if and only if $c_j\ll M^j$ for some $M>0$.

While $\mathcal{M}_g^\infty(K,t)$ is not complete with respect to
its valuation $v_t$, it is nevertheless Henselian. A {\em Henselian}
field is one for which Hensel's lemma holds. To see that $\mathcal{M}_g^\infty(K,t)$ is Henselian, it is enough to examine
the algorithm that proves Hensel's lemma in its simplest incarnation.
Let $f = x^n + a_{n-1}(t) x^{n-1} + \dotsb + a_0(t) x^{n-1}$ be
a polynomial with coefficients in $\mathcal{C}_g^\infty(K,t)$; let
$\bar{f} = x^n + a_{n-1}(0) x^{n-1} + \dotsb + a_0(0) x^{n-1}$ be
its reduction to a polynomial with coefficients in the residue field
$K$ of $\mathcal{C}_g^\infty(K,t)$. If $\bar{f}(0)=0$ and 
$\overline{f'}(0)\ne 0$, the Henselian algorithm produces a root
$x(t)\in K((t))$ of $f(x) = 0$ satisfying $\overline{x(t)} = x(0) = 0$.
We must check that the coefficients of the root $x(t)$ thus produced are majorized by some
$M^j$. Since $K$ is non-archimedean, this follows easily from the fact that the coefficients of $a_0,a_1,\dotsc a_{n-1}$ are majorized by some
$M_0^j, M_1^j, \dotsc M_{n-1}^j$. Hence $x(t)\in \mathcal{M}_g^\infty(K,t)$,
and so $\mathcal{M}_g^\infty(K,t)$ is Henselian.

The Newton-Puiseux method for solving (\ref{eq:fueq}) starts with the
coefficients \[a_{n-1}(t), \dotsc ,a_0(t) \in K((t))\] and
manipulates them to produce (\ref{eq:formsols}). These manipulations are
of four kinds: transforming $t$ linearly,
embedding $K((t))$ in $L((t))$, embedding $K((t))$ in $K((t^{1/l}))$ and
expressing a polynomial 
\[x^n + a_{n-1}(t) x^{n-1} + \dotsb + a_0(t),\;\; a_i\in K((t))\] as a product 
\[
(x^{n_1} + \alpha_{n_1-1}(t) x^{n_1-1} + \dotsb + \alpha_0(t))
(x^{n_2} + \beta_{n_2-1}(t) x^{n_2-1} + \dotsb + \beta_0(t)),\;\;
\alpha_i,\beta_i\in K((t)) \]
by means of Hensel's lemma.
It is clear that the every one of the first three operations takes
a series with a non-trivial radius of convergence to a series
with a non-trivial radius of convergence. That the fourth operation
produces $\alpha_i,\beta_i \in \mathcal{M}_g^\infty(K,t)$ when given
$a_i\in \mathcal{M}_g^\infty(K,t)$ follows from the fact that
$\mathcal{M}_g^\infty(K,t)$ is Henselian.

        Thus the formal solutions (\ref{eq:formsols}) in $L((t^{1/l}))$
 to 
\[x^n + a_{n-1}(t)x^{n-1} + \dotsb +a_0(t) = 0\] constructed by
the Newton-Puiseux method lie in fact in $\mathcal{M}_g^\infty(L,t^{1/l})$,
provided that $a_i(t)\in \mathcal{M}_g^\infty(K,t)$. See \cite{DR}
for explicit expressions for the radii of convergence of (\ref{eq:formsols}).

Thanks to this closure property of $\mathcal{M}_g^\infty(K,t)$,
various matters work out much as for $K((t))$. Any finite Galois extension
of $\mathcal{M}_g^\infty(K,t)$ can be identified with $\mathcal{M}_g^\infty(L,t^{1/l})$
for some finite Galois extension $L$ of $K$ and some positive integer $l$; 
if the extension is unramified, it is of the form
$\mathcal{M}_g^\infty(L,t)$; if it is totally ramified, it is of the form 
$\mathcal{M}^\infty(K,t^{1/n})$, where $n$ is the degree of the extension. Since the closure

of $K\lbrack \lbrack t\rbrack \rbrack$ in $L((t^{1/l}))$ is $L\lbrack \lbrack t^{1/l} \rbrack \rbrack$,
the closure of $\mathcal{C}_g^\infty(K,t)$ in $\mathcal{M}_g^\infty(L,t^{1/l})$ is
$\mathcal{C}_g^\infty(L,t^{1/l})$.

Let $t_0\in K$. Define the {\em specialization map}
$\Sp_{t_0}:\mathcal{M}_g^\infty(K,t)\to K$ taking 
$f\in \mathcal{M}_g^\infty(K,t)$ to $f(t_0)$, if $t_0$ is within
the radius of convergence of $f$, and to $0$ otherwise. If 
$R = \mathcal{M}_g^\infty(L,t^{1/l})$ is a finite Galois extension of
$\mathcal{M}_g^\infty(K,t)$, then $\Sp_{t_0}(R) = L(t_0^{1/l})$ for
every $t_0\in K$. Thus \[t\mapsto \Sp_t(R)\] is a map from $K$ to the
set of finite Galois extensions of $K$.

\begin{lem}\label{lem:plione}
Let $K$ be a $\mathfrak{p}$-adic field.
Let $R$ be a finite Galois extension of $\mathcal{M}_g^\infty(K,t)$.
Then the map \[t\mapsto \Sp_t(R)\] is affinely pliable at $0$.
\end{lem}
\begin{proof} We know that $R$ is of the form $\mathcal{M}_g^\infty(L,t^{1/l})$
for some positive integer $l$ and some finite Galois extension $L$ of $K$.
Let $U = 1 + \pi_K^{2 l + 1} \mathfrak{O}_K$. Suppose $t, t'\in K^*$ belong
to the same coset of $U$. Then $t/t'\in U$, and thus $v_K(t/t'-1)\geq 2 l +1$.
By Hensel's lemma it follows that $x^l = t/t'$ has a root $x_0\in K$.
Choose $l$th roots $t^{1/l}$, ${t'}^{1/l}$ of
$t$ and $t'$ such that $t^{1/l}/{t'}^{1/l} = x_0$. Then
$L(t^{1/l}) = L({t'}^{1/l})$. Therefore the map
\[t\mapsto \Sp_t(R)\]
is affinely pliable at zero.
\end{proof}
\begin{lem}\label{lem:lemf4}
Let $K$ be a $\mathfrak{p}$-adic field. Let $a_0,a_1,\dotsc ,a_{n-1}\in 
\mathcal{M}_g^\infty(K,t)$. Let $\mathcal{M}_g^\infty(L,t^{1/l})$ be
the splitting field of \begin{equation}\label{eq:eqthisone}
x^n + a_{n-1} x^{n-1} + \dotsb + a_0 = 0\end{equation} over
$\mathcal{M}_g^\infty(K,t)$. Let $\eta_1,\eta_2,\dotsc ,\eta_n \in
\mathcal{M}_g^\infty(L,t^{1/l})$ be the roots of (\ref{eq:eqthisone}).
Then there is an $r>0$ such that $\eta_1(t_0),\dotsc ,\eta_n(t_0)$
converge and \[\Sp_{t_0}(\mathcal{M}_g^\infty(L,t^{1/l})) = 
K(\eta_1(t_0),\dotsc,\eta_n(t_0))\] for $t_0\in B_{K,0}(r)-\{0\}$.
\end{lem}
\begin{proof}
Clearly $K(\eta_1(t_0),\dotsc,\eta_n(t_0))\subset \Sp_{t_0} (
\mathcal{M}_g^\infty(L,t^{1/l}))$ for $t$ within the radii of convergence
of $\eta_1,\dotsc \eta_n$. To prove 
$\Sp_t(\mathcal{M}_g^\infty(L,t_0^{1/l})) \subset K(\eta_1(t_0),\dotsc,
\eta_n(t_0))$, it is enough to show that 
\[K(\eta_1(t_0),\dotsc ,\eta_n(t_0))\]
contains a basis of $L$ as a vector space over $K$ as well as an $l$th
root of $t_0$. Let $s_0$ be an $l$th root of $t$ and let $s_1,\dotsc ,s_m$
form a basis of $L$ over $K$. Consider $s_0,\dotsc ,s_m$ as elements of
$\mathcal{M}_g^\infty(L,t^{1/l})$. Since $\mathcal{M}_g^\infty(L,t^{1/l})
= (\mathcal{M}_g^\infty(K,t))(\eta_1,\dotsc ,\eta_n)$, one can reach
$s_i$ after a finite number of additions, substractions, multiplications
and divisions starting from $\eta_1, ...,\eta_n$ and a finite number
of elements of $\mathcal{M}_g^\infty(K,t)$. Each of this operations takes
two series with positive radii of convergence to a series with a positive
radius of convergence. Let $r$ be the minimum of all the radii of convergence
of the finitely many objects appearing in the process. Then, for
$t_0\in B_{K,0}(r)$, each operation $\spadesuit$ takes two series $\rho_1,\rho_2\in
\mathcal{M}_g^\infty(L,t^{1/l})$ to a series $\rho_1\,\spadesuit\, \rho_2 \in
\mathcal{M}_g^\infty(L,t^{1/l})$ taking the value $\rho_1(t_0)\,\spadesuit\,
\rho_2(t_0)$ at $t_0$. Since $\eta_1(t_0),\dotsc,\eta_n(t_0)\in
K(\eta_1(t_0),\dotsc,\eta_n(t_0))$ and $K(\eta_1(t_0),\dotsc,\eta_n(t_0))$
is closed under $\spadesuit = +,-,*,/$, it follows that $K(\eta_1(t_0),\dotsc
,\eta_n(t_0))$ contains $s_0,s_1,\dotsc ,s_n$. Hence 
$\Sp_t(\mathcal{M}_g^\infty(L,t_0^{1/l})) \subset K(\eta_1(t_0),\dotsc,
\eta_n(t_0))$.
\end{proof}

Now let $a_0,a_1,\dotsc,a_{n-1}$ be rational functions on $t$ with coefficients
 in $K$. For every $t_0\in K$, 
\[b_{t_0,0}(t)=a_0(t+t_0), b_{t_0,1}(t)=a_1(t+t_0), \dotsc ,b_{t_0,n-1}(t)=
a_{n-1}(t+t_0)\]
can be seen as elements of $\mathcal{M}_g^\infty(K,t)$. Moreover,
\[b_{\infty,0}(t) = a_0(1/t),\dotsc,b_{\infty,n-1} = a_{n-1}(1/t)\] can be
seen as elements of $\mathcal{M}_g^\infty(K,t)$, as they are
rational functions on $t$.

\begin{prop}\label{prop:propo}
Let $K$ be a $\mathfrak{p}$-adic field. Let $a_0, a_1,\dotsc,a_{n-1} \in K(t)$. Define a function $\mathfrak{S}$ from $K$ to the
set of finite Galois extensions of $K$ as follows: for $t_0\in K$, 
let $\mathfrak{S}(t_0)$ be the splitting field of 
$x^n + a_{n-1}(t_0) x^{n-1} + \dotsb + a_0(t_0) = 0$ over $K$ if
$a_0(t_0),a_1(t_0),\dotsc,a_{n-1}(t_0)$ are finite; let $\mathfrak{S}(t_0)$
be $K$ otherwise. Then $\mathfrak{S}$ is affinely pliable.
\end{prop}
\begin{proof}
Let $t_0\in \mathbb{P}^1(K)$. By Lemma \ref{lem:lemf4}, there are a positive
integer $l$, a finite Galois extension $L$ of $K$ and an open ball $V$
around zero such that, for all $t\in V-\{0\}$,
\[K(\eta_{t_0,1}(t),\dotsc ,\eta_{t_0,n}(t)) = \Sp_t(\mathcal{M}_g^\infty(L,t^{1/l})),\]
where $\eta_{t_0,1}(t),\dotsc ,\eta_{t_0,n}(t)$ are the roots of
\[x^n + b_{t_0,n-1}(t) x^{n-1} + \dotsb + b_{t_0,0} = 0.\]
By Lemma \ref{lem:plione}, $\Sp_t(M_g^\infty(L,t^{1/l}))$ is affinely 
pliable. Therefore
the restriction of $K(\eta_1(t),\dotsc ,\eta_n(t))$ to $V$ is affinely 
pliable at $0$.

It follows from the definition of $b_{t_0,{n-1}},\dotsc ,b_{t_0,0}$ that
\[K(\eta_1(t),\dotsc,\eta_n(t)) = \begin{cases} 
\mathfrak{S}(t+t_0) &\text{if $t_0\ne \infty$}\\
\mathfrak{S}(1/t) &\text{if $t_0=\infty$.}\end{cases}\]
Hence, for every $t_0\ne \infty$ there is an open ball $V_{t_0}$ around $t_0$ such that
$\mathfrak{S}(t)|_{V_{t_0}}$ is affinely pliable at $t_0$.
 Moreover, $\mathfrak{S}(1/t)|_{V_*}$ is affinely pliable at $0$ for some open ball
$V_*$ around $0$. 
This is the same as saying that there is an open subgroup
$U$ of $K$ such that $\mathfrak{S}(1/t)$ depends only on $t U$ for 
$t\in V_*-\{0\}$. Since $U$ is a group, the map $t U\to t^{-1} U$ is 
well-defined and bijective. Hence depending only on $t U$ is the same
as depending only on $(1/t) U$. Therefore we can say that $\mathfrak{S}$
depends only on $(1/t) U$ for $t\in V_\infty$; in other words, 
$\mathfrak{S}(t)$ depends only on $t U$ for $t$ in a neighborhood 
$V_\infty=1/V_*$
of infinity. Thus $\mathfrak{S}(t)$ is affinely
pliable at $0$ when restricted to neighbourhood
$V_\infty$  of infinity.

Since $\mathbb{P}^1(K)$ is compact, it is covered by a finite subcover of
$\{V_{t_0}\}_{t_0\in \mathbb{P}^1(K)}$. Let the subcover be $\{V_s\}_{s\in S}$,
$S$ a finite subset of $\mathbb{P}_1(K)$. By the above $\mathfrak{S}|_{V_s}$ for
every $s\in S$. Since $V_s$ is a ball, its characteristic function 
$t\mapsto \lbrack t\in V_s\rbrack$ is affinely pliable. Hence
\[\mathfrak{S}(t) = \sum_{s\in S} \lbrack t\in V_s\rbrack ((\mathfrak{S}|_{V_s})(t))\]
is affinely pliable.
\end{proof}

Given Proposition \ref{prop:propo} and Lemma \ref{lem:STgood},
 it is a simple matter to show that, given an elliptic curve
$\mathcal{E}$ over $K(t)$,
the map taking an element $t\in K$ to the minimal extension over which
$\mathcal{E}(t)$ acquires good reduction is affinely pliable.

\input{apptrace}

%% file: apptrace.tex
\section{The root number in general}
Let $H_k^*(N)$ be the set of newforms of even positive weight $k$ on
 $\Gamma_0(N)$. Every newform $f\in H_k^*(N)$ has a root number
$\eta_f$. It is a well-known fact that the average of the root numbers
of the elements of $H_2^*(N)$ tends to zero as $N$ goes to infinity.
As some suboptimal bounds
on the error term are labouriously derived in the recent literature, it may
be worthwhile to point out that there is an exact expression for the 
total $\sum_f \eta_f$ of the root numbers of newforms $f\in H_2^*(N)$. 
This expression can be bounded easily from above and below.

Let $W_N$ be the canonical involution for level $N$:
\[W_N: g \mapsto g|_{w_N},\]
where $w_N$ is the matrix 
$\left(\begin{matrix} 0 & -1 \\ N & 0\end{matrix}\right)$.
Every newform $f\in H_k^*(N)$ is an eigenfunction of $W_N$ with 
eigenvalue $\eta_f$.

Let $S_k(N)$ be the space of cusp forms of weight $k$ on $\Gamma_0(N)$. 
For $L M = N$, $f\in H_k^*(M)$,
let $S_k(L;f)$ be the space of linear combinations of $\{f_{|\ell} : \ell | L\}$, where 
\[f_{|\ell}(z) = \ell^{-k/2} f(\ell z) .\]
Since the functions $f_{|\ell}$ for fixed $f$ are linearly independent, 
$\{f_{|\ell} : \ell | L\}$ is actually a basis for $S_k(L;f)$.
By (\cite{AL}, Thm 5) we have
\[S_k(N) = \bigoplus_{L M = N} \bigoplus_{f\in H_k^*(M)} S_k(L; f)\]
as a direct sum of orthogonal Hilbert spaces under the Petersson inner product on $S_k(N)$.

Consider an $f\in H_k^*(L; f)$. For $\ell|L$,
\begin{equation}\begin{aligned}
(W_N f_{|\ell})(z) &= (z \sqrt{N})^{-k} f_{|\ell}\left(\frac{-1}{N z}\right)
= (z \sqrt{N})^{-k} \ell^{k/2} f\left(\frac{-1}{ (M L / \ell) z}\right) \\
&= \left(\frac{L}{\ell}\right)^{k/2} (W_M f) \left(\frac{L}{\ell} z\right)
= \eta_f \left(\frac{L}{\ell}\right)^{k/2} f\left(\frac{L}{\ell} z\right)
= \eta_f f_{|(L/\ell)}(z) .\end{aligned}\end{equation}
Hence the trace of $W_N$ on $S_k(L;f)$ is $\eta_f$ if $L$ is a perfect square
and zero otherwise. Summing over all $f\in H_k^*(M)$ we obtain
\begin{equation}
\Tr(W_N,S_k(N)) = \mathop{\sum_{L M = N}}_{\text{$L$ a square}} 
 \sum_{f\in H_k^*(M)} \eta_f .\end{equation}
By M\"{o}bius inversion
\begin{equation}\label{eq:fromtrace}
 \sum_{f\in H_k^*(N)} \eta_f = \mathop{\sum_{R^2 M = N}} \mu(R) 
 \Tr(W_M,S_k(M)) .\end{equation}

Now consider the curves $\Gamma_0(N)\backslash\mathbb{H}$ and 
$(\Gamma_0(N)\cdot W_N)\backslash\mathbb{H}$, 
where $\Gamma_0(N) * W_N$ is the group
obtained by adjoining $W_N$ to $\Gamma_0(N)$. 
Let $S_k(\Gamma_0(N) * W_N)$ be the set of cusp forms of weight $k$ on
$(\Gamma_0(N) * W_N)\backslash\mathbb{H}$. Write $s_k(\Gamma_0(N))$ and
$s_k(\Gamma_0(N) * W_N)$ for the cardinalities of $S_k(N)$ and 
$S_K(\Gamma_0(N) * W_N)$, respectively. Our goal is to compute
\[\Tr(W_N,S_k(N)) = 2 s_k(\Gamma_0(N) * W_N) - s_k(\Gamma_0(N)) .\]

By Gauss-Bonnet,
\[\frac{1}{2\pi} \Vol(\Gamma_0(N)\backslash\mathbb{H}) = 2 g - 2 + m +
\sum_{i=1}^r (1 - 1/e_i),\]
where $g$ is the genus of $\Gamma_0(N)\backslash\mathbb{H}$, $m$ is the
number of its inequivalent cusps and $e_1$, $e_2$,\dots\: are the orders
of its inequivalent elliptic points. Similarly,
\[\frac{1}{2\pi} \left(\frac{1}{2} \Vol(\Gamma_0(N)\backslash\mathbb{H})\right)
= 
\frac{1}{2\pi} \Vol((\Gamma_0(N) * W_N)\backslash\mathbb{H}) = 
2 g_0 - 2 + m_0 + \sum_{i=1}^{r'} (1 - 1/{e_i'}),\]
where $g_0$ is the genus of $(\Gamma_0(N) * W)\backslash\mathbb{H}$, $m_0$
is the number of its inequivalent cusps and $e_1'$, $e_2'$,\dots\: are the orders
of its inequivalent elliptic points. The relations among $m$, $m_0$,
$e_i$ and $e_i'$ were written out by Fricke (\cite{Fr}, p. 357--367). 
They are as follows. Assume $N>4$.
The involution $W_N$ then 
matches pairs of distinct equivalence classes of cusps
of $\Gamma_0(N)\backslash\mathbb{H}$; therefore, $m = 2 m_0$. 
The equivalence classes of elliptic points of 
$\Gamma_0(N)\backslash\mathbb{H}$ are also paired by $W_N$, which at the
same time introduces $\epsilon_N h(- 4 N)$ new elliptic points, all of order
$2$. Here
\begin{equation}\label{eq:morm}
\epsilon_N = \begin{cases} 2 &\text{if $N\equiv 7 \mo 8$,}\\
4/3 &\text{if $N\equiv 3 \mo 8$,}\\
1 &\text{otherwise,}\end{cases}\end{equation}
and $h(-4 N)$ is the number of equivalence classes of primitive, positive
definite binary quadratic forms of discriminant $-4 N$. Hence
\[\sum_{i=1}^{r'} (1 - 1/e_i') = \frac{1}{2} \sum_{i=1}^r (1-1/e_i) + 
\frac{1}{2} \epsilon_N h(-4 N) .\]
For $k=2$, we have $s_k(\Gamma_0(N)) = g$ and $s_k(\Gamma_0(N) * W) = g_0$.
Hence
\[\begin{aligned}
\Tr(W_N,S_2(N)) = 2 g_0 - g &= 
 \left(\frac{1}{2\pi} (\frac{1}{2} \Vol(\Gamma_0(N)\backslash\mathbb{H})) + 
2 - m_0 - \sum_{i=1}^{r} (1 - 1/e_i)\right) \\ &-
 \frac{1}{2} \left(\frac{1}{2\pi} \Vol(\Gamma_0(N)\backslash\mathbb{H}) +
2 - m - \sum_{i=1}^{r'} (1 - 1/e_i')\right) \\
&= 1 - \frac{1}{2} \epsilon_N h(- 4 N) ,\end{aligned}\]
as was first pointed out by Fricke (op. cit.). For $k>2$, by Riemann-Roch,
\[\begin{aligned}
s_k(\Gamma_0(N)) &= (k-1) (g-1) + \left(\frac{k}{2} - 1\right) m +
\sum_{i=1}^r \lfloor k (e_i-1)/{2 e_i} \rfloor\\
s_k(\Gamma_0(N) * W_N) &= (k-1) (g_0-1) + \left(\frac{k}{2} - 1\right) m_0
+ \sum_{i=1}^{r'} \lfloor k (e_i' - 1)/{2 e_i'} \rfloor \end{aligned}\]
(see, e.g., \cite{Shi}, Thm 2.24). Hence
\[\begin{aligned}
\Tr(W_N,S_k(N)) &= 2 s_k(\Gamma_0(N) * W_N) - s_k(\Gamma_0(N)) \\
              &= (k-1) (2 (g_0 - 1) - (g-1)) + 
\left(\frac{k}{2} - 1 \right) (2 m_0 - m) \\ &+ 
\left(2 \sum_{i=1}^{r'} (1 - 1/e_i') - 
\sum_{i=1}^{r} (1 - 1/e_i) \right) \\
&= (k-1) (2 g_0 - g - 1) + 2 \lbrack k/4\rbrack\, \epsilon_N h(-4 N) \\
&= (k-1) (- \frac{1}{2} \epsilon_N h(- 4 N)) + (k/2) (\epsilon_N h(- 4 N))
- 2 {k/4} (\epsilon_N h( - 4 N))\\
&= \begin{cases} \frac{1}{2} \epsilon_N h(- 4 N) &\text{if $4|k$}\\
-\frac{1}{2} \epsilon_N h(- 4 N) &\text{if $4\nmid k$ .}\end{cases}
\end{aligned}\]
We invoke (\ref{eq:fromtrace}) and conclude that
\[\sum_{f\in H_k^*(N)} \eta_f = \sum_{R^2 M = N} \mu(R) \cdot
\begin{cases} (1 - \frac{1}{2} \epsilon_M h(-4 M)) &\text{if $k=2$,}\\
\frac{1}{2} \epsilon_M h(- 4 M) &\text{if $k>2$, $4|k$,}\\
-\frac{1}{2} \epsilon_M h(- 4 M) &\text{if $k>2$, $4\nmid k$,}\end{cases}\]
provided $N$ is not of the form $R^2$, $2 R^2$, $3 R^2$ or $4 R^2$ for
some square-free integer $R$. Here, as usual,
$\epsilon_N$ is as in (\ref{eq:morm}). 

It is a simple consequence of Dirichlet's formula for the class number that
\[h(d)\ll |d|^{1/2} \log |d| \log \log |d|\]
for any negative $d$ (see, e.g., \cite{N}, p. 254). 
Therefore
\begin{equation}\label{eq:strbound}\begin{aligned}
\left|\sum_{f\in H_k^*(N)} \eta_f\right| 
 &\ll  N^{1/2} \log N \log \log N \prod_{p^2|n} (1 + 1/p)\\
&\ll N^{1/2} \log N (\log \log N)^2 .
\end{aligned}\end{equation}
By Siegel's theorem,
\[h(d)\gg |d|^{1/2-\epsilon} .\]
Hence, for any square-free $N$,
\begin{equation}
\left|\sum_{f\in H_k^*(N)} \eta_f\right| 
 \gg N^{1/2-\epsilon} .\end{equation}

We may finish by commenting on the special cases $N = R^2,\, 2 R^2,\, 3 R^3$,
or, more precisely on the trace $\Tr(W_N,S_k(N))$ for $N = 1, 2, 3$.
For those values of $N$, 
the genera of $\Gamma_0(N)\backslash\mathbb{H}$ and
$(\Gamma_0(N) * W_N)\backslash \mathbb{H}$
are zero. An explicit computation
by means of Riemann-Roch gives
\[\begin{aligned}
\Tr(W_N,S_k(N))  &= \lfloor k/12\rfloor - 1  
&\text{if $N=1$, $k\equiv 2 \mo 12$,}\\
\Tr(W_N,S_k(N))  &= \lfloor k/12\rfloor &\text{if $N=1$, 
$k\not\equiv 2 \mo 12$,}\\
\Tr(W_N,S_k(N))  &= 3 \lfloor k/4 \rfloor - 1 & \text{if $N=2$,}\\
\Tr(W_N,S_k(N))  &= 1 - 3 \{k/3\} & \text{if $N=3$,}
\end{aligned}\]
for $k>2$. (The fact that the genera are zero gives us that
$S_k(N)$ is empty for $k=2$, $N=1,2,3$.)
For $N = R^2,\, 2 R^2$, there is a term of
$\lfloor k/12\rfloor$, resp. $3 \lfloor k/4\rfloor$, which 
dominates
all other terms when $k$ grows more rapidly than $N$. 
For all other $N$, including $N = 3 R^2$,
the bound is (\ref{eq:strbound}), which does not depend on $k$.

%% file: appb.tex
\chapter{Addenda on the parity problem}
\section{The average of $\lambda(x^2 + y^4)$}\label{sec:gebouw}
We prove in this section that the Liouville function averages to zero
over the integers represented by the polynomial $x^2 + y^4$. This is
the same polynomial for which Friedlander and Iwaniec first broke
parity (\cite{FI1}, \cite{FI2}). As $x^2 + y^4$ is not homogeneous,
the results in this section have no apparent bearings on the root numbers 
of elliptic curves. The interest in studying $x^2 + y^4$ resides mainly
in the implied opportunity to test the flexibility of the basic
Friedlander-Iwaniec framework.

As we will see, \cite{FI1} can be used without any modifications; only
\cite{FI2} must be rewritten. We will let $\alpha$ be the Liouville
function or the Moebius function: $\alpha = \lambda$ or $\alpha = \mu$.

\subsection{Notation and identities}
By $n$ we shall always mean a positive integer, and by $p$ a prime.
As in \cite{FI2}, we define 
\[f(n\leq y) = \begin{cases} f(n) &\text{if $n\leq y$}\\
  0 &\text{otherwise,}\end{cases}\]
\[f(n>y) = \begin{cases} f(n) &\text{if $n>y$}\\
  0 &\text{otherwise.}\end{cases}\]
Let \[P(z) = \mathop{\prod_{p\leq z}}_{\text{$p$ prime}} p.\]
For any $n$,
\begin{equation}\label{eq:trot}
f(n>y) = \mathop{\sum_{b c|n}}_{\gcd(n/c,P(z))=1} \mu(b) f(c>y)
.\end{equation}
Write
\[\sum_* \dotsb\;\;\;\;\text{for}\;\; 
\mathop{\sum_{b c|n}}_{\gcd(n/c,P(z))=1} \dotsb\]
Then
\[\begin{aligned}
\sum_* \mu(b) \alpha(c>y) &= \sum_* \mu(b\leq y) \alpha(c>y) +
 \sum_* \mu(b>y) \alpha(c>y) \\
 &= \sum_* \mu(b\leq y) \alpha(c) - \sum_* \mu(b\leq y) \alpha(c\leq y) 
+ \sum_* \mu(b>y) \alpha(c>y) .\end{aligned}\]
Let $w>y$. Proceed:
\[\label{eq:gold}\begin{aligned}
\sum_* \mu(b) \alpha(c>y) &= \sum_* \mu(b\leq y) \alpha(c) 
- \sum_* \mu(b\leq y) \alpha(c\leq y) \\
&+\sum_* \mu(y<b<w) \alpha(c>y) + \sum_* \mu(b>w) \alpha(y<c<w) 
\\ &+ \sum_* \mu(b\geq w) \alpha(c\geq w) .
\end{aligned}\]
We denote the summands on the right side of (\ref{eq:gold}) by
$\beta_1(n)$, $\beta_2(n)$, $\beta_3(n)$, $\beta_4(n)$ and
$\beta_5(n)$.

If $\alpha=\mu$, then, by M\"obius inversion,
\begin{equation}\label{eq:ros1}
\beta_1(n) = 
\sum_* \mu(b\leq y) \alpha(c) = 
\mu(n/\gcd(n,P(z)^{\infty}) \leq y)\, \mu(\gcd(n,P(z))^{\infty}) ,
\end{equation}
whereas, if $\alpha = \lambda$,
\begin{equation}\label{eq:ros2}
\beta_1(n) = 
\sum_* \mu(b\leq y) \alpha(c) = 
\mu(n/\gcd(n,P(z)^{\infty}) \leq y)\, \lambda(\gcd(n,P(z))^{\infty}) .
\end{equation}
Clearly
\[
\beta_2(n) =  \mathop{\sum_{b\leq y}}_{\gcd(b,P(z))=1} 
\sum_{c\leq y}
\mu(b) \alpha(c) \mathop{\mathop{\sum_d}_{b c d = n}}_{\gcd(d, P(z))=1}
1 .\]
If $n < w^2 z$, then
\begin{equation}\label{eq:peas}\beta_5(n) = 
\mathop{\sum_{b c|n}}_{\gcd(n/c,P(z))=1} 
\mu(b\geq w) \alpha(c\geq w) =
\mathop{\sum_{b c = n}}_{\gcd(b,P(z))=1} 
\mu(b\geq w) \alpha(c\geq w) ,\end{equation}
as $\gcd(n/c,P(z))=1$ implies that either $w=1$ or $w>z$, and the latter
possibility is invalidated by $b c d = n$, $b>w$, $c>w$, $n\leq w^2 z$.

Let us be given a sequence $\{a_n\}_{n=1}^{\infty}$ of non-negative
real numbers. For $j=1,\dotsb,5$, we write
\begin{equation}\label{eq:herod}A(x) = \sum_{n=1}^x a_n,\;\;\;
A_d(x) = \mathop{\sum_{1\leq n\leq x}}_{d|x} a_n,\;\;\;
S_j(x) = \sum_{n=1}^x \beta_j(n) a_n .\end{equation}
We will regard $y$, $w$ and $z$ as functions
of $x$ to be set later. For now, we require that $w(x)^2 z(x) > x$.
We have
\[\label{eq:jochaa}
\sum \alpha(n) a_n = 
\sum_{n=1}^x \alpha(n\leq y) a_n + \sum_{j=1}^5 S_j(x) .\]
\subsection{Axioms}
Let $\{a_n\}_{n=1}^\infty$, $a_n$ non-negative, be given. We let
$A(x)$ and $A_d(x)$ be as in (\ref{eq:herod}). We
assume the crude bound
\begin{equation}\label{eq:gargm}
A_d(x) \ll d^{-1} \tau^{c_1}(d) A(x)\end{equation}
uniformly in $d\leq x^{1/3}$, where $c_1$ is a positive constant. We
also assume we can express $A_d$ in the form
\begin{equation}\label{eq:garg0}
A_d(x) = g(d) A(x) + r_d,\end{equation}
where
\[\text{$g:\mathbb{Z}^+\to \mathbb{R}_0^+$ is a multiplicative function,}\]
\begin{equation}\label{eq:garg1}
0\leq g(p)<1,\;\; g(p)\ll p^{-1},
\end{equation} 
\begin{equation}\label{eq:garg2}
\sum_{p\leq x} g(p) = \log \log x + c_2 + O((\log x)^{-1}),
\end{equation} 
\begin{equation}\label{eq:garg3}
\sum_{d\leq D(x)} |r_d(x)| \ll A(x) (\log x)^{-C_1},
\end{equation} 
where 
\begin{equation}
x^{2/3} < D(x) < x,
\end{equation} 
and $C_1$ is a sufficiently large constant ($C_1\leq 65\cdot 2^{c_1} + 4$).
We also assume the following bilinear bound:
\begin{equation}\label{eq:rotstern}
\sum_m \left|\mathop{\mathop{\sum_{N<n<2 N}}_{m n\leq x}}_{
\gcd(n, m P(z)) = 1} \mu(n) a_{m n}\right| \leq A(x) \cdot
(\log x)^{-C_2}\end{equation}
for every $N$ with
\[y(x)<N<w(x) ,\]
where 
\[y(x)\ll D^{1/2}(x) N^{-\epsilon},\;\;\;\log (x^{1/2} w^{-1}(x)) =
o(\log x/C_3 \log \log x),\]
and $C_2$ and $C_3$ are sufficiently large constants.
In \cite{FI2}, conditions (\ref{eq:gargm})--(\ref{eq:garg3}) appear 
(sometimes in stricter forms)
as (1.6), (1.9), (R) and (R1), respectively. Condition (\ref{eq:rotstern})
is a special case of $(\mathrm{B}^*)$ in \cite{FI2} (the case corresponding
to $C=1$, in the notation of the said paper). All of these conditions
are proven for
\[a_n = \{(x,y)\in \mathbb{Z}^2 : x^2 + y^4 = n\}\] 
in \cite{FI1}. Specifically, (\ref{eq:gargm})--(\ref{eq:garg3}) are proven
in \cite{FI1}, section 3, and the rest of \cite{FI1} is devoted to proving
$(\mathrm{B}^*)$. The parameters $D(x)$ and $w(x)$ are given by
\begin{equation}\label{eq:greg}
D\gg x^{2/3-\epsilon},\;\;\;w(x)\gg x^{1/2} (\log x)^{C_4} .\end{equation}
The constants $C_1,\dotsc , C_4$ can be arbitarily large. Notice that
\[x^{3/4} \ll A(x) \ll x^{3/4} .\]
\subsection{Estimates}
We will bound each of $S_j(x)$, $1\leq j\leq 5$. The term 
$S_1(x)$ can be bounded easily as in Lemma \ref{lem:strawb}. Let us bound
$S_2(x)$.
Assume $\log z =  O(\log x/(2 C_3 \log \log x))$. Then
\[z^9\ll D y^{-2},\;\;\;\log D/\log z \gg 2 C_3 \log \log x .\]
It follows that we can use a fundamental lemma (a standard 
formulation of a small sieve). We obtain:
\[\mathop{\sum_d}_{\gcd(d,P(z))=1} a_{b c d} = g(b c) (1 + O((\log x)^{-2 C_3}))
+ O\left(\mathop{\sum_{d\leq D}}_{b c|d} |r_d(x)|\right).\]
Hence, by (\ref{eq:garg2}) and (\ref{eq:garg3}),
\[\begin{aligned}
S_2(n) &= \mathop{\sum_{b\leq y}}_{\gcd(b,P(z))=1} \sum_{c\leq y}
\mu(b) \alpha(c) \mathop{\mathop{\sum_d}_{b c d = n}}_{\gcd(d,P(z))=1} 1\\
&= \mathop{\sum_{b\leq y}}_{\gcd(b,P(z))=1} \sum_{c\leq y} \mu(b) \alpha(c)
 g(b c) (1+ O((\log x)^{- 2 C_3})) A(x) \\
&+ O(\sum_d \tau_3(d) |r_d(x)|)\\
&=  \mathop{\sum_{b\leq y}}_{\gcd(b,P(z))=1} \sum_{c\leq y} \mu(b) \alpha(c)
 g(b c) A(x) + O(A(x) (\log x)^{-C_5}),\end{aligned}\]
where $C_5$ is a large constant. Note that (\ref{eq:garg3}) implies
\[\mathop{\sum_{b\leq y}}_{\gcd(b,P(z))=1} \sum_{c\leq y} \mu(b) \alpha(c)
g(b c) \ll A(x) (\log x)^{-5}.\]
See \cite{FI2}, (2.4). 

To bound $S_3(x)$, a simple application of the bilinear condition
(B.1.9) will suffice:
\[|S_3(x)| = \left|\mathop{\sum_{b, c, d}}_{\gcd(b d,P(z))=1} 
\mu(y<b<w) \alpha(c>y) \right| \leq \sum_m \tau(m) \left|
\mathop{\mathop{\sum_{y<n<w}}_{m n\leq x}}_{\gcd(n,P(z))=1} \mu(n)
a_{m n} \right| .\]
Since $n$ has no small factors, the condition $\gcd(n,m)=1$ may be added
with a total change of at most $O(A(x) (\log x)/z)$. 
The factor $\tau(m)$ may be 
extracted as in \cite{FI2}, p 1047. We obtain
\[S_3 \ll A(x) (\log x)^{-C_6} + A(x)/z .\]
The term $S_4$ can be treated in the same way, with the proviso
that $\alpha$ must be replaced by $\mu$. This replacement induces a total
change of at most $O(A (x) (\log x)/z)$.

All terms up to now have contributed at most $O(A(x) ((\log x)^{-5} +
(\log x)/z)$. One term remains, namely, $S_5$. 
By (\ref{eq:peas}), 
 \[S_5(n) = \mathop{\sum_{b c = n}}_{\gcd(b,P(z))=1} \mu(b\leq w) 
\alpha(c\leq w) .\]
Hence
\[
\mathop{\mathop{\sum_{w\leq x\leq x w^{-1}}}_{\gcd(b,P(z))=1}}_{b c = n}
 1 = \mathop{\sum_{w\leq b\leq x w^{-1}}}_{\gcd(b,P(z))=1} g(b) A(x) +
O\left(\sum_{d\leq x w^{-1}} |r_d(x)|\right) .\]
By (\ref{eq:garg2}) and a fundamental lemma,
\[\mathop{\sum_{w\leq b\leq x w^{-1}}}_{\gcd(b, P(z))=1} g(b) \sim
\frac{1}{\log z} (\log x w^{-1} - \log w) = \frac{\log x w^{-2}}{\log z}
\ll \frac{\log \log x}{\log z} .\]
We are given $w(x)\gg x^{1/2} (\log x)^{-C_4}$; see (\ref{eq:greg}).
Set \[z(x) = e^{\log x/C_3 \log \log x} .\] Then
\[\mathop{\sum_{w\leq b\leq x w^{-1}}}_{\gcd(b, P(z))=1} g(b) \ll
 \frac{(\log \log x)^2}{\log x} A(x) .\]
Hence
\[\sum_n \alpha(n) a_n = \sum_{j=1}^5 S_j(x) + O(A(y)) \ll 
\frac{(\log \log x)^2}{\log x} A(x) ,\]
as was desired. We have proven

\begin{thm} Let $\alpha = \mu$ or $\alpha=\lambda$. Then
\[\mathop{\sum_{a\geq 1} \sum_{b\geq 1}}_{a^2 + b^4 \leq x} \mu(a^2 + b^4)
\ll \left(\mathop{\sum_{a\geq 1} \sum_{b\geq 1}}_{a^2 + b^4 \leq x} 1\right)
\cdot \frac{(\log \log x)^2}{\log x}  \ll
x^{3/4} \frac{(\log \log x)^2}{\log x} .\]
\end{thm}